\documentclass{amsart}

\usepackage{amsmath,amsfonts,amssymb,amsthm,mathtools}

\usepackage{verbatim}
\usepackage{amsmath,amsfonts}
\usepackage[mathscr]{euscript} 
\usepackage{amsthm}
\usepackage{blkarray,bigstrut} 
\usepackage{url}

\usepackage{graphicx} 

\usepackage{enumitem} 

\usepackage{hyperref} 
\hypersetup{colorlinks} 

\usepackage[colorinlistoftodos]{todonotes}

\RequirePackage{cleveref}
\usepackage{hypcap}
\hypersetup{colorlinks=true, citecolor=darkblue, linkcolor=darkblue}
\definecolor{darkblue}{rgb}{0.0,0,0.7}
\newcommand{\darkblue}{\color{darkblue}}

\definecolor{darkred}{rgb}{0.68,0,0}
\newcommand{\darkred}{\color{darkred}}

\definecolor{darkgreen}{rgb}{0,.38,0}
\newcommand{\darkgreen}{\color{darkgreen}}

\newcommand{\defn}[1]{\emph{\darkblue #1}}
\newcommand{\defna}[1]{\emph{\darkred #1}}
\newcommand{\defnb}[1]{\emph{\darkblue #1}}
\newcommand{\defng}[1]{\emph{\darkgreen #1}}

\newcommand{\romb}[1]{{\darkblue #1}}
\newcommand{\romg}[1]{{\darkgreen #1}}

\newcommand{\defngb}[1]{{\textbf{\textit{\darkgreen #1}}}}

\setlist[enumerate]{
	label=\textnormal{({\roman*})},
	ref={\roman*}}

\makeatletter
\def\th@plain{%
	\thm@notefont{}
	\itshape 
}
\def\th@definition{%
	\thm@notefont{}
	\normalfont 
}
\makeatother

\makeatletter
\def\fdsy@scale{1}
\newcommand\fdsy@mweight@normal{Book}
\newcommand\fdsy@mweight@small{Book}
\newcommand\fdsy@bweight@normal{Medium}
\newcommand\fdsy@bweight@small{Medium}

\DeclareFontFamily{U}{FdSymbolB}{}
\DeclareFontShape{U}{FdSymbolB}{m}{n}{
	<-7.1> s * [\fdsy@scale] FdSymbolB-\fdsy@mweight@small
	<7.1-> s * [\fdsy@scale] FdSymbolB-\fdsy@mweight@normal
}{}
\DeclareFontShape{U}{FdSymbolB}{b}{n}{
	<-7.1> s * [\fdsy@scale] FdSymbolB-\fdsy@bweight@small
	<7.1-> s * [\fdsy@scale] FdSymbolB-\fdsy@bweight@normal
}{}

\makeatother

\DeclareSymbolFont{fdrelations}{U}{FdSymbolB}{m}{n}
\SetSymbolFont{fdrelations}{bold}{U}{FdSymbolB}{b}{n}

\DeclareMathSymbol{\lescc}{\mathrel}{fdrelations}{66}

\newtheorem{thm}{Theorem}[section]
\newtheorem{lemma}[thm]{Lemma}

\newtheorem{cor}[thm]{Corollary}
\newtheorem{prop}[thm]{Proposition}

\newtheorem{conv}[thm]{Convention}

\theoremstyle{definition}
\newtheorem{ex}[thm]{Example}

\newtheorem{rem}[thm]{Remark}
\newtheorem{definition}[thm]{Definition}

\numberwithin{figure}{section}
\numberwithin{equation}{section}

\def\wh{\widehat}

\def\emp{\nothing}
\def\sq{\square}

\def\AA{\mathbb A}
\def\nn{\mathbb N}
\def\cc{\mathbb C}

\def\rr{\mathbb R}

\def\rrs{\mathbb R_{>0}}
\def\rrp{\mathbb R_{\ge 0}}

\def\fq{{\mathbb F}_q}

\def\sm{\smallsetminus}

\def\la{\lambda}
\def\ga{\gamma}
\def\si{\sigma}
\def\de{\delta}
\def\ep{\epsilon}
\def\al{\alpha}
\def\be{\beta}
\def\om{\omega}

\def\cB{\mathcal B}

\def\cA{\mathcal A}

\def\cH{\mathcal H}
\def\cI{\mathcal I}

\def\cL{\mathcal L}

\def\cP{\mathcal P}
\def\cQ{\mathcal Q}

\def\ssu{\subset}

\def\<{\langle}
\def\>{\rangle}

\def\width{\text{{\rm width}}}

\def\SL{ {\text {\rm SL} } }
\def\GL{ {\text {\rm GL} } }

\def\rI{ {\text {\sc I} } }
\def\rJ{ {\text {\sc J} } }
\def\ri{ {\text {\rm i} } }

\def\rA{ {{A} } }
\def\rC{ {{C} } }

\def\height{ {\text {\rm height} } }

\def\rk{\textnormal{rk}}
\def\girth{\textnormal{girth}}
\def\polygirth{\textnormal{polygirth}}

\def\rT{{\text {\rm T} } }

\def\rkm{{f}}
\def\rkn{{g}}

\def\Mor{{\Phi}}

\def\0{{\mathbf 0}}

\def\tl{\triangleleft}
\def\LL{{\mathbb{L}}}

\def\nothing{\varnothing}

\def\.{\hskip.06cm}
\def\ts{\hskip.03cm}

\def\sts{\hskip.015cm}

\def\bbe{\textbf{\textit{e}}}
\def\bbc{\textbf{\textit{c}}}

\def\ba{\textbf{\textit{a}}}
\def\bb{\textbf{\textit{b}}}
\def\bd{\textbf{\textit{d}}}

\newcommand{\SYT}{\operatorname{{\rm SYT}}}

\def\nin{\noindent}

\def\SP{{\textsc{\#P}}}
\def\NP{{\textsc{NP}}}
\def\SP{{\textsc{\#{}P}}}
\def\GapP{{\textsc{GapP}}}
\def\poly{{\textsc{P}}}

\def\aD{\textrm{D}}

\def\ah{\textrm{h}}
\def\aI{{ \textrm{I} } }
\def\aIr{\textrm{\em I}}
\def\aJ{\textrm{J}}

\def\aK{\textrm{K}}
\def\aKr{\textrm{\em K}}
\def\aL{\textrm{L}}
\def\aLr{\textrm{\em L}}
\def\aM{\textrm{M}}

\def\aN{\textrm{N}}
\def\aNr{\textrm{\em N}}
\def\an{\textrm{n}}

\def\sn{{{\textsc{n}}}}

\def\aP{\textrm{p}}
\def\aPr{\textrm{\em p}}

\def\ar{d}

\def\av{\textrm{v}}

\def\bA{\textbf{\textrm{A}}\hskip-0.03cm{}}
\def\bAr{\textbf{\textrm{\em A}}\hskip-0.03cm{}}

\def\bB{\textbf{\textrm{B}}\hskip-0.03cm{}}
\def\bBr{\textbf{\textrm{\em B}}\hskip-0.03cm{}}

\def\bC{\textbf{\textrm{C}}\hskip-0.03cm{}}
\def\bD{\textbf{\textrm{D}}\hskip-0.03cm{}}

\def\bM{\textbf{\textrm{M}}\hskip-0.03cm{}}

\def\bMr{\textbf{\textrm{\em M}}\hskip-0.03cm{}}
\def\bNr{\textbf{\textrm{\em N}}\hskip-0.03cm{}}

\def\bN{\textbf{\textrm{N}}\hskip-0.03cm{}}

\def\bK{\textbf{\textrm{K}}\hskip-0.03cm{}}

\def\bT{\mathbf{T}}

\def\qnq{q}
\def\rnr{t}

\def\aAr{\mathrm{A}}
\def\aBr{\mathrm{B}}

\def\aA{\textrm A}
\def\aB{\textrm B}
\def\aC{\textrm C}

\DeclareMathOperator{\Bc}{\mathcal{B}}
\def\aa{\textrm{a}}
\def\ab{\textrm{b}}
\def\abr{\text{\em b}}

\def\ac{c}
\def\acr{c}
\DeclareMathOperator{\Act}{\textnormal{Act}}

\def\ag{\textrm{g}}
\def\af{\textrm{f}}

\def\ap{\om}
\def\apr{\om}

\def\aq{\textrm{q}}
\def\aqr{\textrm{\em q}}
\def\arr{\textrm{\em r}}
\def\as{\textrm{s}}
\def\asr{\textrm{\em s}}

\def\aw{\textrm{w}}

\def\Low{\textrm{\em Low}}

\DeclareMathOperator{\Aunt}{\textnormal{Aunt}}

\DeclareMathOperator{\Cc}{\mathcal{C}}

\DeclareMathOperator{\Cov}{\textnormal{Cov}}
\DeclareMathOperator{\Cnt}{\textnormal{Cont}} 

\DeclareMathOperator{\Des}{\textnormal{Des}}

\DeclareMathOperator{\eb}{\mathbf{e}}
\DeclareMathOperator{\Ec}{\mathcal{E}} 
\DeclareMathOperator{\Fam}{\textnormal{Fam}}
\DeclareMathOperator{\Fb}{\mathbb{F}}

\DeclareMathOperator{\Hc}{\mathcal{H}}

\DeclareMathOperator{\Ic}{\mathcal{I}}
\DeclareMathOperator{\Jc}{\mathcal{J}}
\DeclareMathOperator{\Jcm}{{\mathcal{J}}\hskip-.0572cm}

\DeclareMathOperator{\Lc}{\mathcal{L}}

\DeclareMathOperator{\nb}{\mathbf{n}}

\DeclareMathOperator{\op}{\textnormal{op}} 
\DeclareMathOperator{\Par}{\textnormal{Par}}
\DeclareMathOperator{\Pas}{\textnormal{Pas}}

\DeclareMathOperator{\Rb}{\mathbb{R}}
\newcommand{\salp}{\sim_\alpha} 

\newcommand{\spstar}{\textsf{null}} 
\newcommand{\spstarr}{{\textsf{\em null}}} 
\newcommand{\supp}{\textnormal{supp}}
\newcommand{\Stn}{\textnormal{Stn}} 

\newcommand{\us}{u} 
\newcommand{\ws}{w} 

\DeclareMathOperator{\fb}{\mathbf{f}}
\DeclareMathOperator{\gb}{\mathbf{g}}

\DeclareMathOperator{\hb}{\mathbf{h}}

\newcommand\hbw[1]{{\bT^{\<#1\>}\mathbf{{h}}}}

\DeclareMathOperator{\ub}{\mathbf{u}}
\DeclareMathOperator{\vb}{\mathbf{v}}

\newcommand\vbw[1]{\bT^{\<#1\>}\mathbf{{v}}}

\DeclareMathOperator{\wb}{\mathbf{w}}

\def\Xf{{\wh{X}}}

\DeclareMathOperator{\zb}{\mathbf{z}}
\DeclareMathOperator{\Zb}{\mathbb{Z}}

 \def\precc{\preccurlyeq} 

\DeclareMathOperator{\Ef}{\Theta} 
\DeclareMathOperator{\ef}{\text{\it e}} 
\DeclareMathOperator{\Qf}{ {\Gamma} } 
\DeclareMathOperator{\Vf}{\Omega} 
\DeclareMathOperator{\Vfm}{\Omega^0} 
\DeclareMathOperator{\Vfp}{\Omega^+} 

\DeclareMathOperator{\vf}{{\text{\it \sts v}}} 
\DeclareMathOperator{\vfs}{{\text{\it \sts v}}^\ast} 
\DeclareMathOperator{\wf}{{\text{\it \sts w}}} 

\DeclareMathOperator{\inc}{\textnormal{inc}} 

\def\cD{\Df}

\newcommand{\Df}{\mathscr{D}}

\newcommand{\Gf}{\mathscr{G}}
\newcommand{\Mf}{\mathscr{M}}
\newcommand{\Nf}{\mathscr{N}}

\newcommand{\lcov}{{\leftarrowtail}}

\newcommand{\Pf}{\cP}
\newcommand{\Af}{\mathscr{A}}

\newcommand{\low}{\textnormal{down}}
\newcommand{\up}{\textnormal{up}}
\newcommand{\zf}{z}

\newcommand{\Zf}{{Z}}
\newcommand{\Zlow}{{Z}_{\textnormal{down}}}
\newcommand{\Zup}{{Z}_{\textnormal{up}}}

\title{Log-concave poset inequalities}
\date{\today}

 \author{Swee Hong Chan}
\address[Swee Hong Chan]{Department of Mathematics, Rutgers University,  Piscatway, NJ 08854.}
\email{\texttt{sweehong.chan@rutgers.edu}}

\author[\ts Igor Pak]{Igor Pak}
\address[Igor Pak]{Department of Mathematics, UCLA,  Los Angeles, CA 90095.}
\email{\texttt{pak@math.ucla.edu}}

\begin{document}

\begin{abstract}
We study combinatorial inequalities for various classes of set systems: \ts
matroids, polymatroids, poset antimatroids, and interval greedoids.
We prove log-concave inequalities for counting certain weighted
feasible words, which generalize and extend several previous results
establishing \emph{Mason conjectures} for the numbers of independent sets
of matroids.  Notably, we prove matching equality
conditions for both earlier inequalities and our extensions.

In contrast with much of the previous work, our proofs are
combinatorial and employ nothing but linear algebra.
We use the language formulation of greedoids
which allows a linear algebraic setup, which in turn can be analyzed
recursively.  The underlying non-commutative nature of matrices
associated with greedoids allows us to proceed beyond
polymatroids and prove the equality conditions.  As further
application of our tools, we rederive both \emph{Stanley's
inequality} on the number of certain linear extensions,
and its equality conditions, which we then also
extend to the weighted case.
\end{abstract}

	\maketitle

\section{Introduction}\label{s:intro}

\subsection{Foreword}\label{ss:intro-for}
It is always remarkable and even a little suspicious, when a nontrivial
property can be proved for a large class of objects.  Indeed, this says
that the result is ``global'', i.e.\ the property is a consequence of the
underlying structure rather than individual objects.  Such results are
even more remarkable in combinatorics, where the structures are weak
and the objects are plentiful.  In fact, many reasonable conjectures
in the area fail under experiments, while some are ruled out by
theoretical considerations (cf.~$\S$\ref{ss:hist-pre} and~$\S$\ref{ss:finrem-first}).

This paper is concerned with log-concavity results for counting problems
in the general context of posets, and is motivated by a large body of
amazing recent work in area, see a survey by Huh~\cite{Huh}.  Surprisingly,
these results involve deep algebraic tools which go much beyond previous
work on the subject, see earlier surveys
\cite{Bra,Bre,Bre2,Sta2}.  This leads to several difficult questions,
such as:

\smallskip

\nin
\qquad $\circ$  \ How far do these inequalities generalize?

\smallskip

\nin
\qquad $\circ$  \ How do we extend/develop new algebraic tools to prove these generalizations?

\smallskip

\nin
We aim to answer the first question in as many cases as we can,
both generalizing the inequalities to larger classes of posets
and strengthening these inequalities to match equality conditions
which we also prove. We do this by sidestepping the second question,
or avoiding it completely.  

\medskip

There is a very long and only partially justified tradition in combinatorics
of looking for purely combinatorial proofs of combinatorial results.
Although the very idea of using advanced algebraic tools to prove combinatorial
inequalities is rather mesmerizing, one wonders if these tools are
really necessary.  Are they giving us a \defna{true insight} \ts into
the nature of these inequalities that we were missing for
so long?  Or, perhaps, the absence of purely combinatorial proofs
is a reflection of our continuing \defna{lack of understanding}\ts{}?

We posit that, in fact, all poset inequalities can be obtained
by elementary means (cf.~$\S$\ref{ss:intro-discussion}).
We show how this can be done for a several
large families of inequalities, and intend to continue this work
in the future (see~$\S$\ref{ss:finrem-complexity-GapP}).
There are certain tradeoffs, of course, as we need to introduce a
technical linear algebraic setup (see~$\S$\ref{ss:intro-proof}),
which allows us to quickly reprove both classical and recently established
poset inequalities.  The advantage of our approach is its
flexibility and noncommutative nature, making it amenable to extend
and generalize these inequalities in several directions.

Of course, none of what we did takes anything away from the algebraic
proofs of poset inequalities which remained open for decades ---
the victors keep all the spoils (see Section~\ref{s:hist}).
We do, however, hope the reader will appreciate that our combinatorial
tools are indeed \defna{more powerful} \ts than the algebraic tools,
at least in the cases we consider
(cf.~$\S\S$\ref{ss:finrem-could}--\ref{ss:finrem-AF}).

\medskip

\subsection{What to expect now}\label{ss:intro-about}
A long technical paper deserves a long technical introduction.
Similarly, a friendly and accessible paper deserves a friendly and
accessible introduction.   Naturally, we aim to achieve both
somewhat contradictory goals.

Below we present our main results and applications, all of which require
definitions which are standard in the area, but not a common
knowledge in the rest of mathematics.  We make an effort
to have the introduction thorough yet easily accessible, at the expense
of brevity.\footnote{In an effort to streamline the presentation,
some basic notation is collected in a short Section~\ref{s:def},
which we encourage the reader to consult whenever there is
an apparent misunderstanding or ambiguity.}

In addition, rather than jump to the most general and thus most
involved results, we begin slowly, and take time to introduce
the reader to the world of poset inequalities.  Essentially,
the rest of the introduction can be viewed as an extensive survey
of our own results interspersed with a few examples and some earlier
results directly related to our work.  The reader well versed in
the greedoid literature can speed read a few early subsections.

We say very little about our tools at this stage, even though we consider
them to be our main contribution (see~$\S$\ref{ss:intro-proof}
and~$\S$\ref{ss:intro-discussion}).  These are fully presented in the following
sections, which in turn are followed by the proofs of all the results.
As we mentioned above, our tools are elementary but technical, and are best
enjoyed when the reader is convinced they are worth delving into.

Similarly, in the introduction,
we say the bare minimum about the rich history of the
subject and the previous work on poset inequalities.
This is rather unfair to the many experts in the area whose names and
contributions are mentioned only at the end of the paper.  Our choice
was governed by the effort to keep the introduction from exploding
in size.  We beg forgiveness on this point, and try to mitigate it
by a lengthy historical discussion in Section~\ref{s:hist}, with
quick pointer links sprinkled throughout the introduction.

\medskip

\subsection{Matroids}\label{ss:intro-matroids-old}
A (finite) \defna{matroid}~$\Mf$ \ts is a pair \ts $(X,\Ic)$ \ts of a
\defnb{ground set} \ts $X$, \ts $|X|=n$, and a nonempty collection of \defnb{independent sets}
\ts $\Ic \subseteq 2^X$ \ts that satisfies the following:
\begin{itemize}
	\item (\defng{hereditary property}) \, $S\subset T$, \. $T\in \Ic$ \, $\Rightarrow$ \, $S\in \Ic$\ts, and
	\item (\defng{exchange property}) \,  $S, \ts T\in \Ic$, \.  $|S|<|T|$ \, $\Rightarrow$ \,
$\exists \ts x \in T \setminus S$ \. s.t.\ $S+x \in \Ic$\ts.
\end{itemize}
\defnb{Rank} of a matroid is the maximal size of the independent set: \ts $\rk(\Mf) := \max_{S\in \Ic} \ts |S|$.
A \defnb{basis} of a matroid is an independent set of size \ts $\rk(\Mf)$.  Finally,
let \ts $\cI_k \ts := \ts \bigl\{S\in \cI, \. |S|=k\bigr\}$,
and let \ts $\rI(k)=\bigl|\cI_k\bigr|$ \ts be the \defn{number of independent sets}
in \ts $\Mf$ \ts of size~$k$, \ts $0\le k \le \rk(\Mf)$.

\smallskip

\begin{thm}[{\rm \defng{Log-concavity for matroids}, \cite[Thm~9.9~(3)]{AHK},
formerly \defna{Welsh--Mason conjecture}}{}] \label{t:matroids-AHK}
For a matroid \ts $\Mf=(X,\Ic)$ \ts and integer \ts $1\le k < \rk(\Mf)$, we have:
\begin{equation}\label{eq:matroid-LC}
\rI(k)^2 \. \ge \. \rI(k-1) \ts \cdot \ts \rI(k+1)\ts.
\end{equation}
\end{thm}

\smallskip

See~$\S$\ref{ss:hist-matroid} for the historical background.
The \defnb{log-concavity} in~\eqref{eq:matroid-LC} classically
implies \defnb{unimodality} of the sequence \ts $\{\rI(k)\}$\ts:
$$
\rI(0) \. \le \. \rI(1) \. \le \. \ldots \. \le \. \rI(k) \. \ge \. \rI(k+1)  \. \ge \. \ldots \. \ge \. \rI(m)\., \quad \text{where} \ \
\, m\.=\. \rk(\Mf)\ts.
$$
It was noted in \cite[Lem.~4.2]{Lenz1} that other results in~\cite{AHK}
imply that the inequalities~\eqref{eq:matroid-LC} are \emph{always strict}
(see~$\S$\ref{ss:hist-NBC}).
Further improvements to~\eqref{eq:matroid-LC} have been long conjectured
by Mason~\cite{Mas} and were recently established in quick succession.

\smallskip

\begin{thm}[{\rm  \defng{One-sided ultra-log-concavity for matroids}, \cite[Cor.~9]{HSW}, formerly \defna{weak Mason conjecture}}{}]
\label{t:matroids-HSW}
For a matroid \ts $\Mf=(X,\Ic)$ \ts and integer \ts $1\le k < \rk(\Mf)$, we have:
\begin{equation}\label{eq:matroid-ULC-one-sided}
\rI(k)^2 \, \ge \, \left(1 \. + \, \frac{1}{k}\right)   \,  \rI(k-1)  \   \rI(k+1)\ts.
\end{equation}
\end{thm}

\smallskip

\begin{thm}[{\rm  \defng{Ultra-log-concavity for matroids}, \cite[Thm~1.2]{ALOV} and \cite[Thm~4.14]{BH}, formerly \defna{strong Mason conjecture}}{}]
\label{t:matroids-BH}
For a matroid \ts $\Mf=(X,\Ic)$,  \ts $|X|=n$,  \ts and integer \ts $1\le k < \rk(\Mf)$, we have:
\begin{equation}\label{eq:matroid-ULC}
\rI(k)^2 \, \ge \, \left(1 \. + \, \frac{1}{k}\right)  \left(1 \. + \, \frac{1}{n-k}\right)  \,\.  \rI(k-1)  \   \rI(k+1)\ts.
\end{equation}
\end{thm}

\smallskip

Equation~\eqref{eq:matroid-ULC} is a reformulation of \defnb{ultra-log-concavity} of the sequence \ts $\{\rI(k)\}$\ts:
$$
\ri(k)^2 \. \ge \. \ri(k-1) \ts \cdot \ts \ri(k+1)\,, \quad \text{where} \quad \ri(m)\.:= \. \frac{\rI(m)}{\binom{n}{m}}
$$
can be viewed as the probability that random $m$-subset of~$X$ is independent in~$\Mf$.

\medskip

\subsection{More matroids}\label{ss:intro-matroids}
For an independent set \ts $S\in \cI$ \ts of a matroid \ts $\Mf=(X,\Ic)$, denote by
\begin{equation}\label{eq:matroid-cont-def}
\Cnt(S) \, := \bigl\{x \in X \setminus S~:~S +x \in \cI\bigr\}
\end{equation}
the set of \defnb{continuations} of~$S$. For all \ts $x, y\in \Cnt(S)$, we write \ts
$x\sim_S y$ \ts when \ts $S+x+y \notin\cI$ \ts or when  \ts $x=y$\ts. Note that \ts ``$\sim_S$'' is an equivalence
relations, see Proposition~\ref{p:equiv-matroid}.
We call an equivalence class of the relation \ts $\sim_S$ \ts a \defn{parallel class} of~$S$,
and we denote by $\Par(S)$ the set of parallel classes of~$S$.

For every \ts $0\le k< \rk(\Mf)$, define the \ts \defnb{$k$-continuation number} of a matroid~$\Mf$
as the maximal  number of parallel classes of independent sets of size~$k$\ts:
\begin{equation}\label{eq:matroid-p-def}
\aP(k) \ := \ \max \ts  \bigl\{ \.  \big|\Par(S)\big|  \ : \  S \in \Ic_k \. \bigr\}.
\end{equation}
Clearly, \ts $\aP(k)\le n-k$.

\smallskip

\begin{thm} [{\rm \defng{Refined log-concavity for matroids}}{}]
\label{t:matroids-Par}
For a matroid \ts $\Mf=(X,\Ic)$ \ts and integer \ts $1\le k < \rk(\Mf)$, we have:
\begin{equation}\label{eq:matroid-Par}
\rI(k)^2 \, \ge \, \left(1\. + \. \frac{1}{k}\right)  \left(1\. + \, \frac{1}{\aPr(k-1)-1}\right)  \,\.  \rI(k-1)  \   \rI(k+1)\ts.
\end{equation}
\end{thm}

\smallskip

Clearly, Theorem~\ref{t:matroids-Par} implies Theorem~\ref{t:matroids-BH}.
This is our first result of the long series of generalizations that follow.
Before we proceed, let us illustrate the power of this refinement in a special case.

\smallskip

\begin{ex}[{\rm \defng{Graphical matroids}}{}]\label{ex:intro-graphical}
Let \ts $G=(V,E)$ \ts be a connected graph
with \ts $|V|=\sn$ \ts edges.  The corresponding \defnb{graphical matroid} \ts $\Mf_G=(E,\cI)$
\ts is defined to have independent sets to be all \emph{spanning forests} in~$G$, i.e.\
spanning subgraphs without cycles.  Then \ts $\rI(k)$ \ts is the number of spanning
forests with $k$ edges, bases are \emph{spanning trees} in~$G$, and \ts
$\rk\bigl(\Mf_G\bigr) = \sn-1$.

Let \ts $k=\sn-2$ \ts in Theorem~\ref{t:matroids-Par}.
Observe that \. $\aP(\sn-3) \leq 3$ \. since  \. $T-e-e'$ \. can have
at most three connected components, for every spanning tree $T$ in~$G$ and edges \ts $e,e'\in E$.
Then~\eqref{eq:matroid-Par} gives:
\begin{equation}\label{eq:matroids-graphical-ex-LC}
\frac{\aI(\sn-2)^2}{\aI(\sn-3) \. \cdot \. \aI(\sn-1)}  \  \geq \  \frac{3}{2} \. \bigg(1\. + \. \frac{1}{\sn-2} \bigg)
\ \to \ \frac{3}{2} \quad \text{as} \quad \sn\to\infty\ts.
\end{equation}
This is both numerically and asymptotically better than~\eqref{eq:matroid-ULC},
cf.~$\S$\ref{ss:finrem-graph-32}.
For example, when \. $|E|-\sn \to \infty$, we have:
\[
\frac{\aI(\sn-2)^2}{\aI(\sn-3) \. \cdot \. \aI(\sn-1)}  \  \geq_{\eqref{eq:matroid-ULC}} \  \bigg(1 \. + \. \frac{1}{|E|-\sn +2} \bigg) \.
\bigg(1 \. +\. \frac{1}{\sn-2} \bigg) \ \to \ 1 \quad \text{as} \quad \sn\to\infty\ts.
\]
\end{ex}

\medskip

\subsection{Weighted matroid inequalities}\label{ss:intro-matroid-weighted}
Let \ts $\Mf=(X,\Ic)$ \ts be a matroid, and let \ts $\ap: X \to \rrs$ \ts
be a positive \defnb{weight function} on the ground set~$X$.
We extend the weight function to every independent set \ts $S \in \Ic$ \ts as
follows:
\[
\ap(S)  \, := \,  \prod_{x \in S} \. \ap(x).
\]
For all \ts $1\le k < \rk(\Mf)$, define
\[  \aI_\ap(k) \ :=  \  \sum_{{S \in \Ic_k}} \ap(S)\ts.
\]

\smallskip

\begin{thm}[{\rm \defng{Refined weighted log-concavity for matroids}}{}]
\label{t:matroids-Par-weighted}
Let \ts $\Mf=(X,\Ic)$ \ts be a matroid on \ts $|X|=n$ elements, let
\ts $\apr: X \to \rrs$ \ts be a weight function, and
let \ts $1\le k < \rk(\Mf)$.  Then:
\begin{equation}\label{eq:matroid-Par-weighted}
\rI_\apr(k)^2 \, \ge \, \left(1\. + \. \frac{1}{k}\right)  \left(1\. + \, \frac{1}{\aPr(k-1)-1}\right)  \,\.  \rI_\apr(k-1)  \   \rI_\apr(k+1)\ts.
\end{equation}
\end{thm}

\smallskip

\begin{rem}\label{r:intro-weight}
In this theorem, the setup is more important than the result as it can be easily reduced
to Theorem~\ref{t:matroids-Par}.  Indeed, note that one can take multiple copies of elements
in a matroid~$\Mf$. This implies the result for integer valued~$\ap$.  The full version
follows by homogeneity and continuity.  This natural approach fails for the equality
conditions as strict inequalities are not necessarily preserved in the limit, and for
many generalizations below where we have constraints on the weight function.
See~$\S$\ref{ss:hist-weight} for some background.
\end{rem}

\medskip

\subsection{Equality conditions for matroids}\label{ss:intro-matroid-equality}
For a matroid \ts $\Mf=(X,\Ic)$ \ts on \ts $|X|=n$ \ts elements,
define \. $\girth(\Mf) \ts := \ts \min\bigl\{k\.:\. \rI(k)< \binom{n}{k}\bigr\}$.
By analogy with graph theory, \defn{girth of a matroid}
is the size of the smallest circuit in~$\Mf$.

\smallskip

\begin{thm}[{\rm\defng{Equality for matroids}, \cite[Cor.~1.2]{MNY}}{}] \label{t:matroids-equality}
Let \ts $\Mf=(X,\Ic)$ \ts be a matroid on \ts $|X|=n$ elements,
and let \ts $1\le k < \rk(\Mf)$. Then:
\begin{equation}\label{eq:matroid-ULC-equality}
\rI(k)^2 \, = \, \left(1\. + \. \frac{1}{k}\right)  \left(1\. + \, \frac{1}{n-k}\right)  \,\.  \rI(k-1)  \   \rI(k+1)
\end{equation}
	\underline{if and only if} \, $\girth(\Mf)>(k+1)$.
\end{thm}

\smallskip

See~$\S$\ref{ss:hist-improve} for some background on equality conditions.
The theorem says that in order to have equality~\eqref{eq:matroid-ULC-equality},
we must have probabilities  \ts $\ri(k-1)=\ri(k)=\ri(k+1)=1$.
Now we present a weighted version of Theorem~\ref{t:matroids-equality}.
We say that weight function \ts $\ap: X \to \rrs$ \ts is \defn{uniform} if
\ts $\ap(x)=\ap(y)$ \ts for all $x,y\in X$.

\smallskip

\begin{thm}[{\rm \defng{Weighted equality for matroids}}{}]
\label{t:matroids-equality-weighted}
Let \ts $\Mf=(X,\Ic)$ be a matroid on \ts $|X|=n$ \ts elements,
let \ts $1\le k < \rk(\Mf)$, and let
\ts $\apr: X \to \rrs$ \ts be a weight function. Then:
\begin{equation}\label{eq:matroid-weighted-equality}
\rI_\apr(k)^2 \, = \, \left(1\. + \. \frac{1}{k}\right)  \left(1\. + \, \frac{1}{n-k}\right)  \,\.  \rI_\apr(k-1)  \   \rI_\apr(k+1)
\end{equation}
\underline{if and only if} \, $\girth(\Mf)> (k+1)$, and the weight function~$\apr$ is uniform.
\end{thm}

\smallskip

The uniform condition in the theorem is quite natural for integer weight functions,
as it basically says that in order to have~\eqref{eq:matroid-weighted-equality}
all elements have to be repeated the same number of times.  In other words,
weighted inequalities do not have a substantially larger set of equality cases.

\smallskip

\begin{thm}[{\rm \defng{Refined equality for matroids}}{}]
\label{t:matroids-equality-Par}
Let \ts $\Mf=(X,\Ic)$ be a matroid, \ts $1\le k < \rk(\Mf)$, and let
\ts $\apr: X \to \rrs$ \ts be a weight function. Then:
\begin{equation}\label{eq:matroid-equality-Par}
\rI_\apr(k)^2 \, = \, \left(1\. + \. \frac{1}{k}\right)  \left(1\. + \, \frac{1}{\aPr(k-1)-1}\right)  \,\.  \rI_\apr(k-1)  \   \rI_\apr(k+1)
\end{equation}
\underline{if and only if} \,
there exists \. $\asr(k-1)>0$, such that \ts for every \ts $S\in \cI_{k-1}$ \ts we have:
	\begin{align}
		\big|\Par(S)\big| \ &= \  \aPr(k-1)\.,  \quad \ \text{and} \label{eq:ME1}
\tag{ME1}\\
		  \sum_{x \in \Cc} \.\apr(x) \ &= \ \asr(k-1)  \qquad \text{for every} \ \ \Cc \in \Par(S)\ts. \label{eq:ME2} \tag{ME2}
	\end{align}
\end{thm}

\smallskip

Condition~\eqref{eq:ME1} says that the  \ts
$(k-1)$-continuation number is achieved on all independent sets \ts
$S\in \cI_{k-1}$.  When the weight function is uniform,
condition~\eqref{eq:ME2} is saying that all parallel classes
\ts $\Cc \in \Par(S)$ \ts have the same size.

\medskip

\subsection{Examples of matroids}\label{ss:intro-matroid-examples}
First, we prove that the equality conditions are rarely satisfied
for graphical matroids, see Example~\ref{ex:intro-graphical}.  More precisely,
we prove that the refined log-concavity inequality~\eqref{eq:matroids-graphical-ex-LC}
is an equality only for cycles:

\smallskip

\begin{prop}[{\rm \defng{Equality for graphical matroids}}{}]\label{p:graphical-equality}
Let \ts $G=(V,E)$ \ts be a simple connected graph on \ts $|V|=\sn$ \ts vertices, and let
\ts $\aIr(k)$ \ts be the number of spanning forests with $k$ edges.  Then
\begin{equation}\label{eq:prop-graphical-equality}
\frac{\aIr(\sn-2)^2}{\aIr(\sn-3) \. \cdot \. \aIr(\sn-1)}  \  \ge \  \frac{3}{2} \. \bigg(1\. + \. \frac{1}{\sn-2} \bigg),
\end{equation}
and the equality holds \, \underline{if and only if} \, $G$ \ts is an \ts $\sn$-cycle.
\end{prop}

\smallskip

We now show that the equality conditions in Theorem~\ref{t:matroids-equality-Par}
have a rich family of examples (see~$\S$\ref{ss:hist-matroids-ex} for more on these examples).
The weight function is uniform in all these cases: \ts $\ap(x)=1$ \ts for every \ts $x \in X$.

\smallskip

\begin{ex}[{\rm \defng{Finite field matroids}}{}] \label{ex:intro-finite-field}
Let \ts $\Fb_{\qnq}$ \ts be a finite field with $q$ elements, let \ts $m\geq 1$, and let
\ts $X=\fq^m$.  Let \ts $\cI$ \ts be a set of subsets \ts $S\ssu \fq^m$ \ts which are
linearly independent as vectors.  Finally, let \ts $\Mf(m,q)=(X,\cI)$ \ts be a matroid
of vectors in~$\fq^m$ of rank~$m$.
	
Let \ts $1\le k<m$ \ts and let \ts $S\in \cI_{k-1}$, so we have \ts $\dim_{\fq}\langle S\rangle = k-1$.
For all parallel classes \ts $\Cc \in \Par(S)$ \ts we then have \ts $|\Cc| = q^{k-1}$.  Therefore,
\begin{equation}\label{eq:finite-field-eq}
\bigl|\Par(S)\bigr| \, = \, \frac{q^m \. - \. q^{k-1}}{q^{k-1}} \, = \, q^{m-k+1} \. -\. 1 \ts.
\end{equation}
The conditions~\eqref{eq:ME1} and~\eqref{eq:ME2} are then satisfied with
\ts $\aP(k-1)= q^{m-k+1} \ts -\ts 1$ \ts and \. $\as(k-1)=q^{k-1}$.
We conclude that~\eqref{eq:matroid-Par} is an equality for \ts $\Mf(m,q)$, for all \ts $1\le k<m$.
Curiously, the equality~\eqref{eq:finite-field-eq} is optimal for matroids over \ts $\Fb_{\qnq}$,
and we have the following result (see~$\S$\ref{ss:proof-matroid-cor-field} for the proof).
\end{ex}


\begin{cor}\label{c:intro-field}
Let \. $X \subseteq \Fb_{\qnq}^m$ \. be a set of $n$ vectors which span \ts $\Fb_{\qnq}^m$,
and let \. $\Mf=(X,\Ic)$ \.  be the corresponding matroid of rank \. $m=\rk(\Mf)$.
Then, for all \. $1 \leq k< m$,
we have:	
\begin{equation*}
\rI(k)^2 \, \ge \, \left(1\. + \. \frac{1}{k}\right)  \left(1\. + \, \frac{1}{ q^{m-k+1} \. -\. 2}\right)  \,\.  \rI(k-1)  \   \rI(k+1)\ts.
	\end{equation*}
\end{cor}

\smallskip

\begin{ex}[{\rm \defng{Steiner system matroids}}{}] \label{ex:intro-Steiner}
Fix integers \ts $t<m<n$ \ts and a ground set \ts $X$, with \ts $|X|=n$.  A \defn{Steiner system} \ts $\Stn(t,m,n)$ \ts
is a collection \ts $\cB$ \ts of $m$-subsets \ts $B\ssu X$ called \defn{blocks}, such that each \ts $t$-subset
of~$X$ is contained in exactly one block \ts $B\in \cB$.

Let \ts $\Mf(\cB)=(X,\cI)$ \ts be a matroid with \ts $\rk(\Mf)=\girth(\Mf)=(t+1)$,
where the bases are \ts $(t+1)$-subsets of $X$ that are not contained in any block
of the Steiner system.   It is easy to see that this indeed defines a matroid,
cf.~$\S$\ref{ss:hist-matroids-ex}.
Note that~\eqref{eq:matroid-Par-weighted} is trivially an equality for
all \ts $1\le k<t$.

Let \ts $S\in \cI_{t-1}$ \ts be an independent set of size \ts $(t-1)$.
The parallel classes of $S$ are given by \ts $B_1 \setminus S, \ldots, B_\ell \setminus S$,
where \ts $B_1,\ldots, B_\ell \in \cB$ \ts are blocks of the Steiner system
that contain~$S$,
and \ts $\ell=\frac{n-t+1}{m-t+1}$.  Then we have:
\[
|\Par(S)| \, = \. \ell \,, \quad \text{ and } \quad |\Cc| \. = \. m-t+1
\ \ \ \ \text{for every } \ \Cc \in \Par(S).
\]
Since the choice of $S$ is arbitrary, the conditions \eqref{eq:ME1} and \eqref{eq:ME2}
are satisfied  with \. $\aP(t-1)=\ell$ \.
	  and \. $\as(t-1)=m-t+1$.
We conclude that~\eqref{eq:matroid-Par} is also an equality for \ts $k=t$.
\end{ex}
\medskip

\subsection{Morphism of matroids}\label{ss:intro-morphism}
For a matroid \ts $\Mf = (X,\Ic)$, the \defn{rank function} \ts
\. $\rkm: 2^{X} \to \rrs$  \ts is defined by
\[
\rkm(S) \ := \ \max  \. \big\{ \ts  |A| \,: \, A \subseteq  S, \, A \in \Ic \ts \big\}.
\]
Note that \ts $\rk(\Mf) = \rkm(X)$.
There is an equivalent definition of a matroid in terms of
monotonic submodular rank functions, see e.g.~\cite{Welsh}.

\smallskip

Let \ts $\Mf = (X,\Ic)$ \ts and \ts $\Nf = (Y,\Jc)$ \ts be two
matroids with rank functions~$\rkm$ and~$\rkn$, respectively.
Let \ts $\Mor: X\to Y$ \ts be a function that satisfies
\begin{equation}\label{eq:morphisms-rk}
\rkn\big(\Mor(T)\big) \, - \,  \rkn\big(\Mor(S)\big)  \  \leq \  \rkm(T) \, - \,  \rkm(S)
\quad \text{ for every } \ S \subseteq T \subseteq X\ts.
\end{equation}
In this case we say that \ts $\Phi$ \ts is a \defna{morphism}
of matroids, write \ts $\Phi: \Mf \to \Nf$.
A subset \ts $S \in \Ic$ \ts is said to be a \defn{basis}
of \ts $\Mor$ \ts if \ts $\rkn(\Mor(S))=\rk(\Nf)$.  In other words,
$S$ is contained in a basis of \ts $\Mf$, and \ts $\Mor(S)$ \ts
contains a basis of $\Nf$. Denote by $\Bc$ the set of bases
of \ts $\Mor: \Mf\to \Nf$, and let \ts $\Bc_k:=\Bc\cap\Ic_k$.

\smallskip

Let \. $\ap: X \to \rrs$ \. be a positive weight function on the
ground set~$X$.  As before, for every \ts $0 \le k \leq \rk(\Mf)$, let
\[
\aB_\ap(k) \ :=  \  \sum_{{S \in \Bc_k}} \. \ap(S)\ts, \qquad \text{where}
\qquad \ap(S)  \ := \  \prod_{x \in S} \. \ap(x).
\]

\smallskip

\begin{thm}[{\rm \defng{Log-concavity for morphisms}, {\cite[Thm~1.3]{EH}}}{}]\label{t:morphisms-EH}
Let \ts $\Mf = (X,\Ic)$ \ts and \ts $\Nf = (Y,\Jc)$ \ts be matroids, let \ts $n:=|X|$, and let
\ts $\Mor: \Mf \to \Nf$ be a morphism of matroids.  In addition, let \ts $\apr: X \to \rrs$ \ts
be a positive weight function, and let \ $1\le k <\rk(\Mf)$. Then:
\begin{equation}\label{eq:morphism-EH}
\aBr_\apr(k)^2  \ \geq \
\left(1 \ts + \ts \frac{1}{k}\right)  \left(1 \ts + \. \frac{1}{n-k}\right) \, \aBr_\apr(k-1) \, \aBr_\apr(k+1)\ts.
\end{equation}
\end{thm}

\smallskip

Note that when \ts $Y=\{y\}$ \ts and \ts $\Nf=(Y,\emp)$ \ts is
defined by \ts $\rkn(y)= 0$, we have condition~\eqref{eq:morphisms-rk} holds trivially
and \ts $\Bc=\Ic$.  Thus, the theorem generalizes Theorem~\ref{t:matroids-BH} to the
morphism of matroids setting.  We now give the corresponding
generalization of Theorem~\ref{t:matroids-Par-weighted}.

\medskip

Recall the equivalence relation \. ``$\sim_S$'' \. on the set \.
$\Cnt(S)\subseteq X \setminus S$ \. of continuations of $S \in \Ic$,
see~\eqref{eq:matroid-cont-def}.  Similarly, recall the set \ts $\Par(S)$ \ts
of parallel classes of~$S$, see~\eqref{eq:matroid-p-def}.  For every \. $1\le k \leq \rk(\Mf)$, let
\[
\aP(k) \ := \ \max \. \bigl\{ \ts  \big|\Par(S)\big|  \. : \.  S \in \Bc_k \ts \bigr\},
\]
the maximum of the number of parallel classes of bases of morphism~$\Mor$ of size~$k$.

\smallskip

\begin{thm}[{\rm \defng{Refined log-concavity for morphisms}}{}]\label{t:morphisms-refined}
Let \ts $\Mf = (X,\Ic)$ \ts and \ts $\Nf = (Y,\Jc)$ \ts be matroids, and let
\ts $\Mor: \Mf \to \Nf$ be a morphism of matroids.  In addition, let \ts $\apr: X \to \rrs$ \ts
be a positive weight function, and let \. $1\le k <\rk(\Mf)$.  Then:
\begin{equation}\label{eq:morphism-refined}
\aBr_\apr(k)^2  \ \geq \
\left(1\ts + \ts \frac{1}{k}\right) \left(1\ts + \. \frac{1}{\aPr(k-1)\. -\. 1}\right) \,  \aBr_\apr(k-1) \, \aBr_\apr(k+1)\ts.
\end{equation}
\end{thm}

\smallskip

As before, since \. $\aP(k-1) \leq n-k+1$, the theorem is an extension of
Theorem~\ref{t:morphisms-EH}.

\smallskip

\begin{rem}\label{r:intro-morphism}
The notion of morphism of matroids generalizes many classical notions in combinatorics
such as \emph{graph coloring}, \emph{graph embeddings}, \emph{graph homomorphism},
\emph{matroid quotients}, and are a special case of the \emph{induced matroids}.
We refer to~\cite{EH} for a detailed overview and further references
(see also~$\S$\ref{ss:hist-morphism}).
\end{rem}

\medskip

\subsection{Equality conditions for morphisms of matroids}\label{ss:intro-morphism-equality}
We start with the following characterization of equality in Theorem~\ref{t:morphisms-EH},
which resolves an open problem in \cite[Question~5.7]{MNY}.

\smallskip

\begin{thm}[{\rm \defng{Equality for morphisms}}{}]
\label{t:morphisms-EH-equality}
Let \ts $\Mf = (X,\Ic)$ \ts and \ts $\Nf = (Y,\Jc)$ \ts be matroids, let \ts $n:=|X|$, and let
\ts $\Mor: \Mf \to \Nf$ be a morphism of matroids.  In addition, let \ts $\apr: X \to \rrs$ \ts
be a positive weight function, and let \ $1\le k <\rk(\Mf)$. Suppose \ts $\aBr_\apr(k)>0$.
Then:
\begin{equation}\label{eq:morphism-EH-equality}
\aBr_\apr(k)^2  \ = \
\left(1+ \frac{1}{k}\right)  \left(1+ \frac{1}{n-k}\right) \, \aBr_\apr(k-1) \, \aBr_\apr(k+1)\ts.
\end{equation}
\underline{if and only if} \, $\girth(\Mf)>k+1$, \. weight function $\apr$ is uniform,
and \. $\rkn\big(\Mor(S)\big)=\rk(\Nf)$ \ts for all \ts $S\in \Ic_{k-1}$\ts.
\end{thm}

\smallskip

Our next result  is the following characterization of equality in Theorem~\ref{t:morphisms-refined}.

\smallskip

\begin{thm}[{\rm \defng{Refined equality for morphisms}}{}]
\label{t:morphisms-Par-equality}
Let \ts $\Mf = (X,\Ic)$ \ts and \ts $\Nf = (Y,\Jc)$ \ts be matroids, and let
\ts $\Mor: \Mf \to \Nf$ be a morphism of matroids.  In addition, let \ts $\apr: X \to \rrs$ \ts
be a positive weight function, and let \ $1\le k <\rk(\Mf)$. Suppose \ts $\aBr_\apr(k)>0$.
Then:
\begin{equation}\label{eq:morphism-refined-eq}
\aBr_\apr(k)^2  \ \geq \
\left(1\ts + \ts \frac{1}{k}\right) \left(1\ts + \. \frac{1}{\aPr(k-1)\. -\. 1}\right) \,  \aBr_\apr(k-1) \, \aBr_\apr(k+1)\ts.
\end{equation}	
\underline{if and only if} \, there exists \ts $\asr(k-1)>0$, such that for every \ts $S \in \Ic_{k-1}$ \ts we have:
	\begin{align}
		\big|\Par_{S}\big| \ &= \  \aPr(k-1),   \label{eq:MorEqu1} \tag{MME1}\\
		\sum_{x \in \Cc} \. \apr(x) \ &= \ \asr(k-1)  \qquad \text{for every} \ \ \Cc \in \Par(S),
    \ \text{and}\label{eq:MorEqu2} \tag{MME2}\\
\rkn\bigl(\Mor(S)\bigr) \ts & = \ts \rk(\Nf).  \label{eq:MorEqu3} \tag{MME3}
	\end{align}
	\end{thm}

\medskip

\subsection{Discrete polymatroids}\label{ss:intro-poly}
A \defna{discrete polymatroid}\ts\footnote{\defng{Discrete polymatroids} \ts
are related but should not to be confused with \defng{polymatroids}, which is a family of
convex polytopes, see e.g.~\cite[$\S$44]{Schr} and $\S$\ref{ss:hist-polymatroid}.} \. $\cD$ \ts is a pair \ts $([n],\Jc)$ \ts of a
ground set \ts $[n]:=\{1,\ldots,n\}$ \ts and a nonempty finite collection \ts
$\Jc$  \ts of integer points \ts $\ba = (a_1,\ldots,a_n)\in \nn^n$ \ts
that satisfy the following:
\begin{itemize}
	\item (\defng{hereditary property}) \, $\ba \in \Ic$, \. $\bb \in \nn^{n}$ \. s.t.\ \. $\bb \leqslant \ba$ \ $\Rightarrow$ \  $\bb \in \Ic$\ts, and
	\item (\defng{exchange property}) \, $\ba, \bb \in \Ic$,  \. $|\ba|< |\bb|$ \ $\Rightarrow$ \
$\exists \ts i \in [n]$ \. s.t.\ \. $a_i < b_i$ \. and \. $\ba+\bbe_i \in \Jc$\ts.
\end{itemize}
Here \ts $\bb \leqslant \ba$ \ts is a componentwise inequality, \. $|\ba|:= a_1+\ldots+a_n$, and \ts
$\{\bbe_1,\ldots,\bbe_n\}$ \ts is a standard linear basis in~$\rr^n$.
When \ts $\Jc\subseteq \{0,1\}^n$, discrete polymatroid \ts $\cD$ \ts is a matroid.
One can think of a discrete polymatroid as a set system where multisets are allowed,
so we refer to~$\Jc$ as \defn{independent multisets} and to \ts $|\ba|$ \ts as \defn{size}
of the multiset~$\ba$.

The role of bases in discrete polymatroids
is played by \emph{maximal elements} with respect to the order~``$\leqslant$''; they are called
\defng{M-convex sets} in~\cite[$\S$2]{BH}.  Define \ts $\rk(\cD):= \max \{\ts |\ba|\.:\. \ba\in\Jc\.\}$.
For \ts $0\le k \le \rk(\cD)$,
denote by \ts $\Jcm_k:=\{\ba\in \Jc\.:\.|\ba|=k\}$ \ts the subcollection
of independent multisets of size~$k$, and let \ts $\rJ(k):=\bigl|\Jcm_k\bigr|$.

\smallskip

Let \ts $\ap: \ts [n] \to \rrs$ \ts be a positive \defn{weight function} on~$[n]$.  We extend weight function~$\ap$
to all \ts $\ba \in \Jc$ \ts as follows:
\[  \ap(\ba) \ :=   \. \ap(1)^{a_1} \. \cdots \. \ap(n)^{a_n}\ts.
\]
For every \ts $0\le k \le \rk(\cD)$, define
$$\aJ_\ap(k) \ := \  \sum_{\ba \in  \Jcm_k} \,   \frac{\ap(\ba)}{\ba!} \ , \quad \text{where} \quad
\ba! \, := \, a_1! \, \cdots  \, a_n!
$$

\smallskip

\begin{thm}[{\rm \defng{Log-concavity for polymatroids}, \cite[Thm~3.10~$(4)\Leftrightarrow(7)$]{BH}}{}]
\label{t:polymatroids-BH}
Let \ts $\cD=([n],\Jc)$ \ts be a discrete polymatroid, and
let  \ts $\apr: \ts [n] \to \rrs$ \ts be a positive weight function.
For every \ts $1\le k < \rk(\Mf)$, we have:
\begin{equation}\label{eq:polymatroid-BH}
\rJ_\apr(k)^2 \ \ge  \left(1\. + \. \frac{1}{k}\right)  \,\.   \rJ_\apr(k-1)  \   \rJ_\apr(k+1)\ts.
\end{equation}
\end{thm}

\smallskip

We now give a common generalization of Theorem~\ref{t:matroids-Par-weighted} and
Theorem~\ref{t:polymatroids-BH}.  Fix \ts $\rnr \in [0,1]$, and
let
\[
\pi(\ba)\, := \ \sum_{i=1}^{n} \. \binom{a_i}{2}\ts.
\]
For every \ts $0\le k \le \rk(\cD)$, define
$$\aJ_{\ap,\rnr}(k) \ := \  \sum_{\ba \in  \Jcm_k} \, \rnr^{\pi(\ba)}\,\.  \frac{\ap(\ba)}{\ba!}\,.
$$
Note that \ts $\binom{a}{2}=0$ \ts for $a\in \{0,1\}$, so \ts $\pi(\ba)=0$ \ts for all independent
sets \ts $\ba\in \cI$ \ts in a matroid.

\smallskip

For an independent multiset \ts $\ba\in \Jc$ \ts of a discrete polymatroid \ts $\cD=([n],\Jc)$,
denote by
\begin{equation}\label{eq:multi-Cont-def}
\Cnt(\ba) \, := \bigl\{i \in [n]~:~\ba +\bbe_i \in \Jc\bigr\}.
\end{equation}
the set of \defnb{continuations} of~$\ba$. For all \ts $i, j\in \Cnt(\ba)$, we write \ts
$i\sim_\ba j$ \ts when \ts $\ba+\bbe_i+\bbe_j \notin\Jc$ \ts or \ts $i=j$.  This is an
equivalence relation again, see Proposition~\ref{p:equiv-polymatroid}.
We call an equivalence class of the relation
\ts $\sim_\ba$ \ts a \defn{parallel class} of~$\ba$,
and we denote by $\Par(\ba)$ the set of parallel classes of~$\ba$.

For every \ts $0\le k< \rk(\cD)$, define the \ts \defnb{$k$-continuation number} of a
discrete polymatroid~$\cD$ as
the maximal  number of parallel classes of independent multisets of size~$k$\ts:
\begin{equation}\label{eq:multi-p-def}
\aP(k) \ := \ \max \ts  \bigl\{ \.  \big|\Par(\ba)\big|  \ : \  \ba \in \Jcm_k \. \bigr\}.
\end{equation}
For matroids, this is the same notion as defined above in~$\S$\ref{ss:intro-matroids}.

\smallskip

\begin{thm}[{\rm \defng{Refined log-concavity for polymatroids}}{}]
\label{t:polymatroids-Par}
Let \ts $\cD=([n],\Jc)$ \ts be a discrete polymatroid, and
let  \ts $\apr: \ts [n] \to \rrs$ \ts be a positive weight function.
For every \ts $\rnr \in [0,1]$ \ts and \ts $1\le k < \rk(\Mf)$, we have:
\begin{equation}\label{eq:polymatroid-Par}
\rJ_{\apr,\rnr}(k)^2 \ \ge  \left(1\. + \. \frac{1}{k}\right)  \,
\left(1\. + \, \frac{1-\rnr}{\aPr(k-1) \ts  - \ts 1 \ts +\ts \rnr}\right) \,\.   \rJ_{\apr,\rnr}(k-1)  \   \rJ_{\apr,\rnr}(k+1)\ts.
\end{equation}
\end{thm}

\smallskip

When \ts $\rnr=1$, this gives Theorem~\ref{t:polymatroids-BH}.
When \ts $\cD$ \ts is a matroid and \ts $\rnr=0$,
this gives Theorem~\ref{t:matroids-Par-weighted}.
For general discrete polymatroids~$\cD$ and \ts $0<\rnr<1$, this
is a stronger result.

\smallskip

\begin{ex}[{\rm \defng{Hypergraphical polymatroids}}{}]\label{ex:intro-hyper}
Let \ts $\cH=(V,E)$ \ts be a \defn{hypergraph} on the finite set of
vertices~$V$, with \defn{hyperedges} \. $E=\{e_1,\ldots,e_n\}$,
where \ts $e_i \subseteq V$, \ts $e_i \ne \emp$. Let
\. $W=\{w_1,\ldots,w_n\}$ \. be a collection of subsets of~$V$, such
that \ts $w_i \subseteq e_i$, $w_i \ne \emp$, and every vertex
\ts $v\in V$ \ts belongs to some~$w_i$.  A \defn{hyperpath} is an
alternating sequence \. $v \to w_i \to v' \to w_j \to v'' \to \ldots \to u$,
where \. $v, v'\in w_i$, \ts $v', v''\in w_j$, etc., and
the vertices \. $v,v',v'',\ldots,u \in V$ \. are not repeated.

A \defn{spanning hypertree} in~$\cH$ is a collection~$W$ as above, such
that every two vertices \ts $v,u\in V$ \ts are connected by
exactly one such hyperpath.  Similarly, a  \defn{spanning hyperforest}
in~$\cH$ is a collection~$W$ as above,
such that every two vertices are connected by at most one hyperpath.
In the case all \ts $|e_i|=2$, we get the usual notions of
(undirected) graphs, paths, spanning trees and spanning forests.
We say that \. $\bd = (d_1,\ldots,d_n)$, where \ts $d_i = |w_i|-1 \ge 0$,
is a \defn{degree sequence} of~$W$.  Note that in the
graphical case, we have \ts $d_i \in \{0,1\}$, so a forest is
determined by its degree sequence.  In general hypergraphs this
is no longer true.

Finally, a \defn{hypergraphical polymatroid} corresponding to~$\cH$
is a discrete polymatroid \ts $\Df_\cH = ([n],\Jc)$,
where \ts $\Jc$ \ts is a set of
degree sequences of spanning hyperforests in~$\cH$. Similarly
to graphical matroids (Example~\ref{ex:intro-tree}),
the maximal elements are degree sequences of spanning
hypertrees in~$\cH$.  Therefore, Theorems~\ref{t:polymatroids-BH}
and~\ref{t:polymatroids-Par} give log-concavity for the weighted sum
\ts  $\aJ_{\ap,\rnr}(k)$ \ts over degree sequences with total degree \ts
$d_1+\ldots + d_n=k$.
See~$\S$\ref{ss:hist-polymatroids-ex} for the background of this example.
\end{ex}

\medskip

\subsection{Equality conditions for polymatroids}\label{ss:intro-polymatroid-equality}
A discrete polymatroid \ts $\cD=([n],\Jc)$ \ts is called \defn{nondegenerate} if \.
$\eb_i \in \Jc$ \. for every \ts $i \in [n]$.
Define \.
$\polygirth(\cD) \. := \. \min \ts \bigl\{k\.:\.\rJ(k)<\binom{n+k-1}{k-1}\bigr\}$.
Observe that \ts $\ba\in \Jc$ \ts for all \ts $\ba\in \nn^k$, \. $|\ba|<\polygirth(\cD)$.
Note that the \defn{polygirth} of a discrete polymatroid does not coincide with the girth
of a matroid. In fact, $\polygirth(\cD)=2$ \ts when \ts $\cD$ \ts is a matroid with
more than one element.

\smallskip

To get the equality conditions for~\eqref{eq:polymatroid-Par}, we separate
the cases \ts $\rnr=0$, \. $0<\rnr<1$, and \ts $\rnr=1$. The case \ts $\rnr=0$ \ts coincides
with equality conditions for matroids given in Theorem~\ref{t:matroids-equality-Par}.
Examples in~$\S$\ref{ss:intro-matroid-examples} show that this is a difficult
condition with many nontrivial examples. The other two cases are in fact much
less rich.

\smallskip

\begin{thm}[{\rm \defng{Refined equality for polymatroids}, $\rnr=1$ \ts case}{}]
\label{t:polymatroids-BH-equality}
Let \ts $\cD=([n],\Jc)$ \ts be a nondegenerate discrete polymatroid, let \ts $\apr: \ts [n] \to \rrs$ \ts
be a positive weight function, and let \ts $1\le k < \rk(\Mf)$. Then:
\begin{equation}\label{eq:polymatroid-BH-equality}
\rJ_\apr(k)^2 \ =  \left(1 \. + \. \frac{1}{k}\right)  \,\.   \rJ_\apr(k-1)  \   \rJ_\apr(k+1)\ts.
\end{equation}
\underline{if and only if} \, $\polygirth(\cD) \ts > \ts (k+1)$.
\end{thm}

\smallskip

We are giving the equality condition for~\eqref{eq:polymatroid-BH} in place
of~\eqref{eq:polymatroid-Par}, since \ts $\rJ_{\apr,1}(k) = \rJ_\apr(k)$ \ts for all~$k$.

\smallskip

\begin{thm}[{\rm \defng{Refined equality for polymatroids}, $0<\rnr<1$ \ts case}{}]
\label{t:polymatroids-Par-equality}
Let \ts $\cD=([n],\Jc)$ \ts be a nondegenerate discrete polymatroid, and let \ts $\apr: \ts [n] \to \rrs$ \ts
be a positive weight function.  Fix \ts \ts $1\le k < \rk(\Mf)$ \ts and \ts $0<\rnr<1$.  Then:
\begin{equation}\label{eq:polymatroid-Par-equality}
\rJ_{\apr,\rnr}(k)^2 \ =  \left(1\. + \. \frac{1}{k}\right)  \, \left(1\. + \, \frac{1\ts -\ts \rnr}{\aPr(k-1)\ts -\ts 1 \ts +\ts \rnr}\right) \,\.   \rJ_{\apr,\rnr}(k-1)  \   \rJ_{\apr,\rnr}(k+1)\ts.
\end{equation}
\underline{if and only if} \, $k=1$, \. $\polygirth(\cD) \ts > \ts 2$, and \ts $\apr$ \ts is uniform.
\end{thm}

\smallskip

\begin{rem} \label{r:intro-degenerate} The reason the case \ts $\rnr=0$ \ts
is substantially different, is because the combined weight function
\ts $\rnr^{N(\ba)}\ts \ap(\ba)$ \ts is no longer strictly positive.
Alternatively, one can view the dearth of nontrivial examples in these
theorems as suggesting that the bound in Theorem~\ref{t:polymatroids-Par}
can be further improved for \ts $\rnr>0$.  This is based on the reasoning
that Theorem~\ref{t:matroids-Par} sharply improves over Theorem~\ref{t:matroids-BH}
because there are only trivial equality conditions for the latter
(see Theorem~\ref{t:matroids-equality}), when compared with rich equality
conditions for the former (see Theorem~\ref{t:matroids-equality-weighted}).
\end{rem}

\medskip

\subsection{Poset antimatroids}\label{ss:intro-antimatroid}
Let \ts $X$ \ts be finite set we call \defng{letters}, let \ts $n=|X|$\ts, and let \ts $X^\ast$ \ts be a set of
\defng{finite words} in the alphabet~$X$.  A \defn{language over~$X$} \ts
is a nonempty finite subset \ts $\cL \ssu X^\ast$.
A word is called \defn{simple} if it contains each letter at most once;
we consider only simple words from this point on.  We write \ts $x\in \al$ \ts
if word \ts $\al\in\cL$ \ts contains letter~$x$. Finally, let
\ts $|\al|$ \ts be the length of the word, and denote
\ts $\cL_k := \bigl\{\al\in \cL\.:\. |\al|=k\bigr\}$.

\smallskip

A pair \ts $\Af = (X,\cL)$ \ts is an \defna{antimatroid}, if the language
 \ts $\cL\ssu X^\ast$ \ts satisfies:
\begin{itemize}
\item (\defng{nondegenerate property})  \,
every  \. $x\in X$ \. is contained
in at least one \. $\al\in \cL$\ts,
\item (\defng{normal property}) \, every  \. $\al\in \cL$ \. is simple,
\item (\defng{hereditary property}) \, $\al\be\in \cL$  \ $\Rightarrow$ \ $\al\in \cL$\ts, \. and
\item (\defng{exchange property}) \,  $x\in \al$, \.  $x\notin\beta$,  \. and \.   $\ts\al, \ts \beta \in \cL$  \ $\Rightarrow$ \ $\exists\ts y \in \al$ \. s.t.\ \. $\beta y\in \cL$\ts.
\end{itemize}

\smallskip

Note that for every antimatroid \ts $\Af = (X,\cL)$, it follows from the exchange
property that
$$
\rk(\Af) \. := \. \max\{\ts |\al|\.:\.\al\in\cL\} \. = \. n\ts.
$$
Throughout the paper we use only one class of antimatroids
which we now define (cf.~$\S$\ref{ss:hist-antimatroid}).

Let \ts $\cP=(X,\prec)$ \ts be a poset on \ts $|X|=n$ \ts elements.
A simple word \ts $\al\in X^\ast$ \ts is called \defn{feasible} if \ts $\al$ \ts satisfies:
\begin{itemize}
\item (\defng{poset property}) \, if \ts $\alpha$ contains $x\in X$ and $y \prec x$,
then letter $y$ occurs before letter $x$ in $\alpha$.
\end{itemize}
A \defna{poset antimatroid} \. $\Af_\cP=(X,\cL)$ \ts is defined
by the language \ts $\cL$ \ts of all feasible words in~$X$.  The
exchange property is satisfied because one can always take $y$
to be the minimal letter (w.r.t.\ order~$\prec$) that is not in~$\be$.

\smallskip

Let \ts $\ap: X \to \rrs$ \ts be a positive weight function on~$X$.
Denote by \ts $\Cov(x) := \{y\in X\ts:\ts x\ts \lcov \ts y\}$ \ts the set of elements
which cover~$x$.  
We assume the weight function \ts $\ap$ \ts satisfies the following
\ts (\defng{cover monotonicity property})\ts:
\begin{equation}\label{eq:CM}\tag{CM}
	\ap(x)  \, \geq \,  \sum_{y \ts\in\ts \Cov(x)} \. \ap(y)\., \quad\text{for all \ $x\in X$.}
\end{equation}

\smallskip

\nin
Note that when~\eqref{eq:CM} is equality for all $x\in X$, we have:
\begin{equation}\label{eq:ap Hasse}
	 \ap(x) \, = \, \text{number of maximal chains in~$\cP$ starting at~$x$.}
\end{equation}

\medskip


For all \ts $\al \in \cL$ \ts and \ts $0\le k \le n$, let
$$\aL_\ap(k) \, := \, \sum_{\al\in \cL_k} \. \ap(\al)\,, \quad
\text{where} \quad  \ap(\alpha)  \, := \,  \prod_{x \in \alpha} \. \ap(x)\ts.
$$

\smallskip

\begin{thm}[{\rm \defng{Log-concavity for poset antimatroids}}{}]\label{t:antimatroids}
Let \ts $\cP=(X,\prec)$ \ts be a poset on \ts $|X|=n$ \ts elements,
and let \ts $\cA_\cP=(X,\cL)$ \ts be the corresponding poset antimatroid.
Let \ts $\apr: X \to \rrs$ \ts be a positive  weight function
which satisfies~\eqref{eq:CM}.  Then, for every integer \ts $1\le k < n$,
we have:
\begin{equation}\label{eq:antimatroid-LC}
\aLr_\ap(k)^2 \, \ge \, \aLr_\ap(k-1) \ts \cdot \ts \aLr_\ap(k+1)\ts.
\end{equation}
\end{thm}

\smallskip

\begin{ex}[{\rm \defng{Standard Young tableaux of skew shape}}{}]\label{ex:intro-SYT}
Let \ts $\la=(\la_1,\ldots,\la_\ell)\vdash n$, \ts be a Young diagram, and let \ts $\cP_\la=(\la,\prec)$ \ts
be a poset on squares \. $\big\{(i-1,j-1) \. : \. 1\le i \leq  \la_j, \ts 1\le j \le \ell\big\} \ssu \nn^2$, with \. $(i,j) \preccurlyeq (i',j')$ \. if \ts $i\ge i'$ \ts
and \ts $j\ge j'$.  Following~\eqref{eq:ap Hasse}, let \ts $\ap(i,j) = \binom{i+j}{i}$.
Denote \ts $a_\la(k):= \aL_\ap(k)$, \ts $0\le k \le |\la|$, and we have:
$$
a_\la(k) \, = \, \sum_{\mu\ssu\la,\. |\la/\mu|=k} \. f^{\la/\mu} \. \prod_{(i,j)\in \la/\mu} \. \binom{i+j}{i}\ts,
$$
where \. $f^{\la/\mu} = \bigl|\SYT(\la/\mu)\bigr|$  \. is the number of
\emph{standard Young tableaux} of shape~$\la/\mu$ (see~$\S$\ref{ss:hist-antimatroid-ex}).
Now Theorem~\ref{t:antimatroids}
proves that the sequence \ts $\bigl\{a_\la(k)\}$ \ts is log-concave, for every~$\la$.

This example also shows that the weight function condition~\eqref{eq:CM} is necessary.  Indeed,
let $\la$ be a \ts $m\times m$ \ts square, \ts $n=m^2$, and let \ts $\ap(i,j)=1$.
Then, for all \ts $k\le m$, we have:
$$
b(k) \,:= \, \aL_\om(k) \, = \, \bigl|\cL_k\bigr| \, = \, \sum_{\mu\vdash k} \. f^\mu\..
$$
The sequence \ts
$\{b_k\}$ \ts is the \emph{number of involutions in $S_k$},
see e.g.\ \cite[\href{http://oeis.org/A000085}{A000085}]{OEIS}, which
satisfies \. $\log b_k \, = \,  \frac12 n \log n \ts + \ts O(n)$, and
is actually \emph{log-convex}, see e.g.~\cite[$\S$4.5.2]{Mez}.
\end{ex}

\medskip

\subsection{Equality conditions for poset antimatroids}\label{ss:intro-antimatroid-equality}
Let \ts $\cP=(X,\prec)$ \ts be a poset on \ts $|X|=n$ \ts elements, and let
\ts $\Af_\cP=(X,\cL)$ \ts be the corresponding poset antimatroid.

For a word \ts $\al\in \cL$, denote by
$$\Cnt(\al) \, := \, \{\ts x\in X\.:\. \al x\in \cL\ts\}
$$
the set of \defn{continuations} of the word~$\al$.
Define an equivalence relation \ts ``$\salp$'' \ts on \ts $\Cnt(\alpha)$ \ts
by setting \ts  $x \salp y$ \ts if \ts $\alpha x y \notin \Lc$, see
Proposition~\ref{p:equiv-greedoid}.
We call the equivalence classes of \ts ``$\salp$'' \ts the
\defn{parallel classes} of~$\alpha$, and denote by \ts $\Par(\alpha)$ \ts
the set of these parallel classes.

Let \ts $\alpha \in \Lc$ \ts and \ts $x \in \Cnt(\alpha)$.
We say that $y \in X$ is a \defn{descendent} of $x$ with respect to $\alpha$
if $\alpha x y \in \Lc$ and $\alpha y \notin \Lc$.
Denote by \ts $\Des_\alpha(x)$ \ts the set of descendants of~$x$
with respect to $\alpha$.  We omit~$\al$ when the word is clear
from the context.

\smallskip

\begin{thm}[{\rm \defng{Equality for poset antimatroids}}{}]\label{t:antimatroids-equality}
Let \ts $\cP=(X,\prec)$ \ts be a poset on \ts $|X|=n$ \ts elements,
and let \ts $\cA_\cP=(X,\cL)$ \ts be the corresponding poset antimatroid.
Let \ts $\apr: X \to \rrs$ \ts be a positive  weight function
which satisfies~\eqref{eq:CM}, and fix an integer \ts $1\le k < n$.  Then:
\begin{equation}\label{eq:antimatroid-equality}
\aLr_\ap(k)^2 \ = \ \aLr_\ap(k-1) \. \cdot \. \aLr_\ap(k+1)
\end{equation}
\underline{if and only if} \,
there exists \. $\asr(k-1)>0$, such that \ts for every \ts $\al\in \cL_{k-1}$ \ts and \ts $x \in \Cnt(\alpha)$, we have:
\begin{align}
\sum_{x \ts\in\ts \Cnt(\alpha)} \. \apr(x) \  & = \  \asr(k-1)\ts,   \label{eq:AE1} \tag{AE1}\\
\Des_{\alpha}(x) \ & = \ \Cov(x)\ts,  \quad \text{and} \label{eq:AE2} \tag{AE2} \\
\sum_{y\ts\in \ts \Cov(x)} \. \apr(y)  \ & = \  \apr(x)\ts. \label{eq:AE3} \tag{AE3}
\end{align}
\end{thm}

\smallskip

The following is an example of a poset that satisfies conditions of Theorem~\ref{t:antimatroids}.

\smallskip

\begin{ex}[{\rm \defng{Tree posets}}{}]\label{ex:intro-tree}
Let \ts $\rT=(V,E)$ \ts be a finite \emph{rooted tree} with root at $R\in V$,
and the set of leaves \ts $S\ssu V$.  Suppose further, that all leaves \ts $v\in S$ \ts
are at distance \ts $h$ \ts from~$R$.
Consider a poset \ts $\cP_\rT=(V,\prec)$ \ts with \ts $v\prec v'$ \ts if
the shortest path \ts $v'\to R$ \ts goes through~$v$, for all \ts $v,v'\in V$.
We call \ts $\cP_\rT$ \ts the \defn{tree poset} corresponding to~$T$.
Denote by \ts $S(v) := S \cap \{v'\in V\ts:\ts v'\succcurlyeq v\}$ \ts the
subset of leaves in the order ideal of~$v$.

Let \ts $\ap:X\to \rrs$ \ts be defined by~\eqref{eq:ap Hasse}.  Observe
that \ts $\ap(v)=\bigl|S(v)\bigr|$, since maximal chains in~$\cP_\rT$ are
exactly the shortest paths in~$\rT$ towards one of the leaves, i.e.\
of the form \ts $v\to w$ \ts for some \ts $w\in S$.  Note that
\ts $S(v)\supseteq S(v')$ \ts for all \ts $v\prec v'$,
\ts $S(v)\cap S(v')=\emp$ \ts  for all \ts $v$ and $v'$ \ts
that are incomparable, and \ts $\sum_{x\in \Cov(v)} |S(x)| = |S(v)|$ \ts
for all \ts $v\notin S$.  These imply \eqref{eq:AE1}--\eqref{eq:AE3}
for all \ts $k\le h$, with \ts $\as(k-1)=|S|$.
By Theorem~\ref{t:antimatroids}, we get an equality
\eqref{eq:antimatroid-equality} in this case.
\end{ex}

\smallskip

The following result shows the importance of tree posets
for the equality conditions.

\smallskip

\begin{thm}[{\rm \defng{Total equality for poset antimatroids}}{}]\label{t:antimatroids-equality-total}
Let \ts $\cP=(X,\prec)$ \ts be a poset on \ts $|X|=n$ \ts elements,
and let \ts $\cA_\cP=(X,\cL)$ \ts be the corresponding poset antimatroid.
Let \ts $\apr: X \to \rrs$ \ts be a positive  weight function
which satisfies~\eqref{eq:CM}.  Then:
\begin{equation}\label{eq:antimatroid-equality-total}
\aLr_\ap(k)^2 \ = \ \aLr_\ap(k-1) \. \cdot \. \aLr_\ap(k+1) \quad \text{for all \ $1\le k < \height(\cP)$}
\end{equation}
\underline{if and only if} \, $\cP\cup \wh 0$ \ts is a  tree poset \ts $\cP_\rT$ \ts with a
root at~$\ts\wh 0$, with all leaves at the same distance to the root,
and such that \ts $c\.\apr$ \ts is defined by~\eqref{eq:ap Hasse},
for some constant multiple \ts $c>0$.
\end{thm}

\medskip

\subsection{Interval greedoids}\label{ss:intro-greedoids}
Let \ts $X$ \ts be finite set of letters, and let \ts $\cL\ssu X^\ast$ \ts be a
language over~$X$.  A pair \ts $\Gf=(X,\cL)$ \ts is a \defna{greedoid}, if the
language~$\cL$ satisfies:

\begin{itemize}
\item  (\defng{nondegenerate property})  \, empty word \ts $\varnothing$ \ts is in $\Lc$,
\item (\defng{normal property}) \, every  \ts $\al\in \cL$ \ts is simple,	
\item (\defng{hereditary property}) \, $\al\be\in \cL$  \ $\Rightarrow$ \ $\al\in \cL$, and
\item (\defng{exchange property}) \,  $\al, \be\in \cL$ \. s.t.\  \. $|\alpha| > |\beta|$
 \ $\Rightarrow$ \ $\exists\ts x\in \al$ \, s.t.\  \. $\beta x \in \Lc$\ts.
\end{itemize}

\smallskip

\nin
Let \. $\rk(\Gf) := \max\{\ts |\al|\ts:\ts\al\in \cL\ts\}$ \. be the \defn{rank}
of greedoid~$\Gf$.
Note that every maximal word in~$\cL$ has the same length by the exchange property.
In the literature, greedoids are also defined via \defng{feasible sets} of
letters in $\al\in \cL$, but we restrict ourselves to the language notation.
We use \cite[$\S$8.2.B]{BZ} and  \cite[$\S$V.5]{KLS} as our main references on
interval greedoids; see also~$\S$\ref{ss:hist-greedoid} for some background.

\smallskip

Greedoid \ts $\Gf=(X,\cL)$ \ts is called \defna{interval} if the
language~$\cL$ also satisfies:
\begin{itemize}
\item (\defng{interval property}) \,  $\alpha,\beta,\gamma \in X^*$, \. $x \in X$ \ts \.
s.t.\  \. $\alpha x, \ts \alpha \beta \gamma x \in \Lc$
 \ \, $\Rightarrow$ \ \,$\alpha \beta x \in \Lc$\ts.
\end{itemize}
It is well known and easy to see that antimatroids are interval greedoids.

\smallskip

Let \. $\aq: \cL \to \rrs$ \. be a positive weight function.  Let
$$
\aL_\aq(k) \, = \, \sum_{\al\in \cL_k} \. \aq(\al)\ts.
$$
In the next section, we define the notion of \ts \defna{$k$-admissible} \ts
weight function~$\aq$, see Definition~\ref{d:admissible}.  This notion is
much too technical to state here.  We use it to formulate our first main result:

\smallskip

\begin{thm}[{\rm \defng{Log-concavity for interval greedoids}, {\em \defna{first main theorem}}}{}]
\label{t:greedoid}
Let \ts $\Gf=(X,\cL)$ \ts be an interval greedoid, let \ts $1\le k <\rk(\Gf)$,
and let \. $\aqr:\cL \to \rrs$ \. be a \ts {$k$-admissible} \ts weight function.
Then:
\begin{equation}\label{eq:greedoid-LC}
    \aLr_\aqr(k)^2 \ \geq \   \aLr_\aqr(k-1) \. \cdot \. \aLr_\aqr(k+1)\ts.
\end{equation}
\end{thm}

\smallskip

This is the first main result of the paper, as it implies all previous inequalities for
matroids, polymatroids and antimatroids.

\smallskip

\begin{ex}[{\rm \defng{Directed branching greedoids}}{}]\label{ex:intro-branching}
{\rm Let \ts $G=(V,E)$ \ts be a directed graph
on \ts $|V|=n$ \ts vertices strongly connected towards the root \ts $R\in V$.
An \defn{arborescence} is a tree in~$G$ strongly connected towards the root~$R$.
A word \. $\al=e_1\cdots e_\ell\in E^\ast$ \. is called \defn{pointed} if
every prefix of~$\al$ consists of edges which form an arborescence.
One can think of pointed words as \emph{increasing arborescences} in~$G$
(cf.~$\S$\ref{ss:hist-antimatroid-branching}).

The \defnb{directed branching greedoid} \ts $\Gf_G=(E,\cL)$ \ts is defined
on the ground sets $E$ by the language \ts $\cL\ssu E^\ast$ \ts of pointed
words.  It is well known and easy to see that \ts $\Gf_G$ \ts is an interval
greedoid.  When \ts $G=T$ \ts is a rooted tree, greedoid \ts $\Gf_T$ \ts
is the poset antimatroid corresponding to the tree poset \ts $\cP_P$
(see Example~\ref{ex:intro-tree}).   For general graphs, greedoid \ts $\Gf_G$ \ts
is not necessarily a poset antimatroid.  Theorem~\ref{t:greedoid}
in this case proves log-concavity for the numbers \ts $\aL_\aq(k)$ \ts of
 weighted increasing arborescences, cf.~$\S$\ref{ss:hist-antimatroid-branching}.
}\end{ex}

\medskip

\subsection{Equality conditions for interval greedoids}\label{ss:intro-greedoids-equality}
A word \ts $\be\in X^\ast$ is called a \defn{continuation}
of the word \ts $\al\in \cL$, if \ts $\al\be\in \cL$.
Denote by \ts $\Cnt_k(\al)\ssu X^\ast$ \ts the set of
continuations of the word~$\al$ with \ts $\be\in X^\ast$
\ts of length \ts $|\be|=k$. Note that \ts
$\Cnt(\al)=\Cnt_1(\al)$.
For notational convenience, we define \. $\Cnt(\al)=\varnothing$ \. if $\alpha \notin \cL$.

\smallskip

For every \ts $\alpha \in \Lc$, let
\[
\aL_{\aq,\alpha}(k) \ := \  \sum_{\be\in \Cnt_k(\al)} \, \aq(\alpha\beta)\ts.
\]
Note that \. $\aL_\aq(k)= \aL_{\aq,\varnothing}(k)$ \. and \. $\aL_{\aq,\alpha}(0) = \aq(\al)$.

\smallskip

\begin{thm}[{\rm \defng{Equality for interval greedoids}, \ts cf.~Theorem~\ref{t:greedoid-equality-full}}{}]
\label{t:greedoid-equality}
Let \ts $\Gf=(X,\cL)$ \ts be an interval greedoid, let \ts $1\le k <\rk(\Gf)$,
and let \. $\aqr:\cL \to \rrs$ \. be a \ts {$k$-admissible} \ts weight function. Then:
\[
    \aLr_\aqr(k)^2 \ = \   \aLr_\aqr(k-1) \. \cdot \. \aLr_\aqr(k+1)
\]
\underline{if and only if} \, there is \. $\asr(k-1)>0$,  such that for every \. $\al\in \cL_{k-1}$ \. we have:
$$
\aLr_{\aqr,\alpha}(2)  \ = \  \asr(k-1) \, \aLr_{\aqr,\alpha}(1) \ = \  \asr(k-1)^2 \, \aLr_{\aqr,\alpha}(0)\ts.
$$
\end{thm}

\medskip

This is the second main result of the paper, giving an easy way to check the
equality conditions.  A more detailed and technical condition is given in
Theorem~\ref{t:greedoid-equality-full}, which we use to obtain the
equality conditions for matroids, polymatroids and antimatroids.

\medskip

\subsection{Linear extensions}\label{ss:intro-LE}
Let \ts $\Pf:=(X,\prec)$ \ts be a poset on \ts $n:=|X|$ \ts elements.
A \defna{linear extension} of $\Pf$ is a bijection \. $L: X \to \{1,\ldots, n\}$, such that
\. $L(x) < L(y)$ \. for all \. $x \prec y$.  Fix an element \ts $\zf\in X$.
Denote by \ts $\Ec := \Ec(P)$ \ts the set of linear extensions
of~$\Pf$, let \. $\Ec_k\ts := \ts \{L \in \Ec\.:\.L(\zf)=k\}$, and let \. $e(\cP):=|\Ec|$.
See~$\S$\ref{ss:hist-LE} and~$\S$\ref{ss:hist-LE-AF} for some background.

\smallskip

\begin{thm}[{\rm \defng{Stanley inequality}~\cite[Thm~3.1]{Sta}}{}]\label{t:Sta}
Let \ts  $\cP=(X,\prec)$ \ts be a poset with \ts $|X|=n$ \ts elements, and let \ts $\zf\in X$.
Denote by \ts $\aNr(k):=|\Ec_k|$ \ts the number of linear extensions \ts $L\in \Ec(P)$,
such that \ts $L(\zf)=k$.  Then, for every \. $1< k < n$, we have:
\begin{equation}\label{eq:Sta}
\aNr(k)^2 \,\. \ge \,\. \aNr(k-1) \.\cdot \.  \aNr(k+1)\ts.
\end{equation}
\end{thm}

\smallskip

We now give a weighted generalization of this result.  Let \ts $\ap: X \to \rrs$ \ts
be a positive weight function on~$X$.  We say that \ts $\ap$ \ts is \defn{order-reversing}
if it satisfies
\begin{equation}\label{eq:Rev}\tag{Rev}
x \. \preccurlyeq\. y \quad \Rightarrow 	\quad  \ap(x) \,  \geq \, \ap(y)\ts.
\end{equation}
Fix \ts $\zf\in X$, as above.   Define \. $\ap: \ts\Ec\to \rrs$ \. by
\begin{equation}\label{eq:sta-def-omega-weighted}
\ap(L)  \ := \ \prod_{x \, : \, L(x) < L(\zf)} \. \ap(x)\ts,
\end{equation}
and let
\begin{equation}\label{eq:sta-def-N-weighted}
\aN_\ap(k) \, := \, \sum_{L \in \Ec_k} \. \ap(L)\., \quad \text{for all} \ \ 1\le k \le n\ts.
\end{equation}

\smallskip

\begin{thm}[{\rm \defng{Weighted Stanley inequality}}{}]\label{t:Sta-weighted}
Let \ts  $\cP=(X,\prec)$ \ts be a poset with \ts $|X|=n$ \ts elements, and
let \ts $\ap: X \to \rrs$ \ts be a positive order-reversing weight function.
Fix an element \ts $\zf\in X$.
Then, for every \. $1< k < n$, we have:
\begin{equation}\label{eq:Sta-weighted}
\aNr_\apr(k)^2 \,\. \ge \,\. \aNr_\apr(k-1) \.\cdot \. \aNr_\apr(k+1)\ts,
\end{equation}
where \ts $\aNr_\apr(k)$ \ts is defined by~\eqref{eq:sta-def-N-weighted}.
\end{thm}

\smallskip

\begin{rem}  In~$\S$\ref{ss:LE-belts}, we give further applications
of our approach by extending the set of possible weights in
Theorem~\ref{t:Sta-weighted} to a smaller class of
\emph{posets with belts}.  We postpone this discussion to avoid
cluttering, but the interested reader is encouraged to skip to
that subsection which can be read separately from the rest
of the paper.\footnote{In a followup investigation, we use the
combinatorial atlas technology in \cite{CP1} to prove correlation
inequalities for the numbers of linear extensions of posets.}
\end{rem}

\medskip

\subsection{Two permutation posets examples}\label{ss:intro-LE-examples}
It is not immediately apparent that the numbers of linear extensions
appear widely across mathematics.  Below we present two notable examples
from algebraic and enumerative combinatorics, see~$\S$\ref{ss:hist-LE-ex}
for some background.

\smallskip

\begin{ex}[{\rm \defng{Bruhat orders}}{}]\label{ex:sta-bruhat}
Let \ts $\si\in S_n$ \ts and define the \defn{permutation poset} \ts $\cP_\si = ([n],\prec)$
\ts by letting
$$i \precc j \quad \Leftrightarrow \quad i\le j \ \ \text{and} \ \ \si(i)\le \si(j) \ts.
$$
Fix \ts $\zf\in [n]$.
Viewing \ts $\Ec=\Ec(\cP_\si)$ \ts as a subset of \ts $S_n$,
it is easy to see that \ts $\Ec$ \ts is the lower ideal of  \ts $\si$ \ts
in the (weak) \defn{Bruhat order} \ts $\cB_n=(S_n,\lhd)$.
Thus, \. $\Ec_k\.=\.\big\{\nu\in S_n\.:\. \nu(\zf)=k, \. \nu \unlhd\si\big\}$.

Let \ts $\ap(i)=q^i$, where \. $0< q<1$.  Then \ts $\ap$ \ts is order-reversing.
Now~\eqref{eq:sta-def-omega-weighted} gives \ts
$\ap(\nu) = q^{\be(\nu)}$, where
$$\be(\nu)\, := \, \sum_{i=1}^{\zf-1} \ i \ts\cdot\ts \chi\bigl(k-\nu(i)\bigr) \. \quad \text{and}
\quad \chi(t) \. := \. \left\{\aligned \ 1 & \ \ \text{if} \ \ t>0 \\ \ 0 & \ \ \text{if} \ \ t\le 0\endaligned\right.
$$
Now Theorem~\ref{t:Sta-weighted} gives log-concavity \. $a_q(k)^2 \ge a_q(k-1) \cdot a_q(k+1)$, where \.
$a_q(k):=\aN_\ap(k)\ge 0$ \. is given by
$$
a_q(k) \, = \, \sum_{\nu\ts\in \ts S_n\.: \. \nu \ts\unlhd\ts \si, \. \nu(z)=k} \.  q^{\be(\nu)}.
$$
\end{ex}

\smallskip

\begin{ex}[{\rm \defng{Euler--Bernoulli} and \defng{Entringer numbers}}{}]
\label{ex:intro-EB}
Let \ts $\cQ_m=([2m-1],\prec)$ \ts be a height two poset corresponding to the skew
Young diagram \ts $\de_m/\de_{m-2}$, where \ts $\de_m:=(m,\ldots,2,1)$.
The linear extensions of \ts $\cQ_m$ \ts are in natural bijection with
\defn{alternating permutations} \ts $\si\in S_{2m-1}$ \ts s.t.\ \ts
$\si(1)>\si(2)<\si(3)>\si(4)<\ldots$ \ts  Then the numbers \ts $e(\cQ_m)$
\ts are the \defn{Euler numbers}, which are closely related to the
\defn{Bernoulli numbers}, and have EGF
$$
\sum_{m=1}^{\infty} \. (-1)^{m-1} \. e(\cQ_m) \. \frac{t^{2m-1}}{(2m-1)!} \, = \, \tan(t)\.,
$$
see e.g.\ \cite[\href{http://oeis.org/A000111}{A000111}]{OEIS}.
Fix \ts $\zf=1$.  It is easy to see that triangle of numbers \ts
$a(m,k)=\bigl|\Ec_k(\cQ_m)\bigr|$ \ts are \defn{Entringer numbers}
\cite[\href{http://oeis.org/A008282}{A008282}]{OEIS},
and Stanley's Theorem~\ref{t:Sta} proves their log-concavity:
$$a(m,k)^2 \. \ge \. a(m,k-1) \. a(m,k+1) \quad \text{for} \ \ 1\le k \le 2m-2\ts.
$$
Now, let \ts $\ap(2)=\ap(4)=\ldots = 1$, \ts $\ap(1)=\ap(3)=\ldots = q$, where \ts $0<q<1$.
Similarly to the previous example, we have \. $\ap(\si) \ts = \ts q^{\ga(\si)}$, where
\ts $\ga(\si)$ \ts is the number of permutation entries in the odd positions which are \ts $<k$.
Theorem~\ref{t:Sta-weighted} then proves log-concavity for the corresponding \ts
$q$-deformation of the Entringer numbers.
\end{ex}

\medskip

\subsection{Equality conditions for linear extensions}\label{ss:intro-LE-equality}
Let \ts $\Pf:=(X,\prec)$ \ts be a poset on \ts $|X|=n$ \ts elements.
Denote by \ts $f(x) := \bigl|\{y \in X\,{}:\,{}y\prec x\}\bigr|$ \ts and \ts
$g(x) := \bigl|\{y \in X\,{}:\,{}y\succ x\}\bigr|$ \ts the sizes of lower
and upper ideals of \ts $x\in X$, respectively, excluding the element~$x$.

\smallskip

\begin{thm}[{\rm \defng{Equality condition for Stanley inequality}~\cite[Thm~15.3]{SvH}}{}]\label{t:Sta-equality}
Let \ts  $\Pf=(X,\prec)$ \ts be a poset with \ts $|X|=n$ \ts elements.  Let \ts $\zf\in X$ \ts and
let \ts $\aNr(k)$ \ts be the number of linear extensions \ts $L\in \Ec(P)$,
such that \ts $L(\zf)=k$.  Suppose that \ts $\aNr(k)>0$.
Then \. \underline{the following are equivalent}\ts:
\begin{enumerate}
			[{label=\textnormal{({\alph*})},
		ref=\textnormal{\alph*}}]
\item \ $\aNr(k)^2 \. = \. \aNr(k-1) \. \cdot \. \aNr(k+1)$,
\item \ $\aNr(k+1) \. = \. \aNr(k) \. = \. \aNr(k-1)$,
\item \ we have \. $f(x)>k$ \. for all \. $x\succ \zf$, \. and \. $g(x)>n-k+1$ \. for all \. $x\prec \zf$.
\end{enumerate}
\end{thm}

\smallskip

See~$\S$\ref{ss:hist-LE-equality} for some background.
The weighted version of this theorem is a little more subtle and needs
the following \ts \defng{$(\as,k)$-cohesiveness property}\ts:
\begin{equation}\label{eq:cohesive}\tag{Coh}
\ap\bigl(L^{-1}(k-1)\bigr) \. = \. \ap\bigl(L^{-1}(k+1)\bigr) \. = \. \as\ts, \quad
\text{for all} \ \ L \in \Ec_k\,.
\end{equation}

\smallskip

Note that~\eqref{eq:cohesive} can hold for non-uniform weight functions~$\ap$,
for example for \ts $\Pf = \rA_{k+1} \oplus \rC_{n-k-1}$, i.e.\ the \emph{linear sum} of an
 antichain on which $\ap$ is uniform and a chain on which $\ap$ can be non-uniform.
In fact, if $\zf$ is an element in~$\rA_{k+1}$, we can have \ts $\ap(\zf)$ \ts
different from the rest of the antichain.

\smallskip

\begin{thm}[{\rm \defng{Equality condition for weighted Stanley inequality}}{}]\label{t:Sta-equality-weighted}
Let \ts  $\Pf=(X,\prec)$ \ts be a poset on \ts $|X|=n$ \ts elements, and
let \ts $\ap: X \to \rrs$ \ts be a positive order-reversing weight function.
Fix element \ts $\zf \in X$ \ts and let \ts $\aNr_\om(k)$ \ts be defined as
in~\eqref{eq:sta-def-N-weighted}.  Suppose that \ts $\aNr_\apr(k)>0$.
Then \. \underline{the following are equivalent}\ts:
\begin{enumerate}
			[{label=\textnormal{({\alph*})},
		ref=\textnormal{\alph*}}]
\item \ $\aNr_\apr(k)^2 \. = \. \aNr_\apr(k-1) \. \cdot \.\aNr_\apr(k+1)$,
\item \. there exists \. $\asr=\asr(k,\zf) >0$, s.t. 
$$\aNr_\apr(k+1) \. = \. \asr \. \aNr_\apr(k) \. = \. \asr^2 \. \aNr_\apr(k-1),
$$
\item \. there exists \. $\asr=\asr(k,\zf) >0$, s.t.\
 \. $f(x)>k$ \. for all \. $x\succ \zf$, \. $g(x)>n-k+1$ \. for all \. $x\prec \zf$, and~\eqref{eq:cohesive}.
\end{enumerate}
\end{thm}

\medskip

\subsection{Summary of results and implications} \label{ss:intro-summary}
Here is a chain of matroid results from new to known:
$$
\text{\romb{Thm~\ref{t:matroids-Par-weighted}} \ $\Rightarrow$ \ \romb{Thm~\ref{t:matroids-Par}} \ $\Rightarrow$ \ \romg{Thm~\ref{t:matroids-BH}}
\ $\Rightarrow$ \ \romg{Thm~\ref{t:matroids-HSW}} \ $\Rightarrow$ \ \romg{Thm~\ref{t:matroids-AHK}}}.
$$
The first two of these introduce the \emph{refined log-concave inequalities}, both weighted and unweighted,
and they imply the last three known theorems.  For morphisms of matroids and for polymatroids,
we have two new results which extend two earlier results:
$$
\text{\romb{Thm~\ref{t:morphisms-refined}} \ \ $\Rightarrow$ \ \ \romg{Thm~\ref{t:morphisms-EH}} \qquad
and \qquad \romb{Thm~\ref{t:polymatroids-Par}} \ \ $\Rightarrow$ \ \ \romg{Thm~\ref{t:polymatroids-BH}}.
}
$$

Here is a family of implications of log-concave inequalities across matroid generalizations,
from interval greedoids to polymatroids to matroids, and from interval greedoids to poset antimatroids:
$$
\text{
\romb{Thm~\ref{t:greedoid}}  \ $\Rightarrow_{\text{\S\ref{ss:reduction-polytroid}}}$ \ \romb{Thm~\ref{t:polymatroids-Par}}
\ $\Rightarrow_{\text{\S\ref{ss:intro-poly}}}$ \ \romb{Thm~\ref{t:matroids-Par-weighted}} \qquad
and \qquad \romb{Thm~\ref{t:greedoid}}  \ $\Rightarrow_{\text{\S\ref{ss:reduction-antimatroid}}}$ \ \romb{Thm~\ref{t:antimatroids}}.
}
$$
All these results are new.  Note that both polymatroids and poset antimatroids are different
special cases of interval greedoids, while our results on morphisms of matroids are separate and do not
generalize.

\smallskip

For the equality conditions, we have a similar chain of implications
across matroid generalizations:
$$\aligned &
\text{\romb{Thm~\ref{t:greedoid-equality-full}} \ $\Rightarrow$ \ \romb{Thm~\ref{t:greedoid-equality} }
\ $\Rightarrow$ \ \romb{Thm~\ref{t:polymatroids-Par-equality}} \ts $\cup$ \ts\romb{Thm~\ref{t:polymatroids-BH-equality}} \ $\Rightarrow$ \ \romb{Thm~\ref{t:matroids-equality-Par}}
\ $\Rightarrow$ \ \romb{Thm~\ref{t:matroids-equality-weighted}} \ $\Rightarrow$ \ \romg{Thm~\ref{t:matroids-equality}},
}\\
&
\text{\romb{Thm~\ref{t:morphisms-Par-equality}} \ \ $\Rightarrow$ \ \ \romb{Thm~\ref{t:morphisms-EH-equality}}
\qquad and \qquad
\romb{Thm~\ref{t:greedoid-equality-full}} \ \ $\Rightarrow$ \ \
\romb{Thm~\ref{t:antimatroids-equality}} \ \ $\Rightarrow$ \ \ \romb{Thm~\ref{t:antimatroids-equality-total}}.
}
\endaligned
$$
Of these, only Theorem~\ref{t:matroids-equality} was previously known.  The most general
of these, Theorem~\ref{t:greedoid-equality-full}, is too technical to be stated
in the introduction.  The same holds for Definition~\ref{d:admissible} needed in Theorem~\ref{t:greedoid}.
We postpone both the definition and the general theorem until Section~\ref{s:main}.

\smallskip

Finally, for the Stanley inequality and its equality conditions, we have:
$$
\text{
\romb{Thm~\ref{t:Sta-weighted}} \ \ $\Rightarrow$ \ \ \romg{Thm~\ref{t:Sta}} \qquad
\text{and} \qquad
\romb{Thm~\ref{t:Sta-equality-weighted-more}}
\ \ $\Rightarrow$ \ \ \romb{Thm~\ref{t:Sta-equality-weighted}}
\ \ $\Rightarrow$ \ \ \romg{Thm~\ref{t:Sta-equality}}.
}
$$
In both cases, more general results are new and correspond to the case of weighted linear extensions.

\smallskip

Let us emphasize that while some of these implications are trivial or follow immediately
from definitions, others are more involved and require a critical change of notation
and some effort to verify certain poset and weight function properties.
These implications are discussed in Section~\ref{s:reduction}.

\medskip

\subsection{Proof ideas} \label{ss:intro-proof}
Although we prove multiple results, the proof of each log-concavity
inequality uses the same approach and technology, so we refer to
it as ``the proof''.

At the first level, the proof is an inductive argument proving a stronger
claim about eigenvalues of certain matrices associated with the posets.
The induction is not over posets of smaller size, but over other matrices
which can in fact be larger, but correspond to certain parameters decreasing
as we go along.  The claim then reduces to the base of induction,
which is the only part of the proof requiring a computation.  The latter
involves checking eigenvalues of explicitly written small matrices, making
the proof fully elementary.

Delving a little deeper, we set up a new type of structure which we call a \ts
\defna{combinatorial atlas}.  In the special case of greedoids,
a combinatorial atlas \ts  $\AA$ \ts associated with a
greedoid \ts $\Gf = (X,\cL)$, \ts $|X|=n$, \ts is comprised of:

\smallskip

\nin \hskip1.cm
$\circ$ \ acyclic digraph \. $\Qf_\Gf=(\cL,\Ef)$, \ts with the unique source at the
empty word \ts $\emp\in \cL$, and edges

\nin \hskip1.42cm
corresponding to
multiplications by a letter: \ $
\Ef \ts = \ts \bigl\{(\al,\al x)\.:\. \al,\al x\in \cL, \. x\in X\bigr\}$,

\smallskip

\nin \hskip1.cm
$\circ$ \ each vertex \ts $\al \in \cL$ \ts is associated with a pair \ts
$(\bM_\al,\hb_\al)$, \ts where \ts $\bM_\al=\big(\aM_{i j}\big)$ \ts is a nonnegative

\nin \hskip1.42cm
symmetric \ts $\ar\times \ar$ \ts matrix, \ts $\hb_\al = (\ah_1,\ldots,\ah_\ar)$ \ts is a nonnegative vector,
and \ts $\ar=n+1$,

\smallskip

\nin \hskip1.cm
$\circ$ \ each edge \ts $(\al,\al x)\in \Ef$ \ts is associated with a linear \ts transformation \ts
$\bT_\al^{\langle x\rangle}: \rr^\ar\to \rr^\ar$.

\medskip

\nin
The key technical observation is that under certain conditions on the atlas,
we have every matrix
 \ts $\bM:=\bM_\al$, \ts $\al \in \cL$,  is
\defn{hyperbolic}:
\begin{equation}\label{eq:Hyp-intro}
	\langle \vb, \bM \wb \rangle^2 \  \geq \  \langle \vb, \bM \vb \rangle  \langle \wb, \bM \wb \rangle \quad \text{for every \ \, $\vb, \wb \in \Rb^{\ar}$, \ \  such that \ \ $\langle \wb, \bM \wb \rangle > 0$}.   \tag{Hyp}
\end{equation}
Log-concavity inequalities now follow from~\eqref{eq:Hyp-intro} for the matrix
\ts $\bM_\emp$, by interpreting the inner products as numbers \ts $\aL_\aq(k)$,
\ts $\aL_\aq(k-1)$ \ts and \ts $\aL_\aq(k+1)$, respectively.

We prove~\eqref{eq:Hyp-intro} by induction, reducing the claim for \ts $\bM_\al$ \ts
to that of \ts $\bM_{\al x}$ \ts, for all \ts $x\in \Cnt(\al)$.
Proving~\eqref{eq:Hyp-intro}
for the base of induction required the eigenvalue interlacing argument,
cf.~$\S$\ref{ss:finrem-interlace}.    This is
where our conditions for the weight function~$\ap$ appear in the calculation.
We also need a few other properties
of the atlas.  Notably, we require every matrix \ts $\bM_\al$ \ts to be \defn{irreducible}
with respect to its support, but that is proved by a direct combinatorial argument.

\smallskip

For other log-concavity inequalities in the paper, we consider similar atlas constructions
and similar claims.  For the equalities, we works backwards and observe that we need
equations~\eqref{eq:Hyp-intro} to be equalities.  These imply the local properties
which must hold for certain edges \ts $(\al,\al x)\in \Ef$.  Analyzing these properties
gives the equality conditions we present.

\medskip

\subsection{Discussion} \label{ss:intro-discussion}
Skipping over the history of the subject (see Section~\ref{s:hist}), in recent
years a great deal of progress on the subject was made by Huh and his coauthors.
In fact, until the celebrated Adiprasito--Huh--Katz paper~\cite{AHK}, even the
log-concavity for the number of $k$-forests (Welsh--Mason conjecture for
graphical matroids), remained open.
That paper was partially based on the earlier work \cite{Huh0,Huh1,HK},
and paved a way to a number of further developments, most notably
\cite{ADH,BES,BST,BHMPWa,BHMPWb,HSW,HW}.

From the traditional order theory point of view, the level of algebra
used in these works overwhelms the senses.  The inherent rigidity of
the original algebraic approach required either to extend the algebra
as in the papers above, or to downshift in the technology.
The \emph{Lorentzian polynomials} approach developed by
Br\"and\'en--Huh~\cite{BH18,BH} and by Anari~et.~al~\cite{ALOV}
allowed stronger results such as Theorem~\ref{t:matroids-BH}
and led to further results and applications such as
\cite{ALOV-RW,BLP,HSW,MNY}.  This paper represented the first
major downshift in the technology.

\smallskip

\nin
$(\circ)$ \.
A casual reader can be forgiven in thinking of this paper as a
successful deconstruction of the Lorentzian polynomials into the
terminology of linear algebra.  This is the opposite of what
happens both mathematically and philosophically.  Our approach
does in fact contain much of the Lorentzian polynomials approach
as a special case (cf.~$\S$\ref{ss:finrem-would}).  This can be
made precise, but we postpone that discussion until~\cite{CP}.

However, viewing greedoids and its special
cases as \emph{languages} allows us to reach far beyond what
the Lorentzian polynomials possibly can.\footnote{Lest one think
to use a straightforward generalization to noncommutative polynomials,
try imagining the right notion of a \emph{partial derivative}
which plays a crucial role in \cite{ALOV,BH}.}
To put this precisely, our maps \ts
$\bT_\al^{\langle x\rangle}$ \ts have a complete flexibility in
their definition.  In the world of Lorentzian polynomials, the
corresponding maps are trivial.  We trade the elegance of that
approach to more complexity, flexibility and strength.

\smallskip

\nin
$(\circ)$ \.
The true origin of our ``combinatorial atlas'' technology lies
in our deconstruction of the Stanley inequality~\eqref{eq:Sta}.
This is both one of oldest and the most mysterious results in the
area, and our proof is elementary but highly technical, more so
than our proof of greedoid results.

To understand the conundrum Stanley's inequality represents,
consider the original proof in~\cite{Sta} which is barely a
page long via a simple reduction to the classical
\emph{Alexandrov--Fenchel inequality}.  The latter is
a fundamental result on the subject, with many different
proofs across the fields, all of them difficult (see~$\S$\ref{ss:hist-AF}).
This difficulty represented the main obstacle in obtaining
an elementary proof of Stanley's inequality.

\smallskip

\nin
$(\circ)$ \.
Most recently, the new proof of the Alexandrov--Fenchel inequality
by Shenfeld and van Handel~\cite{SvH19} using ``Bochner formulas'',
renewed our hopes for the elementary proof of Stanley's inequality.
Their proof exploits the finiteness of the set of normals to polytope
facets in a very different way from Alexandrov's original approach
in~\cite{Ale}, see discussion in \cite[$\S$6.1]{SvH19}.
Our next point of inspiration was a most recent paper~\cite{SvH}
by Shenfeld and van Handel, where the authors obtain the
equality conditions for Stanley's inequality
(see Theorem~\ref{t:Sta-equality}) with applications
to Stanley's inequality (cf.~$\S$\ref{ss:finrem-AF}).

Deconstruction of \cite{SvH19,SvH} combined with ideas
from~\cite{BH,Sta} and our earlier work \cite{CPP1,CPP2},
led to our ``combinatorial atlas'' approach.  Both the Stanley
inequality and the conditions for equality followed from our
linear algebra setting and became amenable to generalizations.
Part of the reason for this is the explicit construction of
maps \ts $\bT_\al^{\langle x\rangle}$, which for convex polytopes are shown
in~\cite{SvH19} to exist only indirectly albeit in greater
generality, see also~$\S$\ref{ss:finrem-AF-proof}.

\smallskip

\nin
$(\circ)$ \.
Now, once we climbed the mountain of Stanley's inequality by means
of the new technology, going down to poset antimatroids, polymatroids
and matroids became easier.  Our ultimate extension to interval
greedoids required additional effort, as evidenced in the
technical definitions in Section~\ref{s:main}.  Furthermore,
our approach retained the flexibility of allowing us to match
the results with equality conditions.

\smallskip

\nin
$(\circ)$ \.
In conclusion, let us mention that the ultimate goal we set out
in~\cite{Pak}, remains unresolved.  There, we observed that
the Adiprasito--Huh--Katz inequalities for graphs and Stanley inequalities
for numbers of linear extensions correspond to
nonnegative integer functions in \ts $\GapP=\SP-\SP$.  We
asked whether these functions are themselves in~$\SP$.  This amounts
to finding a \emph{combinatorial interpretation} for the difference
of the LHS and the RHS of these inequalities.  While we use only
elementary tools, the eigenvalue based argument is not direct
enough to imply a positive answer.  See~$\S$\ref{ss:finrem-complexity-GapP}
for more on this problem.

\medskip

\subsection{Paper structure}\label{ss:intro-structure}
We start with basic definitions and notions in Section~\ref{s:def}.
In the next Section~\ref{s:main} we present the main results of the
paper on log-concave inequalities and the matching equality conditions
for interval greedoids.  We follow in Section~\ref{s:reduction} with
a chain of combinatorial reductions explaining how our greedoids
results imply poset antimatroid, polymatroid and matroid results.

In Section~\ref{s:atlas} we introduce the notion of \emph{combinatorial atlas},
which is the main technical structure of this paper.   We then show how
to derive log-concave inequalities in this general setting.  The key
combinatorial properties of the atlases are given in Section~\ref{s:Pull}.
In the next Section~\ref{s:atlas-equality}, we show that under additional
conditions on the atlas, we can characterize the equality conditions.

\smallskip

From this point on, much of the paper occupy proofs of the results:

\smallskip

\nin
\hskip1.cm  $\circ$ \ Thm~\ref{t:greedoid} (interval greedoids inequality) is proved in Section~\ref{s:proof-greedoid},

\smallskip

\nin
\hskip1.cm   $\circ$ \  Thm~\ref{t:greedoid-equality-full} (interval greedoids equality conditions) is proved  in Section~\ref{s:proof-greedoid-equality},

\smallskip

\nin
\hskip1.cm   $\circ$ \  Thm~\ref{t:matroids-Par-weighted}, Thm~\ref{t:matroids-equality-weighted}, Thm~\ref{t:matroids-equality-Par} (matroid inequality and equality conditions) are proved in  Section~\ref{s:proof-matroid};

\nin
\hskip1.42cm  in addition, this section includes proof of Prop.~\ref{p:graphical-equality}, further results on log-concavity for graphs

\nin
\hskip1.42cm
($\S$\ref{ss:proof-matroid-more-graph}), and examples of combinatorial atlases ($\S$\ref{ss:proof-matroid-ex}),

\smallskip

\nin
\hskip1.cm   $\circ$ \  Thm~\ref{t:polymatroids-Par}, Thm~\ref{t:polymatroids-BH-equality} and Thm~\ref{t:polymatroids-Par-equality} (discrete polymatroid inequality and equality conditions)

\nin
\hskip1.42cm
are proved  in Section~\ref{s:proof-polymatroid},

\smallskip

\nin
\hskip1.cm   $\circ$ \  Thm~\ref{t:antimatroids}, Thm~\ref{t:antimatroids-equality} and Thm~\ref{t:antimatroids-equality-total} (poset antimatroid inequality and equality conditions)  are

\nin
\hskip1.42cm
 proved  in Section~\ref{s:proof-antimatroid},

\smallskip

\nin
\hskip1.cm   $\circ$ \  Thm~\ref{t:morphisms-refined}, Thm~\ref{t:morphisms-EH-equality}
and Thm~\ref{t:morphisms-Par-equality} (morphism of matroids inequality and equality conditions)

\nin
\hskip1.42cm
are proved in Section~\ref{s:proof-morphism},

\smallskip

\nin
\hskip1.cm   $\circ$ \ Thm~\ref{t:Sta-weighted} (weighted Stanley's inequality) is proved in Section~\ref{s:proof-Sta}; in addition,
this section includes

\nin
\hskip1.42cm   $\S$\ref{ss:LE-belts} on posets with belts and an example $\S$\ref{ss:LE-example} of a combinatorial atlas in this case,

\smallskip

\nin
\hskip1.cm   $\circ$ \ Thm~\ref{t:Sta-equality-weighted} (equality condition for weighted Stanley's inequality)
is proved in Section~\ref{s:proof-Sta-equality}.

\smallskip

\nin
These last two sections are the most technically involved parts of this paper. Note that
although Sections \ref{s:proof-matroid}--\ref{s:proof-morphism} are somewhat independent,
we do recommend the reader start with the matroid proofs in
Section~\ref{s:proof-matroid} because of the examples and as a starting point of
generalizations, and antimatroid proofs in Section~\ref{s:proof-antimatroid}
because it has the shortest and cleanest reduction to the earlier greedoid results.

\smallskip

We conclude the paper with a lengthy historical Section~\ref{s:hist} which cover
to some degree various background behind results int he introduction.  Since the
material is so vast, we are somewhat biased towards most recent and general results.
We present final remarks and open problems in Section~\ref{s:finrem}.

\bigskip

\section{Definitions and notations}\label{s:def}

\subsection{Basic notation} \label{ss:sed-basic}
We use \ts $[n]=\{1,\ldots,n\}$, \ts $\nn = \{0,1,2,\ldots\}$,
\ts $\Zb_+ = \{1,2,\ldots\}$, \ts $\rrp=\{x\ge 0\}$ \ts and \ts $\rrs=\{x> 0\}$.
For a subset \ts $S\subseteq X$ \ts and element \ts $x\in X$,
we write \ts $S+x:=S\cup \{x\}$ \ts and \ts $S-x:=S\sm\{x\}$.

\medskip

\subsection{Matrices and vectors}
%
Throughout the paper we denote matrices with bold capitalized letter and
the entries by roman capitalized letters: \ts $\bM=(\aM_{ij})$.  We also keep
conventional index notations, so, e.g., \. $\big( \bM^3+ \bM^2 \big)_{ij}$ \. is the $(i,j)$-th
matrix entry of \. $\bM^3+\bM^2$.
We denote vectors by bold small letters, while vector entries by either unbolded uncapitalized letters
or vector components, e.g.\ \ts $\hb = (\ah_1,\ah_2,\ldots)$ \ts and \ts $\ah_i= (\hb)_i$.

A real matrix (resp., a real vector) is \defn{nonnegative} if all its entries are nonnegative real numbers, and
is \defn{strictly positive} if all of its entries are positive real numbers.
The \defn{support} of a real \ts $\ar \times \ar$ \ts symmetric matrix $\bM$ is defined as:
 \[
 \supp\ts (\bM) \ := \ \bigl\{ \.  i\in [\ar] \. : \. \aM_{ij} \neq 0  \ \text{ for some } j \in [\ar] \. \bigr\}.
 \]
In other words,  $\supp(\bM)$ \ts is  the set of indexes for which the corresponding row and
column of $\bM$ are nonzero vectors.
Similarly, the \defn{support} of a real $\ar$-dimensional vector $\hb$ is defined as:
    \[  \supp\ts (\hb) \ := \ \{ \.  i \in [\ar] \. : \. \ah_{i} \neq 0  \. \}. \]
For  vectors \ts $\vb, \wb \in \Rb^{\ar}$, we write \ts $\vb \leqslant \wb$ \ts
to mean the componentwise inequality, i.e.\ \ts $\av_i \leq \aw_i$ \ts for all $i \in [\ar]$.
We write \ts $|\vb| := \av_1+\ldots+\av_{\ar}$.  We also use \ts
$\eb_1,\ldots, \eb_{\ar}$ \ts to denote the standard basis of \ts $\Rb^{\ar}$.

Finally, for a subset \ts $S\subseteq [\ar]$,
   the \defn{characteristic vector} of~$S$  is the vector \. $\vb \in \Rb^{\ar}$ \. such that
  \ts $\av_i=1$ \ts if \ts $i \in S$ and \ts $\av_i=0$ \ts if \ts $i \notin S$.
   We use \ts $\0 \in \Rb^{\ar}$ \ts to denote the zero vector.

 \smallskip

 \subsection{Words}
 For a finite ground set $X$,
 we denote by  $X^*$ the set of all sequences $x_1\cdots x_\ell$ $(\ell \geq 0)$ of elements $x_i \in X$ for $i \in [\ell]$.
 We call an element of $X^*$ a \defn{word} in the \defn{alphabet}~$X$.  By a slight abuse of notation we use \ts
 $x_i$ \ts to also denote the $i$-th letter in the word~$\al$.
 The \defn{length} of a word \ts $\alpha=x_1 \cdots x_\ell$ \ts is the number of letters $\ell$ in the word, and is denoted by $|\alpha|$.
 The \defn{concatenation} \ts $\alpha \beta$ \ts of two words $\alpha$ and $\beta$
 is the string $\alpha$ followed by the string~$\beta$.  In this case \ts $\alpha$ \ts is called a \defn{prefix} of \ts $\alpha \beta$.
 For every \ts $\alpha=x_1\cdots x_\ell \in X^*$, we write \ts $z \in \alpha$ \ts if \ts $x_i=z$ \ts for some \ts $i\in [\ell]$.

\medskip

\subsection{Posets}
    A \defn{poset} \ts $\Pf=(X,\prec)$ \ts is a pair of \defn{ground set}~$X$ and a \defn{partial order} ``$\prec$'' on~$X$.
    For $x,y \in X$, we say that $y$ \defn{covers} $x$ in~$\Pf$, write \ts $x \ts \lcov\ts y$,
    if \ts $x \prec y$, and there exists no \ts $z \in X$ \ts such that \. $x \prec z \prec y$.
    For \ts $x,y \in X$, we write \ts $x \ts || \ts y$ \ts  if $x$ and $y$ are \defn{incomparable} in \ts $\Pf$.
    Denote by \ts $\inc(x)\subset X$ \ts the subset of elements \ts $y\in X$ \ts incomparable with~$x$.

    A \defn{lower ideal} of $\Pf$ is a subset $S \subseteq X$ such that, if $x \in S$ and $y\prec x$, then $y \in S$.
    Similarly, an \defn{upper ideal} of $\Pf$ is a subset $S \subseteq X$ such that, if $x \in S$ and $y \succ x$, then $y \in X$.
    The \defn{Hasse diagram} $\Hc:=\Hc_{\Pf}$ of~$\Pf$ \ts is the acyclic digraph with $X$ as the vertex set,
    and with $(x,y)$ as an edge if \ts $x \ts \lcov \ts y$.

    A \defn{chain} of $\Pf$ is a subset of $X$ that is totally ordered: \. $x_1 \prec x_2 \prec \ldots \prec x_\ell$\ts.
    An antichain is a subset \ts $S\subset X$, such that every two elements in $S$ are incomparable.
    \defn{Height} of a poset \ts $\height(\cP)$ \ts is
    the length of the maximal chain in~$\cP$.  Similarly, \defn{width} of a poset \ts $\width(\cP)$ \ts
    is the size of the maximal antichain in~$\cP$.
Element \ts $x\in X$ \ts is called \defn{minimal} if there is no \ts $y \in X$, s.t.\ $y\prec x$.  Define
\defn{maximal} elements similarly.

%

\bigskip

\section{Combinatorics of interval greedoids} \label{s:main}

\subsection{Preliminaries} \label{ss:main-prelim}
Let \ts $\Gf=(X,\cL)$ \ts be an interval greedoid of rank \ts $m:=\rk(\Gf)$.
Recall the definitions of
\ts $\Par(\alpha)$ \ts and \ts $\Des_\alpha(x)$ \ts given
in~$\S$\ref{ss:intro-antimatroid-equality} above, and note that
\ts ``$\salp$'' \ts remains an equivalence relation, see
Proposition~\ref{p:equiv-greedoid}.

For all \ts $\al\in \cL$ \ts and \ts $x,y\in X$, define \defn{passive} and \defn{active non-continuations} as follows:
\begin{align*}
\Pas_{\alpha}(x,y) \ &:= \ \bigl\{ z \in X \ : \ \alpha z \notin \Lc,  \ \alpha xz, \alpha yz \notin \Lc, \ \alpha xyz \in \Lc  \bigr\},\\
\Act_{\alpha}(x,y) \ &:= \ \bigl\{ z \in X \ :  \ \alpha z \notin \Lc, \   \alpha xz, \alpha yz \in \Lc, \ \alpha xyz \in \Lc  \bigr\}.
\end{align*}

\smallskip

Let \ts $\aq: \Lc \to \rrs$ \ts be a positive \defn{weight function}, which we extend
to \ts $\aq: X^\ast \to \rr$ \ts by setting \ts $\aq(\al)=0$ \ts for all $\al\notin \cL$.
Let \. $\bbc=(\ac_0,\ldots,\ac_m)\in \rrs^{m+1}$, where \ts $m=\rk(\Gf)$,
be a fixed positive sequence, which we call the \defn{scale sequence}.
Consider another weight function \ts $\ap: X^\ast \to \rr$\ts:
\begin{equation}\label{eq:ap-def}
\ap(\al) \, := \, \frac{\aq(\al)}{\ac_{\ell}}\,, \quad \text{where \ $\ell=|\al|$ \ and \ $\al \in X^\ast$,}
\end{equation}
which we call the \defn{scaled weight function}.

\medskip

\subsection{Properties} \label{ss:main-properties}
Fix weight function \ts $\aq: \Lc \to \rrs$ \ts and scale sequence \ts $\bbc\in \rrs^{m+1}$.
For every word \ts $\al\in \cL$ \ts of length \ts $\ell := |\al|$, consider the following
properties.

\bigskip

\nin
{\small \bf 1.} {\defngb{Continuation invariance property:}}
\begin{equation}\label{eq:ContInv} \tag{ContInv}
		\aq(\alpha  xy \beta) \ = \  \aq(\alpha  yx \beta)
\quad \text{for all \  $x,y \in \Cnt(\alpha)$ \, and \ $\beta \in X^*$.}
\end{equation}

\smallskip

\nin
Note that by the exchange property, we have \ts $\alpha xy \beta \in \Lc$ \.  if and only if \. $\alpha yx \beta \in \Lc$.

\bigskip

\nin
{\small \bf 2.} {\defngb{Passive-active monotonicity property:}}
\begin{align}\label{eq:PasAct} \tag{PAMon}
\sum_{z \ts \in \ts \Pas_{\alpha}(x,y)} \, \sum_{\beta \ts \in \ts  \Cnt_k(\alpha xyz)} \, \aq\big(\alpha xyz \beta\big)
 \ \, \geq  \ 	\,
\sum_{z \ts \in \ts  \Act_{\alpha}(x,y)} \, \sum_{\beta \ts \in \ts  \Cnt_k(\alpha xyz)} \, \aq\big(\alpha xyz \beta\big)\ts,
\end{align}
for all distinct \ts $x,y \in \Cnt(\alpha)$, and \ts $k \geq 0$.  We also have a stronger
property stated in terms of~$\cL$\ts.

\bigskip

\nin
{\small \bf 2$'$.}  \defngb{Weak local property:}
\begin{align}\label{eq:LI} \tag{WeakLoc}
x,y,z \in X \ \, \ \text{s.t.} \. \ \ \alpha xz, \ts \alpha yz, \ts \alpha xyz  \in \Lc
\quad \Rightarrow \quad \alpha z \in \Lc.
\end{align}

\smallskip

\nin
Observe that~\eqref{eq:LI} implies that \ts $\Act_{\alpha}(x,y)= \varnothing$ \ts for
all distinct \ts $x,y\in \Cnt(\al)$,
which in turn trivially implies \eqref{eq:PasAct}.  Note also that~\eqref{eq:LI} is a
property of a greedoid rather than the weight function.  Greedoids that
satisfy~\eqref{eq:LI} are called \ts \defna{weak local greedoids}.\footnote{This is
a new class of greedoids which is similar but more general than the
\emph{local poset greedoids}.  See Section~\ref{s:reduction} for the properties
of weak local greedoids, relationships to other classes, and $\S$\ref{ss:finrem-local}
for further background.
}

\bigskip

\nin
{\small \bf 3.} \defngb{Log-modularity property:}
\begin{align}\label{eq:LogMod} \tag{LogMod}
		\ap(\alpha x) \,\. \ap(\alpha y) \ = \  \ap(\alpha) \,\. \ap(\alpha xy) \quad
\text{for  all \ $x,y \in \Cnt(\alpha)$ \  s.t.\ \   $\alpha xy \in \Lc$}.
\end{align}

\bigskip

\nin
{\small \bf 4.} \defngb{Few descendants property:}
\begin{equation}\label{eq:FewDes} \tag{FewDes}
|\Cc| \geq  2 \quad \Rightarrow \quad \Des_{\alpha}(x) \. = \. \varnothing\ts,
\quad \text{ for every } \. x \in \Cc \, \text{ and } \, \Cc \in \Par(\alpha).
\end{equation}

\smallskip

\nin
Note that \eqref{eq:FewDes} is  satisfied if \. $|\Cc|\le 1$, or if \. $\Des_{\alpha}(x) =\varnothing$.

\bigskip

\nin
{\small \bf 5.} \defngb{Syntactic monotonicity property:}
\begin{align}\label{eq:SynMon} \tag{SynMon}
\ap(\alpha x)^2   \ \geq \    \sum_{y \ts\in\ts \Des_{\alpha}(x)}  \ \ap(\alpha) \,\. \ap(\alpha xy)\ts, \quad
\text{for all \ $x \in \Cnt(\alpha)$}.
\end{align}

\medskip

For all \ts $\Cc \in \Par(\alpha)$, define
\begin{align}\label{eq:ab-def}
		\ab_\alpha(\Cc) \ := \
		\begin{cases}
			\displaystyle \ \, \sum_{y \ts\in\ts \Des_{\alpha}(x)}  \ \frac{\ap(\alpha) \,\. \ap(\alpha xy)}{\ap(\alpha x)^2}
			 & \text{ if } \ \Cc  =  \{x\},\\
			 			\quad 0 & \text{ if }  \ |\Cc|\geq 2.
		\end{cases}
\end{align}

\nin
Note that properties~\eqref{eq:FewDes} and~\eqref{eq:SynMon} imply that \. $\ab_{\alpha}(\Cc) \leq 1$ \ts for all
\ts $\Cc \in \Par(\alpha)$.  This sets up our final

\bigskip

\nin
{\small \bf 6.} \defngb{Scale monotonicity property:}
	\begin{align}\label{eq:ScaleMon}\tag{ScaleMon}
\bigg(1\.- \. \frac{\ac_{\ell+1}^2}{\ac_{\ell} \. \ac_{\ell+2}} \bigg)
		\sum_{\Cc \ts\in\ts \Par(\alpha)} \. \frac{1}{1 \. - \. \ab_\alpha(\Cc)}
 \ \leq \ 1  \., \quad
\text{for all
\ $\Cc \in \Par(\alpha)$}.
	\end{align}

\smallskip
\nin
We adopt the convention that \eqref{eq:ScaleMon} is always satisfied whenever \ts $\ac_{\ell+1}^2 \ts \geq \ts \ac_{\ell} \. \ac_{\ell+2}$ \ts (because then the LHS is considered nonpositive), and that \ts $\ab_{\alpha}(\Cc) \ts < \ts 1$ \ts for all  \ts $\Cc \in \Par(\alpha)$ \ts whenever \ts $\ac_{\ell+1}^2 \ts < \ts \ac_{\ell} \ts \ac_{\ell+2}$ \. (as otherwise the LHS is considered to be $\infty$) .
In particular, note that \eqref{eq:ScaleMon} is  satisfied for the uniform scale sequence \ts
$\bbc = (1,\ldots,1)$\ts.

\smallskip

\begin{rem}
The last four properties \eqref{eq:LogMod}, \eqref{eq:FewDes}, \eqref{eq:SynMon}
and \eqref{eq:ScaleMon} have a linear algebraic interpretation as certain
matrix being \emph{hyperbolic}.  We postpone a discussion of this until
the next section.
\end{rem}
\medskip

\subsection{Admissible weight functions}\label{ss:main-results}
We can now give the main definition used in the first main result of the
paper (Theorem~\ref{t:greedoid}).

\smallskip

\begin{definition}[{\rm \defna{$k$-admissible weight functions}}] \label{d:admissible} {
Let \ts $\Gf=(X,\cL)$ \ts be an interval greedoid of rank \ts $m:=\rk(\Gf)$,
and let \ts $1\le k < m$.  Weight function \ts $\apr: \cL\to \rrs$ \ts
is called \. \defng{$k$-admissible}, if there is a scale sequence
\ts $\bbc=(\acr_0,\ldots,\acr_m)\in \rrs^{m+1}$, such that properties \.
\eqref{eq:ContInv}, \. \eqref{eq:PasAct}, \. \eqref{eq:LogMod}, \. \eqref{eq:FewDes}, \. \eqref{eq:SynMon} and \eqref{eq:ScaleMon}
are satisfied for all \ts $\al \in \cL$ \ts of length \ts $|\al| <k$. }
\end{definition}

\smallskip

We can also state our second main result of the paper, which
gives the third equivalent condition in Theorem~\ref{t:greedoid-equality}
that is both more detailed and useful in applications.

\smallskip

\begin{thm}[{\rm \defng{Equality for interval greedoids}, {\em \defna{second main theorem}}}{}]
\label{t:greedoid-equality-full}
Let \ts $\Gf=(X,\cL)$ \ts be an interval greedoid of rank \ts $m:=\rk(\Gf)$,
let \ts $1\le k <m$, and let \.\ts $\aqr:\cL \to \rrs$ \. be a \ts {$k$-admissible}
\ts weight function with a scale sequence \ts $\bbc=(\acr_0,\ldots,\acr_m)\in \rrs^{m+1}$.
Then, \ts \underline{the following are equivalent}\ts:

\medskip

\nin
\qquad $\mathbf{a.}$ \, We have:
\begin{equation}\label{eq:item-equality-greedoid-a}\tag{GE-$a$}
\aLr_\aqr(k)^2 \ = \   \aLr_\aqr(k-1) \. \cdot \. \aLr_\aqr(k+1)\ts.
\end{equation}

\medskip

\nin
\qquad $\mathbf{b.}$ \, There is \. $\asr(k-1)>0$,  such that for every \. $\al\in \cL_{k-1}$ \. we have:
\begin{equation}\label{eq:item-equality-greedoid-b}\tag{GE-$b$}
\aLr_{\aqr,\alpha}(2)  \ = \  \asr(k-1) \, \aLr_{\aqr,\alpha}(1) \ = \  \asr(k-1)^2 \, \aLr_{\aqr,\alpha}(0)\ts.
\end{equation}

\medskip

\nin
\qquad $\mathbf{c.}$ \, There is \. $\asr(k-1)>0$,  such that for every \. $\al\in \cL_{k-1}$ \. we have:
		\begin{align}
			\label{eq:GrEqu1}\tag{GE-$c$\ts$1$}  & \sum_{x \ts \in\ts \Cnt(\alpha)}  \.
\frac{\aqr(\alpha x)}{\aqr(\alpha)} \ = \  \asr(k-1), \ \ \text{and}\\
			\label{eq:GrEqu2}\tag{GE-$c$\ts$2$}
& \bigl(1\ts - \ts \abr_{\alpha}(\Cc)\bigr) \, \sum_{x\ts \in\ts \Cc} \ \frac{\aqr(\alpha x)}{\aqr(\alpha)}
			\ = \ \asr(k-1) \, \bigg( 1 \. - \. \frac{\ac_{k}^2}{\ac_{k-1} \.\ac_{k+1}} \bigg)
			 \quad \ \, \text{ for all \ \ $\Cc \in \Par(\alpha)$,}
		\end{align}
where \ts $\abr_{\alpha}(\Cc)$ \ts is defined in~\eqref{eq:ab-def}.
\end{thm}
\smallskip

Note that~\eqref{eq:GrEqu1} and \eqref{eq:GrEqu2} imply that \eqref{eq:ScaleMon} is always an equality
for \. $\alpha \in \Lc_{k-1}$.

\smallskip

\begin{rem}\label{r:main-nonexist}
Note that the $k$-admissible property of weight functions~$\aq$ is quite constraining
and there are interval greedoid for which there are no such~$\aq$.  Given the abundance
of examples where such weight functions are natural, we do not investigate the structural
properties they constrain (cf.~$\S$\ref{ss:hist-weight}).
\end{rem}

\bigskip


\section{Combinatorial preliminaries} \label{s:reduction}

In this section we present basic properties of matroids, polymatroids,
poset antimatroids, local poset greedoids and interval greedoids.
We include the relations between these classes which will be important
in the proofs.  Most of these are relatively straightforward,
but stated in a different way and often dispersed
across the literature.  We include the short proofs for completeness
and as a way to help the reader get more familiar with the notions.
The reader well versed with greedoids can skip this section
and come back whenever proofs call for the specific results.

\smallskip

\subsection{Equivalence relations}\label{ss:reduction-equiv}
Here we prove that equivalence relations given in the introduction
are well defined.  We include short proofs both for completeness.

\smallskip

\begin{prop}  \label{p:equiv-matroid}
Let \ts $\Mf=(X,\Ic)$ \ts be a matroid, and let \ts $S\in \cI$ \ts be an independent set.
Then the relation \ts ``$\sim_S$'' defined in \ts {\rm $\S$\ref{ss:intro-matroids}}
is an equivalence relation.
\end{prop}

\begin{proof}
Observe that  \. $x \sim_S y$ \. if and only if \ts $x$ \ts and \ts $y$ \ts are \emph{parallel}
in the matroid \ts $\Mf/S$ \ts obtained from~$\Mf$ by contracting over~$S$.
\end{proof}

\smallskip

\begin{prop}  \label{p:equiv-polymatroid}
Let \ts $\cD=([n],\Jc)$ \ts be a discrete polymatroid, and let \ts $\ba\in \Jc$ \ts
be an independent multiset.
Then the relation \ts ``$\sim_\ba$'' defined in \ts {\rm $\S$\ref{ss:intro-poly}}
is an equivalence relation.
\end{prop}

\begin{proof}
It suffices to prove transitivity of ``$\sim_\ba$'', as reflexivity and symmetry
follow immediately from the definition. Let \ts $i\sim_\ba$ \ts and \ts $j \sim_\ba k$.
Suppose to the contrary that \ts $i \not \sim_\ba k$, so  \ts $\ba +\eb_i +\eb_k \in \Jc$.
On the other hand, \ts $\ba +\eb_j \in \Jc$ \ts since \ts $j \in\Cnt(\ba)$.
It then follows from applying the exchange property to \ts $\ba +\eb_j$ \ts
and \ts $\ba +\eb_i +\eb_k$, that either \ts $\ba +\eb_j +\eb_i \in \Jc$ \ts
or \ts $\ba +\eb_j +\eb_k \in \Jc$, both of which give us a contradiction.
\end{proof}

\smallskip

\begin{prop}  \label{p:equiv-greedoid}
Let \ts $\Gf=(X,\cL)$ \ts be an interval greedoid, and
let \ts $\al\in \cL$ \ts be a fixed word.
Then the relation \ts ``$\salp$'' defined in \ts {\rm $\S$\ref{ss:intro-antimatroid}}
is an equivalence relation.
\end{prop}

\begin{proof}
Reflexivity follows immediately from the definition.
For the symmetry, let \ts $x \salp y$ \ts and suppose to the contrary that \ts
$y \not \sim_\alpha x$.  This is equivalent to \ts $\alpha yx \in \Lc$.
On the other hand, \ts $\alpha x \in \Lc$ \ts since \ts $x \in \Cnt(\alpha)$.
It then follows from applying the exchange property to \ts $\alpha x$ \ts
and \ts $\alpha yx$ \ts that \ts $\alpha xy \in \Lc$, which contradicts
\ts $x \salp y$.

For transitivity, let \ts $x\sim_\al y$ \ts and \ts $y\sim_\al z$.
Suppose to the contrary, that \ts $x \not \sim_\al z$, so  \ts $\alpha xz \in \Lc$.
On the other hand, \ts $\alpha y \in \Lc$ \ts since \ts $y \in\Cnt(\al)$.
It then follows from applying the exchange property to \ts $\alpha y$ \ts
and \ts $\alpha xz$, that either \ts $\alpha yx \in \Lc$ \ts or
\ts $\alpha yz \in \Lc$, both of which gives us a contradiction.
\end{proof}


We conclude with another equivalence relation, which will prove important
in~$\S$\ref{ss:proof-morphism-Hyp}.
Let \ts $\Phi: \Mf\to \Nf$ \ts be a morphism of matroids,
let \ts $f$ \ts be the rank function for $\Mf=(X,\Ic)$, and let \ts $g$ \ts
be the rank function for $\Nf=(Y,\Jc)$.  For an independent set \ts $S\in \Ic$,
let \ts $H \subseteq X$ \ts be given by
\begin{equation}\label{eq:reduction-morphisms-H-def}
H \ := \ \bigl\{ x \in X \setminus S \ : \   \. \rkn\big( \Phi( S+x) \big) = \rk(\Nf)-1 \bigr\}.
\end{equation}
Denote by \ts ``$\sim_H$'' \ts  the equivalence relation on~$H$, defined by
\begin{equation}\label{eq:reduction-morphisms-sim-def}
x \ts \sim_H \ts y \quad \Longleftrightarrow
\quad \rkn\big(\Phi(S+x+y)\big) \, = \, \rk(\Nf)-1.
\end{equation}

\begin{prop}\label{p:reduction-equiv-sim-H}
The relation \ts ``$\sim_H$'' \ts defined in~\eqref{eq:reduction-morphisms-sim-def}
is an equivalence relation.
\end{prop}

\begin{proof}
Reflexivity and symmetry follows directly from definition, so it suffices to prove transitivity.
Suppose that \. $x,y,z\in H$ \. are distinct elements, such that \. $x \sim_H y$ \. and
\. $y \sim_H z$.  Assume to the contrary, that \. $x \not\sim_H z$.  This implies
that \. $g\bigl(\Phi(S+x+z)\bigr)=\rk(\Nf)$.  Applying  the exchange property for matroid~$\Nf$
to  \. $\Phi(S+y)$ \. and \. \. $\Phi(S+x+z)$, we have that either \.
$g\bigl(\Phi(S+y+x)\bigr)=\rk(\Nf)$ \. or \. \. $g\bigl(\Phi(S+y+z)\bigr)=\rk(\Nf)$.
This contradicts the assumption, and completes the proof.
\end{proof}

\medskip

\subsection{Antimatroids \. $\ssu$ \. interval greedoids} \label{ss:reduction-antimatroid}
Note that \. (\emph{nondegenerate property}) \. defining the language of a  greedoid is vacuously true for  poset antimatroids.
Also note that two properties defining the language of a greedoid are identical
to those defining antimatroids:  \.
(\emph{normal property}) and \. (\emph{hereditary property}).
  Similarly,
the \. (\emph{exchange property}) \. for antimatroids is more restrictive
than the \. (\emph{exchange property}) \. for greedoids.

It remains to show that the \. (\emph{interval property}) \. holds for antimatroids.
Let \ts $\Af=(X,\cL)$ \ts be an antimatroid.  Suppose \.  $\alpha,\beta,\gamma \in X^*$ \.
and \ts $x \in X$, s.t.\  \. $\alpha x$, \ts $\alpha \beta \gamma x \in \Lc$.
Write \ts $\alpha':=\alpha x$ \ts and $\beta':=\alpha \beta$.
Then note that \ts $x \in \alpha'$ \ts and \ts $x\notin \beta'$, as otherwise \.
$w:=\alpha \beta \gamma x \notin \Lc$ \ts since $w$ is not a simple word,
and \ts $\alpha',\beta' \in \Lc$.  Also note that $x$ is the only letter in
$\alpha'$ that is not contained in~$\beta'$.  It then follows from the
(\emph{exchange property}) for~$\Af$, that \. $\alpha \beta x= \beta' x \in \Lc$,
as desired. \ $\sq$

\smallskip

\begin{prop}  \label{p:reduction-antimatroid-LI}
Let \ts $\cP=(X,\prec)$ \ts be a poset, and let \ts $\Af=(X,\cL)$ \ts be the corresponding antimatroid.
Then $\Af$ \ts satisfies the \emph{(interval property)}, \eqref{eq:FewDes} and  \eqref{eq:LI}.\footnote{Weak
local property does not hold for all antimatroids, but holds for all poset antimatroids.}
\end{prop}

\begin{proof}
The \ts (\emph{interval property}) is proved above for all antimatroids.  For~\eqref{eq:LI},
let \ts $x,y,z \in X$, s.t.\  \. $\alpha xz$, \ts $\alpha yz$, \ts $\alpha xyz  \in \Lc$.
Since $\alpha x z \in \cL$ and $y \notin \alpha xz$,
this implies $z$ is incomparable to $y$ in~$\cP$.
Together with \ts $\alpha yz \in \Lc$, this implies that \ts $\alpha z \in \Lc$, as desired.

For~\eqref{eq:FewDes}, note that \ts $\Af$ \ts satisfies
\begin{equation}\label{eq:poset basic}
	\text{$\alpha x, \. \alpha y \ts \in\Lc$\ts, \ \ $x, \ts y \ts\in X$ \ \ $\Longrightarrow$ \  \ $\alpha xy \in \Lc$\ts.}
\end{equation}
Indeed, this is because \. $\alpha y  \in \Lc$ \. implies that every element in $\Pf$ that is less than $y$ is contained in $\alpha$, so they are also contained in $\alpha x$.
This in turn implies that
$\alpha x y \in \Lc$.
Now note that \eqref{eq:poset basic} implies that $|\Cc|=1$ for every parallel class \ts
$\Cc\in \Par(\al)$ \ts of \ts $\alpha \in \Lc$,
and thus \eqref{eq:FewDes} is satisfied trivially.
\end{proof}

\medskip

\subsection{Matroids \. $\ssu$ \. greedoids}\label{ss:reduction-matroids}
Given a matroid \ts $\Mf = (X,\Ic)$, we construct the corresponding greedoid
\ts $\Gf=(X,\Lc)$, where \ts $\Lc$ \ts is defined as follows:
\[ \alpha \. = \. x_1 \ts\cdots \ts x_\ell \in \Lc \  \ \Longleftrightarrow \ \ \alpha
\ \ \text{is simple and} \ \{x_1,\ldots, x_\ell\} \in \Ic.   \]
Observe that \ts (\emph{nondegenerate property}) for~$\Gf$
follows from matroid~$\Mf$ being nonempty, \.
(\emph{normal property}) follows from definition, \.
(\emph{hereditary property}) for~$\Gf$ follows from the (\emph{hereditary property})
for~$\Mf$, and the (\emph{exchange property}) for~$\Gf$ follows
from (\emph{exchange property}) for~$\Mf$.
%

\smallskip

\begin{prop}\label{p:reduction-matroid-LI}
Given a matroid \ts $\Mf = (X,\Ic)$, the greedoid \ts $\Gf=(X,\Lc)$ \ts
constructed above satisfies the \ts {\em (interval property)},
\eqref{eq:FewDes} and \eqref{eq:LI}.
\end{prop}

\begin{proof}
Now note that, the greedoid \ts $\Gf$ \ts satisfies
\begin{equation}\label{eq:matroid basic}
	\text{$\alpha xy \ts \in\Lc$\ts, \ \ $x, \ts y \ts\in X$ \ \ $\Longrightarrow$ \  \ $\alpha y \in \Lc$\ts.}
\end{equation}
This follows from commutativity of~$\Lc$ and
the (\emph{hereditary property}) of~$\Mf$.
The (\emph{interval property}) for~$\Gf$ follows immediately from~\eqref{eq:matroid basic}.

Now, it follows from  \eqref{eq:matroid basic}  that \ts $\Des_{\alpha}(x) =\varnothing$ \ts
for every \ts $\alpha \in \Lc$ \ts and \ts $x \in X$, and \eqref{eq:FewDes} then follows trivially.
Finally, let \ts $x,y,z \in X$, s.t.\ $\alpha xz$, \ts $\alpha yz$, \ts $\alpha xyz  \in \Lc$.
Applying \eqref{eq:matroid basic} to $\alpha xz \in \Lc$, it then follows that \ts $\alpha z \in \Lc$.
This proves \eqref{eq:LI}, and completes the proof.
\end{proof}

\medskip

\subsection{Discrete polymatroids \. $\ssu$ \. greedoid} \label{ss:reduction-polytroid}
Given a discrete polymatroid \ts $\Df=([n],\Jc)$, we construct the corresponding
greedoid \ts $\Gf=(X,\Lc)$ \ts as follows.
Let \. $ X := \bigl\{x_{i \ts j} \ts :
\ts 1 \leq i,j \leq n \bigr\}$ \. be the alphabet.\footnote{Unlike the rest of the paper, here $|X|=n^2$.}

For every word \ts $\alpha \in X^*$, denote by \. $\ba_\al=(\aa_1,\ldots,\aa_n) \in \nn^n$ \ts
the vector counting the number of occurrences of $x_{i,*}$'s in $\alpha$, i.e.\
\. $\aa_i \ts := \ts \big|\big\{\ts j \in [n] \.:\. x_{i\ts j} \in \alpha \ts  \big\}\big|$.
The word \ts $\al\in X^\ast$ \ts is called \defnb{well-ordered}  if for every letter
\ts $x_{i\ts j}$ \ts in $\alpha$, letter \ts  $x_{i \ts j-1}$ \ts is also
in $\alpha$ \ts \emph{before} \ts $x_{i \ts j}$\ts.

Define \ts $\Lc$ \ts to be the set of simple well-ordered words \ts $\alpha \in X^*$, such that \ts $\ba_\al\in \Jc$.
Note that, each vector $\ba \in \Jc$ corresponds to \. $\binom{|\ba|}{\ts\aa_1\ts, \. \ldots \., \ts \aa_n\ts}$ \. many feasible words \ts
$\alpha \in \Lc$ \ts for which \ts $\ba_\alpha=\ba$.  Namely, these are all permutations of
the word \. $x_{1\ts 1}\cdots x_{1\ts \aa_1} \,
 \cdots \,  x_{n \ts 1} \cdots x_{n\ts \aa_n}$ \. preserving the relative order of letters \.
 $x_{i\ts 1},\ldots, x_{i\ts\aa_i}$.

For the greedoid~$\Gf=(X,\Lc)$, the (\emph{nondegenerate property}) and the (\emph{normal property})
follow from definition.  On the other hand, the (\emph{hereditary property})
and the (\emph{exchange property}) for~$\Gf$ follows from the corresponding properties
for~$\Df$. This completes the proof.  \ $\sq$

\smallskip

\begin{prop}\label{p:reduction-polymatroid-LI}
Given a discrete polymatroid \ts $\Df=([n],\Jc)$, the greedoid \ts $\Gf=(X,\Lc)$ \ts
constructed above satisfies the \ts {\em (interval property)}, \eqref{eq:FewDes}
and \eqref{eq:LI}.
\end{prop}

\begin{proof}
First, let us show that (\emph{interval property}) holds for~$\Gf$.
Let \. $\alpha,\beta,\gamma \in X^*$, and let \. $z=x_{i\ts j} \in X$ \ts \.
s.t.\  \. $\alpha z, \ts \alpha \beta \gamma z \in \Lc$.
Since \. $\alpha \beta \gamma z \in \Lc$, this implies that \. $x_{i \ts j+1}, \ldots, x_{i \ts n} \notin \beta$.
Since \ts $\alpha z \in \Lc$, this implies that $\alpha \beta z$ is well-ordered.
On the other hand, by applying the (\emph{hereditary property}) of~$\Df$ to the word \ts
$\alpha \beta \gamma z$, it then follows that \ts $\ba_{\alpha \beta z} \in \Jc$.
Hence, the word \ts $\alpha \beta z \in \Lc$, which proves the (\emph{interval property}).

Now, note that \ts $\Gf$ \ts satisfies
\begin{equation}\label{eq:polymatroid basic}
	\Des_{\alpha}(x_{ij}) \ \subseteq \ \{x_{i \ts j+1} \} \qquad \text{for every \, $\alpha \in \Lc$ \, and \, $x_{ij} \in \Cnt(\alpha)$.}
\end{equation}
For~\eqref{eq:LI},
let \ts $x,y,z \in X$, s.t.\ \. $\alpha xz$, \ts $\alpha yz$, \ts $\alpha xyz \ts \in \Lc$.
Suppose to the contrary, that \ts $\alpha z \notin \Lc$.
Since \ts $\alpha xz \in \Lc$ \ts and \ts $\alpha yz \in \Lc$,
this implies that  \ts $z \in \Des_{\alpha}(x)$ \ts and \ts $z \in \Des_{\alpha}(y)$.
On the other hand, this intersection is empty by \eqref{eq:polymatroid basic}.
This gives a contradiction, and proves~\eqref{eq:LI}.

For~\eqref{eq:FewDes},
let \ts $\ba=\ba_\al$ \ts where \ts $\alpha \in \Lc$, and let \. $x,y\in \Cnt(\alpha)$ \.
be distinct elements s.t.\  \ts $x \salp y$.
Let \ts $i,j \in [n]$ be such that
\. $\ba_{\alpha x} = \ba+\eb_i$ \.
and \.  $\ba_{\alpha y} = \ba+\eb_j$.
Note that \. $i \neq j$ \. and \. $\ba+\eb_i, \ba+ \eb_j \in \Jc$.
Suppose to the contrary, that \eqref{eq:FewDes} is not satisfied,
so we can assume that \.
$\Des_{\alpha}(x) \neq \varnothing$.
By  \eqref{eq:polymatroid basic}, this implies that \. $\ba+2\eb_i \in \Jc$.
Now, by applying the polymatroid exchange property
to \. $\ba+ \eb_j$ \. and \. $\ba+2\eb_i$,
we then have \. $\ba+ \eb_i + \eb_j \in \Jc$.
This contradicts the assumption that \ts $x \salp y$,
and proves~\eqref{eq:FewDes}.	
\end{proof}

\medskip

\subsection{Exchange property for morphism of matroids}\label{ss:reduction-morphism}
We will also need the following basic result.

\begin{prop}\label{p:reduction-morphism-base-exchange}
Let \ts $\Phi: \Mf \to \Nf$ \ts be a morphism of matroids \ts $\Mf=(X,\Ic)$ \ts
and \ts $\Nf=(Y,\Jc)$.  Let \ts $S,T \ssu X$, \ts $|S|=|T|$ \ts be two distinct bases of~$\Phi$.
Then there exists \ts $z \in S \setminus T$ \ts and \ts $w \in T \setminus S$ \ts such that \ts
$S-z+w$ \ts is also a basis of~$\Phi$.
\end{prop}

\begin{proof}
Fix an arbitrary $z \in S \setminus T$.  We split the proof into two cases.
First,  suppose that \ts $\Phi(S-z)$ \ts contains a basis of~$\Nf$.
Applying the exchange property of $\Mf$ to the independent sets $S-z$ and $T$,
there exists \ts $w \in T \setminus S$ \ts such that \ts $S':=S-z+w$ \ts
is an independent set of~$\Mf$.  Note that \ts $\Phi(S') \supset  \Phi (S-z)$ \ts
contains a basis of $\Nf$ by assumption, so $S'$ is a basis of $\Phi$, as desired.
	
Second, suppose that  \ts $\Phi(S-z)$ \ts does not contain a basis of~$\Nf$.
Applying the exchange property of $\Nf$ to \ts $\Phi(S-z)$ \ts and \ts $\Phi(T)$,
there exists \ts $w \in T \setminus S$ \ts such that \ts $\Phi(S-z+w)$ \ts
contains a basis of~$\Nf$.	Since $\Phi$ is a morphism of matroid, we have
	\begin{align*}
		\rkm(S-z+w) \. -  \. \rkm(S-z) \ \geq \ \rkn\big(\Phi(S-z+w)\big)
\. - \.  \rkn\big(\Phi(S-z)\big) \, = \, 1,
	\end{align*}
where \ts $f$ \ts and \ts $g$ \ts are rank functions in~$\Mf$ and~$\Nf$, respectively.
This implies that \ts $S-z+w$ \ts is an independent set of $\Mf$,
and therefore \ts $S-z+w$ \ts is a basis of the morphism~$\Phi$.
This completes the proof.
\end{proof}

\bigskip

\section{Combinatorial atlases and hyperbolic matrices}
\label{s:atlas}

In this section we introduce combinatorial atlases and present
the local–global principle which allows one to recursively
establish hyperbolicity of vertices. See~$\S$\ref{ss:finrem-main}
for some background.

\subsection{Combinatorial atlas}\label{ss:atlas-def}
Let \ts $\cP=(\Vf,\prec)$ \ts be a locally finite poset of bounded height.\footnote{In our
examples, the poset~$\cP$ can be both finite and infinite.}
Denote by \. $\Qf =(\Vf, \Ef)=\Hc_\cP$ \. be the acyclic digraph given by the
Hasse diagram of~$\cP$.  Let \ts $\Vfm\subseteq\Vf$ \ts be the set of maximal
elements in~$\cP$, so these are \defn{sink vertices} in~$\Qf$.  Similarly, denote by
\ts $\Vfp:=\Vf\sm \Vfm$ \ts the \defn{non-sink vertices}.  We write \ts $\vfs$ \ts
for the set of out-neighbor vertices \ts $v'\in \Vf$, such that \ts $(v,v')\in \Ef$.

\smallskip

\begin{definition}{\rm
A \defna{combinatorial atlas} \ts $\AA=\AA_\cP$ \ts of dimension \ts $\ar$ \ts
is an acyclic digraph   \. $\Qf:=(\Vf, \Ef)=\Hc_\cP$ \. with an additional structure:

\smallskip

$\circ$ \ Each vertex \ts $\vf \in \Vf$ \ts is associated with a  pair \ts $(\bM_{\vf}, \hb_{\vf})$\ts, where \ts
$\bM_{\vf}$ \ts is a nonnegative symmetric \ts

\quad \. $\ar \times \ar$ \ts matrix,
and \ts $\hb_{\vf}\in \rr_{\ge 0}^\ar$ \ts is a nonnegative vector.

\smallskip
	
$\circ$ \  Every vertex $v \in \Omega^+$ has outdegree $d$, and the outgoing edges of each vertex \ts $\vf \in \Vfp$ \ts
are labeled

\quad \. with indices \ts $i\in [\ar]$.
We denote the edge labeled \ts $i$ \ts as  \ts
$\ef^{\<i\>}=(\vf,\vf^{\<i\>})$, where \ts $1 \leq i \leq \ar$.

\smallskip
	
$\circ$ \ Each edge $\ef^{\<i\>}$ is associated to a linear transformation \.
$\bT^{\langle i \rangle}_{\vf}: \ts \Rb^{\ar} \to \Rb^{\ar}$.

\smallskip

\nin
Whenever clear, we drop the subscript~$\vf$ to avoid cluttering.
We call \. $\bM=(\aM_{ij})_{i,j \in [\ar]}$ \. the \defn{associated matrix} of $\vf$,
and \. $\hb = (\ah_i)_{i \in [\ar]}$  \. the \defn{associated vector} of~$\vf$.
In notation above, we have \ts $\vf^{\<i\>}\in \vfs$, for all \ts $1\le i \le \ar$.

}
\end{definition}

\smallskip

\subsection{Local-global principle}\label{ss:atlas-local}
As in the introduction (see $\S$\ref{ss:intro-proof}), matrix $\bM$ is called \defn{hyperbolic}, if
	\begin{equation}\label{eq:Hyp}
\langle \vb, \bM \wb \rangle^2 \  \geq \  \langle \vb, \bM \vb \rangle  \langle \wb, \bM \wb \rangle \quad \text{for every \ \, $\vb, \wb \in \Rb^{\ar}$, \ \  such that \ \ $\langle \wb, \bM \wb \rangle > 0$}. \tag{Hyp}
\end{equation}
For the atlas \ts $\AA$, we say that \ts $\vf\in \Vf$ \ts is \defn{hyperbolic},
if the associated matrix \ts $\bM_{\vf}$ \ts is hyperbolic, i.e.\ satisfies \eqref{eq:Hyp}.
We say that atlas \ts $\AA$ \ts satisfies \ts\defng{hyperbolic property} \ts if every
 \ts $\vf\in \Vf$ \ts is hyperbolic.

\smallskip

Note that property \eqref{eq:Hyp} depends only on the support of $\bM$,
i.e.\ it continues to hold after adding or removing zero rows or columns.
This simple observation will be used repeatedly through the paper.
%

\smallskip

We say that atlas \ts $\AA$ \ts satisfies \ts\defng{inheritance property} \ts if for
every non-sink vertex \ts $\vf\in \Vfp$, we have:
\begin{equation}\label{eq:Inh}\tag{Inh}
\begin{split}
		&	  (\bM \vb )_{i} \ = \ \big \langle   \bT^{\<i\>}\vb, \, \bM^{\<i\>} \,
\bT^{\<i\>}\hb \big \rangle \quad \text{ for every} \ \ \, i \in \supp(\bM) \text{ \ \ and \ \. }  \vb \in \Rb^{\ar},
\end{split}
\end{equation}
where \ts \ts $\bT^{\<i\>}=\bT^{\<i\>}_{\vf}$\., \ts $\hb=\hb_{\vf}$ \. and
\. $\bM^{\<i\>}:=\bM_{{\vf}^{\<i\>}}$ \. is the matrix associated with~$\vf^{\<i\>}$\..

\smallskip

Similarly, we say that atlas \ts $\AA$ \ts satisfies \ts\defng{pullback property} \ts if for
every non-sink vertex \ts $\vf\in \Vfp$, we have:
\begin{equation}\label{eq:Pull}\tag{Pull}
		\sum_{i\ts\in\ts \supp(\bM)} \. \ah_{i} \, \big \langle   \bT^{\<i\>}\vb, \, \bM^{\<i\>} \,  \bT^{\<i\>}\vb  \big \rangle  \ \geq \ \langle \vb, \bM \vb \rangle \qquad \text{ for every } \vb \in \Rb^{\ar}.
\end{equation}

\smallskip

We say that a non-sink vertex \ts $\vf\in \Vfp$ \ts is \defn{regular} if the following positivity conditions are satisfied:
\begin{align}\label{eq:Irr}\tag{Irr}
& \text{The associated matrix \ts $\bM_{\vf}$ \ts restricted to its support is irreducible.} \\
\label{eq:hPos} \tag{h-Pos} &\text{The associated vector \ts $\hb_{\vf}$ \ts restricted to the support of \ts $\bM_{\vf}$ \ts is strictly positive.}
\end{align}
Note that a matrix is \defn{irreducible} if
if it is not similar via a permutation to a block upper triangular matrix that has more than one block of positive size.
\smallskip

We now present the first main result of this section, which is
a local-global principle for \eqref{eq:Hyp}.

\smallskip

\begin{thm}[{\rm \defna{local–global principle}}{}]
\label{t:Hyp}
	Let $\AA$ be a combinatorial atlas that satisfies properties \eqref{eq:Inh} and \eqref{eq:Pull}, and
	let \ts $\vf\in \Vfp$ \ts be a non-sink regular vertex of~$\Qf$.
	Suppose every out-neighbor of \ts $\vf$ \ts is hyperbolic.
	Then \ts $\vf$ \ts is also hyperbolic.
\end{thm}

\smallskip

Theorem~\ref{t:Hyp} reduces checking the property \eqref{eq:Hyp} to sink vertices \ts $v\in \Vfm$.
In our applications, the pullback property \eqref{eq:Pull} is more complicated condition to check
than the inheritance property \eqref{eq:Inh}.  In the next Section~\ref{s:Pull}, we
present conditions implying \eqref{eq:Pull} that are easier to check.

\smallskip

\subsection{Eigenvalue interpretation of hyperbolicity}
The following lemma that gives an equivalent condition to~\eqref{eq:Hyp}
that is often easier to check.  A symmetric matrix $\bM$ satisfies \eqref{eq:OPE} if
\begin{equation}\label{eq:OPE}\tag{OPE}
	\text{$\bM$ \. has at most \. \defng{one positive eigenvalue} \. (counting multiplicity).}
\end{equation}
The equivalence between \eqref{eq:Hyp} and \eqref{eq:OPE}
is well-known in the literature, see e.g.,  \cite{Gre}, \cite[Thm~5.3]{COSW},
\cite[Lem.~2.9]{SvH19} and \cite[Lem.~2.5]{BH}.
We present a short proof for completeness.

\smallskip

	\begin{lemma}\label{l:Hyp is OPE}
	Let \ts {\rm $\bM$} \ts be a self-adjoint  operator on $\Rb^{\ar}$ for an inner product $\langle \cdot, \cdot \rangle$.
	Then  \ts {\rm $\bM$} \ts satisfies \eqref{eq:Hyp} \ts \underline{if and only if}  \. {\rm $\bM$} \ts satisfies \eqref{eq:OPE}.
\end{lemma}

\smallskip

\begin{proof}
For the \ts \eqref{eq:Hyp} \ts $\Rightarrow$ \ts \eqref{eq:OPE} \ts direction,
suppose to the contrary that $\bM$ has eigenvalues \ts $\lambda_1, \lambda_2>0$ (not necessarily distinct).
Let $\vb$ and $\wb$ be  orthonormal eigenvectors of $\bM$  for $\lambda_1$ and $\lambda_2$, respectively.	
It then follows that
	\begin{align*}
		0  \ = \ \langle \vb, \bM \wb \rangle \quad \text{ and } \quad \langle \vb, \bM \vb \rangle \, \langle \wb, \bM \wb \rangle \ = \  \lambda_1  \lambda_2\.,
	\end{align*}
	which contradicts \eqref{eq:Hyp}.
	
For the \ts  \eqref{eq:OPE} \ts $\Rightarrow$ \ts \eqref{eq:Hyp} \ts direction, 	
let \ts $\vb, \wb \in \Rb^{\ar}$ \ts be such that \. $\langle \wb, \bM \wb \rangle > 0$.
	Let $\lambda$ be the largest eigenvalue of $\bM$, and let $\hb$ be a corresponding eigenvector.
	Since 	\. $\langle \wb, \bM \wb \rangle > 0$, this implies that  $\lambda$ is a positive eigenvalue.
	Since $\bM$ has at most one positive eigenvalue (counting multiplicity), it follows that   $\lambda$ is the unique positive eigenvalue of $\bM$, and is a simple eigenvalue.
	In particular, this implies that
		\begin{equation*}
		\langle \wb, \bM \hb \rangle  \, \neq  \, 0,
	\end{equation*}
		as otherwise, we would have $\langle \wb, \bM \wb \rangle \leq 0$.
Let \ts $\zb \in \Rb^{\ar}$ \ts be the vector
	\[ \zb \ := \  \vb  \ - \ \frac{\langle \vb, \bM \hb \rangle}{\langle \wb, \bM \hb \rangle} \, \wb.
\]
It follows that $\langle \zb, \bM \hb \rangle   =0 $.	
	Since $\lambda$ is the only positive eigenvalue of $\bM$,
	we then have
	\begin{equation}\label{eq:Juliet 1}
		\langle \zb, \bM \zb \rangle \ \leq \ 0.
	\end{equation}
	On the other hand, we have
\begin{align*}
\langle \zb, \bM \zb \rangle \ & =  \ \langle \vb, \bM \vb \rangle  \ - \ 2 \. \frac{\langle \vb, \bM \hb \rangle \, \langle \vb, \bM \wb \rangle}{\langle \wb, \bM \hb \rangle}  \ + \  \frac{\langle \vb, \bM \hb \rangle^2 \, \langle \wb, \bM \wb \rangle}{\langle \wb, \bM \hb \rangle^2} \\
& \geq \ \langle \vb, \bM \vb \rangle \ - \  \frac{\langle \vb, \bM \wb \rangle^2}{\langle \wb, \bM \wb \rangle}\.,
\end{align*}	
where the last inequality is due to the AM--GM inequality.
Combining this inequality with \eqref{eq:Juliet 1}, we get
\[
    \langle \vb, \bM \wb \rangle^2 \  \geq \  \langle \vb, \bM \vb \rangle \  \langle \wb, \bM \wb \rangle\ts,
\]
	which proves \eqref{eq:Hyp}.
\end{proof}

\medskip

\subsection{Proof of Theorem~\ref{t:Hyp}}
%

%
Let \ts $\bM:=\bM_{\vf}$ \ts and \ts $\hb:=\hb_{\vf}$ \ts be the associated matrix and the associated vector of $\vf$, respectively.
	Since \eqref{eq:Hyp} is a property that is invariant under restricting to the support of $\bM$,
	it follows from \eqref{eq:Irr} that we can assume that $\bM$ is irreducible.
	
	Let  \ts $\bD:= (\aD_{ij})$  \ts be the \ts $\ar \times \ar$ \ts diagonal matrix given by
\[   \aD_{ii} \ := \   \frac{(\bM \hb)_{i}}{\ah_i} \qquad \text{ for every } \ 1\le i \le \ar\ts.
\]
	Note that $\bD$ is well defined and  	\. $\aD_{i i} >0$,  by \eqref{eq:hPos} and the assumption that $\bM$ is irreducible.
	Define a new inner product \ts $\langle \cdot, \cdot \rangle_{\bD}$ \ts on \ts $\Rb^{\ar}$ \ts by
	\.
	$ \langle \vb, \wb \rangle_{\bD} \ := \  \langle \vb, \bD \wb\rangle$\..
	
	Let \. $\bN := \bD^{-1} \bM$\ts.
	Note that 	\.
	$ \langle \vb, \bN \wb \rangle_{\bD}  =   \langle \vb, \bM \wb\rangle$ \. for every $\vb,\wb \in \Rb^{\ar}$.
	Since $\bM$ is a symmetric matrix,
	this implies that  $\bN$ is  a self-adjoint operator on $\Rb^{\ar}$ for the inner product $\langle \cdot, \cdot \rangle_{\bD}$\..
	A direct calculation shows that	 \ts $\hb$ \ts is an eigenvector of \ts $\bN$ \ts for eigenvalue \ts $\lambda=1$.
	Since $\bM$ is irreducible matrix and $\hb$ is a strictly positive vector,
	it then follows from the  Perron--Frobenius theorem that \ts $\lambda=1$ \ts is the
    largest real eigenvalue of $\bN$, and that it has multiplicity one.	
	
\medskip

\nin
\textbf{Claim:} \ {\em $\lambda=1$ \. is the only positive eigenvalue of \ts $\bNr$ \ts $($counting multiplicity$)$.}

\medskip

	By applying Lemma~\ref{l:Hyp is OPE} to the matrix $\bN$ and the inner product $\langle \cdot, \cdot \rangle_{\bD}$, it then follows that
	\[ \langle \vb, \bN \wb \rangle_{\bD}^2 \  \geq \  \langle \vb, \bN \vb \rangle_{\bD} \  \langle \wb, \bN \wb \rangle_{\bD} \quad   \text{ for every } \vb, \wb \in \Rb^{\ar}. \]
	Since 	\.
	$ \langle \vb, \bN \wb \rangle_{\bD}  =   \langle \vb, \bM \wb\rangle$, this implies \eqref{eq:Hyp} for~$\vf$, and completes the proof
of the theorem.
\qed

\medskip

\begin{proof}[Proof of the Claim]		
		Let $i \in [\ar]$
		 and $\vb \in \Rb^{\ar}$.		
		It follows from \eqref{eq:Inh} that
			\begin{align}\label{eq:general 0}
			\big((\bM \vb)_{i}\big)^2 \ = \  \big \langle    \vbw{i}, \, \bM^{\<i\>} \, \hbw{i}  \big \rangle^2.
		\end{align}
		Since \ts $\bM^{\<i\>}$ \ts satisfies \eqref{eq:Hyp} by the assumption of the theorem, applying \eqref{eq:Hyp} to the RHS of \eqref{eq:general 0} gives:
		\begin{align}\label{eq:general 1}
			\big((\bM \vb)_{i}\big)^2 \   \geq \  \big \langle   \vbw{i}, \, \bM^{\<i\>} \, \vbw{i}  \big \rangle  \,  \big \langle   \hbw{i}, \, \bM^{\<i\>} \, \hbw{i}  \big \rangle,
		\end{align}
Here \eqref{eq:Hyp} can be applied since
		\. $\big \langle   \hbw{i}, \, \bM^{\<i\>} \, \hbw{i}  \big \rangle \. = \. (\bM \hb)_{i} >0 $.
		Now note that
		\begin{align*}
			\big((\bN \vb)_{i}\big)^2 \, \aD_{i i} \ &= \  \big((\bM \vb)_{i}\big)^2  \, \frac{\ah_{i}}{(\bM \hb)_{i}} \ =_{\eqref{eq:Inh}}  \ \big((\bM \vb)_{i}\big)^2  \,  \frac{\ah_{i}}{\big \langle   \hbw{i}, \, \bM^{\<i\>} \, \hbw{i}  \big \rangle } \\
			& \geq_{\eqref{eq:general 1}} \ \ah_{i} \, \big \langle   \vbw{i}, \, \bM^{\<i\>} \, \vbw{i}  \big \rangle.
		\end{align*}
		Summing this inequality over all \ts $i \in [\ar]$, gives:
		\begin{align}\label{eq:general 2}
			\langle \bN \vb, \bN \vb \rangle_{\bD} \ \geq \  \sum_{i=1}^r  \ah_{i} \,   \big \langle   \vbw{i}, \, \bM^{\<i\>} \, \vbw{i}  \big \rangle \ \geq_{\eqref{eq:Pull}}
			\ \langle \vb, \bM \vb \rangle  \ = \  \langle  \vb, \bN \vb \rangle_{\bD}\..
		\end{align}
		Now, let $\lambda$ be an arbitrary eigenvalue of $\bN$,
		and let
		$\gb$ be an eigenvector of  $\lambda$.
		We have:
		\begin{align*}
			\lambda^2  \langle  \gb,  \gb \rangle_{\bD} \ = \  	\langle \bN \gb, \bN \gb \rangle_{\bD}
			\ \geq_{\eqref{eq:general 2}} \ 	 \langle  \gb, \bN \gb \rangle_{\bD}  \ = \ \lambda \. \langle  \gb,  \gb \rangle_{\bD}\..	\end{align*}
		This implies that $\lambda\geq 1$ or $\lambda \leq 0$.
		Since $\lambda =1$ is the largest eigenvalue of \ts $\bN$ \ts and is simple, we obtain the result.
	\end{proof}

\smallskip

\begin{rem}\label{r:Hyp-claim}
In the proof above, neither the Claim nor the proof of the Claim are new,
but a minor revision of Theorem~5.2 in~\cite{SvH19}.  We include the proof
for completeness and to help the reader get through our somewhat cumbersome
notation.  \end{rem}

\bigskip

\section{Pullback property}\label{s:Pull}

In this section we present sufficient conditions for \eqref{eq:Pull} that
are easier to verify, together with a construction of the maps \ts $\bT^{\<i\>}$.

\smallskip

\subsection{Three new properties}\label{ss:Pull-three}
Let $\AA$ be a combinatorial atlas.  
We say that $\AA$ satisfies the \defng{projective property}, if for every non-sink vertex \ts
$\vf\in \Vfp$ \ts and every \ts $i \in \supp(\bM)$, we have:
\begin{equation}\label{eq:Proj}\tag{Proj}
	\begin{split}
		& \big(\vbw{i} \big)_j    \ = \
		\begin{cases}
			\. \vb_{j} & \text{ if } \ \, j \in \supp\big(\bM^{\<i\>}\big) \cap \supp(\bM), \\
			\. \vb_{i} & \text{ if } \ \, j \in  \supp\big(\bM^{\<i\>}\big) \setminus \supp(\bM).
		\end{cases}
	\end{split}
\end{equation}	
We say that $\AA$ satisfies the \defng{transposition-invariant property},  if for every non-sink vertex \ts
$\vf\in \Vfp$, we have:
\begin{equation}\label{eq:TPInv}\tag{T-Inv}
	\aM_{jk}^{\<i\>}	\ = \  \aM_{ki}^{\<j\>}	 \ = \ \aM_{ij}^{\<k\>} \qquad 	\text{ for every distinct $i, j, k \in \supp(\bM)$.}
\end{equation}

\smallskip

Now, let \ts $\vf\in \Vfp$ \ts be a non-sink vertex of~$\Qf$, and let \ts $i \in \supp(\bM)$.
We partition the support of  matrix \ts $\bM^{\<i\>}$ \ts associated with vertex \ts $\vf^{\<i\>}$,
into two parts:
\begin{equation}\label{eq:Fam-def}
\Aunt^{\<i\>} \ := \   \supp\big(\bM^{\<i\>}\big)  \, \cap \,  \big(\supp(\bM) \ts - \ts i\big), \qquad \Fam^{\<i\>} \ := \ \supp\big(\bM^{\<i\>}\big)   \setminus  \big(\supp(\bM) \ts - \ts i\big).
\end{equation}
In other words, \ts $\Aunt^{\<i\>}$ \ts consists of elements  in the support of $\bM$
that do not include~$i$,\footnote{The name ``aunt'' here is referring to the siblings of the parent.}
while \ts $\Fam^{\<i\>}$ \ts consists of $i$ together with elements that initially are not in the support of $\bM$, but is then included in the support of $\bM^{\<i\>}$.\footnote{The name ``family'' here is referring to both the parents and their children.}
%
For every distinct \ts $i,j \in \supp(\bM)$, let
\begin{equation}\label{eq:Kdef}
	 \aK_{ij} \ := \  	\ah_{j} \, \aM_{jj}^{\<i\>}   \ - \  \ah_{j} \sum_{k \ts\in \ts\Fam^{\<j\>} } \aM_{ik}^{\<j\>}.
\end{equation}
Let us emphasize that \ts $\Aunt^{\<i\>}$, \ts $\Fam^{\<i\>}$, and
\ts $\aK_{ij}$ \ts all depend on non-sink vertex \ts $\vf$ \ts of~$\Qf$, even though
\ts $\vf$ \ts does not appear in these notation.

\smallskip
We say that $\AA$ satisfies the \defng{$\aKr$-nonnegative property},
if for every non-sink vertex \ts $\vf\in \Vfp$,
\begin{equation}\label{eq:Knonneg}\tag{K-Non}
	\aK_{ij} \  \geq \ 0 \quad \text{ for every distinct } \ \, i,j \. \in \.\supp(\bM).
\end{equation}

\smallskip

\nin
The main result of this subsection is the following sufficient condition for \eqref{eq:Pull}.

\smallskip	

\begin{thm}\label{t:Pull}
		Let \ts $\AA$ \ts be a combinatorial atlas that satisfies \eqref{eq:Inh}, \eqref{eq:Proj}, \eqref{eq:TPInv} and \eqref{eq:Knonneg}.
		Then  \ts $\AA$ \ts also satisfies  \eqref{eq:Pull}.
\end{thm}

\medskip

\subsection{Symmetry lemma}\label{ss:Pull-lemma}
To prove Theorem~\ref{t:Pull}, we need the following:

\smallskip

\begin{lemma}\label{l:Knonneg symmetric}
	Let \ts $\AA$ \ts be a combinatorial atlas that satisfies  \eqref{eq:Inh}, \eqref{eq:Proj}, and \eqref{eq:TPInv}.
	Then, for every non-sink vertex \ts $\vf\in \Vfp$, we have:
	\[ \aKr_{ij} \ = \ \aKr_{ji} \quad \text{ for every distinct } \ \, i,j \. \in \. \supp(\bMr).    \]
\end{lemma}

\smallskip

\begin{proof}
	Let \. $\eb_1,\ldots, \eb_{\ar}$ \. be the standard basis for \ts $\Rb^{\ar}$.
	It follows from \eqref{eq:Inh} that:
	\begin{align*}
		\aM_{ij} \ &= \ \big(\bM \eb_{j}\big)_{i} \ =_{\eqref{eq:Inh}} \  \big\langle \bT^{\<i\>}\eb_{j}, \bM^{\< i \>} \, \hbw{i}\big\rangle
		\ = \ \sum_{k =1}^{\ar} \. \aM_{jk}^{\<i\>} \, \big(\hbw{i}\big)_k \\
		\ &= \  \sum_{k  \ts\in\ts \Fam^{\<i\>}} \aM_{jk}^{\<i\>} \, \big(\hbw{i}\big)_k \ + \ \sum_{k \ts\in\ts \Aunt^{\<i\>} } \aM_{jk}^{\<i\>} \, \big(\hbw{i}\big)_k\\
		\ &= \ \aM_{jj}^{\<i\>} \, \big(\hbw{i}\big)_j  \ + \ \sum_{k  \ts\in\ts \Fam^{\<i\>}} \aM_{jk}^{\<i\>} \, \big(\hbw{i}\big)_k \ + \ \sum_{k \ts\in\ts \supp(\bM)\setminus \{i,j\}} \aM_{jk}^{\<i\>} \, \big(\hbw{i}\big)_k\..
	\end{align*}

Applying \eqref{eq:Proj} to the equation above, we get:
\begin{align}\label{eq:Ksym 1}
	\aM_{ij} \ &= \ \aM_{jj}^{\<i\>} \, \ah_j  \ \ + \ \sum_{k \ts\in\ts \Fam^{\<i\>}} \aM_{jk}^{\<i\>} \, \ah_i \ \ + \ \sum_{k \ts\in\ts \supp(\bM)\setminus \{i,j\}} \aM_{jk}^{\<i\>} \, \ah_k\..
\end{align}
By the same reasoning, we also get:
\begin{align}\label{eq:Ksym 2}
	\aM_{ji} \ &= \ \aM_{ii}^{\<j\>} \, \ah_i  \ \ + \ \sum_{k  \ts\in\ts \Fam^{\<j\>}} \aM_{ik}^{\<j\>} \, \ah_j \ \ + \ \sum_{k \ts\in\ts \supp(\bM)\setminus \{i,j\}} \aM_{ik}^{\<j\>} \, \ah_k\..
\end{align}
By \eqref{eq:TPInv}, the rightmost sums in \eqref{eq:Ksym 1} and \eqref{eq:Ksym 2} are equal.
On the other hand, the left side of \eqref{eq:Ksym 1} and \eqref{eq:Ksym 2} are equal since $\bM$ is a symmetric matrix.
Equating \eqref{eq:Ksym 1} and \eqref{eq:Ksym 2}, we obtain:
\begin{align*}
	\aM_{jj}^{\<i\>} \, \ah_j  \ + \ \sum_{k  \ts\in\ts \Fam^{\<i\>}} \aM_{jk}^{\<i\>} \, \ah_i  \ = \   \aM_{ii}^{\<j\>} \, \ah_i  \ + \ \sum_{k  \ts\in\ts \Fam^{\<j\>}} \aM_{ik}^{\<j\>} \, \ah_j\.,
\end{align*}
which is equivalent to
\begin{align*}
	\aM_{jj}^{\<i\>} \, \ah_j  \ - \
	\sum_{k  \ts\in\ts \Fam^{\<j\>}} \aM_{ik}^{\<j\>} \, \ah_j \ = \
	\aM_{ii}^{\<j\>} \, \ah_i  \ - \  \sum_{k \ts\in\ts \Fam^{\<i\>}} \aM_{jk}^{\<i\>} \, \ah_i\..
\end{align*}
The lemma now follows by noting that the LHS of the equation above is equal to \ts
$\aK_{ij}$, while the RHS  is equal to \ts $\aK_{ji}$\ts.
\end{proof}

\subsection{Proof of Theorem~\ref{t:Pull}}\label{ss:Pull-proof}
	Let $\vf$ be a non-sink vertex of $\Qf$,
	and let $\vb \in \Rb^{\ar}$.
		The left side of \eqref{eq:Pull} is equal to
		\begin{align}\label{eq:families left}
			\sum_{i \ts\in\ts \supp(\bM)} \ah_{i} \, \big\langle \vbw{i}, \bM^{\<i\>} \vbw{i} \big\rangle  \
= \ \sum_{i \ts\in\ts \supp(\bM)} \ \sum_{j,k \ts\in\ts \supp(\bM^{\<i\>}) } \ah_{i} \, \big(\vbw{i}\big)_j \, \big(\vbw{i}\big)_k \, \aM^{\<i\>}_{jk}\..
		\end{align}

\nin
\underline{First}, this sum can be partitioned into the sum over the following five families:
		\begin{enumerate}
			[{label=\textnormal{({\arabic*})},
				ref=\textnormal{\arabic*}}]
			\item \label{typeDiv 1}
			The  triples $(i,j,k)$,
			where $i \in \supp(\bM)$, and
			$j,k \in \Aunt^{\<i\>}$ are  distinct.
			By \eqref{eq:Proj}, the term in \eqref{eq:families left} is  equal to
			\[ \ah_i \, \av_j  \, \av_k  \, \aM^{\<i\>}_{jk}\..\]
			
			\item \label{typeDiv 2}
						The  triples $(i,j,k)$,
			where $i \in \supp(\bM)$, and
			$j,k \in \Fam^{\<i\>}$ (not necessarily  distinct).
			By \eqref{eq:Proj}, the term in \eqref{eq:families left} is  equal to
			\[ \ah_i \, \av_i^2  \, \aM^{\<i\>}_{jk}\..\]
			
			\item \label{typeDiv 3}
						The  triples $(i,j,k)$,
			where $i \in \supp(\bM)$, \ts
			$j\in \Aunt^{\<i\>}$, and $k \in \Fam^{\<i\>}$.
			By \eqref{eq:Proj}, the term in \eqref{eq:families left} is  equal to
			\[ \ah_i \, \av_i  \, \av_j  \, \aM^{\<i\>}_{jk}\..\]
			\item \label{typeDiv 4}
						The  triples $(i,j,k)$,
			where $i \in \supp(\bM)$, \ts $j \in \Fam^{\<i\>}$, and
			$k \in \Aunt^{\<i\>}$.
			By \eqref{eq:Proj}, the term in \eqref{eq:families left} is  equal to
			\[ \ah_i \, \av_i \, \av_k  \, \aM^{\<i\>}_{jk}\..\]
			
			\item \label{typeDiv 5}
						The  triples $(i,j,k)$,
			where $i \in \supp(\bM)$, and
			$j=k \in \Aunt^{\<i\>}$.
			By \eqref{eq:Proj}, the term in \eqref{eq:families left} is  equal to
			\[ \ah_i \, \av_j^2  \, \aM^{\<i\>}_{jj} \ = \  \frac{\ah_i}{\ah_j} \, \av_j^2 \, \aK_{ij} \ + \    \sum_{k \ts\in\ts \Fam^{\<j\>}}\ah_i \, \av_j^2 \,  \aM^{\<j\>}_{ik}\..\]
			Thus the sum over  this family can be partitioned further into the sum over the following two families:
			\begin{enumerate}
				[{label=\textnormal{(5{\alph*})},
					ref=\textnormal{5\alph*}}]
				\item \label{typeDiv a}	
				The pair $(i,j)$, where $i,j \in \supp(\bM)$ are distinct, with the term
				\[   \frac{\ah_i}{\ah_j} \, \av_j^2 \,  \aK_{ij}\.. \]
				
				\item  \label{typeDiv b}
															The  triples $(i,j,k)$,
									where $i,j  \in \supp(\bM)$ are distinct, \ts
									 and $k \in \Fam^{\<j\>}$,
									with the term
									\[ \ah_i \, \av_j^2 \,  \aM^{\<j\>}_{ik}\.. \]
			\end{enumerate}
		\end{enumerate}
		
		\medskip

\nin
\underline{Second},
		the right side of \eqref{eq:Pull} is equal to
		\begin{align}
			\langle \vb, \bM \vb \rangle  \ &= \  \sum_{i' \ts\in\ts \supp(\bM)} \av_{i'} \, (\bM \vb)_{i'}  \ =_{\eqref{eq:Inh}} \  \sum_{i' \ts\in\ts \supp(\bM)} \av_{i'} \, \big \langle   \vbw{i'}, \, \bM^{\<i'\>} \,   \hbw{i'}  \big \rangle  \notag \\
		\label{eq:families right}	\ &= \ \sum_{i' \ts\in\ts \supp(\bM)} \ \sum_{j', \ts k' \ts\in\ts\supp(\bM^{\<i'\>})}  \av_{i'} \, \big(\vbw{i'}\big)_{j'} \, \big(\hbw{i'}\big)_{k'} \, \aM^{\<i'\>}_{j'k'}.
		\end{align}
		This sum can be partitioned into the sum over the following five families:
		\begin{enumerate}
			[{label=\textnormal{({\arabic*$^\prime$})},
				ref=\textnormal{\arabic*$^\prime$}}]
			\item \label{typeDiv 1'}
						The  triples $(i',j',k')$,
			where $i' \in \supp(\bM)$, and
			$j',k' \in \Aunt^{\<i'\>}$ are  distinct.
			By \eqref{eq:Proj}, the term in \eqref{eq:families right} is  equal to
			\[ \ah_{k'} \, \av_{i'} \, \av_{j'}  \, \aM^{\<i'\>}_{j'k'}\..\]

			\item \label{typeDiv 2'}
									The  triples $(i',j',k')$,
			where $i' \in \supp(\bM)$, and
			$j',k' \in \Fam^{\<i'\>}$ (not necessarily distinct).
			By \eqref{eq:Proj}, the term in \eqref{eq:families right} is  equal to
			\[ \ah_{i'} \, \av_{i'}^2   \, \aM^{\<i'\>}_{j'k'}\..\]

			\item \label{typeDiv 3'}
									The  triples $(i',j',k')$,
			where $i' \in \supp(\bM)$, \ts
			$j' \in \Aunt^{\<i'\>}$, and $k' \in \Fam^{\<i'\>}$.
			By \eqref{eq:Proj}, the term in \eqref{eq:families right} is  equal to
			\[ \ah_{i'} \, \av_{i'} \, \av_{j'}  \, \aM^{\<i'\>}_{j'k'}\..\]
			
			\item \label{typeDiv 4'} 	
			The  triples $(i',j',k')$,
			where $i' \in \supp(\bM)$, \ts
			$j' \in \Fam^{\<i'\>}$, and $k' \in \Aunt^{\<i'\>}$.
			By \eqref{eq:Proj}, the term in \eqref{eq:families right} is  equal to
			\[ \ah_{k'} \, \av_{i'}^2   \, \aM^{\<i'\>}_{j'k'}\..\]
			
			\item \label{typeDiv 5'}
									The  triples $(i',j',k')$,
			where $i' \in \supp(\bM)$, and
			$j'=k' \in \Aunt^{\<i'\>}$.
			By \eqref{eq:Proj}, the term in \eqref{eq:families right} is  equal to
			\[ \ah_{j'} \, \av_{i'} \, \av_{j'}  \, \aM^{\<i'\>}_{j'j'} \ = \    \av_{i'} \, \av_{j'} \, \aK_{i'j'} \ + \    \sum_{k' \ts\in\ts \Fam^{\<j'\>}}\ah_{j'} \, \av_{i'} \, \av_{j'} \,  \aM^{\<j'\>}_{i'k'}\..\]
			Thus the sum over  this family can be partitioned further into the sum over the following two families:
						\begin{enumerate}
			[{label=\textnormal{(5{\alph*$^\prime$})},
	ref=\textnormal{5\alph*$^\prime$}}]
				\item \label{typeDiv a'}	
				The pair $(i',j')$, where $i',j' \in \supp(\bM)$ are distinct, with the term
				\[    \av_{i'} \, \av_{j'} \,  \aK_{i'j'}\.. \]
				
				\item  \label{typeDiv b'}
				The  triples $(i',j',k')$,
				where $i',j' \in \supp(\bM)$ are distinct, and $k' \in \Fam^{\<j'\>}$,
				with the term
				\[ \ah_{j'} \, \av_{i'} \, \av_{j'} \,  \aM^{\<j'\>}_{i'k'}\.. \]
			\end{enumerate}

		\end{enumerate}
		
		\medskip
		
\nin
\underline{Third}, we show that the RHS of \eqref{eq:families left} is at least as large as the RHS of \eqref{eq:families right}.
		We have the following six cases:
		\begin{enumerate}
			\item 		The term in \eqref{typeDiv 1} is equal to that of \eqref{typeDiv 1'} by substituting
			\. $i' \gets j$, $j' \gets k$, $k' \gets i$ (counterclockwise substitution)  to \eqref{typeDiv 1}:
			\[  \ah_i \, \av_j  \, \av_k  \, \aM^{\<i\>}_{jk} \ =_{\eqref{eq:TPInv}} \  \ah_i \, \av_j  \, \av_k  \, \aM^{\<j\>}_{ki}  \ = \  \ah_{k'} \, \av_{i'} \, \av_{j'}  \, \aM^{\<i'\>}_{j'k'}\..  \]
			
		\item 		
		The term in \eqref{typeDiv 2} is equal to that of \eqref{typeDiv 2'}
		by
		substituting \. $i' \gets i$, $j' \gets j$, $k' \gets k$ (identity substitution) to \eqref{typeDiv 2}:
		\[ \ah_i \, \av_i^2  \, \aM^{\<i\>}_{jk} \ = \ \ah_{i'} \, \av_{i'}^2   \, \aM^{\<i'\>}_{j'k'}\..\]
		
		\item 	
		The term in \eqref{typeDiv 3} is equal to that of \eqref{typeDiv 3'}
		by
		substituting \. $i' \gets i$, $j' \gets j$, $k' \gets k$ (identity substitution) to \eqref{typeDiv 3}:
			\[ \ah_i \, \av_i  \, \av_j  \, \aM^{\<i\>}_{jk} \ = \ \ah_{i'} \, \av_{i'} \, \av_{j'}  \, \aM^{\<i'\>}_{j'k'}\..\]
			
		\item The term in \eqref{typeDiv 4} is equal to that of \eqref{typeDiv b'}
		by
				substituting  \. $i' \gets k$, $j' \gets i$, $k' \gets j$  (clockwise substitution) to~\eqref{typeDiv 4}:
							\[ \ah_i \, \av_i \, \av_k  \, \aM^{\<i\>}_{jk} \ = \  \ah_{j'} \, \av_{i'} \, \av_{j'} \,  \aM^{\<j'\>}_{i'k'}\..\]
		\item 		The term in \eqref{typeDiv a} is equal to that of \eqref{typeDiv a'}
		by substituting  \. $i' \gets i$, $j' \gets j$ (identity substitution) to~\eqref{typeDiv a}:
		\begin{align*}
		\frac{\ah_i}{\ah_j} \, \av_j^2 \,  \aK_{ij}  \ + \
		\frac{\ah_j}{\ah_i} \, \av_i^2 \,  \aK_{ji}
	    \ &\geq \  		2 \, \av_i \, \av_j \, \sqrt{\aK_{ij} \, \aK_{ji}}
			\ =_{\text{Lem}~\ref{l:Knonneg symmetric}} \
				\av_i \, \av_j  \, \aK_{ij} \ + \  \av_j \, \av_i  \, \aK_{ji} \\
		\ &= \  \av_{i'} \, \av_{j'} \,  \aK_{i'j'} \ + \  \av_{j'} \, \av_{i'} \,  \aK_{j'i'}\.,
		\end{align*}
	where the first inequality follows from \eqref{eq:Knonneg} and the AM-GM inequality.\footnote{Note that this is only instance of inequality in this proof.}
	
			\item The term in \eqref{typeDiv b} is equal to that of \eqref{typeDiv 4'}
	by
	substituting  \. $i' \gets j$, $j' \gets k$, $k' \gets i$  (clockwise substitution) to \eqref{typeDiv b}:
						\[ \ah_i \, \av_j^2 \,  \aM^{\<j\>}_{ik} \ = \ \ah_{k'} \, \av_{i'}^2   \, \aM^{\<i'\>}_{j'k'}\..\]
		\end{enumerate}		
		This completes the proof of the theorem. \qed
	
\medskip
	
\begin{rem}\label{r:Pull-K-PSD}
The condition \eqref{eq:Knonneg} in Theorem~\ref{t:Pull} can be weakened as follows.
Let \ts $\vf\in \Vfp$ \ts be a non-sink vertex, and let \. $\bK :=(\aK_{ij})_{i,j \in \supp(\bM)}$ be the matrix defined by
\[ \aK_{ij} \ := \
\begin{cases}
				\, \aK_{ij} \ \text{as in \eqref{eq:Kdef}} & \text{ if \, $i,j \ts \in \ts \supp(\bM)$, \, $i \ts\ne\ts j$},\\
				\, \displaystyle -\.\sum_{\ell \ts\in\ts \supp(\bM) \setminus \{i\}} \. \frac{\ah_i}{\ah_\ell} \, \aK_{i,\ell}
	& \text{ if \, $i=j\ts \in \ts \supp(\bM)$.}		
\end{cases}
\]
We claim that the condition \eqref{eq:Knonneg} in Theorem~\ref{t:Pull} can be replaced with
	\begin{equation}\label{eq:Kposdef}\tag{K-PSD}
		\text{The matrix \ $-\bK$ \ \. is positive semidefinite,}
	\end{equation}
for every non-sink vertex $\vf$ of $\Qf$.
This generalization follows from the same proof as Theorem~\ref{t:Pull} by a straightforward modification to step (v).
Note that in this paper we never apply this (slightly more general) version of Theorem~\ref{t:Pull},
as all interesting applications that we found satisfy the stronger condition \eqref{eq:Knonneg},
which is also easier to check.
\end{rem}

\bigskip

\section{Hyperbolic equality for combinatorial atlases}\label{s:atlas-equality}

In this section we characterize when the equality conditions in \eqref{eq:Hyp} hold
for all non-sink vertices in a combinatorial atlas.  For that, we obtain the
equality variation of the local-global principle (Theorem~\ref{t:Hyp}).
See~$\S$\ref{ss:finrem-main} for some background.

\smallskip

\subsection{Statement} \label{ss:atlas-equality-global}
Let \ts $\AA$ \ts be a combinatorial atlas of dimension~$\ar$.
Recall that, for a non-sink vertex \ts $\vf$ \ts of~$\ts\Qf$,
we denote by \ts $\bM=\bM_{\vf}$ \ts the associated matrix of~$\vf$,
by \ts $\hb=\hb_{\vf}$ \ts the associated vector of~$\vf$,
by \ts $\bT^{\<i\>}=\bT^{\<i\>}_{\vf}$ \ts the associated linear
transformation of the edge \ts $\ef^{\<i\>}=(\vf, \vf^{\<i\>})$, and
by \ts $\bM^{\<i\>}$ the associated matrix of the vertex~$\vf^{\<i\>}$.

\smallskip

A \defn{global pair} \ts $\fb,\gb \in \Rb^{\ar}$ \ts
is a pair of nonnegative vectors, such that
 \begin{equation}
 	\label{eq:PosGlob} \tag{Glob-Pos}  \text{$\fb \ts +\ts\gb$ \ is a strictly positive vector}\ts.
 \end{equation}
Here \ts $\fb$ \ts and \ts $\gb$ \ts are global in a sense that they are the same for all
vertices \ts $\vf\in \Vf$.

\smallskip

Fix a number \ts $\as>0$.
We say that a vertex \ts $\vf \in \Vf$ \ts satisfies~\eqref{eq:sEqu}, if
\begin{equation}\label{eq:sEqu}
	\tag{s-Equ}  \langle \fb , \bM \fb \rangle  \ = \  \as \, \langle \gb , \bM \fb \rangle \ = \  \as^2 \, \langle \gb , \bM \gb \rangle,
\end{equation}
where \ts $\bM=\bM_{\vf}$ \ts as above.
Observe that \eqref{eq:sEqu} implies that equality occurs in \eqref{eq:Hyp} for substitutions \ts $\vb\gets \gb$ \ts
and \ts $\wb\gets \fb$, since
\begin{equation}\label{eq:Equ}
  \langle \gb , \bM \fb \rangle^2  \ = \  \as \,\. \langle \gb , \bM \gb \rangle \ \as^{-1} \. \langle \fb , \bM \fb \rangle \ = \ \langle \gb , \bM \gb \rangle  \ \langle \fb , \bM \fb \rangle\..
\end{equation}
We say that the atlas \ts $\AA$ \ts satisfies \defng{$\asr$-equality property}
if~\eqref{eq:sEqu} holds for every \ts $\vf \in \Vf$.

\medskip

We now present the first main result of this section,
which is a \defng{local-global principle} for \eqref{eq:sEqu}.
A vertex \ts $\vf \in \Vfp$ \ts is called \defn{functional source} if
the following conditions are satisfied:
	\begin{align}
	\label{eq:ProjGlob} \tag{Glob-Proj}  &\af_j \ = \  \big(\bT^{\<i\>} \fb \big)_j \quad \text{ and }  \quad  \ag_j \ = \  \big(\bT^{\<i\>} \gb \big)_j  \quad \text{ for every } \ i \in \supp(\bM), \ j \in \supp(\bM^{\<i\>}), \\
			\label{eq:hGlob} \tag{h-Glob}  &\text{$\fb\.=\.\hb_{\vf}$\..} 
\end{align}	
Here condition~\eqref{eq:ProjGlob} means that \ts $\fb,\gb$ \ts are fixed points of the
projection \ts $\bT^{\<i\>}$ \ts when restricted to the support.

We say that an edge \ts $\ef^{\<i\>}=(\vf, \vf^{\<i\>})\in \Ef$ \ts is \defn{functional} if
$\vf$ is a functional source and \. $i \in \supp(\bM) \, \cap \, \supp(\hb)$.
A vertex $\wf \in \Vf$ is a \defn{functional target}
of~$\vf$, if there exists a directed path \. $\vf\to \wf$ \. in $\Qf$ consisting
of only functional edges.  Note that a functional target is not necessarily
a functional source.

\smallskip

\begin{thm}[{\rm \defna{local-global equality principle}}]
\label{t:equality Hyp}
Let \ts $\AA$ \ts be a combinatorial atlas that satisfies properties
\eqref{eq:Inh}, \eqref{eq:Pull}.  Suppose also \ts $\AA$ \ts
satisfies property \eqref{eq:Hyp} for every vertex \ts $v\in \Vf$.
Let \ts $\fb,\gb$ \ts be a global pair of~$\AA$.
Suppose a non-sink vertex \ts $\vf \in \Vfp$ \ts satisfies \eqref{eq:sEqu}
with constant $\asr>0$.  Then every functional target of $\vf$ also
satisfies \eqref{eq:sEqu} with the same constant~$\asr$.
\end{thm}

\medskip

\subsection{Algebraic lemma} \label{ss:atlas-equality-lemma}
We start with the following general algebraic result.  Recall
that a matrix is \emph{hyperbolic} if it satisfies~\eqref{eq:Hyp}.

\smallskip

\begin{lemma}\label{l:equality kernel}
	Let $\bMr$ be a nonnegative symmetric hyperbolic $\arr \times \arr$ matrix.
	Let \ts $\fb,\gb \in \Rb^{\arr}$ \ts be nonnegative vectors, let \ts $\asr>0$,
	and let \ts $\zb := \fb \ts - \ts \asr  \gb$.
	Then \eqref{eq:sEqu} holds \  \underline{if and only if}  \ $\bMr \zb =0$.	
\end{lemma}

\smallskip

\begin{proof}
	The $\Leftarrow$ direction follows from the fact that
	\begin{equation}\label{eq:Ramon 1}
	\begin{split}
		\langle \fb , \bM \fb \rangle  \ - \  \as \, \langle \gb , \bM \fb \rangle \ &= \  \langle  \zb, \bM  \fb \rangle \ = \ \langle  \bM \zb,  \fb \rangle, \quad \text{ and } \quad
		\as \, \langle \gb , \bM \fb \rangle \ - \   \as^2 \, \langle \gb , \bM \gb \rangle \ = \ \as \, \langle \gb, \bM \zb \rangle.
	\end{split}
	\end{equation}
	Thus it suffices to prove the  $\Rightarrow$ direction.
	We will assume  that $\bM$ is nonzero when restricted to the support
of \ts $\gb+\fb$, as otherwise every term in \eqref{eq:sEqu} is equal to 0 and  the lemma follows immediately.
Let \. $\wb :=\gb+\fb$, and  the assumption implies that
\. $  \langle \wb , \bM \wb  \rangle >0$.
By \eqref{eq:Hyp}, we then have that the matrix $\bM$ is negative semidefinite on \ts $(\bM \wb)^\perp$.
Now note that \ts $\zb \in (\bM \wb)^\perp$, since \ts $\langle \zb, \bM \wb\rangle =0$ \ts by \eqref{eq:Ramon 1} and \eqref{eq:sEqu}.
Also note that
\begin{equation}\label{eq:flipflop 1}
	 \langle \zb, \bM \zb\rangle \ = \  \langle \fb , \bM \fb \rangle  \ - \  2 \,  \as \, \langle \gb , \bM \fb \rangle \  +  \ \as^2 \, \langle \gb , \bM \gb \rangle \ =_{\eqref{eq:sEqu}} \ 0.
\end{equation}
 It then follows from these three observations that \ts $\bM \zb=0$, as desired.
\end{proof}

\medskip

\subsection{Proof of Theorem~\ref{t:equality Hyp}}
	By induction, it suffices to show that, for every functional edge \. $(\vf,\vf^{\<i\>})\in \vfs$,
	we have that \ts $\vf^{\<i\>}$ \ts satisfies \eqref{eq:sEqu} with the same constant~$\as>0$.
	
\smallskip

	It follows from \eqref{eq:Inh}, that for every \ts $i \in \supp(\bM)$ \ts we have:
	\begin{align*}
		\big(\bM \gb\big)_{i}  \ = \ \langle \bT^{\<i\>}\gb,  \bM^{\<i\>} \bT^{\<i\>}\hb\rangle
\qquad \text{and} \qquad \big(\bM \hb\big)_{i} \ = \  \langle \bT^{\<i\>} \hb,  \bM^{\<i\>} \, \bT^{\<i\>}\hb\rangle\..
	\end{align*}
	It then follows from \eqref{eq:ProjGlob} and the fact that \. $\fb=\hb =\hb_{\vf}$ \. by \eqref{eq:hGlob} \. that
	\begin{equation}\label{eq:golf 0.5}
			\big(\bM \gb\big)_{i}  \ = \ \langle \gb,  \bM^{\<i\>} \fb\rangle \qquad \text{and} \qquad  \big(\bM \fb\big)_{i} \ = \  \langle  \fb,  \bM^{\<i\>}\fb\rangle.
	\end{equation}

			Let \. $\zb \ts := \,  \fb \ts - \ts s \gb$.
			It then follows from \eqref{eq:sEqu} and \eqref{eq:flipflop 1} that
			\. $\langle \zb, \bM \zb \rangle =0$.
By Lemma~\ref{l:equality kernel}, \eqref{eq:sEqu} implies  that \. $\bM \zb = \0$, which is equivalent to
	\. $\as \,  \bM \gb  =   \bM \fb
	$.
	Together with \eqref{eq:golf 0.5},
	this implies that 	
	\begin{equation}\label{eq:golf 2}
		 \as \, \langle \gb,  \bM^{\<i\>} \fb\rangle \ =_{\eqref{eq:Inh}} \ \as \,  \big(\bM \gb\big)_{i}  \ = \  \big(\bM \fb\big)_{i} \ = \ \langle \fb,  \bM^{\<i\>} \fb\rangle.
	\end{equation}
	On the other hand, we have
	\[  \langle \fb,  \bM^{\<i\>} \fb\rangle \ =_{\eqref{eq:golf 0.5}} \  \big(\bM \fb \big)_{i}  \ =_{\eqref{eq:golf 2}} \  \frac{\as+1}{\as} \, \big(\bM (\fb+\gb) \big)_{i} \ > \, 0, \]
	where the positivity follows by \eqref{eq:PosGlob} and the assumption that \ts $i \in \supp(\bM)$.
	Now note that,
	\begin{equation}\label{eq:golf 3}
		\begin{split}
			& \langle \zb,  \bM^{\<i\>} \zb\rangle
			\ \quad = \    \as^2 \, \langle \gb, \bM^{\<i\>} \gb \rangle \ - \  2 \as \,  \langle \gb, \bM^{\<i\>} \fb \rangle
			\  + \  \langle \fb, \bM^{\<i\>} \fb \rangle \\
			& \quad =_{\eqref{eq:golf 2}} \  \as^2 \, \bigg( \langle \gb, \bM^{\<i\>} \gb \rangle \ - \ \frac{\langle \gb, \bM^{\<i\>} \fb \rangle^2 }{\langle \fb, \bM^{\<i\>} \fb \rangle } \bigg),
		\end{split}
	\end{equation}
	which is nonpositive as $\vf^{\<i\>}$ satisfies \eqref{eq:Hyp}.
	On the other hand, we have
	\begin{align*}
		& \sum_{i\ts\in \ts \supp(\bM)}  \ah_{i} \,  \langle \zb,  \bM^{\<i\>} \zb\rangle \ =_{\eqref{eq:ProjGlob}} \
		\sum_{i\ts\in \ts \supp(\bM)}  \ah_{i} \,  \langle \bT^{\<i\>}\zb,  \bM^{\<i\>} \bT^{\<i\>}\zb\rangle \ \geq_{\eqref{eq:Pull}} \  \langle \zb, \bM \zb \rangle  \ =_{\eqref{eq:flipflop 1}} \ 0.
	\end{align*}	
So the RHS of this inequality is equal to~0, while the LHS
 is a sum of nonpositive terms by~\eqref{eq:golf 3}.
	This implies that
	every term in the first sum is equal to~0, and thus
	 \. $\ah_{i} \langle \zb,  \bM^{\<i\>} \zb\rangle=0$ \. for every $i \in \supp(\bM)$.
	This in turn  implies that \. $\langle \zb,  \bM^{\<i\>} \zb\rangle=0$ \.
whenever \ts $(\vf,\vf^{\<i\>})$ \ts is a functional edge.
	This is equivalent to saying that the left side of \eqref{eq:golf 3}
	is zero, and we have:
	\begin{equation}\label{eq:golf 4}
	  \langle \gb, \bM^{\<i\>} \gb \rangle \ = \ \frac{\langle \gb, \bM^{\<i\>} \fb \rangle^2 }{\langle \fb, \bM^{\<i\>} \fb \rangle } \ =_{\eqref{eq:golf 2}} \ \frac{1}{\as} \, \langle \gb, \bM^{\<i\>} \fb \rangle.
	\end{equation}
	It then follows from \eqref{eq:golf 2} and \eqref{eq:golf 4} that $\vf^{\<i\>}$ satisfies \eqref{eq:sEqu} whenever $(\vf, \vf^{\<i\>})$ is a functional edge, which completes the proof. \qed

\bigskip

\section{Log-concave inequalities for interval greedoids}
\label{s:proof-greedoid}

In this section, we prove Theorem~\ref{t:greedoid} by constructing
a combinatorial atlas corresponding to a greedoid,
and applying both local-global principle in Theorem~\ref{t:Hyp} and sufficient
conditions for hyperbolicity given in Theorem~\ref{t:Pull}.

\smallskip

\subsection{Combinatorial atlas for interval greedoids}
\label{ss:greedoid-def}
Let \ts $\Gf=(X,\cL)$ \ts be an interval greedoid on \ts $|X|=n$ \ts elements,
let \ts $1\le k < \rk(\Gf)$,  and let \ts $\aq:\Lc \to \rrs$ \ts
be the weight function in Theorem~\ref{t:greedoid}.
We define a combinatorial atlas \ts $\AA$ \ts corresponding
to \ts $(\Gf,k,\aq)$ \ts as follows.

\smallskip

Define an acyclic graph \ts $\Qf:=(\Vf, \Ef)$, where the set of vertices \.  $\Vf \ts := \Vfm\cup \Vf^1 \cup \ldots \cup \Vf^{k-1}$ \.
is given by\footnote{Yes, graph~$\Qf$ has uncountably many vertices.}
\begin{align*}
	 \Vf^m \, & := \ \big\{ \ts (\alpha,m,t) \ \mid \  \alpha \in X^* \text{ with }  |\alpha| \leq k-1-m, \ t \in [0,1]  \ts \big\} \ \  \text{ for } \ m\geq 1, \\
	 \Vfm  \ & := \ \big\{ \ts (\alpha,0,1) \ \mid \  \alpha \in X^* \text{ with }  |\alpha| \leq k-1  \ts \big\}.
\end{align*}
Here the restriction \ts $t=1$ \ts in \ts $\Vfm$ \ts is crucial for a technical reason that will be apparent later in the section.

\begin{figure}[hbt]
		\includegraphics[width=8.1cm]{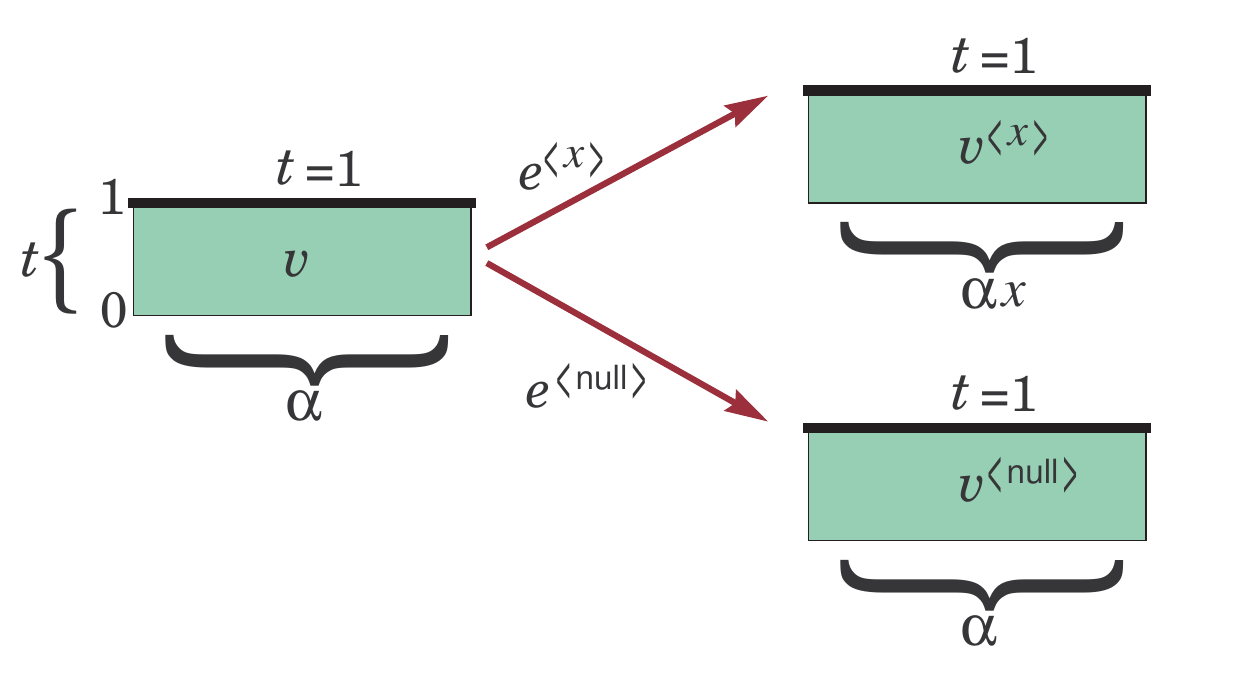}
		\vskip-.25cm
		\caption{Edges of two type: \. $\ef^{\<x\>} = \bigl(\vf,\vf^{\<x\>}\bigr)$, \. $v=(\al,m,t)$, \. $\vf^{\<x\>}=(\al x,m-1,1)$, and \.
$\ef^{\<\spstar\>} = \bigl(\vf,\vf^{\<\spstar\>}\bigr)$, \. $v=(\al,m,t)$, \. $\vf^{\<\spstar\>}=(\al,m-1,1)$. }
		\label{f:matroid-edge}
\end{figure}

Let \. $\Xf:=X\cup \{\spstar\}$ \ts be the set of letters~$X$
with one special element \ts $\spstar$ \ts added.  The reader should think of
element \ts $\spstar$ \ts as the empty letter.  Let \ts $\ar:= |\Xf|=(n+1)$ \ts
be the dimension of the atlas, so each vertex \ts $\vf\in \Vf^m$, \ts $m\geq 1$,  has
exactly \ts $(n+1)$ \ts outgoing edges we label \ts $\bigl(\vf,\vf^{\<x\>}\bigr)\in \Ef$, where
\ts $x\in \Xf$ \ts and \ts $\vf^{\<x\>}\in \Vf^{m-1}$ \ts is defined as follows:
\[ \
\vf^{\<x\>} \ := \
 \begin{cases}
 	(\alpha x,m-1,1) & \text{ if } \ \, x \in  X,\\
 	(\alpha,m-1,1) & \text{ if } \ \, x = \spstar\ts.
 \end{cases}
 \]
Let us emphasize that this is not a typo and we indeed have the last parameter \ts $t=1$,
for all \ts $\vf^{\<x\>}$ (see Figure~\ref{f:matroid-edge}).


\medskip

For every $\alpha \in X^*$ and every \ts $m \in \{1,\ldots, \rk(\Gf)-|\alpha|-1\}$,
we denote by \.
 $\bA(\alpha,m) := (\aA_{x\ts y})_{x,y \in \Xf}$ \. the symmetric \ts $\ar \times \ar$ \ts
 matrix defined as follows:\footnote{When $m=1$ and $\alpha\in \cL$, we have $\Cnt_{m-1}(\alpha)$ consists of exactly one element, namely the empty word.}
$$
\aligned
& \aA_{x\ts y} \, := \ 0 \qquad  \text{for \ \ $x \notin \Cnt(\alpha) \ts + \ts \spstar$ \ \
or \ \  $y \notin \Cnt(\alpha) \ts + \ts \spstar$,} \\
& \aA_{x\ts y} \, := \,  \sum_{\beta \ts\in\ts \Cnt_{m-1}(\alpha xy)}\, \aq(\alpha x y\beta)  \qquad  \text{for \ \ $x\ne y$, \ \ $x,y \in \Cnt(\alpha)$,}\\
& \aA_{x\ts x} \, := \,    \sum_{y\ts\in\ts \Des_{\alpha}(x)} \ \sum_{\beta\ts\in\ts  \Cnt_{m-1}(\alpha xy)} \,
\aq(\alpha x y\beta)\qquad \text{for \ \ $x \in \Cnt(\alpha)$,}\\
& \aA_{x\.\spstar} \, = \, \aA_{\spstar \. x} \, := \, \sum_{\beta \ts\in\ts \Cnt_{m-1}(\alpha x)}\, \aq(\alpha x\beta)
\qquad \text{for \ \ $x \in \Cnt(\alpha)$ \ \ and \ \ $y=\spstar$,}\\
& \aA_{\spstar\. \spstar} \, := \,  \sum_{\beta \ts\in\ts \Cnt_{m-1}(\alpha)} \, \aq(\alpha \beta).
\endaligned
$$

For the second line, note that \eqref{eq:ContInv} implies \. $\aA_{x\ts y} = \aA_{y\ts x}$.
Note also that \ts $\aA_{x\.\spstar} >0$, since by the exchange property
the word \ts $\alpha x \in \Lc$ \ts can be extended to \ts $\alpha x\beta \in \Lc$ \ts
for some \ts $\beta \in X^*$ \ts with \. $|\beta| \leq \rk(\Gf)-|\alpha|-1$.

\smallskip

For each vertex \. $\vf=(\alpha,m,t) \in \Vf$, define the associated matrix as follows:
\[  \bM \, = \, \bM_{(\alpha,m,t)} \ := \  t \, \bA(\alpha,m+1)  \ +  \ (1-t) \, \bA(\alpha,m).
\]
Similarly, define the associated vector \. $\hb = \hb_{(\alpha,m,t)} \in \Rb^{\ar}$ \. with
coordinates
\[  \ah_x \ := \
  \begin{cases}
  	\, t & \text{ if } \ x \in X,\\
  	\, 1-t & \text{ if } \ x =\spstar.
  \end{cases}    \]
Finally, define
the linear transformation \. $\bT^{\<x\>}: \Rb^{\ar} \to \Rb^{\ar}$ \.
associated to the edge \ts $(\vf,\vf^{\<x\>})$, as follows:
\begin{align*}
	  \big( \bT^{\<x\>}\vb\big)_y \ := \
	\begin{cases}
		\, \av_y & \text{ if } \  \, y  \. \in \. \supp(\bM), \\
		\, \av_x & \text{ if } \ \, y \.\in \. \Xf  \setminus  \supp(\bM).
	\end{cases}
\end{align*}

%
%

\medskip

\subsection{Properties of the atlas} \label{ss:greedoid-prop}
We now show that our combinatorial atlas \ts $\AA$ \ts
satisfies the conditions in Theorem~\ref{t:Hyp}, in the following series of lemmas.

\smallskip

\begin{lemma}\label{l:greedoid Irr}
	For every vertex \. $\vf\ts =\ts (\alpha,m,t)\in \Vf$, we have:
	\begin{enumerate}
		\item \label{item:Irr 1} the support of the associated matrix \ts $\bMr_{\vf}$ \ts is given by
		 $$\supp(\bMr_{\vf} \, ) \ = \  \supp(\bAr \, (\alpha,m+1)) \ = \ \supp(\bAr \, (\alpha,m)) \ = \
		 \begin{cases}
		 	\Cnt(\alpha) \ts + \ts \spstarr  & \text{ if }\  \alpha \in \Lc\ts,\\
		 	\varnothing & \text{ if } \ \alpha \notin \Lc.
		 \end{cases}	 $$
		\item \label{item:Irr 2} vertex \. $\vf$ \. satisfies \eqref{eq:Irr}, and
		\item \label{item:Irr 3} vertex \. $\vf$ \. satisfies \eqref{eq:hPos} \. for \. $t \in (0,1)$.
	\end{enumerate}
\end{lemma}

\smallskip

\begin{proof}
	Part \eqref{item:Irr 1} follows directly from the definition
	of matrices \. $\bM$, \. $\bA(\alpha,m+1)$, and \. $\bA(\alpha,m)$.
	Part \eqref{item:Irr 3} follows from the fact that \ts $\hb_{\vf}$ \ts is a strictly positive vector when \ts $t\in (0,1)$.
	
	We now prove part \eqref{item:Irr 2}.
	If \ts $\alpha \notin \Lc$, then $\bM$ is a zero matrix and
	\ts $\vf$ \ts trivially satisfies \eqref{eq:Irr}.
	If \ts $\alpha \in \Lc$,
	then it
	 follows from the definition of \. $\bM = \bigl(\aM_{x\ts y}\bigr)$,
 that \ts $\aM_{x\. \spstar} >0$ \ts for every \ts $x \in \Cnt(\alpha)$.
	Since the support of~$\ts\bM$ \ts is \ts $\Cnt(\alpha) + \spstar$,
	this proves \eqref{eq:Irr}, as desired.
\end{proof}

\smallskip

\begin{lemma}\label{l:greedoid Proj}
	For every greedoid \ts $\Gf=(X,\cL)$, the atlas \ts $\AA$ \ts satisfies \eqref{eq:Proj}.
\end{lemma}

\begin{proof}
		Let $\vf=(\alpha,m,t)\in \Vf^m$, $m\ge 1$, be a non-sink vertex of~$\Qf$.
	The condition \eqref{eq:Proj} follows directly from the definition of \ts $\bT^{\<x\>}$.	
\end{proof}
	
\smallskip

\begin{lemma}\label{l:greedoid Inh}
	For every greedoid \ts $\Gf=(X,\cL)$, the atlas \ts $\AA$ \ts satisfies  \eqref{eq:Inh}.
\end{lemma}

\begin{proof}
Let $\vf=(\alpha,m,t)\in \Vf^m$, $m\ge 1$, be a non-sink vertex of~$\Qf$.
Let \. $x \in \supp(\bM) =\Cnt(\alpha) \ts \cup \ts \{\spstar\}$.
By the linearity of $\bT^{\<x\>}$, it suffices to show that for every \.
$y \in \Cnt(\alpha) \ts \cup \ts \{\spstar\}$, we have:
\[  \aM_{xy} \ = \ \big\< \bT^{\<x\>} \eb_y\ts, \ts\bM^{\<x\>} \,  \bT^{\<x\>} \hb \big\>,
\]	
where \. $\big\{\eb_{y}, \. y \in \Xf\big\}$ \. is the standard basis for \ts $\Rb^{\ar}$.
We present only the proof for the case \ts $x,y \in \Cnt(\alpha)$, as the proof
of the other cases are analogous.	
		
First suppose that \ts $x,y \in \Cnt(\alpha)$ \ts are distinct.
Then:
	\begin{equation*}
	\begin{split}
		& \big\< \bT^{\<x\>} \eb_y\ts, \ts\bM^{\<x\>} \,  \bT^{\<x\>} \hb \big\> \quad = \ \
		\sum_{z \ts\in\ts \supp(\bM^{\<x\>})}  \aM^{\<x \>}_{yz} \, \big(\bT^{\<x\>} \hb \big)_z \\
		& \quad  = \  \sum_{z \ts\in\ts \supp(\bM^{\<x\>}), \. z\ne \spstar}
\, t\. \bA(\alpha x, m)_{yz}  \  + \ (1-t)\. \bA(\alpha x, m)_{y\ts \spstar}   \\
		& \quad  = \  \sum_{z \ts\in\ts X}  \, \sum_{\beta \in \Cnt_{m-1}(\alpha xyz)} \,  t\, \aq(\alpha x yz\beta)
 \ + \  \sum_{\beta \ts\in\ts \Cnt_{m-1}(\alpha xy)} \,(1-t) \,\aq(\alpha x y\beta)
\\
	& \quad  = \    \sum_{\gamma \ts\in\ts \Cnt_{m}(\alpha xy)} \, t \, \aq(\alpha x y\gamma) \ + \  \sum_{\beta \in \Cnt_{m-1}(\alpha xy)} \,(1-t) \,\aq(\alpha x y\beta),
	\end{split}
	\end{equation*}
	where we substitute \. $\gamma \gets z \beta$ \. in the first term of the last equality.
This implies that
	\begin{align*}
	\big\< \bT^{\<x\>} \eb_y\ts, \ts\bM^{\<x\>} \,  \bT^{\<x\>} \hb \big\>  \  = \ t \, \bA(\alpha, m+1)_{xy} \ + \
		(1-t) \,  \bA(\alpha, m)_{xy} \ = \   \aM_{xy}\.,
	\end{align*}
		which proves \eqref{eq:Inh} for this case.
		
		\smallskip
		
Now suppose that \. $x=y \in \Cnt(\alpha)$.
Then:
	\begin{align*}
	& \big\< \bT^{\<x\>} \eb_x, \bM^{\<x\>} \,  \bT^{\<x\>} \hb \big\> \ = \
	  \sum_{y  \ts\in\ts \supp(\bM^{\<x\>}) \setminus \supp(\bM)}   \  \sum_{z \ts\in\ts \supp(\bM^{\<x\>})}  \,   \aM^{\<x \>}_{yz} \, \big(\bT^{\<x\>} \hb \big)_z\..
	\end{align*}
	By the same argument as above,
 this equation becomes
	\begin{align*}
	 & \sum_{y  \ts\in\ts \supp(\bM^{\<x\>}) \setminus \supp(\bM)} \ \sum_{\gamma \ts\in\ts \Cnt_{m}(\alpha xy)} \.  t \, \aq(\alpha xy \gamma)
 \ + \  \sum_{\beta \ts\in\ts \Cnt_{m-1}(\alpha xy)} \.  (1-t) \,\aq(\alpha xy \beta)
 \\
	 &
\qquad = \   \sum_{y \ts\in\ts \Des_{\alpha}(x)} \ \sum_{\gamma \ts\in\ts\Cnt_{\alpha xy}(m)} \, t \, \aq(\alpha xy \gamma)
	 \ + \ \sum_{y \ts\in\ts \Des_{\alpha}(x)} \ \sum_{\beta \ts\in\ts\Cnt_{\alpha x y}(m-1)} \, (1-t) \, \aq(\alpha x y\beta) \\
	& \qquad  = \ t \,  \bA(\alpha, m+1)_{xx} \ \, + \ \,
	(1-t) \,  \bA(\alpha, m)_{xx} \ \ = \ \.  \aM_{xx}\,,
\end{align*}
which proves \eqref{eq:Inh} for this case.  This completes the proof.
\end{proof}

\medskip

\begin{lemma}\label{l:greedoid TPInv}
		Let \ts $\Gf=(X,\Lc)$ \ts be an interval greedoid,  and suppose the weight function \ts $\aqr:\Lc \to \rrs$ \ts
satisfies \eqref{eq:ContInv}.
		Then the atlas \ts $\AA$ \ts satisfies~\eqref{eq:TPInv}.
\end{lemma}

\smallskip

\begin{proof}
Let $\vf=(\alpha,m,t)\in \Vf^m$, $m\ge 1$, be a non-sink vertex of~$\Qf$, and
let \ts $x,y,z$ \ts be distinct elements of \. $\supp(\bM) = \Cnt(\alpha) + \spstar$.
We present only the proof for the case when
\. $x,y,z \in \Cnt(\alpha)$,  as other cases follow analogously.

First suppose that \. $\alpha \ts x'y'z' \ts \notin \ts\Lc$ \.
for every permutation \. $(x',y',z')$ \. of \. $\{x,y,z\}$.
Then
\[  \aM_{yz}^{\<x\>}	\ = \  \aM_{zx}^{\<y\>}	 \ = \ \aM_{xy}^{\<z\>} \ = \ 0,  \]
and \eqref{eq:TPInv} is satisfied.
So, without loss of generality, we assume that \. $\alpha x y z \in \Lc$.
 It then follows from the interval exchange property that
  $\alpha \ts x' y' z' \in \Lc$ for every permutation \. $(x',y',z')$ \. of \. $\{x,y,z\}$.
 This allows us to apply \eqref{eq:ContInv} for $\alpha \in \Lc$ and any two elements from $\{x,y,z\}$.

We now have
\begin{align*}
	\aM_{yz}^{\<x\>} \ \, &= \  \, \bA(\alpha x,m)_{yz} \,
	\ = \ \, \sum_{\beta \ts\in\ts \Cnt_{m-1}(\alpha xyz)} \, \aq(\alpha  xyz \beta) \\
	\ &=_{\eqref{eq:ContInv}} \  \sum_{\beta \ts\in\ts \Cnt_{m-1}(\alpha yxz)} \, \aq(\alpha  yxz \beta)
 \ \, =  \  \, \bA(\alpha y,m)_{xz} \ \, =  \ \, \aM^{\<y\>}_{xz}\,.
\end{align*}
By an analogous argument, it follows that
\ts $\aM_{yz}^{\<x\>} = \aM_{xy}^{\<z\>}$, and thus
 \eqref{eq:TPInv} is satisfied, as desired.
\end{proof}

\smallskip

\begin{lemma}\label{l:greedoid Knonneg}
			Let \ts $\Gf=(X,\Lc)$ \ts be a greedoid, let \ts $1\le k < \rk(\Gf)$,
and suppose the weight function  \ts $\aqr:\Lc \to \rrs$ \ts  satisfies \eqref{eq:ContInv} and \eqref{eq:PasAct}.
	Then the atlas \ts $\AA$ \ts satisfies \eqref{eq:Knonneg}.
\end{lemma}

\smallskip

\begin{proof}
Let $\vf=(\alpha,m,t)\in \Vf^m$, $m\ge 1$, be a non-sink vertex of~$\Qf$.
		We need to check the condition
		\eqref{eq:Knonneg} for distinct \. $x,y \in \supp(\bM) = \Cnt(\alpha) + \spstar$.

	First suppose that $x,y$ are distinct elements of  $\Cnt(\alpha)$.
	We have:
	\begin{align*}
			&\aM_{yy}^{\<x\>}  \ = \  \bA(\alpha x,m)_{yy}  \ \, = \ \,
\sum_{z \ts\in\ts \Des_{\alpha x}(y)}  \ \sum_{\beta \ts\in\ts \Cnt_{m-1}(\alpha xyz)} \, \aq(\alpha xyz \beta).
	\end{align*}
 	Note that $z \in X$ in the equation above is summed over the set
	\begin{align*}
		 	\big\{ z \in X \, : \,  \alpha xz \notin \Lc, \. \alpha x y z \in \Lc  \big\}.
	\end{align*}
By  the interval exchange property, every element $z$ in the set above also satisfies $\alpha z \notin \Lc$.
	We can then partition the set above into
	\begin{align*}
		& \big\{z  \, : \,  \alpha z \notin \Lc, \. \alpha xz \notin \Lc, \. \alpha yz \notin \Lc, \. \alpha xyz \in \Lc  \big\} \  \cup \ \big\{z  \, : \,  \alpha z \notin \Lc, \. \alpha xz \notin \Lc, \. \alpha yz \in \Lc, \. \alpha xyz \in \Lc  \big\} \\
		& \qquad \ \ = \ \Pas_{\alpha}(x,y) \  \cup \  \big\{z  \, : \,  \alpha z \notin \Lc, \. \alpha xz \notin \Lc, \. \alpha yz \in \Lc, \. \alpha xyz \in \Lc  \big\}.
	\end{align*}
	On the other hand, we have:
	\begin{align*}
		&\sum_{z \ts \in \ts \Fam^{\<y\>}} \. \aM_{xz}^{\<y\>}  \ = \
		 \sum_{z \ts \in \ts \Fam^{\<y\>}} \. \bA(\alpha y,m)_{xz}
		 \ = \ \sum_{z \in \Des_{\alpha }(y)} \ \sum_{\beta \in \Cnt_{m-1}(\alpha yxz)} \, \aq(\alpha yxz \beta) \\
& \qquad =_{\eqref{eq:ContInv}} \ \sum_{z \in \Des_{\alpha }(y)} \ \sum_{\beta \in \Cnt_{m-1}(\alpha xyz)} \, \aq(\alpha xyz \beta),
		\end{align*}
where in the last equality we apply \eqref{eq:ContInv} to swap $x$ and~$y$.
 	Note that $z \in X$ in the equation above is summed over the set
		\begin{align*}
	\{ z \in X \, : \,  \alpha z \notin \Lc, \. \alpha yz \in \Lc, \. \alpha x y z \in \Lc  \},
	\end{align*}
		which can be partitioned into
	\begin{align*}
	& \big\{z  \, : \,  \alpha z \notin \Lc, \. \alpha xz \in  \Lc, \. \alpha yz \in \Lc, \. \alpha xyz \in \Lc  \big\} \  \cup \ \big\{z  \,: \,  \alpha z \notin \Lc, \. \alpha xz \notin \Lc, \. \alpha yz \in \Lc, \. \alpha xyz \in \Lc  \big\} \\
			& \qquad = \ \Act_{\alpha}(x,y) \  \cup \  \big\{z  \, :\,  \alpha z \notin \Lc, \. \alpha xz \notin \Lc, \. \alpha yz \in \Lc, \. \alpha xyz \in \Lc  \big\}.
\end{align*}
It follows from the calculations above that
\begin{align*}
& \aM_{yy}^{\<x\>}  	\  \ -  \, \sum_{z \ts \in \ts \Fam^{\<y\>}} \. \aM_{xz}^{\<y\>}  \ \ = \
 \sum_{z \ts\in\ts \Pas_{\alpha}(x,y)} \
\sum_{\beta \ts\in\ts \Cnt_{m-1}(\alpha xyz)} \, \aq(\alpha xyz \beta)  \\  & \hskip2.7cm -  \
\sum_{z \ts\in\ts \Act_{\alpha}(x,y)} \
\sum_{\beta \ts\in\ts \Cnt_{m-1}(\alpha xyz)} \, \aq(\alpha xyz \beta),
	\end{align*}
	which is nonnegative by \eqref{eq:PasAct}.
	This proves \eqref{eq:Knonneg} in this case.
	
	\smallskip

	Now suppose that \ts $x \in \Cnt(\alpha)$ \ts and \ts $y=\spstar$.
	Without loss of generality, we assume that \ts $\alpha \in \Lc$,
	as otherwise \ts
	$\aM^{\<x\>}_{yz}$ \ts is always equal to the zero matrix and
 \eqref{eq:Knonneg} is  trivially satisfied.
	Then we have:
\begin{align*}
\aM_{\spstar\. \spstar}^{\<x\>} \ = \ \bA(\alpha x,m)_{\spstar\.\spstar}  \ =
\sum_{\beta \in \Cnt_{m-1}(\alpha x)} \, \aq(\alpha x \beta).
\end{align*}
	On the other hand, we have \. $\supp(\bM) \ts = \ts \supp\big(\bM^{\<\spstar\>}\big) \ts = \ts \Cnt(\alpha) + \spstar$,
	which implies that
	\[  \Fam^{\<\spstar\>} \ = \ \supp\big(\bM^{\<\spstar\>}\big)   \setminus  \big(\supp(\bM) \ts - \ts \spstar\big)  \ = \ \{\spstar\}.  \]
Therefore, we have:	
\begin{align*}
		\sum_{z \ts \in \ts \Fam^{\<\spstar\>}} \aM_{xz}^{\<\spstar\>}
		\ &= \  \aM^{\<\spstar\>}_{x \ts \spstar}  \ = \ \bA(\alpha ,m)_{x\ts\spstar}  \ = \
\sum_{\beta \in \Cnt_{m-1}(\alpha x)} \, \aq(\alpha x\beta).
	\end{align*}
	It thus follows from the calculations above
	that \. $\aK_{x\ts \spstar}= 0$.
	This completes the proof of  \eqref{eq:Knonneg}.
\end{proof}

\medskip

\subsection{Basic hyperbolicity} \label{ss:greedoid-basic}
To prove hyperbolicity of vertices in $\Vfm$, we need the following straightforward
linear algebra lemma.  We include the proof for completeness.

\smallskip

\begin{lemma}\label{l:Hyp N}
	Let \ts $\bNr=\big(\aNr_{i\ts j}\big)$ \ts be a nonnegative symmetric \ts $(n+1)\times (n+1)$  matrix,
such that its nondiagonal entries are equal to~$1$. Suppose that
$$(\ast) \qquad
\aNr_{1\ts 1},\ldots, \aNr_{n\ts n} \.\leq\. 1 \qquad \text{and} \qquad
\aNr_{n+1\ts n+1}  \ \geq \  \sum_{i=1}^n  \, \frac{\aNr_{n+1\ts n+1}-1}{1-\aNr_{i\ts i}} \, \ \text{ if } \ \aNr_{ii}<1  \ \text{for all } \ i \in [n].
$$
 	Then $\bNr$ satisfies \eqref{eq:Hyp}.
\end{lemma}

\begin{proof}  Fix \ts $\ep>0$.
Substituting \ts $\aN_{i\ts i}\gets \aN_{i\ts i}-\ep$ \ts
for every \ts $1\le i \le n$ \ts if necessary,
we can assume that all inequalities in~$(\ast)$ are strict.
Note that \eqref{eq:Hyp} is preserved under taking the limit \ts
$\ep \to 0$, so it suffices to prove the result in this case.

We prove that $\bN$ satisfies \eqref{eq:OPE} by induction on~$n$.
By Lemma~\ref{l:Hyp is OPE} this implies~\eqref{eq:Hyp}.
The base case $n=0$ is trivial.  Assume that the claim is true for \ts $(n-1)$.
Let \. $\lambda_1 \geq \ldots \geq \lambda_{n+1}$ \. be the eigenvalues of~$\bN$,
and let \. $\lambda_1' \geq \ldots \geq \lambda_n'$\.  be the eigenvalues
of the matrix obtained by removing the first row and column of~$\bN$.
By the \emph{Cauchy interlacing theorem}, we have
\[
    \lambda_1 \ \geq \  \lambda_1' \ \geq \ \lambda_2  \ \geq \ \ldots \ \geq \ \lambda_n' \ \geq \  \lambda_{n+1}\..
\]

Note that \. $\lambda_2', \ldots, \lambda_n'$ \. are nonpositive by induction.
It then follows that \. $\lambda_3,\ldots, \lambda_{n+1}$ \. are nonpositive.
By the Perron--Frobenius theorem, we also have \ts $\lambda_1 >0$.
It thus suffices to show that \ts $\lambda_2 \leq 0$, which will follow from
showing that \ts $\det(\bN)$ \ts has sign \ts $(-1)^n$.
Observe that \. $\det(\bN)$ \. is equal to
{\small
	\begin{align*}
		\begin{vmatrix}
			\aN_{1\ts 1} - 1  &  0  &  \cdots & 0  &    1-\aN_{n+1\ts n+1} \\[2 pt] 		
			0  & 	 \aN_{2\ts 2} -1 &  &   0  &   1-\aN_{n+1\ts n+1}  \\[2 pt]
			\vdots &  & \ddots &   \vdots  & \vdots  \\[3 pt]
			0  &  0 & \cdots  \,  &
			\aN_{n\ts n} -1   & 1 - \aN_{n+1\ts n+1}  \\[3 pt]
			1 & 1 & \ldots & 1 & \aN_{n+1\ts n+1}
		\end{vmatrix} \ = \
\begin{vmatrix}
			\aN_{1\ts 1} - 1  &  0  &  \cdots & 0  &    1-\aN_{n+1\ts n+1} \\[2 pt] 		
			0  & 	 \aN_{2\ts 2} -1 &  &   0  &   1-\aN_{n+1\ts n+1}  \\[2 pt]
			\vdots &  & \ddots &   \vdots  & \vdots  \\[3 pt]
			0  &  0 & \cdots  \,  &
			\aN_{n\ts n} -1   & 1 - \aN_{n+1\ts n+1}  \\[3 pt]
			0 & 0 & \ldots & 0 & \aJ
		\end{vmatrix},
	\end{align*}	
}
\[ \text{where} \qquad \aJ \ := \ \aN_{n+1\ts +1}  \ - \
\sum_{i=1}^{n} \, \frac{\aN_{n+1\ts n+1}\. -\. 1}{1\. -\. \aN_{i\ts i}} \, > \, 0\,,
\]
by the assumption~$(\ast)$.  Therefore, we have
	\begin{align*}
		\det(\bN)  \ = \ \aJ \.\cdot \. \prod_{i=1}^n \.	\big(\aN_{i\ts i}\. - \. 1 \big)\ts,
	\end{align*}
and by the assumptions on $\ts \aN_{i\ts i}$ \ts this determinant has sign \ts $(-1)^n$.
This completes the proof.
\end{proof}

\bigskip

\subsection{Proof of Theorem~\ref{t:greedoid}} \label{ss:greedoid-proof}
We first show that every sink vertex in the combinatorial atlas \ts $\AA$ \ts is hyperbolic.

\smallskip

\begin{lemma}\label{l:greedoid BHyp}
Let \ts $\Gf=(X,\cL)$ \ts be an interval greedoid on \ts $|X|=n$ elements,
let \ts $1\le k <\rk(\Gf)$, and let  \ts
$\aqr:\Lc \to \rrs$ \ts be a $k$-admissible weight function.
Then every vertex in \ts $\Vfm$ \ts satisfies \eqref{eq:Hyp}.
\end{lemma}

\begin{proof}
	Let \ts $\vf=(\alpha, 0,1)\in \Vfm$ \ts be a sink vertex.
	It suffices to show that \ts $\bA(\alpha,1)$ \ts
	satisfies \eqref{eq:Hyp}. First note that if \ts
$\alpha \notin \Lc$, then \ts $\bA(\alpha,1)$ \ts
	is a zero matrix, and \eqref{eq:Hyp} is trivially true.
	Thus, we can assume that \ts $\alpha \in \Lc$.
	We write \. $\aA_{x,y} \ts := \ts \bA(\alpha,1)_{xy}$ \.
for every \ts $x,y \in X$.
	
\smallskip

Let \ts $\Cc\in \Par(\alpha)$ \ts be a parallel class.
Suppose that \ts $|\Cc|\geq 2$, and let \ts $x,y$ \ts
be distinct elements of~$\Cc$.
	
\smallskip

\nin
\textbf{Claim:} \. {\em For every \ts $z \in \Xf$, we have \,
$\apr(\alpha y) \.  \aAr_{xz} \, = \,   \apr(\alpha x)  \. \aAr_{yz}$\ts.}

	\begin{proof}[Proof of Claim]
		First suppose that $z \in \{x,y\}$.
		It then follows from \eqref{eq:FewDes} and the fact that $\alpha x y \notin \Lc$ that
		\. $\aA_{x,z} \ = \  \aA_{y,z} \ = \ 0$,
		which implies the claim in this case.
		
		Now suppose that \ts $z \in X \setminus \{x,y\}$.
		It follows from the exchange property that \ts $\alpha xz \in \Lc$ \ts if and only if \ts $\alpha yz \in \Lc$.
		There are now two cases. If \ts $\alpha xz \notin \Lc$ \ts and \ts $\alpha yz \notin \Lc$,
		then again we have \. $\aA_{x,z}  =   \aA_{y,z}  =  0$, which implies the claim.
		If \ts $\alpha xz \in \Lc$ \ts and \ts $\alpha yz \in \Lc$, we then have:
		\begin{align*}
			\aA_{xz} \, &= \,  \aq(\alpha xz)  \, =_{\eqref{eq:LogMod}} \, \ac_{\ell+2} \.\frac{\ap(\alpha x) \, \ap(\alpha z) }{ \ap(\alpha)}\,,  \qquad
			\aA_{yz} \, = \,  \aq(\alpha yz)  \, =_{\eqref{eq:LogMod}} \, \ac_{\ell+2} \.\frac{\ap(\alpha y) \, \ap(\alpha z) }{ \ap(\alpha)}\,,
		\end{align*}
	where \ts $\ell :=|\alpha|$.
	This implies the claim in this case.
	Finally, let \ts $z=\spstar$. Then we have \.  $\aA_{xz}  =  \ac_{\ell+1}  \. \ap(\alpha x)$ \. and \.
$\aA_{yz}  =  \ac_{\ell+1} \. \ap(\alpha y)$, which implies the claim.
	\end{proof}
	
\smallskip

Deduct the $y$-row and $y$-column of \ts $\bA(\alpha,1)$ \ts
by \. $\frac{\ap(\alpha y)}{\ap(\alpha x)}$ \. of the $x$-row and $x$-column of \ts $\bA(\alpha,1)$.
It then follows from the claim that the resulting matrix has $y$-row and $y$-column is equal to zero.
	Also, note that \eqref{eq:Hyp} is preserved under this transformation.
	Applying this linear transformation repeatedly,
	and by restricting to the support of resulting matrix which preserves \eqref{eq:Hyp},
	without loss of generality we can
	assume that \ts $|\Cc|=1$ \ts for every parallel class \ts $\Cc\in \Par(\alpha)$.
	Then the matrix \. $\bA(\alpha,1)$ \. is equal to
		\begin{align*}
 \begin{pmatrix}
			\ac_{\ell+2} \, \ab_{\alpha}(\Cc_1) \, \frac{\ap(\alpha x_1)^2}{\ap(\alpha)} \,    &    \aq(\alpha x_1 x_2) & \hspace{3 pt} \cdots \hspace{3 pt}
			& \aq(\alpha x_1 x_n) &    \aq(\alpha x_1)  \\[2 pt] 		
			\aq(\alpha x_2 x_1) & \ac_{\ell+2} \, \ab_{\alpha}(\Cc_2) \,  \frac{\ap(\alpha x_1)^2}{\ap(\alpha)}
			&  &   \vdots  &   \vdots  \\[2 pt]
			\vdots &  & \ddots &   \aq(\alpha x_{n-1} x_n) & \aq(\alpha x_{n-1})   \\[3 pt]
		\aq(\alpha x_n x_1) & \cdots & \aq(\alpha x_n x_{n-1}) &
			\ac_{\ell+2} \, \ab_{\alpha}(\Cc_\ar) \,  \frac{\ap(\alpha x_1)^2}{\ap(\alpha)} & \aq(\alpha x_n) \\[3 pt]
		\aq(\alpha x_1) & \ldots & \aq(\alpha x_{n-1}) & \aq(\alpha x_{n}) &  \aq(\alpha)
		\end{pmatrix},
	\end{align*}
	where \ts $\Cc_i = \{x_i\}$ \ts for $i \in [n]$,  the rows and columns
are indexed by \. $\Xf=\big\{x_1,\ldots, x_n, \spstar\big\}$, and  $\ab_{\alpha}(\Cc)$ is as defined in  \eqref{eq:ab-def}.
	We now rescale  the $x_i$-row and $x_i$-column
	by $\frac{\sqrt{\ap(\alpha)}}{\sqrt{\ac_{\ell+2}} \ts  \ap( \alpha x_i)}$,
	and the $\spstar$-row and $\spstar$-column by
	$\frac{\sqrt{\ac_{\ell+2}}}{\ac_{\ell+1} \, \sqrt{\ap(\alpha)}}$.
Again, note that \eqref{eq:Hyp} is preserved under this transformation.
	It then follows from \eqref{eq:LogMod} that the matrix becomes
	  		\begin{align*}
	  	\begin{pmatrix}
	  	\ab_{\alpha}(\Cc_1) &   1 & \hspace{3 pt} \cdots \hspace{3 pt}
	  		& 1 &    1  \\[2 pt] 		
	  		1 &  \ab_{\alpha}(\Cc_2)
	  		&  &   \vdots  &   \vdots  \\[2 pt]
	  		\vdots &  & \ddots &  1 & 1  \\[3 pt]
	  	1 & \cdots & 1 &
	  		 \ab_{\alpha}(\Cc_n)  & 1 \\[3 pt]
	  		1 & \ldots & 1 & 1 &  \frac{\ac_{\ell+2} \, \ac_{\ell}}{ \ac_{\ell+1}^2}
	  	\end{pmatrix}.
	  \end{align*}
	It follows from \eqref{eq:SynMon} and \eqref{eq:ScaleMon} that this matrix satisfies conditions~$(\ast)$
in Lemma~\ref{l:Hyp N}.  Hence, by the lemma, this matrix satisfies \eqref{eq:Hyp}.  We conclude
that  $\vf$ satisfies  \eqref{eq:Hyp}, as desired.	
\end{proof}

\smallskip

We can now prove that every vertex in $\Qf$ is hyperbolic.

\smallskip

\begin{prop}\label{p:greedoid}
Let \ts $\Gf=(X,\cL)$ \ts be an interval greedoid on \ts $|X|=n$ elements,
let \ts $1\le k <\rk(\Gf)$, and let  \ts
$\aqr:\Lc \to \rrs$ \ts be a $k$-admissible weight function.
Then every vertex in \ts $\Vf$ \ts satisfies \eqref{eq:Hyp}.\end{prop}

\begin{proof}
	We will show that every vertex in \ts $\Vf^m$ \ts for \ts $m \leq k-1$ \ts
satisfies \eqref{eq:Hyp} by induction on $m$.
	The claim is true for \ts $m=0$ \ts by Lemma~\ref{l:greedoid BHyp}.
	Suppose that the claim is true for \ts $\Vf^{m-1}$.
Note that the atlas \ts $\AA$ \ts satisfies the assumptions of Theorem~\ref{t:Hyp}
by Lemmas~\ref{l:greedoid Proj}, \ref{l:greedoid Inh}, \ref{l:greedoid TPInv}, and~\ref{l:greedoid Knonneg}.
It then follows that every regular vertex in \ts $\Vf^m$ \ts satisfies \eqref{eq:Hyp}.

	On the other hand, by Lemma~\ref{l:greedoid Irr}, the regular vertices of \ts $\Vf^m$ \ts
are those of the form \ts $\vf=(\alpha,m,t)$ \ts with \ts $t \in (0,1)$.
Since \eqref{eq:Hyp} is preserved under taking the limits \ts $t\to 0$ \ts and \ts $t\to 1$,
it then follows that every vertex in $\Vf^m$ satisfies \eqref{eq:Hyp}, and the proof is complete.
\end{proof}

\medskip

\begin{proof}[Proof of Theorem~\ref{t:greedoid}]
	Let \ts $\bM=\bM_{\vf}$ \ts be the matrix associated with the vertex \ts $\vf=(\varnothing,k-1,1)$.	
Let \ts $\vb$ \ts and~$\ts \wb$ \ts be the characteristic vectors of \ts $X$ \ts and \ts $\{\spstar\}$, respectively.
	Then:
\begin{equation}\label{eq:greedoid-last}
\aL_\aq(k+1) \ = \ \< \vb, \bM \vb \>\ts, \qquad \aL_\aq(k) \ = \ \< \vb, \bM \wb \>\ts, \qquad \aL_\aq(k-1) \ = \ \< \wb, \bM \wb \>\ts.
\end{equation}	
	By Proposition~\ref{p:greedoid}, vertex \ts $\vf$ \ts satisfies \eqref{eq:Hyp}.
Substituting~\eqref{eq:greedoid-last} into \eqref{eq:Hyp}, gives the
log-concave inequality \eqref{eq:greedoid-LC} in the theorem.  \end{proof}

\bigskip

\section{Proof of equality conditions for interval greedoids}
\label{s:proof-greedoid-equality}

In this section we prove Theorem~\ref{t:greedoid-equality-full}.
The implication \. $\eqref{eq:item-equality-greedoid-b} \Rightarrow \eqref{eq:item-equality-greedoid-a}$ \.
is obvious. We now prove the other implications.

\subsection{Proof of \. {\rm \eqref{eq:item-equality-greedoid-a} $\Rightarrow$  \eqref{eq:GrEqu1} $\&$ \eqref{eq:GrEqu2}}}
\label{ss:proof equality greedoid a implies c}

Let \ts $\AA$ \ts be the combinatorial atlas defined in~$\S$\ref{ss:greedoid-def},
that corresponds to \ts $(\Gf,k, \aq)$.
Recall that every vertex of $\Qf$ satisfies \eqref{eq:Hyp}
by Proposition~\ref{p:greedoid}.

\smallskip

As at the end of previous section, let \ts $\vb, \wb \in \Rb^{\ar}$ \ts be
the characteristic vectors of \ts $X$\ts and \ts $\{\spstar\}$, respectively.
It is straightforward to verify that $\vb,\wb$ is a global pair of $\Qf$, i.e.\
they satisfy \eqref{eq:PosGlob}.

Let \ts $\vf=(\varnothing,k-1,1)\in \Vf$ \ts and let \ts $\bM=\bM_{\vf}$ \ts
be the matrix associated with~$\vf$.  Recall that \ts $\bM=\bA(\varnothing,k)$ \ts
and we have equalities~\eqref{eq:greedoid-last} again:
	\begin{equation}\label{eq:kilo 3}
\aL_\aq(k+1) \ = \ \< \vb, \bM \vb \>\ts, \qquad \aL_\aq(k) \ = \ \< \vb, \bM \wb \>\ts, \qquad \aL_\aq(k-1) \ = \ \< \wb, \bM \wb \>\ts.
	\end{equation}
Note also that \. $\aL_\aq(k+1), \aL_\aq(k), \aL_\aq(k-1)>0$ \. since \. $k <  \rk(\Gf)$.
It then follows from \eqref{eq:item-equality-greedoid-a}, that $\vf$ satisfies \eqref{eq:sEqu} for some $\as>0$.

Let us show that, for every \ts $\alpha \in \Lc$ \ts of length $(k-1)$, we have:
	\begin{equation}\label{eq:lima 4}
		\langle \vb, \bA(\alpha,1) \vb \rangle  \ = \  \as \, \langle \wb, \bA(\alpha,1) \vb \rangle \ = \  \as^2 \, \langle \wb , \bA(\alpha,1) \wb \rangle.
	\end{equation}
	First, suppose that \ts $k=1$.
	It then follows that \ts $\alpha =\varnothing$ \ts and \ts $\vf=(\varnothing,0,1)$.
	Equation \eqref{eq:lima 4} now follows from the fact that $\vf$ satisfies \eqref{eq:sEqu}.

Now suppose that \ts $k>1$. Then it is straightforward to verify that
$\vf$ is a functional source, i.e.\ satisfies \eqref{eq:ProjGlob} and \eqref{eq:hGlob}, where we apply the substitution $\fb \gets \vb$ for \eqref{eq:hGlob}.
By Theorem~\ref{t:equality Hyp}, every functional target of $\vf$ also
satisfies \eqref{eq:sEqu} with the same \ts $\as >0$.
On the other hand, it is straightforward to verify that the functional
targets of $\vf$ are  those of the form \. $(\alpha, 0,1)$.
Combining these observations, we conclude \eqref{eq:lima 4}.

\smallskip

	Let \ts $\zb :=  \vb - \as \wb$.
It follows from
\eqref{eq:lima 4}
that \. $\langle \zb, \bA(\alpha,1) \zb \rangle =0$.
It then follows from  Lemma~\ref{l:equality kernel}
that \. $\bA(\alpha,1) \zb = \0$, which is equivalent to
	\. $ \as \, \bA(\alpha,1)  \wb  = \bA(\alpha,1)  \vb$.
	This implies that
		\begin{align*}
		\as \, \aq(\alpha)  \ = \ \as \, \big(\bA(\alpha,1)  \wb\big)_{\spstar} \ = \ \big(\bA(\alpha,1)   \vb\big)_{\spstar}
		\ = \   \sum_{ x \ts\in\ts \Cnt(\alpha)} \.\aq(\alpha x),
	\end{align*}	
	which proves \eqref{eq:GrEqu1} for $\as(k-1)=\as$.

\smallskip
	
	Let \ts $x\in \Cnt(\alpha)$ \ts be an arbitrary continuation.
	By the same reasoning as above,  we have:
$$\as \,  \aq(\alpha x)   \ = \
\as \, \big(\bA(\alpha,1)  \wb\big)_{x}  \ =  \ \big(\bA(\alpha,1)  \vb\big)_{x}\,.
$$
On the other hand, we also have:
	\[
	\big(\bA(\alpha,1)  \vb\big)_{x} \ = \ \sum_{y \ts\in\ts\Des_\alpha(x)} \aq(\alpha xy )  \ + \      \sum_{\substack{y \ts\in\ts \Cnt(\alpha) \\ y \ts\not \salp\ts x } } \aq(\alpha xy )\,.
	\]
	It then follows that:
\begin{align}\label{eq:papa 4}
		\as \,  \aq(\alpha x) \ &= \         \sum_{y \ts\in\ts \Des_\alpha(x)} \aq(\alpha xy )  \ + \      \sum_{\substack{y \ts\in\ts \Cnt(\alpha) \\ y \ts\not \salp\ts x } } \. \aq(\alpha xy )\ts.
\end{align}
	Let $\Cc$ be the parallel class in \ts $\Par_\alpha$ \ts containing~$x$.
	We now show that \eqref{eq:papa 4} is equivalent to \eqref{eq:GrEqu2}.
	
\smallskip

Applying \eqref{eq:LogMod} to \eqref{eq:papa 4} and dividing both sides by \ts $\ap(\alpha)$,  we get:
	\begin{align}\label{eq:papa 1}
		\as \,  \ac_{k} \,  \frac{\ap(\alpha x)}{\ap(\alpha)}    \ &= \
		\sum_{y \ts\in\ts\Des_\alpha(x)} \ac_{k+1} \, \frac{\ap(\alpha xy)}{\ap(\alpha)}  \ + \      \sum_{\substack{y \ts\in\ts \Cnt(\alpha) \\ y \ts \not \salp \ts x} } \ac_{k+1} \, \frac{\ap(\alpha x )\. \ap(\alpha y)}{\ap(\alpha)^2} \,.
	\end{align}

	Now note that, \eqref{eq:GrEqu1} gives:
\begin{align}\label{eq:papa 3}
	\sum_{\substack{y \ts\in\ts \Cnt(\alpha) \\ y \ts\not\ts \salp x} }  \frac{\ap(\alpha y)}{\ap(\alpha)} \ = \ \as \, \frac{\ac_{k-1}}{\ac_{k}} \ - \  \aa_{\alpha}(\Cc),
\end{align}
where
\[\aa_{\alpha}(\Cc) \ :=  \ \sum_{y \ts\in\ts \Cc} \. \frac{\ap(\alpha y)}{\ap(\alpha)}.\]

Now note that,  when  $|\Cc|\geq 2$,
\begin{align*}
	\sum_{y \ts\in\ts \Des_\alpha(x)}  \frac{\ap(\alpha xy)}{\ap(\alpha)} \ =_{\eqref{eq:FewDes}}  \ 0 \ = \   \frac{\ap(\alpha x) }{\ap(\alpha)}  \, \aa_{\alpha}(\Cc) \, \ab_{\alpha}(\Cc),
\end{align*}
where the last equality is because \. $\ab_{\alpha}(\Cc) =0$ \. when $|\Cc| \geq 2$.
On the other hand, when $\Cc=\{x\}$,
\begin{align*}
	\sum_{y \ts\in\ts \Des_\alpha(x)}  \frac{\ap(\alpha xy)}{\ap(\alpha)}  \ = \
	\frac{\ap(\alpha x)^2 }{\ap(\alpha)^2}  \ \ab_{\alpha}(\Cc) \ = \ \frac{\ap(\alpha x) }{\ap(\alpha)}  \ \aa_{\alpha}(\Cc) \ \ab_{\alpha}(\Cc),
\end{align*}
where the last equality is because \. $\aa_{\alpha}(\Cc)= \frac{\ap( \alpha x)}{\ap(\alpha)}$ \.  when $\Cc=\{x\}$.
This allows us to conclude that
\begin{align}\label{eq:papa 2}
	\sum_{y \ts\in\ts \Des_\alpha(x)}  \frac{\ap(\alpha xy)}{\ap(\alpha)}  \ = \   \frac{\ap(\alpha x) }{\ap(\alpha)}  \ \aa_{\alpha}(\Cc) \ \ab_{\alpha}(\Cc),
\end{align}

Substituting \eqref{eq:papa 3} and \eqref{eq:papa 2} into \eqref{eq:papa 1}, we obtain:
\begin{align*}
	\as \,  \ac_{k} \,  \frac{\ap(\alpha x)}{\ap(\alpha)}    \ &= \
	\ac_{k+1} \, \frac{\ap(\alpha x)}{\ap(\alpha)} \,  \aa_{\alpha}(\Cc) \, \ab_{\alpha}(\Cc)   \ + \       \ac_{k+1} \, \frac{\ap(\alpha x )}{\ap(\alpha)} \,  \bigg(  \as \, \frac{\ac_{k-1}}{\ac_{k}} \ - \  \aa_{\alpha}(\Cc) \bigg).
\end{align*}

		This is equivalent to
	\[ \as \bigg( \frac{\ac_{k-1}}{\ac_{k}} - \frac{\ac_{k}}{\ac_{k+1}} \bigg) \ = \ \aa_{\alpha}(\Cc) \, \big(1-\ab_{\alpha}(\Cc) \big),  \]
	which is  \eqref{eq:GrEqu2}.  This completes the proof.  \qed

		\medskip
		
\subsection{Proof of \. {\rm \eqref{eq:GrEqu1} \& \eqref{eq:GrEqu2} $\Rightarrow$ \eqref{eq:item-equality-greedoid-b}}}
		Write $\as := \as(k-1)$.
			We have \. $\aL_{\aq,\alpha}(0) = \aq(\alpha)$ \. by definition, and
			\begin{align*}
				\aL_{\aq,\alpha}(1) \ = \   \sum_{x \ts\in\ts \Cnt(\alpha)}   \aq(\alpha x) \ =_{\eqref{eq:GrEqu1}}  \  \as \, \aq(\alpha)  \ = \ \as \, \aL_{\aq,\alpha}(0),
			\end{align*}
			which proves the first part of \eqref{eq:item-equality-greedoid-b}.
			For the second part of \eqref{eq:item-equality-greedoid-b}, we have:
			\begin{align}\label{eq:lima 1}
				\aL_{\aq,\alpha}(2) \ = \     \sum_{x \in \Cnt(\alpha)} \left( \sum_{y \ts\in\ts \Des_\alpha(x)} \aq(\alpha xy )  \ + \    \sum_{\substack{y \ts\in\ts \Cnt(\alpha) \\ y \ts\not\ts \salp x} } \aq(\alpha xy ) \right).
			\end{align}
On the other hand, we showed in the proof above (see~\S\ref{ss:proof equality greedoid a implies c}),
that \eqref{eq:GrEqu2} is equivalent to~\eqref{eq:papa 4}.  Therefore, for every \ts $x \in \Cnt(\alpha)$, we have:
			\begin{align*}
			\as \,  \aq(\alpha x) \ &= \     \sum_{y \ts\in\ts \Des_\alpha(x)} \aq(\alpha xy )  \ + \      \sum_{\substack{y \ts\in\ts \Cnt(\alpha) \\ y \ts\not\ts \salp x } } \aq(\alpha xy ).
		\end{align*}
		Substituting this equation into \eqref{eq:lima 1},  we conclude:
		\begin{align*}
			\aL_{\aq,\alpha}(2) \ = \     \sum_{x \ts\in\ts \Cnt(\alpha)}
			\as \, \aq(\alpha x)	\ =_{\eqref{eq:GrEqu1}} \ 	 \as^2 \, \aq(\alpha) \ = \ \as^2 \, \aL_{\aq,\alpha}(0).
		\end{align*}
			This proves the second part of \eqref{eq:item-equality-greedoid-b}, and completes the proof. \qed


\bigskip

\section{Proof of matroid inequalities and equality conditions}\label{s:proof-matroid}

\smallskip

In this section we give proofs of Theorem~\ref{t:matroids-Par-weighted},
Theorem~\ref{t:matroids-equality-weighted}, Theorem~\ref{t:matroids-equality-Par},
Proposition~\ref{p:graphical-equality} and give further extension of
graphical matroid results.  We conclude with two explicit examples
of small combinatorial atlases.

\smallskip

\subsection{Proof of Theorem~\ref{t:matroids-Par-weighted}}\label{ss:proof-matroid-Par-weighted}
We deduce the result from Theorem~\ref{t:greedoid}.
Let \. $\Gf=(X,\Lc)$ \. be the interval greedoid
constructed in~$\S$\ref{ss:reduction-matroids}, and
corresponding to matroid \. $\Mf=(X,\Ic)$.
Let \. $1\le k < \rk(\Mf)$ \. and let \. $\ap:X \to \rrs$ \. be
\. as in the theorem.
Define the weight function  \. $\aq:\Lc \to \rrs$ \. by the product formula:
\begin{equation}\label{eq:matroid-proof-q-product}
	  \aq(\alpha) \ := \ \ac_{\ell} \, \prod_{x \in \alpha} \. \ap(x),
\end{equation}
where \ts $\ell:=|\al|$, and \ts $\ac_\ell$ \ts is given by
\begin{equation}\label{eq:matroid-proof-c-ell}
\ac_{\ell} \ := \  \left\{
\aligned
& \. 1 \quad \text{ for } \ \ell \neq k+1\ts, \\
& \. 1\, + \, \frac{1}{\aP(k-1)-1} \ \quad \text{ for } \ \ \ell = k+1\ts.
\endaligned\right.
\end{equation}
Since every permutation of an independent set gives rise to a feasible word, we then have:
\begin{align*}
	& \aL_\aq(k-1) \ = \ (k-1)! \, \cdot \,   \aI_\ap(k-1),  \qquad
	\aL_\aq(k) \ = \ k!  \, \cdot \, \aI_\ap(k), \quad \text{and}\\
	& \qquad \aL_\aq(k+1) \ = \ (k+1)! \, \bigg(1+\frac{1}{\aP(k-1)-1}\bigg) \,\cdot \,  \aI_\ap(k+1)\ts.
\end{align*}
This reduces \eqref{eq:matroid-Par-weighted} to~\eqref{eq:greedoid-LC}.

\smallskip

By Theorem~\ref{t:greedoid}, it remains to show that $\aq$ is a $k$-admissible weight function.
First note that the weight function \ts $\aq$ \ts is multiplicative and thus
satisfies \eqref{eq:ContInv} and \eqref{eq:LogMod}.  By Proposition~\ref{p:reduction-matroid-LI},
greedoid \ts $\Gf$ \ts satisfies~\eqref{eq:LI},  which in turn implies~\eqref{eq:PasAct}.
By the same proposition, greedoid \ts $\Gf$ \ts is interval and satisfies~\eqref{eq:FewDes}.
Further, property~\eqref{eq:matroid basic} implies that \ts $\Des_{\alpha}(x) =\varnothing$ \ts
for every $\alpha \in \Lc$ and $x \in X$, which in turn trivially implies \eqref{eq:SynMon}.

To verify \eqref{eq:ScaleMon}, first suppose that \ts $\ell <k-1$.
Then \. $\ac_{\ell} =\ac_{\ell+1}=\ac_{\ell+2}=1$, which implies that
the LHS of \eqref{eq:ScaleMon} is equal to~0 while the RHS of \eqref{eq:ScaleMon}
is equal to~1, as desired. 
Now suppose that \ts $\ell=k-1$.
Note that \ts $\ab_{\alpha}(\Cc)=0$ \ts for every $\alpha \in \Lc$ \ts and \ts $\Cc \in \Par(\alpha)$, since \ts $\Des_{\alpha}(x)=\varnothing$.
Then, for every \ts $\alpha \in \Lc$ \ts of length \ts $k-1$, we have:
\begin{align*}
\bigg(1\.- \. \frac{\ac_{k}^2}{\ac_{k-1} \. \ac_{k+1}} \bigg)  \sum_{\Cc \in \Par(\alpha)}  \frac{1}{1-\ab_\alpha(\Cc)} \ = \ \bigg(1\.- \. \frac{\ac_{k}^2}{\ac_{k-1} \. \ac_{k+1}} \bigg) \, \big|\Par_\alpha \big|  \ 
= \ \frac{\big|\Par_\alpha\big|}{\aP(k-1)} \ \leq  \ 1.
\end{align*}
This finishes the proof of \eqref{eq:ScaleMon}.

In summary, greedoid \ts $\Gf=(X,\Lc)$ \ts satisfies
\eqref{eq:ContInv}, \. \eqref{eq:PasAct}, \. \eqref{eq:LogMod},
\. \eqref{eq:FewDes}, \. \eqref{eq:SynMon} and \eqref{eq:ScaleMon}.
By Definition~\ref{d:admissible},
we conclude that weight function \ts $\aq$ \ts is $k$-admissible,
which completes the proof of the theorem. \qed

\medskip

\subsection{Proof of Theorem~\ref{t:matroids-equality-Par}}\label{ss:proof-matroid-eq-Par}
We will prove the theorem as a consequence of Theorem~\ref{t:greedoid-equality-full}.
From the proof of Theorem~\ref{t:matroids-Par-weighted} given above, it suffices
to show that \eqref{eq:GrEqu1} and \eqref{eq:GrEqu2} are equivalent to
\eqref{eq:ME1} and \eqref{eq:ME2} for the greedoid~$\Gf$.
	
Let \ts $\alpha \in \Lc$ \ts of length \. $|\al|=k-1$.
We denote by $\as_{\Mf}(k-1)$ the constant that appears in \eqref{eq:ME2},
and $\as_{\Gf}(k-1)$ the constant that appears in \eqref{eq:GrEqu2}.
Recall that \ts
$\ab_{\alpha}(\Cc) \ = \ 0$ \ts for every $\alpha \in \Lc$ \ts and \ts $\Cc \in \Par(\alpha)$.
Note that
	\[  \sum_{x \ts\in\ts \Cc} \frac{\aq(\alpha x)}{\aq(\alpha)} \ =  \ \sum_{x \ts\in\ts \Cc} \. \ap(x)
\quad \text{and}\quad  1 \. - \. \frac{\ac_{k}^2}{\ac_{k-1} \.\ac_{k+1}}  \ = \  \frac{1}{\aP(k-1)}\,,
\]
	where the first equality follows from the product formula~\eqref{eq:matroid-proof-q-product}
	and the second equality is because of  the choice of constants \ts $\ac_\ell$ \ts in~\eqref{eq:matroid-proof-c-ell}.
	This implies that \eqref{eq:GrEqu2} and \eqref{eq:ME2} are equivalent under the substitution \. $\as_{\Mf}(k-1) \ts := \ts \as_{\Gf}(k-1)/\aP(k-1)$.
	
	Now, let \. $S=\{x_1,\ldots, x_{k-1}\}$ \. be an arbitrary independent set of size \ts $(k-1)$,
	and let \ts $\alpha:=x_1\ts\cdots \ts x_{k-1}$\ts.
	We have:
	\begin{align*}
		\sum_{x \ts\in\ts \Cnt(\alpha)}  \. \frac{\aq(\alpha x)}{\aq(\alpha)}  \ = \ \sum_{x \ts\in\ts \Cnt(\alpha)} \. \ap(x) \ =_{\eqref{eq:ME2}}  \    \big|\Par_S\big| \,\cdot \, \as_{\Mf}(k-1).
	\end{align*}
	This implies  that \eqref{eq:GrEqu1} and \eqref{eq:ME1} are equivalent, and
	completes the proof of the theorem. \qed

\medskip

\subsection{Proof of Theorem~\ref{t:matroids-equality-weighted}}\label{ss:proof-matroid-eq-weighted}
The direction \. $\Leftarrow$ \. is trivial, so it suffices
to prove the \. $\Rightarrow$ \. direction.

Let \ts $S$ \ts be an arbitrary independent set of size \ts $k-1$.
Recall that \ts $\aP(k-1)\le n-k+1$.
From the equality \eqref{eq:matroid-weighted-equality} and inequality
\eqref{eq:matroid-Par-weighted}, it follows that \. $\aP(k-1)=n-k+1$.
	On the other hand, it follows from equation~\eqref{eq:ME1} in Theorem~\ref{t:matroids-equality-Par}, that
	\. $|\Par_S|=\aP(k-1)$.
	Combining these two observations, we obtain:
	\begin{equation}\label{eq:hotel 3}
		S \cup \{x , y \} \  \text{ is an independent set for every distinct $x,y \in X \setminus S$.}
	\end{equation}

Let us show that every $(k+1)$-subset of $X$ is independent.
Fix an independent set \ts $U$ \ts of size $k+1$, and take an arbitrary
$(k+1)$-subset \ts $T$ \ts of~$X$.  If \ts $T=U$ \ts then we are done,
so suppose that \ts $T \neq U$.  Let \ts $x \in T \setminus U$ \ts
and let \ts $y \in U \setminus T$.  	
	Let $U'$ be the $(k+1)$-subset given by \. $U' \ts := \ts U +x-y$.
It follows from \eqref{eq:hotel 3} that \ts $U'$ \ts is an independent set.
Observe that the size of the intersection has increased: \ts $|T \cap U'| > |T \cap U|$,
Letting \ts $U\gets U'$, we can iterate this argument until we eventually get \ts
$U'=T$, as desired.
	
We can now prove that the weight function  \ts $\ap:X \to \rrs$ \ts is uniform.
Let \ts $x,y\in X$ \ts be  distinct  elements.  Let $S$ be a $(k-1)$-subset
of~$X$ that  contains neither $x$ nor~$y$.  It follows from the argument in the previous paragraph that $S$
is an independent set of the matroid~$\Mf$, and every parallel class
of~$S$ has cardinality~1.  By applying \eqref{eq:ME2} to the parallel
class \ts $\Cc_1=\{x\}$ \ts and \ts $\Cc_2=\{y\}$,
we conclude that \ts $\ap(x)=\ap(y)$.  This completes the proof.
\qed

\medskip

\subsection{Proof of Proposition~\ref{p:graphical-equality}}\label{ss:proof-matroid-prop-graph}
The inequality~\eqref{eq:prop-graphical-equality} in the proposition is a restatement
of~\eqref{eq:matroids-graphical-ex-LC}.  Thus, we need to show that equality in
\eqref{eq:matroids-graphical-ex-LC}  holds
if and only if $G$ is an $\sn$-cycle.  The \. $\Leftarrow$ \.
direction follows from a direct calculation, so it suffices
to prove the \. $\Rightarrow$ \. direction.

We first show that \ts $\deg(v)\geq 2$ \ts for every $v \in V$.
Suppose to the contrary, that there exists $v \in V$ such that $\deg(v)=1$.
Let $e$ be the unique edge adjacent to~$v$, and let  $S\ssu E$ be a forest with \ts
$\sn-3$ \ts edges not containing~$e$.  Then $v$ is a leaf vertex in the
contraction graph~$G/S$.  On the other hand, the graph \ts $G/S$ \ts
is the complete graph $K_3$ by \eqref{eq:ME1}, a contradiction.

We now show that \ts $\deg(v) \leq 2$ \ts for every $v \in V$.
Suppose to the contrary, that \ts $\deg(v) \geq 3$ \ts for some \ts $v \in V$.
Let \ts $e,f,g\in E$ \ts be three distinct edges adjacent  to~$v$.
Then there exists a spanning tree $T$ in $G$ that contains \ts $e,f,g$.
Let \. $S= T -\{e,f\}$, and let $x$ and $y$ be the other endpoint of $e$ and $f$, respectively.
Note that $S$ is a forest with \ts $\sn-3$ \ts edges, so it follows from \eqref{eq:ME2} that
there exists \ts $\as(\sn-3)$ \ts many edges connecting the component of $G/S$ containing~$x$,
and the component of $G/S$ containing~$y$.  Now let \ts $U:=T-\{f,g\}= S + e -g$,
which is another forest with \ts $\sn-3$ \ts edges.  Note that there are at least \ts
$\as(\sn-3)+1$ \ts edges connecting the component of $G/U$ containing $\{v,x\}$,
and the component of $G/S$ containing~$y$, namely the edge $f$ and the other \ts
$\as(\sn-3)$ \ts many edges connecting the component of $G/S$ containing $x$ and
the component of $G/S$ containing~$y$. This contradicts~\eqref{eq:ME2}.

Finally, observe that the $\sn$-cycle is  the only connected graph for which
every vertex has degree two.  This completes the proof.  \ $\sq$

\medskip

\subsection{More graphical matroids}\label{ss:proof-matroid-more-graph}
The following result is a counterpart to the Proposition~\ref{p:graphical-equality} proved above.

\smallskip

\begin{thm}\label{t:graphical-strict}
Let \ts $G=(V,E)$ \ts be a simple connected graph on \ts $|V|=\sn$ \ts vertices, and let
\ts $\aIr(k)$ \ts be the number of spanning forests with $k$ edges.  Then
\begin{equation}\label{eq:graphical-strict}
	\frac{\aIr(k)^2}{\aIr(k+1) \. \cdot \. \aIr(k-1)}  \  \geq  \
\bigg(1\. + \. \frac{1}{k} \bigg) \bigg(1\. + \. \frac{1}{\binom{\sn-k+1}{2} -1} \bigg),
\end{equation}
and the inequality is always strict if \. $1< k < \sn-2$.
\end{thm}

\begin{proof}
The inequality~\eqref{eq:graphical-strict} follows immediately from \eqref{eq:matroid-Par}
and the fact that \ts $\aP(k-1) \leq \binom{\sn-k+1}{2}$.

For the second part, suppose to the contrary, that we have equality in~\eqref{eq:graphical-strict}
for some simple connected graph~$G$.
	It then follows from \eqref{eq:ME1} and \eqref{eq:ME2}, that there exists $\as>0$ such that $G$ satisfies the following
\defng{clique-partition property}:

\smallskip

\begin{center}\begin{minipage}{13cm}%
{Let \. $A_1,\ldots, A_{\sn-k+1}$ \. be a partition of $V$, such that
each $A_i\ssu V$ spans a connected subgraph of $G$.  Then the graph
obtained by contracting each $A_i$ to one vertex (loops are removed but multiple edges remain)
is the complete graph \ts $K_{\sn-k+1}$\ts, with the multiplicity of every edge equal to~$\as$.
}
\end{minipage}\end{center}

\smallskip

\nin	
Now, start with an arbitrary partition \. $A_1,\ldots, A_{\sn-k+1}$ \. of $V$ such that
each \ts $A_i\ssu V$ \ts is nonempty and spans a connected subgraph of~$G$.  We get our contradiction
if this partition does not satisfy the clique-partition property above.  Since $k>1$, without loss of
generality, we assume that $A_1$ has at least two vertices.

Let $x$ be a vertex in $A_1$ that is adjacent to a vertex in $A_2$.
If $x$ is adjacent to any other vertex in $A_i$, $i\ge 3$, then by moving $x$ to
$A_2$ we create a new  partition $A_1',\ldots, A_{\sn-k+1}'$ and note that now
there are $\as+1$  edges connecting $A_2'$ and $A_i'$, contradicting
the clique-partition property.  Thus, $x$ is not adjacent to \. $A_3,\ldots, A_{\sn-k+1}$,
and we  can then move  $x$ to $A_2$ to create a new partition.
By iteratively moving elements to $A_2$ until only one element $y$ remains in $A_1$,
and applying the clique-partition property to the resulting partition, we conclude that~$y$
is adjacent to $A_3, \ldots, A_{\sn-k+1}$.

We now return to the original partition \. $A_1,\ldots, A_{\sn -k+1}$,
and we move $y$ to $A_3$  to obtain a new partition \. $A_1'',\ldots A_{\sn-k+1}''$.
In this new partition, there are $\as+1$ edges connecting \ts $A_3''$ \ts and \ts $A_4''$,
a contradiction. Note that here that part \ts $A_4''$ \ts is nonempty since \ts $k <\sn-2$.
This completes the proof.
\end{proof}

\smallskip

\begin{rem}\label{r:proof-matroid-strict}
The inequality~\eqref{eq:graphical-strict} is incomparable with~\eqref{eq:matroid-ULC},
and is stronger only for very dense graphs:
\[
|E| \ \geq \ \binom{\sn-k+1}{2} \. + \. k \. - \. 1.
\]
To explain this, note that \ts $\aP(k-1)$ \ts is usually smaller than the binomial coefficient above.
This is why the inequality~\eqref{eq:graphical-strict} is strict for \. $1< k < \sn-2$.  This also
underscores the power of our main matroid inequality~\eqref{eq:matroid-Par}.
\end{rem}

\medskip

\subsection{Proof of Corollary~\ref{c:intro-field}}\label{ss:proof-matroid-cor-field}
The result follows immediately from Theorem~\ref{t:matroids-Par} and the fact that
$$\Par(S) \. \leq  \. q^{m-k+1}-1 \quad\text{for every \ \  $S \in \Ic_{k-1}$\ts.}
$$
This is because  the contraction \ts $\Mf/S$ \ts with parallel elements
removed is a realizable matroid over \ts $\Fb_{\qnq}$ \ts of rank~$m-k+1$,
which can have at most \. $q^{m-k+1}-1$ \. nonzero vectors. \qed
\medskip

\subsection{Examples of combinatorial atlases}  \label{ss:proof-matroid-ex}
The numbers of independent sets can grow rather large, so we give two rather
small matroid examples to help the reader navigate the definitions.
We assume that the weight function \ts $\ap$ \ts is uniform in
both examples.

\smallskip

\begin{ex}[{\rm Free matroid}]  \label{ex:free-atlas}
Let \ts $\Mf=(X,\Ic)$ \ts be a \emph{free matroid} on $n=4$ elements: \. $X=\{x_1,\ldots x_4\}$ \.
and \. $\Ic = 2^X$.  In this case, we have \ts $\rI(1)=\rI(3)=4$, \ts $\rI(2)=6$, and the
inequality \eqref{eq:matroid-ULC} is an equality.

Following $\S$\ref{ss:reduction-matroids} and the proof of Theorem~\ref{t:matroids-Par-weighted},
the corresponding greedoid \ts $\Gf=(X,\Lc)$ \ts has all simple words
in~$X^\ast$.  Let \ts $k=2$ and \ts $\alpha=\varnothing$.
Then \ts $\bA(\alpha,k-1)$ \ts and \ts $\bA(\alpha,k)$ \ts are \ts $(n+1)\times (n+1)$ \ts
matrices given by
{\small
\begin{align*}
	\bA(\emp,1) \ = \
	\begin{pmatrix}
		0 & 1 & 1 & 1 & 1 \\
		1 & 0 & 1 & 1 & 1 \\
		1 & 1 & 0 & 1 & 1 \\
		1 & 1 & 1 & 0 & 1 \\
		1 & 1 & 1 & 1 & 1 	
	\end{pmatrix} \qquad \text{and} \qquad \bA(\emp,2) \ = \
	\begin{pmatrix}
	0 & 3 & 3 & 3 & 3\\
	3 & 0 & 3 & 3 & 3\\
	3 & 3 & 0 & 3 & 3\\
	3 & 3 & 3 & 0 & 3\\
	3 & 3 & 3 & 3 & 4\\		
\end{pmatrix},
\end{align*}}
where the rows and columns are labeled by \. $\{x_1,x_2,x_3, x_4, \spstar\}$.
Recall that  each entry of the matrices is counting the number of certain
feasible words, and only words of length \ts $k+1=3$ \ts are weighted by \.
$1+\frac{1}{\aP(k-1)-1}=\frac{3}{2}$\ts.

As in the proof of Theorem~\ref{t:greedoid} (see~$\S$\ref{ss:greedoid-proof}),
let \. $\vb, \wb\in \rr^5$ \. be the vectors given by
\[ \vb \, := \, (1,1,1,1,0)^\intercal  \quad \text{and} \quad \wb \, := \, (0,0,0,0,1)^\intercal.
\]
Inequality~\eqref{eq:matroid-ULC} in this case is equivalent to~\eqref{eq:greedoid-LC},
which in turn can be rewritten as:
\begin{equation}\label{eq:hyp specific}
	  \langle \vb, \bA(\alpha,k) \wb \rangle^2   \  \geq \  \langle \vb, \bA(\alpha,k) \vb \rangle \ \langle \wb, \bA(\alpha,k) \wb \rangle\ts.
\end{equation}
In this case the equality holds, since
\[ \langle \vb, \bA(\emp,2) \wb \rangle \, = \, 12, \qquad \langle \vb, \bA(\emp,2) \vb \rangle \, = \, 36, \qquad  \langle \vb, \bA(\emp,2) \vb \rangle \, =  \, 4,\]
as implied by Theorem~\ref{t:matroids-equality}.
\end{ex}

\smallskip

\begin{ex}[{\rm Graphical matroid}]  \label{ex:graphical-atlas}
Let \ts $G=(V,E)$ \ts be a graph as in the figure below, where
\. $\sn:=|V|=4$ \. and \. $E=\{a,b,c,d,e\}$.  Let \ts $\Mf=(E,\cL)$ \ts
be the corresponding graphical matroid (see Example~\ref{ex:intro-graphical}).
In this case \ts $n=|E|=5$ \ts and \ts $\rk(\Mf) = \sn-1=3$.

\begin{figure}[hbt]
	\centering
	\includegraphics[width = 0.12\linewidth]{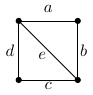}
\end{figure}

Let \ts $\al=\emp$ \ts and \ts $k=2$.
Then \ts $\bA(\alpha,k-1)$ \ts and \ts $\bA(\alpha,k)$ \ts are \ts $(n+1)\times (n+1)$ \ts
matrices given by
{\small
\begin{align*}
	\bA(\emp,1) \ = \
	\begin{pmatrix}
		0 & 1 & 1 & 1 & 1 & 1\\
		1 & 0 & 1 & 1 & 1 & 1\\
		1 & 1 & 0 & 1 & 1 & 1\\
		1 & 1 & 1 & 0 & 1 & 1\\
		1 & 1 & 1 & 1 & 0 & 1\\		
		1 & 1 & 1 & 1 & 1 & 1\\				
	\end{pmatrix} \qquad \text{and} \qquad \bA(\emp,2) \ = \
\begin{pmatrix}
	0 & 3 & 4.5  & 4.5  & 3 & 4\\
	3 & 0 & 4.5  & 4.5  & 3 & 4\\
	4.5  &  4.5  & 0& 3 & 3 & 4\\
	4.5  &  4.5  &  3 & 0 & 3 & 4\\
	3 &  3 &  3 &  3 & 0 & 4\\
	4 & 4 & 4 & 4 &4 & 5
\end{pmatrix},
\end{align*}}
where the rows and columns are labeled by \. $\{a, b, c, d, e,  \spstar\}$.
As in the previous example, each entry of the matrices is counting the number
of certain feasible words, and only words of length $k+1=3$ is weighted by \.
$1+\frac{1}{\aP(k-1)-1}=\frac{3}{2}$\ts.

As above,
let \. $\vb, \wb\in \rr^6$ \. be the vectors given by
\[ \vb \, := \, (1,1,1,1,1,0)^\intercal  \quad \text{and} \quad \wb \, := \, (0,0,0,0,0,1)^\intercal.
\]
Inequality~\eqref{eq:matroid-ULC} in this case is equivalent to~\eqref{eq:greedoid-LC},
which in turn can be rewritten as~\eqref{eq:hyp specific}.  Note that in this case we have
\[
\langle \vb, \bA(\alpha,k) \wb \rangle \ = \ 72, \qquad \langle \vb, \bA(\alpha,k) \vb \rangle \ = \ 20, \qquad
\langle \vb, \bA(\alpha,k) \vb \rangle \ =  \ 5,\]
and indeed we have a strict inequality \. $20^2 > 72 \times 5$, as implied
by Theorem~\ref{t:matroids-equality}.
\end{ex}

\bigskip


\section{Proof of discrete polymatroid inequalities and equality conditions}
\label{s:proof-polymatroid}

\smallskip

In this section we give proofs of Theorem~\ref{t:polymatroids-Par}, Theorem~\ref{t:polymatroids-BH-equality} and
Theorem~\ref{t:polymatroids-Par-equality}.

\smallskip

\subsection{Proof of Theorem~\ref{t:polymatroids-Par}}
We deduce the result from Theorem~\ref{t:greedoid}.
This proof is similar to the argument in the proof of Theorem~\ref{t:matroids-Par-weighted}
in the previous section, so we will emphasize the differences.

Let \. $\Gf=(X,\Lc)$ \. be the interval greedoid
constructed in~$\S$\ref{ss:reduction-polytroid}, and
corresponding to discrete polymatroid \. $\Df=([n],\Jc)$.
Let \. $1\le k < \rk(\Df)$, let \. $0 < t \le 1$, \. and let \. $\ap:X \to \rrs$ \. be
\. as in the theorem.

Let \ts $\ba_\al=(\aa_1,\ldots,\aa_n) \in \nn^n$ \ts be the vector corresponding
to the word \ts $\al\in \Lc$.
We define the weight function  \ts $\aq:\Lc \to \rrs$ \ts by the product formula
\begin{equation*}
	\aq(\alpha) \ := \ \ac_{\ell} \,\. t^{\pi(\ba_\al)} \. \ap(\ba_\al),
\end{equation*}
where \ts $\ell:=|\al|=|\ab_\al|$, and \ts $\ac_\ell$ \ts is given by
\begin{equation}\label{eq:polymatroid-proof-c-ell}
\ac_{\ell} \ := \  \left\{
\aligned
& \. 1 \quad \text{ for } \ \ell \neq k+1\ts, \\
& \. 1\, + \, \frac{1-\rnr}{\aP(k-1)\ts -\ts 1 \ts +\ts \rnr} \ \quad \text{ for } \ \ \ell = k+1\ts.
\endaligned\right.
\end{equation}
Using this weight function, we obtain:
\begin{align*}
	\aL_\aq(k) \ &= \  \sum_{\alpha \in \Lc_k} \, t^{\pi(\ba_\al)} \. \ap(\ba_\alpha)  \
	= \ \sum_{\ab \in \Jcm_k } \,  t^{\pi(\ba_\al)} \.\ap(\ba) \,\.  \frac{k!}{\aa_1! \ts\cdots\ts \aa_n!} \ = \ k! \,\cdot\, \aJ_{\ap,t}(k),
\end{align*}
where the third equality follows from every permutation of a feasible word that
is well-ordered is again a feasible word.  By the same calculation, we have
\begin{align*}
	\aL_\aq(k-1) \ = \ (k-1)! \,\cdot \,  \aJ_{\ap,t}(k-1) \ ,   \quad
	\aL(k+1) \ = \ (k+1)! \, \bigg(1+\frac{1-\rnr}{\aP(k-1)-1+\rnr}\bigg) \, \cdot \, \aJ_{\ap,t}(k+1).
\end{align*}
This reduces \eqref{eq:polymatroid-Par} to~\eqref{eq:greedoid-LC}.

\smallskip

By Theorem~\ref{t:greedoid}, it remains to show that $\aq$ is a $k$-admissible weight function.
First note that the weight function \ts $\aq$ \ts is multiplicative and thus
satisfies \eqref{eq:ContInv} and \eqref{eq:LogMod}.  By Proposition~\ref{p:reduction-polymatroid-LI},
greedoid \ts $\Gf$ \ts satisfies~\eqref{eq:LI},  which in turn implies~\eqref{eq:PasAct}.
By the same proposition, greedoid \ts $\Gf$ \ts is interval and satisfies~\eqref{eq:FewDes}.

We now verify \eqref{eq:SynMon}, which is no longer similar to the matroid case.
It follows from \eqref{eq:polymatroid basic} that the right side of \eqref{eq:SynMon} is~$0$,
unless \ts $x=x_{ij}$ \ts and \ts $\Des_{\alpha}(x)=x_{i \ts j+1}$.  In the latter case, we have
\begin{equation}\label{eq:polymatroid Sinmon}
	 \frac{\ap(\alpha x)}{\ap(\alpha)}   \  = \ \rnr^{j-1} \, \ap(i) \, , \quad  \sum_{y \in \Des_{\alpha}(x)} \. \frac{\ap(\alpha xy)}{\ap(\alpha x)}  \ = \ \rnr^{j} \. \ap(i),
\end{equation}
and \eqref{eq:SynMon} follows from \ts $\rnr \leq 1$.

For \eqref{eq:ScaleMon}, the same argument as for matroids works for \ts $\ell <k-1$.
Now suppose that $\ell=k-1$.  Note that, for every $\alpha \in \Lc$ and $\Cc \in \Par(\alpha)$,
we have \. $\ab_{\alpha}(\Cc) \leq \rnr$ \.
by \eqref{eq:polymatroid Sinmon}.
Hence, for every $\alpha \in \Lc_{k-1}$\ts, we have:
\begin{align*}
	& \bigg(1\.- \. \frac{\ac_{k}^2}{\ac_{k-1} \. \ac_{k+1}} \bigg)  \sum_{\Cc \in \Par(\alpha)}  \, \frac{1}{1-\ab_\alpha(\Cc)} \ \leq  \
		\bigg(1\.- \. \frac{\ac_{k}^2}{\ac_{k-1} \. \ac_{k+1}} \bigg) \, \sum_{\Cc \in \Par(\alpha)}  \frac{1}{1-\rnr}
	\\
	\quad  &= \  \bigg(1\.- \. \frac{\ac_{k}^2}{\ac_{k-1} \. \ac_{k+1}} \bigg)  \, \frac{\big|\Par_\alpha \big|}{1-\rnr}  \ = \ \frac{\big|\Par_\alpha\big|}{\aP(k-1)} \ \leq  \ 1\,,
\end{align*}
which proves \eqref{eq:ScaleMon}.

In summary, greedoid \ts $\Gf=(X,\Lc)$ \ts satisfies
\eqref{eq:ContInv}, \. \eqref{eq:PasAct}, \. \eqref{eq:LogMod},
\. \eqref{eq:FewDes}, \. \eqref{eq:SynMon} and \eqref{eq:ScaleMon}.
By Definition~\ref{d:admissible},
we conclude that weight function \ts $\aq$ \ts is $k$-admissible,
which completes the proof of the theorem. \qed

\medskip

\subsection{Proof of Theorem~\ref{t:polymatroids-BH-equality}}
\label{ss:proof-polymatroid-eq-BH}
We deduce the result from Theorem~\ref{t:greedoid-equality-full}.
The proof below only assumes  that \ts $0< t\le 1$.
For the \. $\Rightarrow$ \. direction,
let \ts $\alpha \in \Lc$ \ts with \ts $|\alpha | =k-1$.
Note that, since $\ac_k=\ac_{k-1}=1$, we have
\[ \sum_{x \in \Cc} \frac{\aq(\alpha x)}{\aq(\alpha)} \ = \  \sum_{x \in \Cc} \frac{\ap(\alpha x)}{\ap(\alpha)} \ = \ \aa_\alpha(\Cc). \]
By \eqref{eq:GrEqu2}, there exists \ts $\as>0$, s.t.\ for every \ts
$\Cc \in \Par(\alpha)$ \ts we have:
\begin{align}\label{eq:oscar 1}
	\aa_\alpha(\Cc) \, \big(1-\ab_{\alpha}(\Cc) \big) \,
	\ &= \ \as \, \bigg( 1 \. - \. \frac{\ac_{k}^2}{\ac_{k-1} \.\ac_{k+1}} \bigg) \ = \  \as \, \frac{1-\rnr}{\aP(k-1)}\,.
\end{align}
Summing over all $\Cc \in \Par(\alpha)$, we get:
\begin{align*}
	   \sum_{\Cc \in \Par_\alpha} \aa_{\alpha}(\Cc) \, \big(1-\ab_{\alpha}(\Cc) \big)
	\ = \  \as \, \. (1-\rnr) \, \frac{\big|\Par_\alpha \big|}{\aP(k-1)}\,.
\end{align*}
On the other hand, the equality \eqref{eq:GrEqu1} gives:
$$\as \, = \, \sum_{\Cc \in \Par_\alpha} \. \aa_{\alpha}(\Cc).
$$
Combining these equations, we obtain:
\begin{align*}
	 \sum_{\Cc \in \Par_\alpha} \aa_{\alpha}(\Cc) \, \big(1-\ab_{\alpha}(\Cc) \big)
	\ = \  \sum_{\Cc \in \Par_\alpha}  \aa_\alpha(\Cc) \,  (1-\rnr) \, \frac{\big|\Par_\alpha \big|}{\aP(k-1)}\,,
\end{align*}
which is equivalent to
\begin{align*}
	\sum_{\Cc \in \Par_\alpha} \aa_{\alpha}(\Cc) \, \bigg(1 \. - \. \ab_{\alpha}(\Cc) \. -  \.  (1-\rnr) \. \frac{\big|\Par_\alpha \big|}{\aP(k-1)} \bigg)
	\ = \  0.
\end{align*}
Now note that, the LHS of the equation above is always nonnegative since  \. $\ab_{\alpha}(\Cc) \leq \rnr$ \. by \eqref{eq:polymatroid Sinmon}, and \.
$\big|\Par_\alpha \big| \leq \aP(k-1)$ \. by definition.
Therefore, the equality hold for both inequalities, so in particular we have:
\begin{align*}
	   \ab_{\alpha}(\Cc) \. = \. \rnr \quad \text{ for every  \, $\alpha \in \Lc_{k-1}$ \, and \, $\Cc \in \Par(\alpha)$}.
\end{align*}
Since \ts $\rnr >0$ \ts by assumption,
it follows from \eqref{eq:FewDes} that
\begin{align}\label{eq:oscar 3}
	   |\Cc| \. = \. 1 \  \text{ and } \  \ab_{\alpha}(\Cc) =\rnr >0  \quad \text{ for every  \, $\alpha \in \Lc$ \, with \, $|\alpha|=k-1$ \, and \, $\Cc \in \Par(\alpha)$}.
\end{align}
Restating this equation in the language of polymatroids, we conclude: \. for every \ts $\ba \in \Jcm_{k-1}$,
and every  \ts $i,j \in [n]$ \ts (not necessarily distinct), we have:
	\begin{equation}\label{eq:oscar 4}
	\text{$\ba+\eb_i, \ba+\eb_j \in \Jc$ \quad $\Longrightarrow$ \quad $\ba+\eb_i+\eb_j\in \Jc$.}	
\end{equation}

\smallskip

We can now show that every \ts $\nb=(\an_1,\ldots,\an_n) \in \nn^{n}$
with \ts $|\nb|=k+1$ \ts is contained in~$\Jc$.  We follow the corresponding
argument int the matroid case.  Let \ts $\ba \in \Jc$ \ts with \ts $|\ba|=k+1$.
If $\ba=\nb$,  we are done, so suppose that \ts $\ba \neq \nb$.
Then there exists \ts $i,j \in [n]$, such that \ts
$\aa_i > \an_i$ \ts and \ts $\aa_j < \an_j$.
By the polymatroid hereditary property, we have \ts $\ba-\eb_i \in \Jc$.
Since \ts $\eb_j \in \Jc$ \ts by the assumption that the polymatroid is normal,
we can then apply the exchange property to $\eb_j$ and \ts $\ba-\eb_i$ \ts
to conclude that \. $\ba-\eb_i+ \eb_j-\eb_h \in \Jc$  \. for some \ts $h \in [n]$.
Let \ts $\ub:=\ba -\eb_i-\eb_h$.
Note that \ts $\ub \in \Jcm_{k-1}$ by hereditary property,
and
\[ \text{$\ub+\eb_j \. = \. \ba-\eb_i+\eb_j-\eb_h \ts\in\ts \Jc$, \quad and \quad $ \ub+ \eb_h \. = \. \ba-\eb_i \in \Jc$.}   \]
It then follows from \eqref{eq:oscar 4} that
\[  \ba -\eb_i+ \eb_j \, = \,  \ub+\eb_j+\eb_h \in \Jc.  \]
Make substitution
$\ba \gets \ba-\eb_i+\eb_j$ \. and \. iterate this argument until
eventually \ts $\nb=\ba$, as desired.
This proves the \. $\Rightarrow$ \. direction.

\smallskip

For the \. $\Leftarrow$ \. direction, assume now that \ts $\rnr=1$.
The equality now follows from a direct calculation, since
 that
 \begin{align*}
 1 \. + \. \frac{1-\rnr}{\aP(k-1)-1 +\rnr} \, = \, 1, \quad \text{and} \quad	\aJ_\ap(\ell) \, = \, \frac{\big( \ap(1)\ts + \.\ldots \.+ \ts \ap(n) \big)^\ell}{\ell !} \quad \text{for every \. $\ell \leq k+1$}.
 \end{align*}
This completes the proof.  \qed

\medskip

\subsection{Proof of Theorem~\ref{t:polymatroids-Par-equality}}
\label{ss:proof-polymatroid-eq-Par}
Assume now that \ts $0<\rnr<1$.
From the proof above, it remains to show that \ts $k=1$ \ts and that the
weight function  \ts $\ap$ \ts is uniform.

Let \ts $i,j\in [n]$ \ts  be  distinct  elements,
let  \. $\alpha:=x_{i \ts 1} \cdots x_{i \ts k-1}$, let \.
$x:= x_{i\ts k}$, and let \. $y:=x_{j \ts 1}$.
By  \eqref{eq:oscar 3}, we have \ts $\Cc_1=\{x\}$ \ts and \ts $\Cc_2=\{y\}$ \ts
are both parallel classes of~$\alpha$.
 It then follows from \eqref{eq:oscar 1} and \eqref{eq:oscar 3}, that
 \begin{align*}
 		\aa_{\alpha}(\Cc_1)
 	\ &= \ \aa_{\alpha}(\Cc_2).
 \end{align*}
 On the other hand, we have
 \[  \aa_{\alpha}(\Cc_1) \ = \ \rnr^{k} \, \ap(i), \qquad
 \aa_{\alpha}(\Cc_2) \ = \ \rnr \, \ap(j),
 \]
 so \ts $\rnr^{k-1}= \ap(j)/\ap(i)$.
Since the choice of $i$ and $j$ was arbitrary, we can switch $i$ and $j$ \ts to obtain
  \ts $\rnr^{k-1}= \ap(j)/\ap(i)$.  This implies that \ts $\ap(i) =\ap(j)$ \ts and \ts $k=1$,
which  proves the \. $\Rightarrow$ \. direction.

\smallskip

For the \. $\Leftarrow$ \. direction, assume now that \ts $k=1$.
From the proof above, $\ap(i)=C$ for every $i \in [n]$ and some~$C>0$.
It then follows from a direct calculation that
  $$
  1+ \frac{1-\rnr}{\aP(k-1)-1 +\rnr} \ = \ \frac{n}{n-1+\rnr}\,,
  $$
and
  \begin{align*}
   \aJ_{\ap,t}(0) \ = \ 1\,,  \qquad \aJ_{\ap,t}(1) \ = \  C\ts \rnr \ts n\,, \qquad \aJ_{\ap,t}(2) \ = \  C^2\left( \frac{\rnr^3 \, n}{2} \, + \, \frac{\rnr^2 \, n \, (n-1)}{2}\right).
  \end{align*}
Thus, the equality~\eqref{eq:polymatroid-Par-equality} holds in this.  This completes the proof.
 \qed

 \bigskip


 \section{Proof of poset antimatroid inequalities and equality conditions}
\label{s:proof-antimatroid}

\smallskip

In this section we prove  Theorem~\ref{t:antimatroids} and Theorem~\ref{t:antimatroids-equality}.

\smallskip

\subsection{Proof of Theorem~\ref{t:antimatroids}}
As in the previous sections, we deduce the result from Theorem~\ref{t:greedoid}.
Let \ts $\cP=(X,\prec)$ be a poset on \ts $|X|=n$ \ts elements and let
\ts $\Af=(X,\cL)$ \ts be the corresponding poset antimatroid which is an
interval greedoid by the argument in $\S$\ref{ss:reduction-antimatroid}.
In the notation of Section~\ref{s:main}, let \ts $c_\ell=1$ \ts for all \ts $\ell\ge 1$,
and let \ts $\aq(\al):=\ap(\al)$, so that \eqref{eq:antimatroid-LC} coincides with
\eqref{eq:greedoid-LC} in this case.

It remains to show that $\ap$ is a $k$-admissible weight function.
First note that $\ap$ satisfies \eqref{eq:ContInv} and \eqref{eq:LogMod}
since the weight function is multiplicative. The condition \eqref{eq:ScaleMon}
is also trivially satisfied.  By Proposition~\ref{p:reduction-antimatroid-LI},
both~\eqref{eq:LI} and \eqref{eq:FewDes} are satisfied, and the former implies
\eqref{eq:PasAct}.

For \eqref{eq:SynMon},
mote that the poset ideal greedoid $\Gf$ satisfies
\begin{equation}\label{eq:poset 2}
	\Des_{\alpha}(x)  \ \subseteq \ \bigl\{ \ts  y \in X \, : \,  x \.\lcov \.y \ts  \bigr\},
\end{equation}
for every \ts $\alpha \in \Lc$ \ts and every \ts $x \in \Cnt(\alpha)$.
It then follows that
\begin{align}\label{eq:poset 3}
	\frac{\ap(\alpha x)}{\ap(\alpha)} \ = \  \ap(x)  \ \geq_{\eqref{eq:CM}} \
\sum_{y \, : \, x \ts \lcov\ts y} \. \ap(y) \ \geq_{\eqref{eq:poset 2}} \
	\sum_{y \in \Des_{\alpha}(x)}    \ap(y) \ = \
	\sum_{y \in \Des_{\alpha}(x)}    \frac{\ap(\alpha xy)}{ \ap(\alpha x)},
\end{align}
which proves \eqref{eq:SynMon}.
Hence $\aq$ is indeed a $k$-admissible weight function, which completes the proof. \qed

\medskip

\subsection{Proof of Theorem~\ref{t:antimatroids-equality}}

We deduce the result from Theorem~\ref{t:greedoid-equality-full}.
From the argument above, it suffices to show that \eqref{eq:GrEqu1}
and \eqref{eq:GrEqu2} are equivalent to properties \eqref{eq:AE1}--\eqref{eq:AE3}.
First note that,
\[  \sum_{x \in \Cnt(\alpha)} \frac{\aq(\alpha x)}{\aq(\alpha)} \ = \ \sum_{x \in \Cnt(\alpha)} \ap(x),
\]
for every \ts $\alpha \in \Lc$.
This implies that \eqref{eq:GrEqu1} is equivalent to~\eqref{eq:AE1}.

Let \ts $\alpha \in \Lc$, let \ts $x \in \Cnt(\alpha)$, and
let \ts $\Cc = \{x\}$ \ts be the parallel class in \ts $\Par(\alpha)$ \ts
containing $x$.  Since \. $\ac_{k+1}=\ac_{k} = \ac_{k-1}=1$,
it then follows that the RHS of \eqref{eq:GrEqu2} is equal to~0,
so  \eqref{eq:GrEqu2} is equivalent to
\begin{equation*}
	  \sum_{y \in \Des_{\alpha}(x) } \. \frac{\ap(\alpha) \, \ap(\alpha xy)}{\ap(\alpha x)^2}
\ = \ \ab_{\alpha}(\Cc) \ = \  1.
\end{equation*}
This implies that \eqref{eq:GrEqu2} is equivalent to equality in \eqref{eq:poset 3},
which in turn is equivalent to equality in both \eqref{eq:CM} and \eqref{eq:poset 2}.
The latter is equivalent to \eqref{eq:AE2} and~\eqref{eq:AE3}, which completes the proof. \qed

\bigskip

\section{Proof of morphism of matroids inequalities and equality conditions} \label{s:proof-morphism}

In this section we give proofs of Theorem~\ref{t:morphisms-refined},
Theorem~\ref{t:morphisms-EH-equality} and Theorem~\ref{t:morphisms-Par-equality}.

\smallskip

\subsection{Combinatorial atlas construction}\label{ss:proof-morphism-construction}
Let \ts $\Mf=(X,\Ic)$ \ts and \ts $\Nf=(Y,\Jc)$ \ts be two matroids, and let
\ts $\Phi: \Mf \to \Nf$ \ts be a morphism of matroids.  Let \ts $1\le k < \rk(\Mf)$ \ts
and let \ts $\ap:X \to \rrs$ \ts be the weight function as in
Theorem~\ref{t:morphisms-refined}.  We now define a combinatorial atlas \ts
$\AA$ \ts that corresponds to \ts $(\Phi,k,\ap)$.

Let \ts $\Gf=(X,\Lc)$ \ts be the greedoid which corresponds to matroid~$\Mf$,
see~$\S$\ref{ss:reduction-matroids}.  We extend~$\ts\ap\ts$ to a \emph{nonnegative weight
function} \ts $\aq : \Lc_{\Mf} \to \Rb_{\geq 0}$ \ts by the product formula:
\begin{equation*}
	\aq(x_1 \ts \cdots \ts x_\ell) \ := \
	\begin{cases}
	\displaystyle \. \ac_{\ell} \, \ap(x_1) \. \cdots \.  \ap(x_\ell) \ & \text{ if } \ \{x_1,\ldots,x_\ell\} \in \Bc_\ell\.,\\
	\. 0 \ & \ \text{otherwise}\ts,
	\end{cases}
\end{equation*}
where \ts $c_\ell$ \ts is defined in~\eqref{eq:matroid-proof-c-ell}.
Let
\ts $\Qf:= (\Vf,\Ef)$ \ts be the acyclic graph and let \ts $\AA$ \ts be the
combinatorial atlas defined in \S\ref{ss:proof-matroid-Par-weighted} that
corresponds to the greedoid~$\Gf$.  Note that \ts $\Qf$ \ts depends only
on the matroid~$\Mf$, but \ts $\AA$ \ts depends also on the morphism~$\Phi$.
Note also that, unlike the weight function in \S\ref{ss:intro-greedoids} and~$\S$\ref{ss:proof-matroid-Par-weighted},
here  the weight $\aq$  is \emph{not} strictly positive, so Theorem~\ref{t:greedoid} does not apply in this case.

In this section, we rework the combined proofs of Theorem~\ref{t:greedoid} and
Theorem~\ref{t:matroids-Par-weighted} to apply for morphisms of matroids.
Recall the properties we need to establish as summarized in the key results:
$$
\text{Theorem~\ref{t:Hyp}:} \ \quad \left\{\,
\aligned
& \. \text{\eqref{eq:Inh}, \eqref{eq:Pull} hold for~$\AA$} \\
& \text{$\vf\in \Vfp$ \ts satisfies \eqref{eq:Irr}, \eqref{eq:hPos}}\\
& \text{ \eqref{eq:Hyp} \, holds for all \. $\vf^{\<i\>}\in \vfs$}
\endaligned \,
\right\}
 \quad  \text{$\Rightarrow$ \quad
 \eqref{eq:Hyp} holds for \ts $\vf$.}
$$
$$
\text{Theorem~\ref{t:Pull}: \ \ \ \. $\Bigl\{$\eqref{eq:Inh}, \eqref{eq:Proj}, \eqref{eq:TPInv}, \eqref{eq:Knonneg}$\Bigr\}$
\quad $\Rightarrow$ \quad \eqref{eq:Pull}.\hskip1.78cm}
$$

Now, observe that \eqref{eq:Inh}, \eqref{eq:Proj}, \eqref{eq:TPInv}, and \eqref{eq:Knonneg}
are closed properties, i.e.\ preserved under taking limits.  Thus, they follow from the
arguments in~$\S$\ref{ss:greedoid-prop}.  On the other hand, properties \eqref{eq:Irr}
and \eqref{eq:hPos} need to be verified separately, a the arguments in $\S$\ref{ss:greedoid-prop}
use the strict positivity of~$\aq$.

\smallskip

\begin{lemma}\label{l:morphism Irr}
Let \ts $\vf=(\alpha,m,t)\in \Vf^m$ \ts be a non-sink vertex of the
acyclic graph \ts $\Qf$ \ts defined above,
where \. $\alpha \in X^\ast$, \. $|\al| \le k-1-m$, \. $0< m \leq k-1$, and \ts $0<t<1$.
Then \ts $\vf$ \ts satisfies \eqref{eq:Irr} and  \eqref{eq:hPos}.
\end{lemma}

\smallskip

\begin{proof}
For the second part, it follows from the definition of \ts $\hb_{\vf}$ \ts
that the vector is strictly positive for all \ts $t \in (0,1)$.
Thus, vertex \ts $\vf$ \ts satisfies \eqref{eq:hPos},
for all \ts $t \in (0,1)$.

For the first part, let \ts $\bM_{\vf}$ \ts be the associated matrix
of~$\vf$.
Without loss of generality, we can assume that \ $\alpha \in \Lc$,
as otherwise \ts $\bM_{\vf}=0$ \ts and \eqref{eq:Irr} holds trivially.

\medskip

\noindent \textbf{Claim}: \.
{\em Every \. $x,y \in X$ \. in the support of \ts $\bMr_{\vf}$
belong to the same irreducible component of \ts $\bMr_{\vf}$.}

\medskip

When \ts $\spstar$ \ts is not in the support of~$\bM_{\vf}$,
property \eqref{eq:Irr} follows from the claim.  Now assume
that \ts $\spstar$ \ts is  in the support of~$\bM_{\vf}$\ts.
By the Claim, it remains to show that \ts $\spstar$ \ts belong to the same
irreducible component of \emph{some} \ts $x \in X$ \ts in the support
of $\bM_{\vf}$\ts.	

Let \. $\al = x_1\ts\cdots\ts x_\ell$, where \ts $\ell\le k-m-1$.  By the assumption,
there exits a subset \ts $S \in \Bc_\Phi$ \ts
such that \ts $\{x_1,\ldots, x_\ell\} \ssu S$ \ts  and \ts $|S| \in \{\ell+m-1,  \ts \ell+m, \ts \ell+m+1\}$.
To see this, observe that if  \ts $\spstar$ \ts is  in the support of~$\bM_{\vf}$,
then either \ts $\aA(\al,m)_{\spstar \ts \spstar}\ne 0$, or
\ts $\aA(\al,m)_{\spstar \ts x}\ne 0$, or
\ts $\aA(\al,m+1)_{\spstar \ts x}\ne 0$, for some \ts $x\in X$.

By adding extra elements to \ts $S$ \ts if necessary,
without loss of generality we can assume that \ts
$|S|=\ell+m+1$.
	Let $x \in S \setminus \{x_1,\ldots,x_\ell\}$. Then we have
\ts \ts $\aA(\al,m+1)_{\spstar \ts x}\ne 0$.  This implies that \ts
$x$ and $\spstar$ belong to the same irreducible component of  $\bA(\alpha,m+1)$.
Since \ts $0<t<1$, this implies that \ts
$x$ and $\spstar$ belong to the same irreducible component of $\bM_{\vf}$\ts ,
and completes the proof of the lemma.
\end{proof}

\begin{proof}[Proof of Claim]
Let \ts $\ell=|\alpha|$, as above.
Since $x$ is contained in the support of \ts $\bM_{\vf}$,
this implies that there exits \ts $S \in \Bc_\Phi$ \ts
such that \ts $\{x_1,\ldots, x_\ell,x\} \ssu S$ \ts  and
		\ts $|S| \in \{\ell+m, \ts \ell+m+1, \ts \ell+m+2\}$.
		Similarly, there exists \ts
		$T \in \Bc_\Phi$ \ts
		such that \ts $\{x_1,\ldots, x_\ell,y\} \ssu T$ \ts  and
		\ts $|T| \in\{\ell+m, \ts \ell+m+1, \ts \ell+m+2\}$.
		By adding extra elements to $S$ and $T$ if necessary,
		without loss of generality we can  assume that
		$|S|=|T|=\ell+m+2$.
		
For \ts $S=T$, the claim follows immediately from the definition of \ts $\bM_{\vf}$ \ts and \ts $\aq$,
since \. $\aA(\al,m+1)_{x\ts y}\ne 0$ \. in this case.   So assume that	\ts $S \neq T$.
By the exchange property for morphism of matroids (Proposition~\ref{p:reduction-morphism-base-exchange}),
there exists \. $z \in S \setminus T$ \. and \. $w \in T \setminus S$ \.
such that \. $S-z+w \ts \in \ts \Bc_\Phi$.

Let \ts $S':=S-z+w$.
Note that \. $|S' \setminus \{x_1,\ldots, x_\ell, x, w\}| = m \geq 1$,  and let \ts $x' \in S' \setminus \{x_1,\ldots, x_\ell, x, w\}$.
Note that \. $x' \in S \setminus \{x_1,\ldots, x_{\ell}\}$, which implies that
\. $\aA(\al,m+1)_{x\ts x'}\ne 0$ \. in this case. Therefore, elements $x$ and~$x'$ belongs to
the same irreducible component of \ts $\bA(\alpha,m+1)$,
and thus the same irreducible component of \ts $\bM_{\vf}$ \ts since \ts $0<t<1$.
Note also that we have \. $|S' \cap T| > |S \cap T|$ \. by the construction of~$S'$.
Substitute \ts $x \gets x'$ \ts and \ts $S \gets S'$, and iteratively apply the same argument, until
the set $S$ will eventually becomes~$T$.  This implies that $x$ and $y$ are contained in
the same irreducible component of \ts $\bM_{\vf}$\ts, as desired. 	\end{proof}

\medskip

\subsection{All atlas vertices are hyperbolic} \label{ss:proof-morphism-Hyp}
We first  show that every sink vertex in $\AA$ satisfies \eqref{eq:Hyp}.
We then use Theorem~\ref{t:Hyp} to obtain the result.

\smallskip

Let \ts $\Gf=(X,\Lc)$ \ts be the greedoid corresponding to matroid \ts $\Mf=(X,\Ic)$.
Let \ts $\alpha=x_1\.\cdots\. x_\ell \in \Lc$ \ts of length \ts $\ell:= |\alpha| \leq k-1$,
let \. $S:= \{x_1,\ldots, x_\ell\}$,
and let \ts $\bA(\alpha,1)$ \ts be the matrix defined in~$\S$\ref{ss:greedoid-def}
for~$\Gf$.  Recall that
$$\ap(S)  \ = \ \ap(x_1) \.\cdots\. \ap(x_\ell) \ = \ \frac{\aq(x_1 \ldots x_\ell)}{\ac_\ell}.
$$
For each \ts $x \in X$,
divide  the $x$-row and $x$-column of $\bA(\alpha, 1)$
by  \ts $\sqrt{\ac_{\ell+2} \, \.\ap(S)} \, \ap(x)$.
Multiply the \ts $\spstar$-row \ts and the \ts $\spstar$-column \ts by \.
$\frac{1}{\ac_{\ell+1}} \ts \sqrt{\frac{\ac_{\ell+2}}{\ap(S)}}$\ts.
Denote by \ts $\bB$ \ts the resulting matrix.
Note that \eqref{eq:Hyp} is preserved under this transformation,
so it suffices to show that \ts $\bB$ \ts satisfies \eqref{eq:Hyp}.
Observe that \. $\bB=\bigl(\aB_{xy}\bigr)_{\ts x,\. y \in \Xf}$ \. is given by
$$
\aligned
& \aB_{x \ts y} \ = \ \left\{ {\aligned \. 1 \quad & \text{ if } \ S+x+y \in \Bc \\
\. 0 \quad & \text{ if } \ S +x+y \notin \Bc \endaligned}\right.
\qquad \text{for distinct \, $x,y \in X$,}\\
& \aB_{x \ts \spstar} \ = \ \left\{ {\aligned \. 1 \quad & \text{ if } \ S +x \in \Bc\\
 0 \quad & \text{ if } \ S +x \notin \Bc \endaligned}\right.
 \qquad \text{for \, $x \in X$} \,, \\
& \aB_{\spstar \ts \spstar} \ = \ \left\{ {\aligned \.
 \frac{\ac_\ell \, \ac_{\ell+2}}{\ac_{\ell+1}^2} \quad &\text{ if } \ S  \in \Bc \\
0 \quad &\text{ if } \ S  \notin \Bc\endaligned}\right.\\
&\aB_{x \ts x} \ = \ 0 \qquad \text{for \, $x \in X$.}
\endaligned
$$

\smallskip

We now split the proof into three cases, each discussed as a separate lemma.

\smallskip

\begin{lemma}\label{l:morphism BHyp 1}
	Suppose that \.  $\rkn\big(\Phi(S)\big) = \rk(\Nf)$.
	Then the matrix \ts $\bBr$ \ts satisfies \eqref{eq:Hyp}.
\end{lemma}

\begin{proof}
By the assumption of the lemma, every independent set of $\Mf$ containing \.
$\{x_1,\ldots, x_\ell\}$ \. is also a basis of~$\Phi$.  It then follows that \.
$\bB=\big(\aB_{x y}\big)$ \. is equal to
	\begin{align*}
		\aB_{xy} \ = \ \aB_{yx} \ = \
		\begin{cases}
			1 \ & \text{ if }  \ x,y \in \Cnt(S)  \ \text{ and } \ x \not \sim_S y, \\
			1  \ & \text{ if } \  x \in \Cnt(S) \ \text{ and } \ y=\spstar,\\
			\frac{\ac_\ell \, \ac_{\ell+2}}{\ac_{\ell+1}^2} \ & \text{ if } \ x = y=\spstar,\\
			0 \ & \ \text{otherwise.}
		\end{cases}
	\end{align*}
	In particular,
	if $x \salp y$,
	then $x$-row and $x$-column of $\bB$
	is equal to $y$-row and $y$-column of $\bB$.
	Now, choose a representative element $x_i$ for each equivalence class $\Cc_i$ in $\Par(\alpha)$.
	For every other $y$ in $\Cc_i$,
	we subtract from the $y$-row and $y$-column of $\bA$
	the $x_i$-row and $x_i$-column of $\bA$, respectively.
Note that \eqref{eq:Hyp} is preserved under these transformations.
Restricting to the support, we obtain:	
{\small \[
\bN \ := \
	\begin{pmatrix}
		0 & 1  &\ldots & 1 & 1  \\
		1  & 0 & \ldots &1  & 1 \\
		\vdots & \vdots & \ddots & \vdots  & \vdots \\
		1 & 1 & \ldots & 0 & 1  \\
		1 & 1 & \ldots & 1  & \frac{\ac_{\ell} \ac_{\ell+2}}{\ac_{\ell+1}^2}
	\end{pmatrix},
\]}
	where the rows and columns are indexed by \. $\{x_1,\ldots, x_m, \spstar\}$, with \ts $m:=|\Par(\alpha)|$.
	Now note that
$$
\frac{\ac_{\ell} \ac_{\ell+2}}{\ac_{\ell+1}^2} \ = \ \left\{
\aligned
& \. 1   & \text{for} \ \, \ell<k-1, \\
& \. 1\. + \. \frac{1}{\aP(k-1)\. - \.1} \quad  & \text{for} \ \, \ell=k-1\ts.
\endaligned\right.
$$
In both cases, we have:
	\begin{equation}\label{eq:november 1}
		1  \ \leq \ \frac{\ac_{\ell} \ac_{\ell+2}}{\ac_{\ell+1}^2} \ \leq \  1 + \frac{1}{|\Par(\alpha)|-1} \  = \ 1+\frac{1}{m-1}\..
	\end{equation}
This implies that matrix \ts $\bN$ \ts satisfies the conditions in Lemma~\ref{l:Hyp N}.
Thus \ts $\bN$ \ts is hyperbolic, as desired.
\end{proof}

\smallskip

\begin{lemma}\label{l:morphism BHyp 2}
	Suppose that \.  $\rkn\big(\Phi(S)\big) = \rk(\Nf)-1$.
	Then the matrix \ts $\bBr$ \ts satisfies \eqref{eq:Hyp}.
\end{lemma}
\smallskip

\begin{proof}
By  assumptions of the lemma, we can partition \. $\Cnt(\alpha):= X_1 \ts \cup \ts X_2$ \. into two subsets:
	\begin{align*}
		X_1 \ &:= \  \big\{x \in \Cnt(\alpha) \, : \, \rkn\big(\Phi(S\ts + \ts x)  \big) \ = \ \rk(\Nf)  \big\}, \\
		X_2 \ & := \  \big\{x \in \Cnt(\alpha) \, : \. \rkn\big(\Phi(S \ts + \ts x)  \big) \ = \ \rk(\Nf)-1   \big\}.
	\end{align*}
We now make the observations in three possible cases of \ts $x,y\in X$:

\medskip

\nin
$(1)$ \.
 For every \ts $x,y \in X_1$, we have \. $S+x+y\in \Bc$ \. if and only if \. $x \not \sim_S y$.
This is because  \. $\Phi(S +x+y) \supseteq \Phi(S +x)$, which implies that \. $\Phi(S +x+y)$ \.
contains a basis of~$\Nf$, and because  \. $S +x+y\in \Ic$ \. if and only if \. $x \not \sim_S y$.
\smallskip

\nin
$(2)$ \.  For every \ts  $x \in X_1$ \ts and \ts $y \in X_2$,
we have \. $S +x+y$ \. is a basis of~$\Phi$.
This is because  \. $\Phi(S +x+y) \supseteq \Phi(S+x)$, which implies that
\. $\Phi(S +x+y)$ \. contains a basis of~$\Nf$, and because
			\begin{align*}
			\rkm(S +x+y) \. - \. \rkm(S+y) \, \geq \,
			\rkn\bigl(\Phi(S +x+y)\bigr) \. - \. \rkn\bigl(\Phi(S+y)\bigr) \, =  \, \rk(\Nf) - \bigl(\rk(\Nf)-1\bigl) \, = \,  1,
		\end{align*}
		which implies that \. $S+x+y\in \Ic$.
		\smallskip

\nin
$(3)$ \.  	For every \. $x,y \in X_2$,
		we have \. $S+x+y$ \. is not a basis of~$\Phi$.
		This is because  \. $\rkn\bigl(\Phi(S+x)\bigr) \.= \.\rkn\bigl(\Phi(S+y)\bigr) \. = \. \rkn\bigl(\Phi(S)\bigr) \. = \. \rk(\Nf)-1$, which implies that
		$\rkn\bigl(\Phi(S+x+y)\bigr)\. = \.\rk(\Nf)-1$.

\medskip

It  follows from the observations above that
	\begin{align*}
		\aB_{x\ts y} \ = \ \aB_{y\ts x} \ = \
		\begin{cases}
			1 & \ \text{ if }  \ x,y \in X_1  \, \text{ and } \, x \not \sim_S y\ts, \\
			1 & \ \text{ if } \  x \in X_1 \, \text{ and } \, y\in X_2\ts,\\
			1 & \ \text{ if } \  x \in X_1 \, \text{ and } \, y=\spstar,\\
			0 & \ \text{ otherwise. }
		\end{cases}
	\end{align*}
In particular, for \. $x,y \in X_1$ \ts and \ts $x \sim_S y$,
we have \ts $x$-row ($x$-column) of~$\bB$
equal to \ts $y$-row ($y$-column) of~$\bB$.
Similarly, for \. $x,y \in X_2$,
we have \ts $x$-row  ($x$-column) of $\bB$
is equal to \ts $y$-row ($y$-column) of~$\bB$.

Now let \. $x_1,\ldots, x_m$ \. be  representatives of the equivalence classes under
the relation \. ``$\sim_S$'' \. on $X_1$, and let $y$ be a representative element of~$X_2$.
For every other \ts $z \in X_1$ \ts in the same equivalence class of~$x_i$,
we subtract from the $z$-row ($z$-column) of $\bB$ the $x_i$-row ($x_i$-column) of~$\bB$.
For every other $w \in X_2$, subtract from the $w$-row ($w$-column) of~$\bB$ the
$y$-row ($y$-column) of~$\bB$.  Recall that \eqref{eq:Hyp} is preserved under these
transformations.

By applying these transformations and restricting to the support,
	we obtain the following matrix:{\small
	\[  \begin{pmatrix}
		0 & 1 &\ldots & 1 & 1  \\
		1  & 0 & \ddots &\vdots  & \vdots \\
		\vdots & \ddots & \ddots & 1  & 1 \\
		1 & \ldots & 1 & 0 &   0 \\
		1 & \ldots & 1 & 0  & 0
	\end{pmatrix}, \]}
where the rows and columns are indexed by $\{x_1,\ldots, x_m, y,  \spstar\}$.
The eigenvalues of this matrix are \. $\lambda_1=m$, $\lambda_2=0$, $\lambda_3=\ldots = \lambda_{m+2} =-1$.
This implies that the matrix satisfies  \eqref{eq:OPE}, and by Lemma~\ref{l:Hyp is OPE} also \eqref{eq:Hyp},
as desired.		
\end{proof}

\smallskip

\begin{lemma}\label{l:morphism BHyp 3}
	Suppose that \.  $\rkn\big(\Phi(S)\big) = \rk(\Nf)-2$.
	Then the matrix \ts $\bBr$ \ts satisfies \eqref{eq:Hyp}.
\end{lemma}
\smallskip

\begin{proof}
Let \ts $H \subseteq X$ \ts be given in~\eqref{eq:reduction-morphisms-H-def}, and let \ts
``$\sim_H$'' \ts be an equivalence relation defined by~\eqref{eq:reduction-morphisms-sim-def}.
Let us show that for every \. $x,y \in H$, we have:
\begin{equation}\label{eq:proof-morphism-BHyp3-claim}
\text{$S +x+y\in \Bc$ \quad $\Longleftrightarrow$ \quad
 $x \not \sim_H y$.}
\end{equation}
The \. $\Rightarrow$ \. direction is clear, so it suffices to prove the \. $\Leftarrow$ \. direction.
Let \. $x,y \in H$ \. such that \. $x \not \sim_H y$.
Then  we have:
	\[
\rkm(S+x+y) \. - \. \rkm(S) \ \geq  \  \rkn\big(\Phi(S+x+y) \big) - \rkn\big(\Phi(S)\big) \ = \  \rk(\Nf) - \bigl(\rk(\Nf)-2\bigr) \ = \ 2,
\]
which implies that \. $S+x+y\in \Ic$. Since $\Phi(S+x+y)$ is a basis of $\Nf$ by assumption,
it then follows that \. $S+x+y$ is a basis of $\Phi$, as desired.
	
	It then follows from the claim above that
	\begin{align*}
		\aB_{xy} \, = \, \aB_{yx} \, = \,
		\begin{cases}
			\. 1 & \ \text{ if }   \ x, y \in H \ \text{ and } \ x \not \sim_H y,  \\
			\. 0 & \ \text{ otherwise}.
		\end{cases}
	\end{align*}	
Note that, if \. $x,y \in H$ \. and \. $x \sim_H y$, then \ts $x$-row ($x$-column) of~$\bB$
is equal to $y$-row ($y$-column) of~$\bB$.  Also note that, the support of $\bB$ is contained
in~$H$.
	
Let \. $x_1, \ldots, x_m$ \. be the representatives of the equivalence classes \.
$\Cc_1,\ldots, \Cc_m$ \. of the relation \. ``$\sim_H$''.
For every other \ts $y \in \Cc_i$, we subtract from the $y$-row ($y$-column)
of~$\bB$ the $x_i$-row ($x_i$-column) of~$\bB$.  By applying this transformation
and restricting to the support of the resulting matrix, we obtain the following matrix:
\[
\begin{pmatrix}
		0 & 1  &\ldots & 1  \\
		1  & 0 & \ddots &\vdots   \\
		\vdots & \ddots & \ddots & 1 \\
		1 & \ldots & 1 & 0 		
\end{pmatrix},
\]
where the rows and columns are indexed by \. $\{x_1,\ldots, x_m\}$.
	The eigenvalues of this matrix are
	$\lambda_1=m-1$, $\lambda_2=\ldots=\lambda_m=-1$.
		This implies that the matrix satisfies  \eqref{eq:OPE}, and
by Lemma~\ref{l:Hyp is OPE} also \eqref{eq:Hyp}, as desired.		
\end{proof}

\smallskip

\begin{lemma}\label{l:morphism BHyp}
Let \ts $\Phi: \Mf \to \Nf$ be a morphism of matroid, let \ts $1\le k <\rk(\Mf)$,
and let \ts $\apr: X \to \rrs$ \ts be a positive weight function.  Let \ts $\AA$ \ts
be a combinatorial atlas corresponding to~$\Phi$.  Then every sink vertex
\. $\vf=(\alpha, 0,1)\in \Vfm$ \. satisfies \eqref{eq:Hyp}.
\end{lemma}

\smallskip

\begin{proof}
Let \. $\al=x_1\ts\cdots\ts x_\ell$ \. and \. $S=\{x_1,\ldots,x_\ell\}$.
It suffices to show that \ts $\bA(\alpha,1)$ \ts satisfies \eqref{eq:Hyp}. 	
If \. $\alpha \notin \Lc$ \. or \. $\rkn(S)  < \rk(\Nf)-2$,
then \ts $\bA(\alpha,1)$ \ts is equal to a zero matrix,
so \eqref{eq:Hyp} is trivially satisfied.
Now suppose that \. $\alpha \in \Lc_{\Mf}$ \. and \. $\rkn(S)  \geq  \rk(\Nf)-2$.
Then it follows from Lemma~\ref{l:morphism BHyp 1},
Lemma~\ref{l:morphism BHyp 2} and Lemma~\ref{l:morphism BHyp 3},
that \ts $\bA(\alpha,1)$ \ts satisfies \eqref{eq:Hyp}.
\end{proof}

\smallskip

\begin{lemma}
\label{l:proof-morphism-Hyp-main}
Let \ts $\Phi: \Mf \to \Nf$ be a morphism of matroid, let \ts $1\le k <\rk(\Mf)$,
and let \ts $\apr: X \to \rrs$ \ts be a positive weight function.  Let \ts $\AA$ \ts
be a combinatorial atlas corresponding to~$\Phi$.  Then every vertex
\. $\vf\in \Vf$ \. satisfies \eqref{eq:Hyp}.
\end{lemma}

\smallskip

\begin{proof}
Let \ts $\vf=(\alpha,m,t)\in \Vf^m$. We prove that \ts $\vf$ \ts satisfies \eqref{eq:Hyp}
by induction on~$m$.  The claim is true for \ts $m=0$ \ts by Lemma~\ref{l:morphism BHyp}.
Suppose that the claim is true for $\Vf^{m-1}$.
It then follows from Theorem~\ref{t:Hyp} that every regular vertex in $\Vf^m$
satisfies \eqref{eq:Hyp}.  On the other hand, by Lemma~\ref{l:morphism Irr},
the regular vertices of $\Vf^m$ contain  those of the form \ts $\vf=(\alpha,m,t)$,
where \ts $t \in (0,1)$.  Since \eqref{eq:Hyp} is a property that is preserved
under taking limits \ts $t\to 0$ \ts and \ts $t\to 1$, we conclude that
every vertex in $\Vf^m$ satisfies \eqref{eq:Hyp}. This completes the proof.
\end{proof}

\medskip

\subsection{Proof of Theorem~\ref{t:morphisms-refined}} \label{ss:proof-morphism}
Let \ts $\bM_{\vf}$ \ts be the associated matrix of the vertex \ts $\vf:=(\varnothing,k-1,1)\in \Vf$.
Let $\vb$ and $\wb$ be the characteristic vector of $X$ and $\{\spstar\}$, respectively.
Then
\begin{equation}\label{eq:proof-morphisms-vectors}
\aligned \< \wb, \bM \wb \> \  & = \  (k-1)! \.\cdot\. \aB_\ap(k-1), \qquad \< \vb, \bM \wb \> \ = \  k! \.\cdot\. \aB_\ap(k),  \\
 \text{and} \quad &
\< \vb, \bM \vb \> \  = \ (k+1)! \, \left(1\.+\. \frac{1}{\aP(k-1)-1}\right)\.\cdot\. \aB_\ap(k+1).
\endaligned
\end{equation}
Since \ts $\vf$ \ts satisfies \eqref{eq:Hyp} by Lemma~\ref{l:proof-morphism-Hyp-main},
it then follows from the equations above that
\[
\aB_{\ap}(k)^2  \ \geq \  \left(1+ \frac{1}{k}\right) \, \left(1\.+ \.\frac{1}{\aP(k-1)-1}\right)
\.\cdot\. \aB_{\ap}(k+1) \,\ts \aB_{\ap}(k-1), \]
as desired.\qed

\medskip

\subsection{Proof of Theorem~\ref{t:morphisms-Par-equality}}

We first prove the \. $\Leftarrow$ \. direction.
It follows from \eqref{eq:MorEqu3} that
\begin{align}\label{eq:morphism 1}
	\aB_\ap(k+1) \, = \, \aI_\ap(k+1),  \quad  \aB_\ap(k) \, = \, \aI_\ap(k), \quad \text{and}\quad
\aB_\ap(k-1) \, = \, \aI_\ap(k-1),
\end{align}
where \.
$$\aI_{\ap}(r) \, := \, \sum_{S\ts\in\ts\Ic_r} \. \ap(S).
$$
Similarly, \ts $\aP(k-1)$ \ts coincide for $\Phi$ and~$\Mf$.
Thus \eqref{eq:MorEqu1} is equivalent to \eqref{eq:ME1} for $\Mf$,
and \eqref{eq:MorEqu2} is equivalent to \eqref{eq:ME2} for $\Mf$.
It then follows from Theorem~\ref{t:matroids-Par-weighted} that
\[ 		 \aI_\ap(k)^2  \ = \  \left(1+ \frac{1}{k}\right)
\left(1+ \frac{1}{\aP(k-1)-1}\right) \, \aI_{\ap}(k+1) \, \aI_{\ap}(k-1),  \]
which together with \eqref{eq:morphism 1} proves the \. $\Leftarrow$ \. direction.

\smallskip

We now prove the \. $\Rightarrow$ \. direction.
It follows from the same argument as in the \. $\Leftarrow$ \. direction,
that it suffices to show that \eqref{eq:MorEqu3} is satisfied.
Let \ts $\AA$ \ts be the combinatorial atlas that corresponds to
\ts $(\Phi,k, \ap)$ \ts from  \S\ref{ss:proof-morphism-construction}.
In particular, every vertex of $\Qf$ satisfies \eqref{eq:Hyp} by
Lemma~\ref{l:proof-morphism-Hyp-main}.

As in \S\ref{ss:proof-morphism},
let \ts $\vb, \wb \in \Rb^{\ar}$ \ts be the characteristic vector of
 \ts $X$ \ts and \ts $\{\spstar\}$, respectively.
It is straightforward to verify that \ts $\vb, \wb$ \ts is a global pair
for $\Qf$, i.e.\ they satisfy~\eqref{eq:PosGlob}.

Let \. $\vf=(\varnothing,k-1,1)\in \Vf$ \. and let \ts $\bM=\bM_{\vf}$ \ts
be the associated matrix.
Note that  \. $\aB_\ap(k+1)$, \. $\aB_\ap(k)$ \. and \. $\aB_\ap(k-1)>0$ \. by the assumption of the theorem.
It then follows from \eqref{eq:morphism 1} that $\vf$ satisfies \eqref{eq:sEqu} for some \ts $\as>0$.

We now show that, for every \ts $\alpha \in \Lc_{k-1}$ \ts such that \. $\langle \vb, \bA(\alpha,1) \vb \rangle >0$, we have:
\begin{equation}\label{eq:morphism 3}
	\langle \vb, \bA(\alpha,1) \vb \rangle  \ = \  \as \, \langle \wb, \bA(\alpha,1) \vb \rangle \ = \  \as^2 \, \langle \wb , \bA(\alpha,1) \wb \rangle \ >  \ 0.
\end{equation}
	First suppose that $k=1$.
This implies that \ts $\alpha =\varnothing$ \ts and \ts $\vf=(\varnothing,0,1)$.
Thus, \eqref{eq:morphism 3} follows from the fact that
$\vf$ satisfies \eqref{eq:sEqu}.

Suppose now that \ts $k>1$.  It is easy to see  that
$\vf$ is a functional source in this case, i.e.\ it satisfies \eqref{eq:ProjGlob}
and \eqref{eq:hGlob},  where we apply the substitution $\fb \gets \vb$ for \eqref{eq:hGlob}.  By Theorem~\ref{t:equality Hyp},
every functional target of $\vf$ in $\Qf$ also
satisfies \eqref{eq:sEqu} with the same \ts $\as >0$.
On the other hand, observe that the functional targets of $\vf$ in $\Vfm$
contain those  of the form \. $(\alpha, 0,1)$, with \ts $\alpha \in \Lc_{k-1}$ \ts
satisfying \. $\langle \vb, \bA(\alpha,1) \vb \rangle >0$.
Combining these two observations, we obtain \eqref{eq:morphism 3}.

\vskip.5cm

\nin
\textbf{Claim:} \. {\em For every \ts $T \in \Ic_{k-1}$,
 we have \. $\rkn\big(\Phi(T)\big) \. \neq \.  \rk(\Nf)-1$.}

\smallskip

\begin{proof}
Let \ts $T=\{y_1,\ldots,y_{k-1}\}$ \ts and let \ts $\be = y_1 \. \cdots \. y_{k-1} \in \Lc$.
For every \ts $\ell \geq 0$, let
$$
\aB_{\ap,T}(\ell) \, := \, \sum_{\substack{S\in \Bc_{|T|+\ell} \\ S\ts\supseteq\ts T}} \, \ap(S)\ts.
$$
Then:
\begin{equation}\label{eq:morphism 2}
	\begin{split}
\< \wb, \bA(\beta,1) \wb \> \ & = \  \aB_{\ap,T}(0), \qquad	\< \wb, \bA(\beta,1) \vb \> \ = \  \aB_{\ap,T}(1),\\
		\< \vb, \bA(\beta,1) \vb \> \ &= \ 2 \, \left(1 \. + \. \frac{1}{\aP(0)\ts - \ts 1}\right)\, \aB_{\ap,T}(2).
	\end{split}
\end{equation}
Now suppose to the contrary that  \. $\rkn\big(\Phi(T)\big) =  \rk(\Nf)-1$.
Since \ts $\aB_\ap(k+1)>0$, there is a basis \ts $S\in \Bc_{k+1}$.
Applying the exchange property for $\Phi$, it follows that there exist \. $x,y \in S \setminus T$,
such that  \. $T \ts \cup \{x,y\}\in \Bc_{k+1}$.  This implies that  \. $\aB_{\ap,T}(2) >0$,
which in turn implies that \. $\< \vb, \bA(\beta,1) \vb \> >0$ \. by~\eqref{eq:morphism 2}.
Hence \eqref{eq:morphism 3} applies to~$\beta$, which implies that
\. $\< \wb, \bA(\beta,1) \wb \> >0$.  Again, by \eqref{eq:morphism 2} we conclude
that $\aB_{\ap,T}(0)>0$.
This contradicts  the assumption that \. $\rkn\big(\Phi(T)\big) =  \rk(\Nf)-1$.
This completes the proof of the claim. \end{proof}

\smallskip

It remains to prove~\eqref{eq:MorEqu3}, i.e.\ that every \ts $T \in \Ic_{k-1}$ \ts
satisfies \. $\rkn\big(\Phi(T)\big)  =  \rk(\Nf)$.
Suppose to the contrary that \. $\rkn\big(\Phi(T)\big) <   \rk(\Nf)$.
Since \ts $\aB_\ap(k-1)>0$, there is at least one basis \ts $S\in \Bc_{k-1}$.
By the exchange property of the matroid~$\Mf$ the basis exchange graph is connected,
i.e.\ there exist a sequence of bases \. $T_1,\ldots, T_m \in \Ic_{k-1}$\ts, such that \.
$\bigl|T_{i+1} \ts \sm \ts T_{i}\bigr| =1$, \. $T_1=T$, and \. $T_m=S$.
Since  \. $\rkn\big(\Phi(T) \big) < \rk(\Nf)$ \.
and \. $\rkn\big(\Phi(S) \big) = \rk(\Nf)$,
there exists \ts $i \in [m]$ \ts such that \.
$\rkn\big(\Phi(T_i) \big) = \rk(\Nf)-1$.  This contradicts the claim above,
and completes the proof of~\eqref{eq:MorEqu3}. \qed

\medskip

\subsection{Proof of Theorem~\ref{t:morphisms-EH-equality}}
The \. $\Leftarrow$ \. direction is straightforward.
For the \. $\Rightarrow$ \. direction,
it follows from \eqref{eq:MorEqu3} in Theorem~\ref{t:morphisms-Par-equality},
that for every \ts $S \subseteq X$, \. $|S|=k-1$, the image \ts $\Phi(S)$ \ts
contains a basis of~$\Nf$.  This implies that \eqref{eq:morphism 1} holds.
It then follows from Theorem~\ref{t:matroids-equality},
that every subset of $X$ of size~$k+1$ is independent,
and the weight \ts $\apr:X \to \rrs$ \ts is uniform.
This completes the proof. \qed

\bigskip


\section{Proof of log-concavity for linear extensions}
\label{s:proof-Sta}

In this section we give proofs of Theorem~\ref{t:Sta-weighted} and some
variations of the results for posets with belts ($\S$\ref{ss:LE-belts}).
We also give an example of a combinatorial atlas in this case
($\S$\ref{ss:LE-example}).

\smallskip

\subsection{New notation}
\label{ss:LE-notation}
In the next two sections we fix a ground set \ts $X$ \ts and an element \ts $z \in X$\ts.
Let \. $\Pf=(X_P,\prec)$ \. be a poset for which the ground set \ts $X_P$ \ts is a subset of  \ts $X$\ts.
Let $k \in \{2,\ldots,|X_P|-1\}$,
and let \. $\ap: X \to \rrs$ \. be the order-reversing weight function,
see~$\S$\ref{ss:intro-LE}.
 We define a combinatorial atlas \ts $\AA(\Pf,k):=\AA(\Pf,z,k,\ap)$ \ts as follows.

Recall that \ts $\Ec(\Pf)$ \ts denotes the set of linear extensions of~$\Pf$.
By a slight abuse of notation, in the next two sections a linear extension $\alpha$ of $\Pf$ is a simple word \ts $\al:=x_1\ldots x_{|\al|} \in X_P^*$ \ts of length $|X_P|$ such that  \. $x_i \prec x_j$ \. in $\Pf$ implies that $i\leq j$.
Denote by \ts $\Ec(P,k)$ \ts the set of linear extensions \ts $\alpha \in \Ec(\Pf)$ \ts  such that  \ts $\alpha_k=z$\ts.

For a simple word \ts $\al\in X^*$, we write \ts $x\tl z$ \ts if
\ts $x$ \ts appears to the left of \ts $z$ in~$\al$.
Following~\eqref{eq:sta-def-omega-weighted}, for a word \ts $\al \in X^*$,
let
$$
\ap(\al)  \ := \ \prod_{x \ts\tl\ts \zf} \. \ap(x)\ts,
$$
and \. $\ap(S) :=  \sum_{\al \in S} \ap(\al)$ \. for every $S \subseteq X^*$.

%
%


Let \. $\Zlow:= X -\zf$ \. and denote every element in \ts $\Zlow$ \ts
as \ts $x_{\low}$ \ts instead of $x$.
Similarly, let \. $\Zup:=X -\zf$, and denote every element in \ts $\Zup$  as \ts  $x_{\up}$ \ts instead of $x$.
Since \ts $\Zlow$ \ts \ts $\Zup$ \ts are two copies of the same set, labels
``$\low$'' and ``$\up$''  are used to distinguish them.  We write \.
$\Zf := \Zlow \ts \cup  \ts \Zup$.  Note that  \. $\ar:=|\Zf|=2n-2$ \. since $\Zlow$
and $\Zup$ do not intersect because of the labeling.
We will sometimes drop the ``$\low$'' and ``$\up$'' labels from
$x_{\low}$ and $x_{\up}$ \. when the labels are either clear from the context
or are irrelevant to the discussion.
We denote by \ts $\min(\Pf,\low) \subseteq \Zf_{\low}$ \ts the set of elements of $\Zf_{\low}$ that correspond to minimal elements of $\Pf$,
and by \ts $\max(\Pf,\up) \subseteq \Zf_{\up}$ \ts the set of elements of $\Zf_{\up}$ that correspond to maximal elements of $\Pf$.
More generally, for a subset \ts $S \subseteq X-z$\ts, we denote by \ts $S_{\low}\subseteq \Zf_\low$ \ts the subset in $\Zf_\low$ that corresponds to $S$,
and by \ts $S_{\up}\subseteq \Zf_\up$ \ts the subset in $\Zf_\up$ that corresponds to $S$.

\smallskip

Let \. $\Pf^{\op}:=(X, \prec^{\op})$ \.
denote the \defn{opposite poset} of $\Pf$,
defined by \. $x \prec^{\op} y$ \.
if and only if \. $y \prec x$.\footnote{Sometimes, \. $\Pf^{\op}$ \. is
	also called \defng{dual} or \defng{reverse poset}. }
For every \. $\alpha=x_1\ldots x_{\ell} \in X^*$,
we denote by \. $\alpha^{\op}:=x_\ell \ldots x_1$.
Let \. $\Ec^{\op}$ \. denote the set of linear extensions of \. $\Pf^{\op}$,
and note that \. $|\Ec^{\op}| = |\Ec|=e(\cP)$.
Denote  by \. $\ap^{\op}: X \to \rrs$ \.
the weight function defined by \. $\ap^{\op}(x) := \ap(x)^{-1}$.
Note that \. $\ap^{\op}$ \. is an order-reversing weight function for \. $\Pf^{\op}$.
It then follows that
\begin{equation}\label{eq:tango 4}
\aN_{\ap}(\Pf,k) \ = \   	\aN_{\ap^{\op}}(\Pf^{\op},|X_P|-k+1) \. \prod_{x \in X-z} \ap(x).
\end{equation}
In the subsequent two sections, we shall frequently utilize this technique of interchanging between $\Pf$ and $\Pf^{\op}$ to streamline certain parts of the proofs.

\medskip

\subsection{Combinatorial atlas construction}
\label{ss:LE-atlas}
%


We denote by  \. $\bC(\Pf,k):=\bC(\Pf,k,\ap):= \bigl(\aC_{x\ts y}\bigr)_{x,y \in \Zf}$ \.
the symmetric \ts $\ar \times\ar$ \ts
matrix where,\footnote{Here \. $\omega(\Ec(\Pf-x-y,k-1))$ \. is the sum of $\omega$-weight of all linear extensions $\al$ of $\Pf-x-y$ for which $\al_{k-1}=z$.}
\begin{equation}\label{eq:C-1}\tag{DefC-1}
\begin{split}
\aC_{x \ts y} \ &:= \ \begin{cases}
	\ap(x) \, \omega(\Ec(\Pf-x-y,k-1)) & \text{if  $x \in \min(\Pf,\low)$,  $y \in  \max(\Pf,\up)$},\\
		\ap(x) \, \ap(y) \, \omega(\Ec(\Pf-x-y,k-1)) & \text{if  $x,y \in \min(\Pf,\low)$,  $x\neq y$},\\
	 \omega(\Ec(\Pf-x-y,k-1)) & \text{if  $x,y \in \max(\Pf,\up)$,  $x\neq y$},\\
0 & \text{otherwise},
\end{cases} \\
\aC_{x \ts x} \ &:= \ \sum_{y \in \min(\Pf-x,\low), \. y \succ x} \ap(x) \ap(y) \, \omega(\Ec(\Pf-x-y,k-1)) \quad \text{ for $x \in \min(\Pf,\low)$},\\
\aC_{x \ts x} \ &:= \ \sum_{y \in \max(\Pf-x,\up), \. y \prec x}\omega(\Ec(\Pf-x-y,k-1)) \quad \text{ for $x \in \max(\Pf,\up)$},\\
\aC_{x \ts x} \ &:= \ 0 \quad \text{ for $x \notin \min(\Pf,\low) \cup \max(\Pf,\up)$.}
\end{split}
\end{equation}
Equivalently, \. $\bC(\Pf,k)$ \. is given by\footnote{Here $\{ x\beta y  \in \Ec(P,k) \}$ is the set of linear extensions \. $\alpha \in \Ec(P,k)$ \. such that \ts $\alpha_1=x$ \ts and \ts $\alpha_{|X_{\Pf}|}=y$\ts. The word $\beta \in X^*$ here denotes $\alpha_2 \cdots \alpha_{|X_{\Pf}|-1}$. }
\begin{equation}\label{eq:C-2}\tag{DefC-2}
\begin{split}
	\aC_{x \ts y} \ &:= \ \aC_{y \ts x}  \ := \
	\begin{cases}
	\ap\big(\{ x\beta y  \in \Ec(P,k) \}\big) & \text{if \  $x \in \min(\Pf,\low)$, \. $y \in \min(\Pf,\up)$,}\\
		\ap\big(\{ xy\beta  \in \Ec(P,k+1) \}\big) & \text{if \  $x,y \in \min(\Pf,\low), \. x\neq y$,}\\
			\ap\big(\{\beta xy  \in \Ec(P,k-1) \}\big) & \text{if \  $x,y \in \max(\Pf,\up), \. x\neq y$,}\\
			0 & \text{ if $\{x,y\} \nsubseteq \min(\Pf,\low) \cup \max(\Pf,\up)$},
	\end{cases}\\
\aC_{x \ts x} \ &:= \
\begin{cases}
 \ap\big(\{ xy \beta \in \Ec(\Pf,k+1) \mid  y \succ x \}\big) & \text{if \ $x \in \min(\Pf,\low)$},\\
  \ap\big(\{\beta yx \in \Ec(\Pf,k-1) \mid  y \prec x \}\big) & \text{if \ $x \in \max(\Pf,\up)$}.
\end{cases}
\end{split}
\end{equation}
Note that both definitions will be frequently employed throughout the next two sections, chosen based on their suitability.
Also note that it follows from the definition that \ts $\bC$ \ts is a nonnegative symmetric matrix.

\smallskip

Note that, it follows from    \eqref{eq:C-2} that, for every \. $x \in \min(\Pf,\low)$,
\begin{equation}\label{eq:C-sum-1}
\begin{split}
		\sum_{y \in \Zf_{\low}} \aC_{x \ts y} \ &= \  \ap\big(\{x\beta \mid \Ec(\Pf,k+1) \}\big) \ = \   \ap(x) \, \ap\big(\Pf-x,k\big),\\
				\sum_{y \in \Zf_{\up}} \aC_{x \ts y} \ &= \  \ap\big(\{x\beta \mid \Ec(\Pf,k) \}\big) \ = \ \ap(x)\,  \ap\big(\Pf-x,k-1\big).
\end{split}
\end{equation}
 Similarly, for every \. $x \in \max(\Pf,\up)$\.,
 \begin{equation}\label{eq:C-sum-2}
 	\begin{split}
 		\sum_{y \in \Zf_{\low}} \aC_{x \ts y} \ &= \  \ap\big(\{\beta \. x \mid \Ec(\Pf,k) \}\big) \ = \  \ap\big(\Pf-x,k\big),\\
 		\sum_{y \in \Zf_{\up}} \aC_{x \ts y} \ &= \  \ap\big(\{\beta \. x \mid \Ec(\Pf,k-1) \}\big) \ = \  \ap\big(\Pf-x,k-1\big).
 	\end{split}
 \end{equation}

Let \. $\fb,\gb \in \rr^{\ar}$ \. be the indicator vector of $\Zf_{\low}$ and $\Zf_{\up}$, respectively.
It follows from \eqref{eq:C-sum-1} and  \eqref{eq:C-sum-2}  that
\begin{equation}\label{eq:Cfg}\tag{Cfg}
	\begin{split}
		& \<\fb, \bC(\Pf,k) \gb \> \ = \  \aN_{\ap}(\Pf,k), \qquad
		\<\fb, \bC(\Pf,k) \fb \> \ = \   \aN_{\ap}(\Pf,k+1),\\
		&\<\gb, \bC(\Pf,k) \gb \> \ = \   \aN_{\ap}(\Pf,k-1),
	\end{split}
\end{equation}
where recall that \. $\aN_{\ap}(\Pf,k)$ \. is the sum of $\ap$-weight of linear extensions of $\Pf$ such that $z$ is the $k$-th smallest element.

\smallskip

Let \. $\Qf:=\Qf(\Pf,k)  :=  (\Vf,\Ef)$ \. be the acyclic graph with \.  $\Vf=   \Vf^0 \. \cup \. \Vf^1$,
where
\begin{align*}
	\Vf^1 \, & := \ \{t \in \rr \mid 0 \leq t \leq 1  \}, \qquad \Vfm  \  := \ Z.
\end{align*}
For a non-sink vertex  \ts $v=t \in \Omega^1$ \ts and $x \in Z$, the corresponding outneighbor in $\Omega^0$ is  \ts $v^{\<x\>}:= x$\ts.

Define the combinatorial atlas \ts $\AA(\Pf,k)$ \ts of dimension \. $\ar$ \.
corresponding to \ts poset $\Pf$, and $k \in \{3,\ldots,|X_P|-1\}$
  by the acyclic graph~$\Qf$ and the linear algebraic
data defined as follows.
For each vertex \ts $\vf=x \in \Vf^0$\ts, the associated matrix is
\[ \bM_v \ := \
\begin{cases}
\ap(x) \, \bC(\Pf-x,k-1) \quad  \text{ if  $x \in \min(\Pf,\low)$},\\
 \bC(\Pf-x,k-1) \quad  \text{ if  $x \in \max(\Pf,\up)$},
\end{cases}
\]
and is equal to the zero matrix otherwise.
For each vertex \ts $\vf=t \in \Vf^1$\ts,  the associated matrix is
\[  \bM \ := \ \bM_{\vf} \ := \  t \, \bC(\Pf,k)  \ +  \ (1-t) \, \bC(\Pf,k-1), \]
and the  associated vector \. $\hb := \hb_{\vf} \in \Rb^{\ar}$ \. is defined to have coordinates
\[  \ah_x \ := \
\begin{cases}
	t & \text{ if } \ \  x \in \Zf_{\low}\.,\\
	1-t & \text{ if } \ \ x \in \Zf_{\up}\..
\end{cases} \]
Finally, let
the linear transformation \. $\bT^{\<x\>}: \Rb^{\ar} \to \Rb^{\ar}$ \.
associated to the edge \ts $(\vf,\vf^{\<x\>})$,  be
\begin{align*}
	\big( \bT^{\<x\>}\vb\big)_y \ := \
	\begin{cases}
		\av_y & \text{ if } \  y  \in \supp(\bM), \\
		\av_x & \text{ if } \ y \in \Zf  \setminus  \supp(\bM).
	\end{cases}
\end{align*}

\medskip

\subsection{Properties of the matrix $\bC(\Pf,k)$}

In this subsection we gather properties of the matrix \ts $\bC(\Pf,k)$ \ts that will be used in this paper.

\smallskip

\begin{lemma}\label{lem:C-supp}
	Let $\Pf$ be a poset, and let
	and let $k\in \{2,\ldots,|X_P|-1\}$ such that $\aNr(\Pf,k)>0$.
	Then
	\begin{itemize}
		\item The support of $\bC\.(\Pf,k)$ is equal to $\min(\Pf,\low) \cup \max(\Pf,\up)$, and
		\item The matrix $\bC\.(\Pf,k)$ is irreducible when restricted to the support.
	\end{itemize}
\end{lemma}

\smallskip

\begin{proof}
	Let $n:=|X_P|$.
	It follows from \eqref{eq:C-1} that the support of \ts $\bC\.(\Pf,k)$ \ts  is a subset of \ts $\min(\Pf,\low) \cup \max(\Pf,\up)$\ts.
	Now note that,	since
	$\aN(\Pf,k)>0$,
	there exists a linear extension \ts $\alpha=x_1\cdots x_n \in \Ec(P,k)$\ts, and
	note that \ts $x_1={(x_1)}_{\low} \in \min(\Pf,\low)$  \ts and \ts $x_n={(x_n)}_{\up} \in \max(\Pf,\up)$\ts.	
	Now, let $y$ be an arbitrary element  of \ts $\min(\Pf,\low) \cup \max(\Pf,\up)$\ts.
	For the first claim it suffices to show that \. $y \in \supp \. (\bC(\Pf,k))$\., and for the second claim it suffices to show that  that $y$ is contained in the same irreducible component (of the matrix \ts $\bC(\Pf,k)$\ts) as $(x_1)_{\low}$ and $(x_n)_{\up}$.
	
	By switching to the dual poset in \eqref{eq:tango 4} if necessary,
	we will without loss of generality assume that  \ts $y=y_\low \in \min(\Pf,\low)$\ts.
	Let $\alpha'$ be the linear extension obtained from $\alpha$ by \emph{demoting} $y$ to be the smallest element,
	i.e.
	\[ \alpha':= y \, x_1\cdots x_{i-1} x_{i+1} \cdots x_n, \qquad \text{where} \quad \alpha=: x_1\cdots x_{i-1} \, y \, x_{i+1} \cdots x_n.  \]
	Note that $\alpha'$ is still a linear extension of $\Pf$ since $y$ is a minimal element of $\Pf$.
	Now note that either \ts $\alpha' \in \Ec(\Pf,k)$ \ts or \ts $\alpha' \in \Ec(\Pf,k+1)$\ts.
	In the first case we then have   \ts $(\bC\.(\Pf,k))_{y \ts x_n} >0$\ts, so $y$ is contained in the support of \ts $\bC\.(\Pf,k)$ \ts and is in the same irreducible component as $x_n$.
	In the second case we then have  \ts $(\bC\.(\Pf,k))_{y \ts x_1} >0$\ts, so $y$ is contained in the support of \ts $\bC\.(\Pf,k)$ \ts and is in the same irreducible component at $x_1$.
	This completes the proof.
\end{proof}

\smallskip

\smallskip

\begin{lemma}\label{lem:C-center}
	Let $\Pf$ be a poset,
	and let \. $k \in \{3,\ldots,|X_P|-1\}$ \. such that \. $\aNr(\Pf,k)>0$ \. and \. $\aNr(\Pf,k-1)>0$.
	Then, for every \. $x \in \min(\Pf,\low) \cup \max(\Pf,\up)$\.,
	\[ \aNr(\Pf-x,k-1)>0. \]
\end{lemma}

\smallskip

\begin{proof}	By switching to the dual poset in \eqref{eq:tango 4} if necessary,
	we will without loss of generality assume that \ts $x=x_\low \in \min(\Pf,\low)$\ts.
	By assumption there exists linear extensions \ts $\alpha \in \Ec(\Pf,k-1)$ \ts and
	\ts $\beta \in \Ec(\Pf,k)$\ts.
	Let $\alpha'$ and $\beta'$ be the linear extension of $\Pf$ obtained from $\alpha$ and $\beta$ by demoting $x$ to be the smallest element, respectively.
	It then follows that \. $\alpha' \in \Ec(\Pf,k-1) \cup \Ec(\Pf,k)$ \. and \. $\beta' \in \Ec(\Pf,k) \cup \Ec(\Pf,k+1)$\..
	If \. $\alpha' \in \Ec(\Pf,k)$ \. then we are done, as removing the smallest element from $\alpha'$ (which is $x$) will give us a linear extension in $\Ec(\Pf-x,k-1)$.
		If \. $\beta' \in \Ec(\Pf,k)$ \. then we are also done, as removing the smallest element from $\beta'$ (which is $x$) will give us a linear extension in $\Ec(\Pf-x,k-1)$.
		So we assume that \. $\alpha' \in \Ec(\Pf,k-1)$ \. and \. $\beta' \in \Ec(\Pf,k+1)$\..
		
		This assumption implies that there exists \ts $y \in X_P$ \ts which appears to the right of $z$  in $\alpha'$, but appears to the left of $z$ in $\beta'$.
		This in turn implies that $y$ is incomparable to $z$ in $\Pf$.
		Now, let \. $j$  \. be the smallest integer in the set
\[
        \big\{  i \, : \, \ x_i' \. \| \. z \text{ \ in $\Pf$}, \, k \le i\le |X_P| \. \big\},
\]
		where \ts $x_1'\cdots x_n':=\alpha'$\ts.
		Note that this set is non-empty by the preceding argument.
		Let $\gamma$ be the linear extension of $\Pf$ obtained from $\alpha'$ by demoting \. $x_j'$ \. to the $k-1$-th position.
		Then \. $\gamma \in \Ec(\Pf,k)$ \. and furthermore $x$ is the smallest element in $\gamma$.
		Then,  removing the smallest element of $\gamma$ gives us a linear extension in $\Ec(\Pf-x,k-1)$, and the proof is complete.
\end{proof}

\smallskip

\begin{rem}\label{r:LE-connected}
	The arguments  in Lemma~\ref{lem:C-supp} and Lemma~\ref{lem:C-center} are
	variations of the maximality argument that appears in the proof of
	Thm~8.9 in~\cite{CPP2}.  We refer to~\cite[$\S$12.3, $\S$14.2]{CP-survey}
    for a detailed survey.
\end{rem}

\medskip

\subsection{Properties of the combinatorial atlas}
\label{ss:LE-prop}
We now show that the atlas \ts $\AA(\Pf,k)$ \ts defined above,
satisfies all four conditions in Theorem~\ref{t:Pull},
namely properties \eqref{eq:Inh}, \eqref{eq:Proj},
\eqref{eq:TPInv} and \eqref{eq:Knonneg}.   We prove these
properties one by one, in the following series of lemmas.
For every lemma in this subsection we assume that $\Pf=(X_P,\prec)$ is a poset, and  $k \in \{3,\ldots, |X_P|-1\}$ such that \ts $\aN(\Pf,k)>0$ \ts and \ts $\aN(\Pf,k-1)>0$\ts.

%
%
%

\begin{lemma}\label{l:Stanley Inh}
	The atlas \ts $\AA(\Pf,k)$ \ts satisfies \eqref{eq:Inh} and \eqref{eq:Proj}.
\end{lemma}

\begin{proof}
	Let \ts $\vf=t \in \Vfp$ \ts be a non-sink vertex of $\Qf$.
		The property \eqref{eq:Proj} follows directly from the definition of \ts $\bT^{\<x\>}$.	
		For \eqref{eq:Inh}, let \ts $x \in \supp(\bM)$.
	By linearity of \ts $\bT^{\<x\>}$, it suffices to show that, for every \ts $y \in \Zf$, we have:
	\begin{equation}\label{eq:Inh Stanley 1}
		\aM_{xy} \ = \ \big\< \bT^{\<x\>} \eb_y, \bM^{\<x\>} \,  \bT^{\<x\>} \hb \big\>,
	\end{equation}	
	where \. $(\eb_y)_{y \in \Zf}$ \. is the standard basis for \ts $\Rb^{\ar}$.
	Note that we can assume $y \in \supp(\bM)$,
			as otherwise \.   $\bM \eb_y \. = \. \bT^{\<x\>} \eb_y \. = \.   \0$,
	and \eqref{eq:Inh Stanley 1} then follows trivially.
	It then follows from Lemma~\ref{lem:C-supp} that
	\. $x,y \in  \min(\Pf,\low) \cup \max(\Pf,\up)$\..
	Without loss of generality, assume that \ts $x = x_{\low} \in \Zlow$ \ts and \ts $y = y_{\low} \in \Zlow$,
	as the proofs of the other cases are analogous.
	
We split the proof of~\eqref{eq:Inh Stanley 1} into two  cases.	
First suppose that $x$ and $y$ are distinct.
It then follows that  \.   $\bT^{\<x\>} \eb_y \. = \.   \eb_y$, and
\begin{equation*}
	\begin{split}
		&\big\< \bT^{\<x\>} \eb_y, \bM^{\<x\>} \,  \bT^{\<x\>} \hb \big\> \ = \
		\sum_{\us \in \Zf }  \aM^{\<x \>}_{\us y} \, \big(\bT^{\<x\>} \hb \big)_\us  \\
		&\quad   = \  \sum_{\us \in  \Zlow}  \big(\bC(\Pf-x,k-1)\big)_{\us y} \,  t
		\  + \ \sum_{\us \in   \Zup}  \big(\bC(\Pf-x,k-1)\big)_{\us y} \,  (1-t) \\
		&\quad =_{\eqref{eq:C-sum-1}} \  \ap(x) \. \ap(y) \. \ap(\Ec(\Pf-x-y,k-1)) \. t \ + \
	\ap(x) \. \ap(y) \. 	\ap(\Ec(\Pf-x-y,k-2)) \. (1-t)   \\
		&\quad =_{\eqref{eq:C-1}} (\bC(\Pf,k))_{x \ts y} \. t \ + \  (\bC(\Pf,k-1))_{x \ts y} \. (1-t) \quad = \quad   \aM_{x \ts y},
	\end{split}
\end{equation*}
	as desired.

Now suppose that $x=y$.
Then
	\begin{equation*}
	\begin{split}
		&\big\< \bT^{\<x\>} \eb_x, \bM^{\<x\>} \,  \bT^{\<x\>} \hb \big\> \ = \
		\sum_{\substack{w \ts\in\ts \Zlow \\ w \ts\succ\ts x}} \. \sum_{\us \in \Zf } \. \aM^{\<x \>}_{\us w} \, \big(\bT^{\<x\>} \hb \big)_\us  \\
		&\quad   = \  \sum_{\substack{w \in \Zlow \\ w \succ x}}  \left( \sum_{\us \in  \Zlow} \.  \big(\bC(\Pf-x,k-1\big))_{\us w} \. t \,  + \,   \sum_{\us \in   \Zup}  \. \big(\bC(\Pf-x,k-1)\big)_{\us w}  \. (1-t) \right) \\
		&\quad =_{\eqref{eq:C-sum-1}}  \
	 \sum_{\substack{w \in \Zlow \\ w \succ x}} \ap(x) \. \ap(w) \. \ap(\Ec(\Pf-x-w,k-1)) \. t \ + \ \ap(x) \. \ap(w) \. \ap(\Ec(\Pf-x-w,k-2)) \. (1-t)  \\
			& \quad =_{\eqref{eq:C-1}} \
			 \big(\bC(\Pf, k\big))_{xx} \. t  \, + \,
		 \big(\bC(\Pf, k-1)\big)_{x,x} \. (1-t) \quad = \quad   \aM_{xx}\.,
	\end{split}
\end{equation*}
which completes the proof.
\end{proof}

\smallskip

\begin{lemma}\label{l:Stanley TPInv}
	The atlas \ts $\AA(\Pf,k)$ \ts satisfies \eqref{eq:TPInv}.
\end{lemma}

\begin{proof}
	Let \. $\vf=t \in \Vfp$ \. be a non-sink vertex, and let \. $a,b,c$ \.
be distinct elements of $\supp(\bM)$.
It follows from Lemma~\ref{lem:C-supp} that  \. $a,b,c \in \min(\Pf,\low) \cup \max(\Pf,\up)$\..
Without loss of generaltiy assume that \. $a,b,c \in \min(\Pf,low)$.
It then follows from \eqref{eq:C-1} that
\[  \aM^{\<a\>}_{b \ts c} \ = \  \aM^{\<b\>}_{c \ts a} \ = \ \aM^{\<c\>}_{a \ts b} \ = \ \ap(a) \. \ap(b) \. \ap(c) \. \ap(\Ec(\Pf-a-b-c,k-2)), \]
and the lemma follows.
\end{proof}

\smallskip

\begin{lemma}\label{l:Stanley Knonneg}
	The atlas \ts $\AA(\Pf,k)$ \ts  satisfies \eqref{eq:Knonneg}.
\end{lemma}

\begin{proof}
	Let \. $\vf=t\in \Vfp$ \. be a non-sink vertex.
	We need to check the condition
	\eqref{eq:Knonneg} for distinct \. $x,y \in \supp(\bM)$.
	It follows from Lemma~\ref{lem:C-supp} that 	\. $x,y \in  \min(\Pf,\low) \cup \max(\Pf,\up)$\..
	We will without loss of generality assume that \ts $x=x_{\low} \in \min(\Pf,\low)$ \ts and \ts $y=y_{\low} \in \min(\Pf,\low)$\ts, as the proof of  other cases are analogous.
	
	It follows from \eqref{eq:C-1} that
	\begin{align*}
	&\aM_{yy}^{\<x\>}  \ = \ \sum_{w \in \min(\Pf-x-y,\low), \. w \succ y} \ap(x) \. \ap(y) \. \ap(w) \. \omega(\Ec(\Pf-x-y-w,k-2)).
\end{align*}
Now note that, it follows from Lemma~\ref{lem:C-supp} that
\begin{align*}
	\supp(\bM) \ts - \ts y \ &= \   (\min(\Pf,\low)-y) \cup \max(\Pf,\up),\\
	\supp\big(\bM^{\<y\>}\big)  \ &= \  \min(\Pf-y,\low) \cup \max(\Pf-y,\up),
\end{align*}
Then,  the set \. $\Fam^{\<y\>}$ \. defined in~\eqref{eq:Fam-def}, in this case is  equal to
	\begin{align*}
	\Fam^{\<y\>} \ &= \
\supp\big(\bM^{\<y\>}\big)   \setminus  \big(\supp(\bM) \ts - \ts y\big) \ = \ \min(\Pf-y,\low) \setminus \min(\Pf,\low) \\
&= \ \{ w \in \min(\Pf-y,\low) \. \mid \. w \succ y \}.
	\end{align*}
	This implies that
	\begin{align*}
		\sum_{\ws \ts \in \ts \Fam^{\<y\>}} \. \aM_{x\ws}^{\<y\>}  \ &= \
		\sum_{\ws \ts \in \ts \Fam^{\<y\>}}  \ap(x) \. \ap(y) \. \ap(w) \. \ap(\Ec(\Pf-x--yw,k-2)) \\
		\ &= \
		\sum_{w \in \min(\Pf-y,\low), \. w \succ y} \ap(x) \. \ap(y) \. \ap(w) \omega(\Ec(\Pf-x-y-w,k-2)).
	\end{align*}
	Taking the difference of the two equations above,
	we get:
	\begin{align*}
		\aM_{yy}^{\<x\>}   -
		\sum_{\ws \ts \in \ts \Fam^{\<y\>}} \aM_{x\ws}^{\<y\>}
		\ = \ \sum_{w \in \min(\Pf-x-y,\low), \. w \succ y, \. w \succ x} \ap(x) \. \ap(y) \. \ap(w) \. \omega(\Ec(\Pf-x-y-w,k-2)).
	\end{align*}
	This is clearly nonnegative, and thus  \eqref{eq:Knonneg} holds, as desired.
\end{proof}

\smallskip

\begin{lemma}\label{l:Stanley Irr}
	Every \. $\vf\in \Vfp$ \. satisfies \eqref{eq:Irr}.
	Furthermore, every $\vf=t\in \Vfp$ \. satisfies \eqref{eq:hPos},
	for all \. $0<t<1$.
\end{lemma}

\smallskip

\begin{proof}
	Property \eqref{eq:Irr} follows directly from Lemma~\ref{lem:C-supp},
	and Property \eqref{eq:hPos} follows from the observation that $\hb_{\vf}$ is a positive vector when $t \in (0,1)$.
\end{proof}

\medskip

\subsection{Sink vertices are hyperbolic}\label{ss:LE-sink}
Before we can apply the local-global principle, we need
the following result:

\smallskip

\begin{lemma}\label{l:Stanley BHyp}
Let \. $\Pf=(X_P,\prec)$ \. be a finite poset with $|X_P|=3$,  let \. $\apr:X \to \rrs$ \.
be an order-reversing weight function.
Then the matrix \. $\bC\.(\Pf,2)$ \.
satisfies \eqref{eq:Hyp}.
\end{lemma}

\smallskip

\begin{proof}
Let \. $\{x,y,z\}:= X_P$.
We index the rows and columns of
 \ts $\bC(\Pf,2)$ \ts with  \. $\{x_\low, y_\low, x_\up, y_\up\}$\..
%

We now split the proof of the lemma into seven cases,
	depending on the relative order of \ts $\{x,y,\zf\}$.
First, suppose that \. $x,y,\zf$ \. are incomparable to each other.
Then
{\small	\begin{equation}\label{BaseA1}
		\tag{C1}
		\bC(\Pf,2) \quad = \quad
		\begin{pmatrix}
			0 &  \ap(x) \, \ap(y) & 0 &  \, \ap(x)   \\
			 \, \ap(x) \, \ap(y) & 0  &  \,  \ap(y) & 0   \\
			0 &  \, \ap(y)   & 0 & 1   \\
			 \, \ap(x) & 0  &  1 & 0
		\end{pmatrix}.
	\end{equation}
}	
	We now divide $x_\low$-row and $x_\low$-column by $\ap(x)$,
	and the \ts $y_\low$-row \ts and the \ts $y_\low$-column by
	\ts $\ap(y)$. Recall that \eqref{eq:Hyp} is preserved under this transformation.
	Then the matrix becomes
{\small	  	\begin{equation*}
	  	\begin{pmatrix}
	  		0 &  1 & 0 &   1  \\
	  		1  & 0  &  1   & 0   \\
	  		0 &  1   & 0 &  1  \\
	  		1  & 0  &  1 & 0
	  	\end{pmatrix}.
	  \end{equation*}
}
The eigenvalues of this matrix are \. $\{2 \, , 0,0, -2\, \}$.
This implies that the matrix satisfies \eqref{eq:OPE}.  By Lemma~\ref{l:Hyp is OPE}
we also have \eqref{eq:Hyp}, as desired.

\smallskip

Second, suppose that $x \prec y$,  and $\zf$ are incomparable to both elements.
Then
{\small
		\begin{equation}\label{BaseA2}
		\tag{C2}
		\bC(\Pf,2) \quad = \quad
		\begin{pmatrix}
			 \, \ap(x) \, \ap(y)  & 0 & 0 &  \, \ap(x)   \\
			0 & 0  &   0 & 0   \\
			0 & 0   & 0 & 0  \\
			  \, \ap(x) & 0  & 0  & 1
		\end{pmatrix}.
	\end{equation}
}	Restricting the rows and columns to the support \. $\big\{x_\low, y_\up\big\}$,
	we get
{\small
\[  \begin{pmatrix}
		 \, \ap(x) \, \ap(y) &  \, \ap(x) \\
		 \, \ap(x)  & 1
	\end{pmatrix}.
	\]
}	
This matrix has determinant
\[
    \ap(x) \,  \big( \ap(y) - \ap(x) \big) \ \leq \ 0,
\]
	where the inequality follows from $\ap$ being order-reversing.
	This implies that the matrix satisfies \eqref{eq:OPE}, and thus also \eqref{eq:Hyp}, as desired.
	
\smallskip

	In the remaining cases, element $\zf$ is comparable to either~$x$ or $y$, or both.
	By symmetry, without loss of generality, we assume that \ts $x \prec \zf$.
	Third, suppose that \.
	$x \prec \zf$, \. $x\prec y$, and \. $y \. || \. \zf$.
	Then:
{\small			
        \begin{equation}\label{BaseA3}
		\tag{C3}
		\bC(\Pf,2) \quad = \quad
		\begin{pmatrix}
		  \, \ap(x) \, \ap(y)   &0  & 0 &  \, \ap(x)   \\
			0 & 0  &   0 & 0   \\
			0 & 0   & 0 & 0  \\
			  \, \ap(x) \,  & 0  & 0  & 0
		\end{pmatrix}.
	\end{equation}
}		Restricting the rows and columns to the support \ts $\{x_\low, y_\up\}$,
	we get
{\small 	\[  \begin{pmatrix}
		 \, \ap(x) \, \ap(y) &  \, \ap(x) \\
		 \, \ap(x)  & 0
	\end{pmatrix}.
	\]}
	This matrix has a negative determinant, so it satisfies \eqref{eq:OPE}.
Thus, it also satisfies \eqref{eq:Hyp}, as desired.
	
\smallskip
	
Fourth,	suppose that \. $x \prec \zf$, \. $x \. || \. y$, and \. $y \. || \. \zf$.
	Then:
{\small
        \begin{equation}\label{BaseA4}
		\tag{C4}
		\bC(\Pf,2) \quad = \quad
		\begin{pmatrix}
			0   &  \, \ap(x) \, \ap(y)  & 0 &  \, \ap(x)   \\
			 \, \ap(x) \, \ap(y) & 0  &   0 & 0   \\
			0 & 0   & 0 & 0  \\
			 \, \ap(x) \,  & 0  & 0  & 0
		\end{pmatrix}.
	\end{equation}
}	By restricting the rows and columns to the support \. $\{x_\low, y_\low, y_\up\}$, followed by
	dividing the $x_\low$-row and $x_\low$-column by $\ap(x)$,
	and the $y_\low$-row and the $y_\low$-column by
	$\ap(y)$,
	we get
{\small 	\[  \begin{pmatrix}
		0 &  1 & 1   \\
		 1  & 0 & 0 \\
		 1  & 0 & 0
	\end{pmatrix}.
	\]}
The eigenvalues of this matrix are \. $\{\sqrt{2} \, , 0, -\sqrt{2} \, \}$,
	 so it satisfies \eqref{eq:OPE}. Thus it also satisfies \eqref{eq:Hyp}, as desired.

Fifth, suppose that \.
$x \prec \zf$, \. $y \prec \zf$, and \. $x \. || \. y$.
Then:
{\small \begin{equation}\label{BaseA5}
	\tag{C5}
	\bC(\Pf,2) \quad = \quad
	\begin{pmatrix}
		0   &  \, \ap(x) \, \ap(y)  & 0 & 0   \\
		 \, \ap(x) \, \ap(y) & 0  &   0 & 0   \\
		0 & 0   & 0 & 0  \\
		0  & 0  & 0  & 0
	\end{pmatrix}.
\end{equation}
}
The eigenvalues of this matrix are \. $\{ \, \ap(x) \, \ap(y), 0, 0, -  \, \ap(x) \, \ap(y)\}$,
so it satisfies \eqref{eq:OPE}. Thus, it also satisfies \eqref{eq:Hyp}, as desired.

For the sixth  case, suppose that \. $x \prec \zf \prec y$.  Then:
{\small 	\begin{equation}\label{BaseA6}
		\tag{C6}
		\bC(\Pf,2) \quad = \quad
		\begin{pmatrix}
			0   &0  &   \, \ap(x)  & 0   \\
			0 & 0  &   0 & 0   \\
			  \, \ap(x)  & 0   & 0 & 0  \\
			0  & 0  & 0  & 0
		\end{pmatrix}.
	\end{equation}
}
	The eigenvalues of this matrix are \. $\{ \, \ap(x), 0, 0, -  \, \ap(x)\}$,
	so it satisfies \eqref{eq:OPE}. Thus it also satisfies \eqref{eq:Hyp}, as desired.

Seventh and final case, suppose that \. $x \prec y \prec \zf$.
Then:
{\small
\begin{equation}\label{BaseA7}
	\tag{C7}
	\bC(\Pf,2) \quad = \quad
	\begin{pmatrix}
		 \, \ap(x) \, \ap(y)   &0  & 0  & 0   \\
		0 & 0  &   0 & 0   \\
		0  & 0   & 0 & 0  \\
		0  & 0  & 0  & 0
	\end{pmatrix}.
\end{equation}
}
The eigenvalues of this matrix are \. $\{ \, \ap(x) \, \ap(y), 0, 0, 0\}$,
so it satisfies \eqref{eq:OPE}. Thus it also satisfies \eqref{eq:Hyp}, as desired.
This completes the proof.
\end{proof}

\medskip

\subsection{Proof of Theorem~\ref{t:Sta}}\label{ss:LE-Sta-proof}
We can now prove that the matrix \. $\bC(\Pf,k)$ \. is always hyperbolic.

\smallskip

\begin{prop}\label{p:Stanley}
	Let \ts $\Pf=(X_P,\prec)$ \ts be a finite poset, let $k \in \{2,\ldots,|X_P|-1\}$, and let \ts
$\apr:X \to \rrs$ \ts be an order-reversing weight function.  Then  the matrix \. $\bC\.(\Pf,k)$ \.
satisfies \eqref{eq:Hyp}.
\end{prop}

\smallskip

\begin{proof}
We will prove the proposition  by induction on $|X_P|$.
The base case $|X_P|=3$ follows from Lemma~\ref{l:Stanley BHyp}.
Suppose that the claim is true for $|X_P|-1$.

First note that, if  \ts $\aN(\Pf,k-1) = \aN(\Pf,k)=\aN(\Pf,k+1)=0$\ts,
then  \ts $\bC(\Pf,k)$ \ts is the zero matrix, and \eqref{eq:Hyp} immediately follows.
So we will assume that either one of  \ts $\aN(\Pf,k-1), \aN(\Pf,k),\aN(\Pf,k+1)$ \ts is nonzero.

We split the proof into case (1) and case (2):
For case(1), suppose that at least two of the three numbers are nonzero.
Since the sequence \ts $\aN(\Pf,k-1), \aN(\Pf,k),\aN(\Pf,k+1)$ \ts cannot have internal zeroes (this follows from the demotion argument in the proof of Lemma~\ref{lem:C-center}),
this reduces to either \ts $\aN(\Pf,k-1), \aN(\Pf,k)>0$ \ts or  \ts $\aN(\Pf,k+1), \aN(\Pf,k)>0$\ts.
	By switching to the dual poset in \eqref{eq:tango 4} if necessary, we can without loss of generality assume that
\ts $\aN(\Pf,k-1), \aN(\Pf,k)>0$\ts.
We split the proof further into case (1a), case (1b), and case (1c).

For case (1a), assume that $k\geq 3$.
Let \ts $\AA(\Pf,k)$ \ts be the atlas defined in \S\ref{ss:LE-atlas}.
It follows from Lemma~\ref{l:Stanley Inh}, \ref{l:Stanley TPInv}, \ref{l:Stanley Knonneg}
that this atlas satisfies the assumptions of Theorem~\ref{t:Hyp} (note that these lemmas require $k\geq 3$).
Also note that every sink vertex in $\Vf^0$ satisfies \eqref{eq:Hyp} by the induction assumption, as they correspond to posets with cardinality $|X_P|-1$.
It then follows from Theorem~\ref{t:Hyp} that every regular vertex in $\Vf^1$ satisfies \eqref{eq:Hyp}.
On the other hand, it follows 	from Lemma~\ref{l:Stanley Irr} that every   $\vf=t\in \Vfp$ \.
with  \. $0<t<1$ \. is a regular vertex.
This implies that, for \. $0<t<1$\., the matrix \. $t \, \bC(\Pf,k)  \. +  \. (1-t) \, \bC(\Pf,k-1)$ \. satisfies \eqref{eq:Hyp}.
By taking the limit $t\to 0$ and $t\to 1$, we then conclude that both
\. $\bC(\Pf,k)$ \. and \. $\bC(\Pf,k-1)$ \. satisfies \eqref{eq:Hyp}, as desired.

For case (1b), assume that $k=2$ and $\aN(\Pf,k+1)>0$.
Then by applying the same argument as in case (1a) to the atlas \. $\AA(\Pf,k+1)$\.,
it follows that both
\. $\bC(\Pf,k+1)$ \. and \. $\bC(\Pf,k)$ \. satisfies \eqref{eq:Hyp}, as desired.

For case (1c), assume that \ts $k=2$ \ts  and \ts $\aN(\Pf,k+1)=0$\ts.
 The assumptions imply that $X_P$ can be partitioned into \. $\{x\} \. \cup \. \{z\} \. \cup \{T\}$\., where $x$ is the only element in $X_P$ incomparable to $z$, and $T$ is the upper ideal of $z$ in $\Pf$.
 Also note that the support of  \. $\bC(\Pf,k)$ \. is contained in \. $\{x_\low\} \cup T_{\up}$\..
Now suppose that there exists \ts $y \in \min(T)$ \ts such that \ts $x \. || \. y$.
Let  \ts $\Pf':=(X_P,\prec')$ \ts  be the poset  with the same ground set as $\Pf$ and with $\prec'$ being obtained from $\prec$ by removing the relation \ts $z\prec y$\ts.
Now note that  \ts $\aN(\Pf,k+1)>0$ \ts by construction, so it follows from case (1b) that \ts $\bC(\Pf',k)$ \ts satisfies \eqref{eq:Hyp}.
On the other hand,
it follows from the construction that \ts $\bC(\Pf,k)$ \ts is equal to
\ts $\bC(\Pf',k)$ \ts when restricted to rows and columns indexed by \. $\{x_\low\} \cup T_{\up}$\..
Since \eqref{eq:Hyp} is preserved under restricting to principal submatrices, we have
\ts $\bC(\Pf,k)$ \ts also satisfies \eqref{eq:Hyp}, as desired.
Now suppose that every element $T$ is ordered to be greater than $x$ by $\prec$.
This implies that, for every \. $y \in \max(\Pf,up) \. \cap \. T_{\up}$\.,
\[   \ap(\Ec(\Pf-x-y,k-1)) \ = \  \ap(\Ec(\Pf-y,k-1)), \]
because $x$ is the second smallest element in every linear extension counted by \. $\Ec(\Pf-y,k-1)=\Ec(\Pf-y,1)$.
This implies that
\begin{align*}
(\bC(\Pf,k))_{x_\low, y_\up} \ &=_{\eqref{eq:C-1}} \   \ap(x) \. \ap(\Ec(\Pf-x-y,k-1)) \ = \ \ap(x) \. \ap(\Ec(\Pf-y,k-1)) \\
&=_{\eqref{eq:C-sum-2}} \ap(x) \. \sum_{w \in Z_{\up}} (\bC(\Pf,k))_{w, y_\up}.
\end{align*}
Let \ts $\bD$ \ts be the matrix obtained by deducting the $x_{\low}$-row by the $\ap(x)$ times the sum of the other rows, followed by deducting the the $x_{\low}$-column by the $\ap(x)$ times sum of the other columns (note that this operation preserves \eqref{eq:Hyp}).
Then it follows from the previous equation that
the entries of \ts $\bD$ \ts are given by
\begin{align*}
	(\bD)_{u \ts v} \ = \
	\begin{cases}
		(\bC(\Pf,k))_{u \ts v} & \text{ if } u,v \in Z- x,\\
		0& \text{ if } u \in Z-x \text{ and } v=x,\\
		- \ap(x) \. \ap(\Ec(\Pf-x,k-1)) & \text{ if } u =v=x.
	\end{cases}
\end{align*}
It then follows that $\bC(\Pf,k)$ satisfies \eqref{eq:Hyp} if and only if, the restriction of $\bD$ to rows and columns indexed by $T_{\up}$, satisfies \eqref{eq:Hyp}.
Now let \. $\Pf'$ \. be the induced subposet of $\Pf$ on \. $X_P-x$\..
Note that $\bC(\Pf',k)$ satisfies \eqref{eq:Hyp} by the induction assumption.
Also note that \ts $\bC(\Pf',k)$ \ts is equal to
\ts $\bD$ \ts when restricted to rows and columns indexed by \. $T_{\up}$\..
Since \eqref{eq:Hyp} is preserved under restricting to principal submatrices, it follows that
this submatrix of $\bD$ also satisfies \eqref{eq:Hyp}.
Combined with previous observations, we then conclude that $\bC(\Pf,k)$ satisfies \eqref{eq:Hyp}, as desired.
This completes the proof of case (1).

For case (2), suppose that exactly one of \ts $\aN(\Pf,k-1), \aN(\Pf,k),\aN(\Pf,k+1)$ \ts is nonzero.
The proof then splits into three subcases.
For case (2a), let \ts $\aN(\Pf,k)=0$ \ts while
\ts $\aN(\Pf,k-1) = \aN(\Pf,k+1)=0$\ts.
This implies that \. $X_P$ \. can be partitioned into \ts $S \cup T \cup \{x\}$\ts,
where  \ts $S$ \ts is the lower ideal of $z$ in $\Pf$ and \ts $|S|=k-1$\ts, and  \ts $T$ \ts is the upper ideal of $x$ in $\Pf$.
Then the entries of \ts $\bC(\Pf,k)$ \ts is given by
\begin{align*}
	(\bC(\Pf,k))_{x,y} \ = \
	\begin{cases}
		\ap(x) & \text{ if } x \in \min(\Pf,\low) \cap S_{\low}, \.  y \in \max(\Pf,\up) \cap T_{\up},\\
		0 & \text{ otherwise}.
	\end{cases}
\end{align*}
By a direct computation, the eigenvalues of this matrix are \. $\lambda,-\lambda,0,\ldots,0$\.,
where
\[ \lambda \ := \  |\max(\Pf,\up) \cap T_{\up}| \. \sum_{x \in \min(\Pf,\low) \cap S_{\low}} \ap(x). \]
Thus this matrix satisfies \eqref{eq:OPE}, and so it satisfies \eqref{eq:Hyp}.

For case (2b), let  \. $\aN(\Pf,k+1)>0$ \. while $\aN(\Pf,k-1)=\aN(\Pf,k)=0$.
Let $S$ be the lower ideal of $z$ in $\Pf$. This implies that \ts $|S|=k$\ts,
and the support of \ts $\bC(\Pf,k)$ \ts is contained in \ts $S_{\low}$\ts.
Now let \ts $\Pf':=(X_{\Pf'},\prec')$ \ts  be the poset with ground set \. $X_{\Pf'} := S \cup \{z\} \subseteq X$\., and with relations $\prec'$ given by
\begin{alignat*}{2}
	&\forall \ x,y \in S, \qquad &&x \. \prec' \. y  \quad \Longleftrightarrow \quad  x \. \prec \. y,\\
	&\forall \ x \in S, \qquad &&x \. || \. z.
\end{alignat*}
It  follows from the construction that, for all  $x,y\in S_\low$,
\[  \big(\bC(\Pf,k)\big)_{x_\low,y_\low}  \. =  \. \Ec(\Pf-S-x) \big(\bC(\Pf',k)\big)_{x_\low,y_\low}\..  \]
Since \eqref{eq:Hyp} is a property that is preserved by restricting to principal submatrices,
we have that  \. $\bC(\Pf,k)$ \. satisfies \eqref{eq:Hyp} if \. $\bC(\Pf',k)$ \. also satisfies
\eqref{eq:Hyp}.
Also note that \ts $\aN(\Pf',k-1)>0$\ts, \ts $\aN(\Pf',k)>0$.
Thus by the same argument as in case (1), it follows that \ts $\bC(\Pf',k)$ \ts satisfies \eqref{eq:Hyp}, which in turn implies that \ts $\bC(\Pf,k)$ \ts satisfies \eqref{eq:Hyp}.

Finally, for case (2c), let \. $\aN(\Pf,k-1)>0$ \. while $\aN(\Pf,k+1)=\aN(\Pf,k)=0$.
This case follows by applying the same argument as in case (2b) to the dual poset $\Pf^{\op}$ in \eqref{eq:tango 4}.
This completes the proof.
\end{proof}

\smallskip

\begin{proof}[Proof of Theorem~\ref{t:Sta-weighted}]
It follows from Proposition~\ref{p:Stanley} that the matrix \. $\bC(\Pf,k)$ \. satisfies \eqref{eq:Hyp}, and it follows from \eqref{eq:Cfg} that the theorem is a special case of
 \. $\bC(\Pf,k)$ \.  satisfying
\eqref{eq:Hyp}. This completes the proof.
\end{proof}

\medskip

\subsection{Example of a combinatorial atlas}\label{ss:LE-example}
%
%
Let $\Pf$ be the poset on $X=\{a,b,c,d,z\}$, with the order given by \.
$a \prec b \prec c$, \. $a \prec z$, \. $d \prec c$.
Fix \. $z \in \Pf$ \. as in Stanley's inequality,
and with uniform weight on all linear extensions.
Let  \ts $k=3$.
Then the matrices $\bC(\Pf,k)$ and $\bC(\Pf, k+1)$ are given by

{\small
\begin{align*}
	\bC(\Pf, 3) \ = \
	\begin{pmatrix}
		1 & 0 & 0 & 1 & 0 & 0 & 2 & 0\\
		0 & 0 & 0 & 0 & 0 & 0 & 0 & 0\\
		0 & 0 & 0 & 0 & 0 & 0 & 0 & 0\\
		1 & 0 & 0 & 0 & 0 & 0 & 1 & 0\\
		0 & 0 & 0 & 0 & 0 & 0 & 0 & 0\\
		0 & 0 & 0 & 0 & 0 & 0 & 0 & 0\\
		2 & 0 & 1 & 0 & 0 & 0 & 2 & 0\\
		0 & 0 & 0 & 0 & 0 & 0 & 0 & 0	
	\end{pmatrix},
	\qquad
	\bC(\Pf, 4) \ = \
	\begin{pmatrix}
		1 & 1 & 0 & 0 & 0 & 0 & 2 & 0\\
		1 & 0 & 0 & 0 & 0 & 0 & 0 & 0\\
		0 & 0 & 0 & 0 & 0 & 0 & 0 & 0\\
		0 & 0 & 0 & 0 & 0 & 0 & 1 & 0\\
		0 & 0 & 0 & 0 & 0 & 0 & 0 & 0\\
		0 & 0 & 0 & 0 & 0 & 0 & 0 & 0\\
		2 & 0 & 0 & 1 & 0 & 0 & 3 & 0\\
		0 & 0 & 0 & 0 & 0 & 0 & 0 & 0	
	\end{pmatrix},
\end{align*}
}

\nin
where the rows and columns are labeled by \.
$\bigl\{a_{\min}, b_{\min}, c_{\min}, d_{\min}, a_{\max}, b_{\max}, c_{\max}, d_{\max}\bigr\}$.
In this notation, we have:
\[
\fb \ = \ (1,1,1,1,0,0,0,0)^\intercal, \qquad  \gb \ = \ (0,0,0,0,1,1,1,1)^\intercal.
\]
Recall that the inner products of these two vectors with the matrix \. $\bC(\Pf, k)$ \. has the following combinatorial interpretation by \eqref{eq:Cfg}:
\[
\langle \fb, \bC(\Pf, 3) \fb \rangle \ = \    \aN(4) , \qquad   \langle \fb, \bC(\Pf, 3) \gb \rangle \ = \  \, \aN(3), \qquad  \langle \gb, \bC(\Pf, 3) \gb \rangle \ = \   \aN(2).
\]
Stanley's inequality~\eqref{eq:Sta} is equivalent to
\begin{equation}\label{eqhyp graph}
	\langle \fb,  \bC(\Pf,k) \gb \rangle^2   \  \geq \  \langle \fb,  \bC(\Pf,k) \fb \rangle \.\cdot\. \langle \gb,  \bC(\Pf,k) \gb \rangle.
\end{equation}
In this example, we have:
\[ \aN(4) \ = \ 3, \qquad \aN(3) \ = \ 3, \qquad \aN(2) \ = \ 2, \]
and the log-concavity in Stanley's inequality holds: \. $3^2 \geq 3 \times 2$.

\medskip

\subsection{Posets with belts}\label{ss:LE-belts}
Let \ts $\cP=(X, \prec)$ \ts be a poset with $|X|=n$ elements.
We say that \ts $\cP$ \ts has \defn{belt} at \ts $\zf\in X$ \ts if \ts $\inc(\zf)$ \ts
is either empty or a chain in~$\cP$.  Note that \ts $\width(\cP)=2$ \ts if and only if
$\cP$ has a belt at every element \ts $\zf \in X$.
Below we show how to strengthen Theorem~\ref{t:Sta-weighted} for posets
with belts.

Let \ts $\ap:X\to \rrs$ \ts be an order-reversing weight function defined by~\eqref{eq:Rev},
and fix element \ts $\zf\in X$.
Rather than use multiplicative formula \eqref{eq:sta-def-omega-weighted} to extend $\ap$
to~$\Ec$, we define \. $\ap: \Ec\to \rrs$ \. by the \defn{tropical formula}\ts:
\begin{equation}\label{eq:sta-def-omega-tropical}
\aq(L) \, := \, \max \ts \bigl\{\ts\ap(x) \. : \. L(x)< L(\zf)\ts\bigr\}.
\end{equation}

\smallskip

\begin{thm}[{\rm \defng{Tropical Stanley inequality for posets with belts}}{}]\label{t:Sta-belt}
Let \ts  $\cP=(X,\prec)$ \ts be a poset with \ts $|X|=n$ \ts elements, and suppose
\ts $\cP$ \ts has a belt at \ts $\zf\in X$.  Let \ts $\apr: X \to \rrs$ \ts be a
positive order-reversing weight function.  Define \. $\aqr: \Ec\to \rrs$ \. by the
tropical formula~\eqref{eq:sta-def-omega-tropical}.  Then, for every \. $1< k < n$,
we have:
\begin{equation}\label{eq:Sta-weighted-tropical}
\aNr_\aqr(k)^2 \,\. \ge \,\. \aNr_\aqr(k-1) \.\cdot \. \aNr_\aqr(k+1)\ts,
\end{equation}
where \ts $\aNr_\aqr(k)$ \ts is defined by~\eqref{eq:sta-def-N-weighted}.
\end{thm}


More generally, let \ts $\ap: \text{Low}(\cP) \to \rrs$  \ts be a weight
function on the set of \emph{lower ideals} of the poset~$\cP$.  Suppose
\ts $\ap$ \ts satisfies the following (\defng{submodular property})\ts:
\begin{equation}\label{eq:sta-omega-submod}\tag{Submod}
	\ap\bigl(S+x+y\bigr) \.\cdot \.  \ap(S)  \  \leq \ \omega\big(S+x\big)^2,
\end{equation}
for all \ts $x,y \in \inc(\zf)$, \ts $x\prec y$, and for all \ts $S\ssu X$ \ts such that
\. $S, \ts S+x, \ts S+x+y \in \ts \text{Low}(\cP)$.  We can then
define
\begin{equation}\label{eq:sta-def-omega-sub}
\aq(L) \, := \, \ap(A), \quad \text{where} \quad A := \bigl\{x\in X\.:\.L(x)\prec L(\zf)\bigr\}.
\end{equation}

\smallskip

\begin{thm}[{\rm \defng{Submodular Stanley inequality for posets with belts}}{}]\label{t:Sta-belt-submodular}
Let \ts  $\cP=(X,\prec)$ \ts be a poset with \ts $|X|=n$ \ts elements, and suppose
\ts $\cP$ \ts has a belt at \ts $\zf\in X$.  Let \. $\ap: \Low(\cP) \to \rrs$
\ts be a positive weight function on the set of lower ideals of~$\cP$ which
satisfies~\eqref{eq:sta-omega-submod}.  Define \. $\aqr: \Ec\to \rrs$ \.
by~\eqref{eq:sta-def-omega-sub}.  Then, for every \. $1< k < n$,
we have:
\begin{equation}\label{eq:Sta-weighted-submodular}
\aNr_\aqr(k)^2 \,\. \ge \,\. \aNr_\aqr(k-1) \.\cdot \. \aNr_\aqr(k+1)\ts, \quad
\text{where} \quad
\aNr_\aqr(m) \, := \, \sum_{L \in \Ec_m} \. \aqr(L)\., \quad \text{for all} \ \ 1\le m \le n\ts.
\end{equation}
\end{thm}


\begin{proof}[Proof of Theorem~\ref{t:Sta-belt-submodular}]
	The result follows the same argument as Theorem~\ref{t:Sta-weighted}
	with two changes.  First, in the proof of Lemma~\ref{l:Stanley BHyp}, the case \eqref{BaseA1}
	does not need to be verified since~$\Pf$ has a belt.  Second, the case \eqref{BaseA2} is
	instead verified  through \eqref{eq:sta-omega-submod}.  We omit the details.
\end{proof}
\smallskip

\begin{proof}[Proof of Theorem~\ref{t:Sta-belt}]
The result is a direct consequence of Theorem~\ref{t:Sta-belt-submodular},
as the tropical weight function in \eqref{eq:sta-def-omega-tropical} clearly
satisfies~\eqref{eq:sta-omega-submod}.  The details are straightforward.
\end{proof}

\bigskip


\section{Proof of equality conditions for linear extensions}
\label{s:proof-Sta-equality}

In this section we extend and prove Theorem~\ref{t:Sta-equality-weighted},
see also~$\S$\ref{ss:hist-LE-equality}.

\smallskip

\subsection{More equality}\label{ss:proof-Sta-equality-more}
Let \ts  $\Pf=(X_P,\prec)$ \ts be a poset on \ts $|X_P|=n$ \ts elements,
and let \ts $\ap:X\to\rrs$ \ts be the order-reversing weight.
Recall that in the notation of this section, \. $\aN_\ap(\Pf,k)=\aN_{\ap}(k)$, \. with the latter as defined
in~\eqref{eq:sta-def-N-weighted}.

We add two more items to Theorem~\ref{t:Sta-equality-weighted} and reformulate
it in terms of words, to prove a stronger result:

\smallskip

\begin{thm}[{\rm cf.\ Theorem~\ref{t:Sta-equality-weighted}}{}]\label{t:Sta-equality-weighted-more}
Let \ts  $\Pf=(X_P,\prec)$ \ts be a poset with \ts $|X_P|=n$ \ts elements, and
let \ts $\ap: X_P \to \rrs$ \ts be a positive order-reversing weight function.
Fix element \ts $\zf \in X_P$\ts.  Suppose that \ts $\aNr_\apr(\Pf,k)>0$.
Then \. \underline{the following are equivalent}\ts:
\begin{enumerate}
			[{label=\textnormal{({\alph*})},
		ref=\textnormal{\alph*}}]
\item \label{item:equality Stanley a}
 \ $\aNr_\apr(\Pf,k)^2 \. = \. \aNr_\apr(\Pf,k-1) \. \cdot \.\aNr_\apr(\Pf,k+1)$,
\item \label{item:equality Stanley b}
\. there exists \. $\asr=\asr(k,\zf) >0$, s.t.\ 
$$\aNr_\apr(\Pf,k+1) \. = \. \asr \. \aNr_\apr(\Pf,k) \. = \. \asr^2 \. \aNr_\apr(\Pf,k-1),
$$
\item \label{item:equality Stanley c}
\. there exists  \. $\asr=\asr(k,\zf) >0$, s.t.
		 \[ \aNr_{\ap}(\Pf-S-T,3) \ = \  \asr \, \aNr_{\ap}(\Pf-S-T,2) \ = \ \asr^2 \, \aNr_{\ap}(\Pf-S-T,1) \]
		 for every lower set $S$ and upper set $T$ of $\Pf-z$  satisfying \. $|S|=k-2$,  \. $|T|=n-k-1$.
\item \label{item:equality Stanley d}
		there exists \. $\asr=\asr(k,\zf) >0$, s.t.\ \. $\apr(x_{k-1}) = \apr(x_{k+1}) =  \asr$, and
for every \. $\gamma= x_1\.\cdots\.x_n\in \Ec_k$, we have \.
$\zf \. || \. x_{k-1}$\ts, \. $\zf \.|| \. x_{k+1}$\ts.
\item \label{item:equality Stanley e}
\. there exists \. $\asr=\asr(k,\zf) >0$, s.t.\ \. $\apr(x_{k-1}) = \apr(x_{k+1}) =  \asr$,
\. $f(x)>k$ \. for all \. $x\succ \zf$,
and \. $g(x)>n-k+1$ \. for all \. $x\prec \zf$,
\end{enumerate}
\end{thm}

\smallskip

The direction \. $\eqref{item:equality Stanley b} \Rightarrow \eqref{item:equality Stanley a}$ \.
is  trivial.
For the direction \. $\eqref{item:equality Stanley c} \Rightarrow \eqref{item:equality Stanley b}$,
note that we have
\begin{equation}\label{eq:LE-equality-double-sum}
\aN_\ap(k) \, = \, \sum_{S,T}  \ap(S) \, |\Ec(S)| \. |\Ec(T)| \.  \aN_\ap(\Pf-S-T,2),
\end{equation}
summed over all lower sets $S$ and upper sets $T$ of \ts $\Pf-z$ \ts satisfying \. $|S|=k-2$,  \. $|T|=n-k-1$.
Note that the analogous formulas also hold for \. $\aN_\ap(k\pm 1)$.
Together with \eqref{item:equality Stanley c}, this implies that
\begin{equation}\label{eq:extra 1}
	\aN_{\ap}(k) \ \leq \ \as  \, \aN_{\ap}(k-1), \qquad  \aN_{\ap}(k) \ \leq \ \frac{1}{\as}  \, \aN_{\ap}(k+1).
\end{equation}
This in turn implies that \. $\aN_{\ap}(k)^2 \. \leq \. \aN_{\ap}(k-1) \, \aN_{\ap}(k+1)$.
On the other hand, by Theorem~\ref{t:Sta-weighted} we already know the inequality in the opposite direction: \.
$\aN_{\ap}(k)^2 \. \geq \. \aN_{\ap}(k-1) \, \aN_{\ap}(k+1)$.
This implies the equality in \eqref{eq:extra 1}, which in turn implies  \eqref{item:equality Stanley b}, as desired.

Below we prove
$\eqref{item:equality Stanley a} \Rightarrow \eqref{item:equality Stanley c}$,
$\eqref{item:equality Stanley c} \Leftrightarrow \eqref{item:equality Stanley d}$,
and $\eqref{item:equality Stanley d} \Leftrightarrow \eqref{item:equality Stanley e}$,  thus completing the proof of Theorem~\ref{t:Sta-equality-weighted-more}.

\medskip

\subsection{Proof of \. $\eqref{item:equality Stanley a} \Rightarrow \eqref{item:equality Stanley c}$}
We start with the following preliminary result.

\smallskip

\begin{lemma}\label{l:equality Stanley 1}
Let \ts $\Pf=(X_P,\prec)$ \ts be a finite poset, and let \ts $\apr:X \to \rrs$ \ts
be an order-reversing  weight function, and let $k \in \{3,\ldots, |X_P|-1\}$.
Suppose that there exists $\asr>0$, such that
\[
\aNr_{\ap}(\Pf,k+1) \ = \  \asr \, \aNr_{\ap}(\Pf,k) \ = \ \asr^2 \, \aNr_{\ap}(\Pf,k-1) \  > \ 0.
\]
Then, for every \. $x \in \min(\Pf)$,
\begin{equation}\label{eq:LE-Sta-eq-lemma}
\aNr_{\ap}(\Pf-x,k) \ = \  \asr \, \aNr_{\ap}(\Pf-x,k-1) \ = \ \asr^2 \, \aNr_{\ap}(\Pf-x,k-2) \ > \ 0.
\end{equation}
\end{lemma}

\smallskip

\begin{proof}
	Let \ts $\AA(\Pf,k)$ \ts be the combinatorial atlas
	defined in~\S\ref{ss:LE-atlas}.
	 Let \. $\fb, \gb \in \Rb^{\ar}$ \. be
	the characteristic vectors of \ts $\Zlow$ \ts and \ts $\Zup$, respectively.
	Clearly, \. $\fb,\gb$ \. is a global pair for \ts $\AA$, i.e., they satisfy
	\eqref{eq:PosGlob}.  This allows us to apply Theorem~\ref{t:equality Hyp}
	in the reductions below.
	
	Let $\vf=t=1 \in \Vf^1$. It then follows from the assumptions of the lemma and \eqref{eq:Cfg} that the vertex~$\vf$ satisfies \eqref{eq:sEqu}.
	Also note that $\vf$ satisfies $\eqref{eq:Hyp}$ by Proposition~\ref{p:Stanley}.
	On the other hand, it is straightforward to verify that $\vf$ is a functional vertex of~$\Qf$, i.e.\
	it satisfies \eqref{eq:ProjGlob} and \eqref{eq:hGlob}.  By Theorem~\ref{t:equality Hyp},
	every functional target of~$\vf$ also
	satisfies \eqref{eq:sEqu} with the same \ts $\as >0$.
	On the other hand, it is easy to see that the functional targets of~$\vf$
	include  vertices of the form  \. $x \in \Vf^0$,
	where $x =x_{\low} \in \min(\Pf,\low)$.
	Hence $x$ satisfies \eqref{eq:sEqu}, which implies that
	\begin{equation*}
	\langle \fb, \bC(\Pf-x, k-1) \fb \rangle  \ = \  \as \, \langle \fb,\bC(\Pf-x, k-1) \gb \rangle \ = \  \as^2 \, \langle \gb , \bC((\Pf-x, k-1) \gb \rangle.
\end{equation*}
	It now follows from \eqref{eq:Cfg} that
	\begin{equation*}
		\aN_{\ap}(\Pf-x,k) \ = \  \asr \, \aN_{\ap}(\Pf-x,k-1) \ = \ \asr^2 \, \aN_{\ap}(\Pf-x,k-2).
	\end{equation*}
Also note that \. $\aN_{\ap}(\Pf-x,k-1)>0$ \. by Lemma~\ref{lem:C-supp}. The proof is now complete.
\end{proof}

\smallskip

We can now prove\. $\eqref{item:equality Stanley a} \Rightarrow \eqref{item:equality Stanley c}$.
Since \. $\aN_\ap(\Pf,k) >0$,  it follows from \eqref{item:equality Stanley a} that
\[
\aN_\ap(\Pf,k+1) \ = \  \as \, \aN_\ap(\Pf,k) \ = \ \as^2 \, \aN_\ap(\Pf,k-1) \ > \ 0 \qquad
\text{for} \qquad \as\,:=\, \frac{\aN_\ap(k+1)}{\aN_\ap(k)}\, >\. 0\ts.
\]
		 By applying  Lemma~\ref{l:equality Stanley 1}
		 for $k-2$ many times, we have
		 \begin{equation}\label{eq:tango 3}
		 	 \aN_{\ap}(\Pf-S,3) \ = \  \as \, \aN_{\ap}(\Pf-S,2) \ = \ \as^2 \, \aN_{\ap}(\Pf-S,1) \ > \ 0,
		 \end{equation}
	 for every lower set $S$ of \ts $\Pf-z$ \ts satisfing $|S|=k-2$.
	 Recall that $n:=|X_P|$.
	 Now note that, by applying \eqref{eq:tango 4}, we have that \eqref{eq:tango 3} is equivalent to
	  		 \begin{equation}\label{eq:tango-new-1}
	  	\aN_{\ap^{\op}}(\Pf^{\op}-S,n-k) \ = \  \as \, \aN_{\ap^{\op}}(\Pf^{\op}-S,n-k+1)  \ = \ \as^2 \, \aN_{\ap^{\op}}(\Pf^{\op}-S,n-k+2) \ > \ 0.
	  \end{equation}
	  By applying   Lemma~\ref{l:equality Stanley 1}
	  for another $n-k-1$ many times,
	  \eqref{eq:tango-new-1} is equivalent to
	  	  		 \begin{equation}\label{eq:tango-new-2}
	  	\aN_{\ap^{\op}}(\Pf^{\op}-S-T,1) \ = \  \as \, \aN_{\ap^{\op}}(\Pf^{\op}-S-T,2)  \ = \ \as^2 \, \aN_{\ap^{\op}}(\Pf^{\op}-S-T,3) \ > \ 0,
	  \end{equation}
  	 for every upper set $T$ of \ts $\Pf-z$ \ts satisfing $|T|=n-k-1$.
  	 Finally, by applying \eqref{eq:tango 4} again,
  	 it follows that \eqref{eq:tango-new-2} is equivalent to
  	 	  	  		 \begin{equation*}
  	 	\aN_{\ap}(\Pf-S-T,3) \ = \  \as \, \aN_{\ap}(\Pf-S-T,2)  \ = \ \as^2 \, \aN_{\ap}(\Pf-S-T,1) \ > \ 0,
  	 \end{equation*}
   and the proof is now complete.\qed

\medskip

\subsection{Proof of \. $\eqref{item:equality Stanley c} \Rightarrow \eqref{item:equality Stanley d}$}
Let \. $S:=\{x_1,\ldots, x_{k-2}\}$ \. and \. $T:=\{x_{k+2},\ldots,x_n\}$\..
Let \. $\fb, \gb \in \Rb^{\ar}$ \. be
the characteristic vectors of \ts $\Zlow$ \ts and \ts $\Zup$, respectively.
It follows from \eqref{eq:Cfg}
and \eqref{item:equality Stanley c} that
\begin{equation}\label{eq:tango 6}
	  \langle \fb, \bC(\Pf-S-T,2) \vb \rangle  \ = \ \as \,   \langle \vb, \bC(\Pf-S-T,2) \vb \rangle \ = \ \as^2 \,  \langle \wb, \bC(\Pf-S-T,2) \wb \rangle.
\end{equation}
for some \. $\as>0$.


Let \. $\zb \ts :=\ts \fb - \as \gb$.
It follows from \eqref{eq:tango 6}  that
\. $\< \zb, \bC(\Pf-S-T,2) \zb \> = 0$.
By  Lemma~\ref{l:equality kernel}, this implies \. $\bC(\Pf-S-T,2) \. \zb \. = \. \0$.
On the other hand, the matrix \. $\bC(\Pf-S-T,2)$ \. is one of the seven matrices
in \eqref{BaseA1}--\eqref{BaseA7} because $\Pf-S-T$ is a poset with three elements.
From the seven matrices, only \eqref{BaseA1} and \eqref{BaseA2}
can have \. $\bC(\Pf-S-T,2) \zb =\0$.
Now note that in both cases we have \. $x\.||\.\zf$ \. and $y\.||\.\zf$.
By a direct calculation,  in both cases we have:
$$
\bC(\Pf-S-T,2) \zb =\0 \quad \Longleftrightarrow \quad \as=\ap(x)=\ap(y).
$$
This proves \eqref{item:equality Stanley d}, as desired. \qed

\medskip

\subsection{Proof of \. $\eqref{item:equality Stanley d} \Rightarrow \eqref{item:equality Stanley c}$}
It follows from  \eqref{item:equality Stanley d}, that, for every lower set $S$ and upper set $T$ of $\Pf-z$  satisfying \. $|S|=k-2$,  \. $|T|=n-k-1$.,
we have:
	\begin{equation}\label{eq:insomnia 1}
		\aN_{\ap}(\Pf-S-T,2) \ \leq \ \as \, \aN_{\ap}(\Pf-S-T,1) \qquad \text{and} \qquad \aN_{\ap}(\Pf-S-T,2) \ \leq \ \frac{1}{\as} \, \aN_{\ap}(\Pf-S-T,3),
	\end{equation}
where \ts $\as>0$ \ts is given in~\eqref{item:equality Stanley d}.
Summing over all such \ts $S,T$ \ts as in~\eqref{eq:LE-equality-double-sum},
we obtain:
\begin{equation}\label{eq:insomnia 2}
	\aN_{\ap}(k) \ \leq \ \as  \, \aN_{\ap}(k-1), \qquad  \aN_{\ap}(k) \ \leq \ \frac{1}{\as}  \, \aN_{\ap}(k+1).
\end{equation}
This implies that \. $\aN_{\ap}(k)^2 \. \leq \. \aN_{\ap}(k-1) \, \aN_{\ap}(k+1)$.
On the other hand, by Theorem~\ref{t:Sta-weighted} we already know the inequality in the opposite direction: \.
$\aN_{\ap}(k)^2 \. \geq \. \aN_{\ap}(k-1) \, \aN_{\ap}(k+1)$.
This implies the equality in~\eqref{eq:insomnia 2}, which in turn implies the
equality in~\eqref{eq:insomnia 1}, as desired. \qed

\medskip

\subsection{Proof of \. $\eqref{item:equality Stanley d} \Leftrightarrow \eqref{item:equality Stanley e}$}
Note that both items have the same weight function assumption, which reduces the claim to the following
lemma of independent interest.

\smallskip

\begin{lemma}\label{l:Sta-equality-comb}
	Let \ts  $\Pf=(X,\prec)$ \ts be a poset with \ts $|X|=n$ \ts elements.
	Fix element \ts $\zf \in X$ \ts and suppose that \ts $\aNr(k)>0$.
	Then \. \underline{the following are equivalent}:
	\begin{enumerate}			
		\item \label{item:equality Stanley new 1}
		 $f(y)>k$ \. for all \. $y\succ \zf$, \. and \. $g(y)>n-k+1$ \. for all \. $y\prec \zf$.
		\item \label{item:equality Stanley new 2}
for every \. $\gamma= x_1\.\cdots\.x_n\in \Ec_k$\ts, we have \.
$\zf \. || \. x_{k-1}$ \. and \. $\zf \.|| \. x_{k+1}$\ts.
	\end{enumerate}
\end{lemma}

\begin{proof}
We first prove 	\. \eqref{item:equality Stanley new 1} $\Rightarrow$ \eqref{item:equality Stanley new 2}.
Suppose to the contrary that \ts $\zf$ \ts is comparable to \. $y:=x_{k+1}$.
Then \ts $\zf \prec y$, and it follows from \eqref{item:equality Stanley new 1} that \. $f(y) > k$.
This implies that there are at least \ts $(k+1)$ \ts elements in \ts $\gamma$ \ts that appear before~$y$,
contradicting the assumption that \ts $y=x_{k+1}$.  An analogous argument shows that $\zf$ is incomparable
to \ts $x_{k-1}$\ts.


We now prove \. \eqref{item:equality Stanley new 2} $\Rightarrow$ \eqref{item:equality Stanley new 1}.
Let $y \in X$ be such that $\zf \prec y$,
and suppose to the contrary that $f(y) \leq k$.
Let $Q \subseteq X$ be given by
\[Q \ := \ \{ \ts  x \in X  \. : \.  x \prec \zf, \.  x \prec y \ts \}. \]
Note that \. $|Q| \leq f(y)-1 \leq k-1$, and that $Q$ is a lower ideal of~$\Pf$.
Let \. $R \subseteq X$ \. be given by
\[R \ := \ \big\{ \ts  x \in X  \. : \.  x \prec z \. \text{ or \. $x\. || \.z$}  \ts \big\}. \]
Note that $R$ is a lower ideal of $\Pf$, that \. $z, y\notin R$, and that \. $Q \subseteq R$.
Also note $|R| = n-g(\zf)-1$.
Since \. $g(\zf) \leq n-k$ \. by the assumption that \ts $\aN(k)>0$, it  follows that \.
$|R| \geq k-1$.

We conclude that there exists a lower ideal $U$ of $\Pf$ such that \. $Q \subseteq U \subseteq R$ and $|U|=k-1$.
This in turn implies that
there exists a linear extension \. $\gamma=x_1\.\cdots\. x_n \in \Ec$, such that
\[ U\, = \, \{x_1,\ldots, x_{k-1} \}, \quad  x_k \ = \ \zf, \quad x_{k+1} \, = \, y.    \]
It then follows from \eqref{item:equality Stanley new 2} that $\zf$ and $y$ are incomparable, and we get a contradiction.
The same argument shows that \. $g(y)>n-k+1$ \. for all \. $y\prec \zf$.
This completes the proof of the lemma.
\end{proof}

\bigskip


{\small

\section{Historical remarks} \label{s:hist}

\subsection{}\label{ss:hist-pre}
Unimodality is surprisingly difficult to establish even in some classical cases.
For example, Sylvester in 1878 famously resolved Cayley's 1856 conjecture on unimodality
of $q$-binomial coefficients \ts $\binom{n}{k}_q$ \ts using representations
of \ts $\SL(2,\cc)$, see~\cite{Syl}.
In 1982, a linear algebraic deconstruction was obtained by Proctor~\cite{Pro}.
The first purely combinatorial proof was obtained O'Hara's~\cite{OH} only in 1990,
while the strict unimodality for \ts $k,n-k\ge 8$ \ts
was proved in~2013, by the second author and Panova~\cite{PP}.

Log-concavity is an even harder property to establish.  Over the years, a number
of tools and techniques for log-concavity were found, across many areas of
mathematics and applications, from elementary combinatorial
to analytic, from Lie theoretic to topological.
As Huh points out in~\cite{Huh}, sometimes there is only one known approach
to the problem.
We refer to surveys
\cite{Bre,Bre2,Sta2} for an overview of classical unimodality and
log-concavity results, to~\cite{Bra} for a more recent overview
emphasizing enumerative results and analytic methods, and to~\cite{SW}
for a survey on the role of log-concavity in analysis and probability.

\subsection{}\label{ss:hist-main}
Mason's matroid log-concavity conjectures were stated in~\cite{Mas},
motivated by the earlier work and conjectures in graph theory and
combinatorial geometry.  Many more similar and related conjectures
were stated over the years.  Some of them became famous quickly, and some
were proved quickly, see e.g.\ a celebrated paper by Heilmann and Lieb \cite{HL}
on log-concavity of the \emph{matching polynomial} for a graph.
On the other hand, \defn{Rota's unimodality conjecture} was mentioned in passing
in~\cite{Rota}, reiterated in~\cite[p.~209]{RH}, generalized to log-concavity
by Mason and Welsh, and proved only recently (Theorem~\ref{t:matroids-AHK}).
We refer to \cite[$\S$14.2]{Ox} for a detailed overview of the early
work on the subject.

\subsection{}\label{ss:hist-modern}
In modern times, the algebraic approach was pioneered by Stanley,
who used the \emph{hard Lefschetz theorem} to establish the
\emph{Sperner property} of certain families of posets~\cite{Sta-Lef}.
This easily implied the \emph{Erd\H{o}s--Moser conjecture} and laid
ground for many recent developments.  In fact, Stanley's approach was
itself a rethinking of Sylvester's proof we mentioned above,
see \cite{Sta-Lie}, and it was later deconstructed in~\cite{Pro}.

In the past decade, Huh and coauthors pushed the algebraic approach
to resolve several conjectures which remained open for decades.
They established the hard Lefschetz theorem and the Hodge--Riemann
relations in a number of algebraic settings, which imply the
log-concavity results.  We will not attempt to review this work
largely because it is thoroughly surveyed in Huh's ICM survey~\cite{Huh}.
Below is a quick recap of results used directly in this paper.

\subsection{} \label{ss:hist-matroid-basics}
Matroids are often associated with several important sequences, including
the $f$-vector whose components are the numbers \ts $\rI(k)$, and the $h$-vector,
which can be computed by a certain linear transformation of the $f$-vector.
Both are coefficients of specializations of the \emph{Tutte polynomial}
associated with the matroid.  We refer to \cite{Bry,BO} for
the introduction and further references.

\subsection{}\label{ss:hist-matroid}
In their celebrated paper~\cite{AHK}, Adiprasito, Huh and Katz proved the log-concavity
of the \emph{characteristic polynomial} of a matroid, which is a generalization
of the graph chromatic polynomial, and a specialization of the Tutte polynomial.
They deduce the Welsh--Mason Conjecture~\eqref{eq:matroid-LC} indirectly, via
an observation by Brylawski~\cite{Bry1} (see also~\cite{Lenz2}).
This culminated a series of previous papers \cite{Huh0,Huh1,HK} on the
subject (see also~\cite{AS-Whitney}).

The inequality~\eqref{eq:matroid-ULC} is the strongest of the \emph{Mason's
conjectures}~\cite{Mas}.  This inequality was recently proved independently
by Br\"and\'en and Huh~\cite{BH18,BH}, and by Anari~et.~al~\cite{ALOV} in
the third paper of the series.  These papers use interrelated ideas, and
avoid much of the algebraic technology in~\cite{AHK}.  Let us mention
a notable application in~\cite{ALOV-RW} which proved that the base exchange
random walk mixes in polynomial time. This was yet another long standing
open problem in the area~\cite{FM}.

\subsection{} \label{ss:hist-NBC}
Brylawski \cite[$\S$6]{Bry} and Dawson \cite[Conj.~2.5]{Daw}
conjectured that matroid $h$-vectors are log-concave.  This was
resolved in \cite{ADH} and~\cite{BST}.
The latter paper proves
a stronger version of log-concavity, while the former
proves further results for the \emph{no broken circuit} (NBC)
complex, another popular matroid construction, see~\cite{Bry1}.

For how log-concavity of $h$-vectors implies log-concavity
of $f$-vectors, see e.g.\ \cite[Cor.~8.4]{Bre2},
\cite[Prop.~6.13]{Bry}, and \cite[Prop.~2.7]{Daw}\footnote{There
is an unfortunate typo in the statement of the Dawson's proposition.}.
As we mentioned in the Introduction (see~$\S$\ref{ss:intro-matroids}),
Lenz \cite{Lenz1} showed that log-concavity of the $h$-vector
implies \emph{strict log-concavity} of the $f$-vector.
See also~\cite{DKK} for many low-dimensional examples.

\subsection{} \label{ss:hist-matroids-ex}
The matroid in Example~\ref{ex:intro-finite-field} is a special case of a
matroid \emph{realizable} over~$\fq$, see e.g.~\cite[$\S$6.5]{Ox}.
In Example~\ref{ex:intro-Steiner}, we consider a subclass of \defn{paving matroids}
defined as matroids \ts $\Mf$ \ts with
\ts $\girth(\Mf)=\rk(\Mf)$, see~\cite{Welsh}.
Our construction of Steiner matroids follows~\cite{Jer,Kahn}.
Notably, Jerrum considers matroid corresponding to \ts $\Stn(5, 8, 24)$.
We refer to~\cite{Dembo} for more on finite geometries arising in this example.

\subsection{}\label{ss:hist-morphism}
Theorem~\ref{t:morphisms-EH} for morphism of matroids is proved by Eur and
Huh in~\cite{EH}, by extending the approach in~\cite{BH}.  The notion of
the morphisms is quite elegant, and follows a long series of combinatorial
papers of Las~Vergnas on the subject, which includes a definition of the
Tutte polynomial in this case.  We refer to~\cite{EH} for an overview and
many references, and to~\cite{Chm} for the extensive survey of generalizations
of the Tutte polynomial to general topological embeddings.

\subsection{}\label{ss:hist-polymatroid}
Discrete polymatroids are also called \emph{integral polymatroids} in~\cite{Edm},
and appear in the context of discrete convex sets~\cite{Mur} and
integral generalized permutohedra~\cite{Post}.
We refer to~\cite{HH} for their history and algebraic motivation.
Note that discrete polymatroids are explicitly treated in~\cite[$\S$4.1]{Mur} and~\cite{BH} under
the equivalent formulation of \emph{M-convex sets}.  They are a part the
definition of Lorentzian polynomials, so in fact weighted polymatroids and
Lorentzian polynomials are closely related notions.\footnote{The \emph{completely
log-concave polynomials} considered in~\cite{ALOV} are not necessarily
homogeneous and thus more general; they coincide with Lorentzian polynomials
in the homogeneous case, see~\cite[p.~826]{BH}.}
Although Theorem~\ref{t:polymatroids-BH} is not stated in this
form, it follows easily from the results in~\cite{BH}.
Indeed, we need Theorem 3.10 combined with taking derivatives and
limits in proof of Theorem~4.14, and where Theorem~2.10 is
substituted with Theorem~2.30 (all in~\cite{BH}).  The details are straightforward.

\subsection{} \label{ss:hist-polymatroids-ex}
We refer to \cite[$\S$14]{BKP} and~\cite[$\S$12]{Post}, for the background
on hypergraphical polymatroids in Example~\ref{ex:intro-hyper},
and further references.  Note that there
are many notions of ``hypertree'' and ``hyperforest'' available in the literature.
We refer to \cite[$\S$10.2]{GP} for a quick overview, and
to~\cite{Ber} for background on hypergraphs and more traditional definitions.

\subsection{}\label{ss:hist-weight}
The notion of \emph{weight function} originates in statistical physics
and is now standard in probability and graph theory.  In the context of
graph polynomials it comes up in connection to the \emph{Potts model}
which is equivalent to the \emph{random cluster model}.  We refer to
\cite{Sok} for an extensive introduction, and to \cite{Gri} for a
thorough treatment.

\subsection{}\label{ss:hist-improve}
The equality conditions have long emerged an
important counterpart to the inequalities, see e.g.\ \cite{BB,HLP}.
They serve as a key check on the inequality: if the equality occurs
rarely or never, perhaps there is a way to sharpen the inequality
either directly or by introducing additional parameters.
\emph{Strict log-concavity} inequalities are especially
suggestive of possible quantitative results.

For example, in his pioneer paper~\cite{Huh0}, Huh proved the
log-concavity of the \emph{chromatic polynomial} of a graph,
establishing several conjectures going back to~\cite{Read}.
In a followup paper~\cite{Huh1}, Huh proved a strict log-concavity
conjecture of Hoggar~\cite{Hog}.  There are no explicit stronger
bounds implying strict log-concavity in the style of
Theorem~\ref{t:morphisms-refined} and~\cite{BST}.

In the opposite direction, when there are many special cases when the
inequality becomes an equality, the equality conditions are unlikely
to be very precise.  It seems, this is the case of our equality conditions
for matroid log-concavity given in~$\S$\ref{ss:intro-morphism-equality}
(see also~$\S$\ref{ss:intro-antimatroid-equality}).
In the context of this paper, the only nontrivial equality condition
known prior to this work for matroid inequalities is
Theorem~\ref{t:matroids-equality} proved by Murai, Nagaoka and Yazawa
in~\cite{MNY} using an algebraic argument built on~\cite{BH}.

\subsection{}\label{ss:hist-greedoid}
Greedoids were defined and heavily studied by Korte and Lov\'{a}sz
as set systems on which the \emph{greedy algorithm} provably works, thus the name.
They generalize matroids, which in turn generalize graphs, where the greedy
algorithm is classically defined to compute the \emph{minimal spanning tree} (MST).
For general greedoids, the reader should think of the (greedy) \emph{Prim's algorithm}
for the MST in undirected graphs,
rather than \emph{Kruskal's algorithm}, as a starting point of the generalization.
The approach to greedoids in terms of languages goes back to original papers.
We refer to \cite{KLS} for a foundational monograph on the subject, and
to~\cite{BZ} for a relatively short and digestible survey.

\subsection{}\label{ss:hist-antimatroid}
Antimatroids is a subclass of greedoids named after the \emph{anti-exchange property},
which is a key axiom in their definition via set systems \cite[$\S$3.1]{KLS}.
There are many examples of antimatroids coming from graph theory
(e.g.~\emph{branching process}) and
discrete geometry (e.g.~\emph{shelling process}), although
poset antimatroids have a combinatorial nature they also have some
geometric aspects (see e.g.~\cite{KL}).  Much of the terminology in
the area is rather unfortunate and can be somewhat confusing, so we refer
the reader to the top of page 335 in~\cite{BZ}, which defines classes
of greedoids in terms of properties of the corresponding lattices of
feasible sets.  See Figure~\ref{f:diagram} below for the diagram
of relationships between main greedoid classes (see also~\cite[p.~301]{KLS}
for a larger diagram).

\subsection{}\label{ss:hist-antimatroid-ex}
Standard Young tableaux (see Example~\ref{ex:intro-SYT})
are fundamental in algebraic combinatorics.
They play a key role in representation theory of \ts $S_n$ \ts and \ts $\GL(N,\cc)$,
and the geometry of the Grassmannian, see e.g.~\cite{Ful} and \cite[$\S$7]{EC}.
Numbers \ts $f^{\la/\mu}=|\SYT(\la/\mu)|$ \ts have an elegant \emph{Aitken--Feit
determinant formula} \cite{Ait,Feit}, see also \cite[Cor.~7.16.3]{EC}.
For the sequence \ts $\{b_k\}$ \ts in Example~\ref{ex:intro-SYT}, see e.g.~\cite[Ex.~VIII.5]{FS}.

\subsection{}\label{ss:hist-antimatroid-branching}
Enumeration of increasing arborescences (also called \emph{branchings} and
\emph{search trees}) in Example~\ref{ex:intro-branching}
without graphical constraints is common in enumerative combinatorics, see e.g.~\cite{BBL,FS}.
Maximal arborescences (also called \emph{directed spanning trees})
also appear in connection to the \emph{reachability problem} in
network theory, see \cite{BP,GJ}, and can be sampled by the
\emph{loop-erased random walk} and its relatives, see \cite{GP,Wil}.

\subsection{}\label{ss:hist-LE}
%
Linear extensions of a finite poset \ts $\cP$ \ts are in obvious bijection with
maximal chains in the lattice \ts $\LL(\cP)$ \ts of lower order ideals of~$\cP$.
Lattice \ts $\LL(\cP)$ \ts is always distributive, and by
\emph{Birkhoff's representation theorem} (see e.g.\ \cite[Thm~3.4.1]{EC}),
every finite distributive lattice can be obtained that way.  We refer to
\cite{BrW,Tro} for definitions and standard results on posets and linear
extensions.

\subsection{}\label{ss:hist-LE-AF}
Stanley's inequality~\eqref{eq:Sta} was originally conjectured by
Chung, Fishburn and Graham in~\cite{CFG}, extending an earlier unimodality
conjecture by R.~Rivest (unpublished).  The proof in~\cite{Sta} is a simple
application of the \emph{Alexandrov--Fenchel inequality}.
Until now, no direct combinatorial proof of Stanley's inequality was known
in full generality, although \cite{CFG} gives a simple proof for posets
of width two (see also~\cite{CPP2}).  Most recently, the authors and
Panova obtained a $q$- and multivariate analogues of Stanley's
inequality for posets of width two~\cite{CPP2}.  These notions are specific to
the width two case and are incompatible with the weighted analogue
(Theorem~\ref{t:Sta-weighted}) nor the case of posets with belts
(Theorem~\ref{t:Sta-belt-submodular}).

\subsection{}\label{ss:hist-LE-ex}
The connection between linear extensions of two dimensional posets and lower order
ideals of Bruhat order used in Example~\ref{ex:sta-bruhat} has been discovered a number
of times in varying degree of generality, see \cite{BW,FW} (see also~\cite{DP}).
Statistics \ts $\beta : S_n \to \nn$ \ts in that example seems different from
other permutation statistics which appear in the context of log-concavity,
see e.g.\ \cite{Bra,Bre}.

Statistic \ts $\ga$ \ts on the alternating permutations in Example~\ref{ex:intro-EB}
is more classical.  Note, however, a major difference: while much of the
literature studies permutation statistics as polynomials in \ts $\nn[q]$ \ts whose
coefficients can sometimes form a log-concave sequence, we study values of these
polynomials at fixed \ts $q\in \rr$.
For more on the Euler and Bernoulli numbers and the connection between them,
see e.g.~\cite[$\S$IV.6.1]{FS}.  For log-concavity of Entringer numbers
and their generalizations, see~\cite{B+,G+}.

\subsection{}\label{ss:hist-AF}
The Alexandrov--Fenchel inequality is a classical result in convex
geometry which remains mysterious despite a number of different proofs,
see e.g.~\cite[$\S$20]{BZ-book} and~\cite[$\S$7.3]{Sch}.  It generalizes the
\emph{Brunn--Minkowski inequality} to \emph{mixed volumes}, and
has remarkable applications to the \emph{van der Waerden conjecture},
see e.g.~\cite{vL}.  Let us single out one of the original
proofs by Alexandrov using polytopes~\cite{Ale}, the inspirational
(independent) proofs by Khovanskii and Teissier using Hodge theory, see
\cite[$\S$27]{BZ-book} and~\cite{Tei}, a recent
concise analytic proof by Cordero-Erausquin et al.~\cite{CKMS},
and the proof by Shenfeld and van Handel~\cite{SvH19}, which partly
inspired this paper.

\subsection{}\label{ss:hist-equality}
For geometric inequalities such as the \emph{isoperimetric inequalities},
the equality conditions are classical problems going back to antiquity
(see e.g.~\cite{BZ-book,Sch}).  In many cases, the equality conditions
are equally important and are substantially harder to prove than
the original inequalities. For example, in the
\emph{Brunn--Minkowski inequality}, the equality conditions
are crucially used in the proof of the \emph{Minkowski theorem}
on existence of a polytope with given normals
and facet volumes (see e.g.\ $\S$7.7 and $\S$36.1 in~\cite{Pak}).
For poset inequalities, the equality conditions
are surveyed in~\cite{Win}.

\subsection{}\label{ss:hist-LE-equality}
The equality conditions for Stanley's inequality for the case when $\ap$ is uniform
(Theorem~\ref{t:Sta-equality}), were recently obtained by Shenfeld and van Handel in
\cite[Thm~15.3]{SvH}.  They used a sophisticated geometric analysis to prove
equality conditions of the Alexandrov--Fenchel inequality for convex polytopes.
%
%
We should mention that part~\eqref{item:equality Stanley d}
of Theorem~\ref{t:Sta-equality-weighted-more} is inspired by our results in
\cite[$\S$8]{CPP2} for the Kahn--Saks inequality, which in some sense are
more general.  Finally, the $q$-analogue of the equality condition
was obtained in the same paper \cite[Thm~1.5]{CPP2} for posets of width two.
}

\bigskip

{\small

\section{Final remarks and open problems}\label{s:finrem}

\subsection{}\label{ss:finrem-first}
Unimodality is so natural, sooner or later combinatorialists start seeing
it everywhere, generating a flood of conjectures.  In the spirit of the
``strong law of small numbers'' \cite{Guy}, many such conjectures do in
fact hold in small examples but fail in larger cases.  Sometimes, it takes
years of real or CPU time until large counterexamples are found
(see e.g.~\cite{RR}), in which case they are published.
Notable unimodality disproofs can be found in \cite{Bjo,Stanton,Ste},
all related to poset inequalities in some way.

Log-concavity is a stronger property than unimodality, but is also
more natural.  Indeed, in the absence of symmetry there is no
natural location of the mode (maximum) of the sequence.  While the
mode location is critical in establishing unimodality, it is irrelevant
for log-concavity.  Moreover, as was pointed out in \cite[p.~38]{RT},
log-concavity of  polynomial coefficients is preserved under
multiplication of polynomials, an important property of
poset polynomials.  Similarly, it was shown in~\cite{Lig}
(see also~\cite{Gur}), that ultra-log-concavity is preserved
under \emph{convolution}, yet another property of some poset polynomials.

\subsection{}\label{ss:finrem-Rota}
In his discussion of influence of Rota on matroid theory,
Kung writes that Rota was motivated in his unimodality conjecture
(see $\S$\ref{ss:hist-main}) in part by the mixed volumes
which are ``somewhat analogous'' to the Whitney numbers, see \cite[$\S$3.1]{Kung}.
This seems extremely prescient from the point of view of this
paper, as we prove matroid log-concavity with a technology that
originates in the ``right'' proof of the Alexandrov--Fenchel
inequality.  One could argue that we inadvertently fulfilled
Rota's unstated prediction (cf.~\cite{AS-Whitney}).

\subsection{}\label{ss:finrem-posets}
As we mentioned in the introduction, traditionally matroids are viewed as
a subclass of lattices, see e.g.~\cite{Ox,Welsh}. Similarly, greedoids are usually
defined by their feasible sets a more general subclass of posets
(cf.~$\S$\ref{ss:hist-greedoid}).  Thus, the title of the paper.

\subsection{}\label{ss:finrem-main}
Our proofs in Section~\ref{s:atlas} borrows heavily from~\cite{SvH19},
although they are written in a very different language (see also
Remark~\ref{r:Hyp-claim}). According to the authors, the idea of this proof
can be traced back to the work of Lichnerowicz~\cite{Lic}, see \cite[$\S$6.3]{SvH19}
for a further discussion.

The proof of Theorem~\ref{t:equality Hyp} is a modification on the argument
in~\cite{Ale}, which in turn is based on~\cite{Weyl}.  In the draft of
the paper, we were not aware of the connection and used a similar but
longer argument.  This simplification was kindly proposed to us by
Ramon van Handel (personal communication).

\subsection{}\label{ss:finrem-interlace}
In the proof of Theorem~\ref{t:greedoid} given in $\S$\ref{ss:greedoid-proof},
at a critical step (in the base of induction), we employed \emph{Cauchy's interlacing theorem}.
In fact, interlacing of eigenvalues is surprisingly powerful,
see e.g.\ \cite{Hua,MSS} for notable recent applications.

\subsection{}\label{ss:finrem-AF-proof}
As we mentioned earlier, our proof of Theorem~\ref{t:Sta-weighted} is inspired
by the approach of Shenfeld and van Handel~\cite{SvH19}.  Indeed, the mixed
volumes in Alexandrov--Fenchel inequality can be converted into inner products
in~\eqref{eq:Hyp}, where the vectors are given by the support functions of
the polytopes.  We present this proof in~\cite{CP}.
Technically, one can object that we assumed the diagonal
entries of \ts $\bM$ \ts are assumed to be nonnegative.  In fact, this assumption
is made for convenience as nonnegativity holds in our examples, but allowing
\ts $\aM_{ii}$ \ts to be negative does not change the proof.

Now, it is shown in \cite[$\S$5]{SvH19}, that the corresponding
matrices and vectors for simple, strongly-isomorphic polytopes
satisfy all conditions of Theorem~\ref{t:Hyp}.  Note that in that setting~\eqref{eq:Pull}
is always an equality, see \cite[Eq.~(1.2)]{SvH19} and \cite[Eq.~(5.23)]{SvH} for the proof.
On the other hand, the inequality~\eqref{eq:Pull} can  be strict in our setting.
The comparison between our proof Theorem~\ref{t:Sta-equality-weighted-more} and
the proof of Theorem~\ref{t:Sta-equality} in~\cite{SvH} is also curious
and we don't fully understand it.  We should mention the crucial use of
the opposite poset \ts $\cP^{\op}$, which does not seem to show up
in this context.  It would be interesting to find further applications
in this ``duality'' approach (cf.~$\S$\ref{ss:finrem-reversible}).

\subsection{}\label{ss:finrem-matroid-equality}
Although Theorem~\ref{t:matroids-equality} says that there are no interesting examples
of equality of log-concavity for matroids, the examples in~$\S$\ref{ss:intro-matroid-examples}
suggest that the family of matroids with equality in Theorem~\ref{t:matroids-Par-weighted}
is rather rich.  While our Theorem~\ref{t:matroids-equality-weighted} gives some natural
necessary and sufficient conditions, it would be interesting to see if this
description can be used to obtain a full classification of such matroids
in terms of known classes of set systems.

\subsection{}\label{ss:finrem-could}
Our work is completely independent of the algebraic approach
in~\cite{AHK}, yet some glimpses of similarity to more recent developments
are noticeable if one squints hard enough.  For instance, we need the element
\ts $\spstar$ \ts in the proof of Theorem~\ref{t:greedoid},
for roughly the same technical reason that papers \cite{BHMPWa,BHMPWb} need to
use the \emph{augmented Bergman fan} in place of the (usual) \emph{Bergman fan}
employed in~\cite{AHK}.

\subsection{}\label{ss:finrem-would}
The connection between our proof and Lorentzian polynomial approach is
somewhat indirect to make any formal conclusions.  On the one hand,
we can use combinatorial atlases to emulate everything Lorentzian
polynomial do \cite{CP}.  On the other hand, the atlas we construct for matroids
and polymatroids is sufficiently flexible to allow our refined inequalities.
On a technical level, in notation of Section~\ref{s:Pull}, the matrix
\ts $\bK=(\aK_{ij})$ \ts
which arises when we emulate Lorentzian polynomials, is always zero
(cf.~Remark~\ref{r:Pull-K-PSD}).
Thus, it would be interesting to see if the tools in \cite{ALOV} and \cite{BH} can
be modified to yield our Theorems~\ref{t:matroids-Par-weighted} and~\ref{t:polymatroids-Par}.

\subsection{}\label{ss:finrem-BL}
Most recently, Br\"and\'en and Leake showed in~\cite{BL} how to obtain the log-concavity
of the characteristic polynomial of a matroid using a purely Lorentzian polynomial
approach, avoiding the use of algebra altogether.  While it is too early to say,
we intend to see if the combinatorial atlas technology can be combined with
that approach.

\subsection{}\label{ss:finrem-AF}
It would be interesting to see if one can derive Theorems~\ref{t:Sta-weighted}
from the Alexandrov--Fenchel inequality.  If this is possible, do the tools
in~\cite{SvH} extend to prove Theorem~\ref{t:Sta-equality-weighted}?

\subsection{} \label{ss:finrem-graph-32}
In the Example~\ref{ex:intro-graphical}, the asymptotic constant \ts $3/2$ \ts
is probably far from tight for \emph{dense graphs}, say with \ts $\Omega(\sn^2)$ \ts
edges.  What's the right constant then?

\subsection{} \label{ss:finrem-greedoid}
When it comes to interval greedoids, there are more questions than answers.
For example, since there is a Tutte polynomial for greedoids defined in~\cite{GM},
does it make sense to define an NBC complex?  Are there
any log-concavity results for characteristic polynomials in some
special cases?  Can one define morphism of antimatroids or
interval greedoids?  Are there any other interesting classes
of interval greedoids whose log-concavity is worth studying?

\subsection{}\label{ss:finrem-local}
\emph{Weak local greedoids} introduced in~$\S$\ref{ss:main-properties}
by the weak local property~\eqref{eq:LI}, is a new class of
greedoids. It contains poset antimatroids,
matroids, discrete polymatroids, and \emph{local poset greedoids},
see Figure~\ref{f:diagram}.
We do not consider the latter in this paper, but they play an
important role in greedoid theory, see~\cite[Ch.~VII]{KLS}.
To understand the relationship between weak local greedoids
and local poset greedoids, note the excluded minor
characterization of local poset greedoids in \cite[Cor.~VII.3.2]{KLS}.
By contrast, weak local greedoids exclude the same minor under
contraction, but not necessarily deletion.

\begin{figure}[hbt]
		\includegraphics[height=6.75cm]{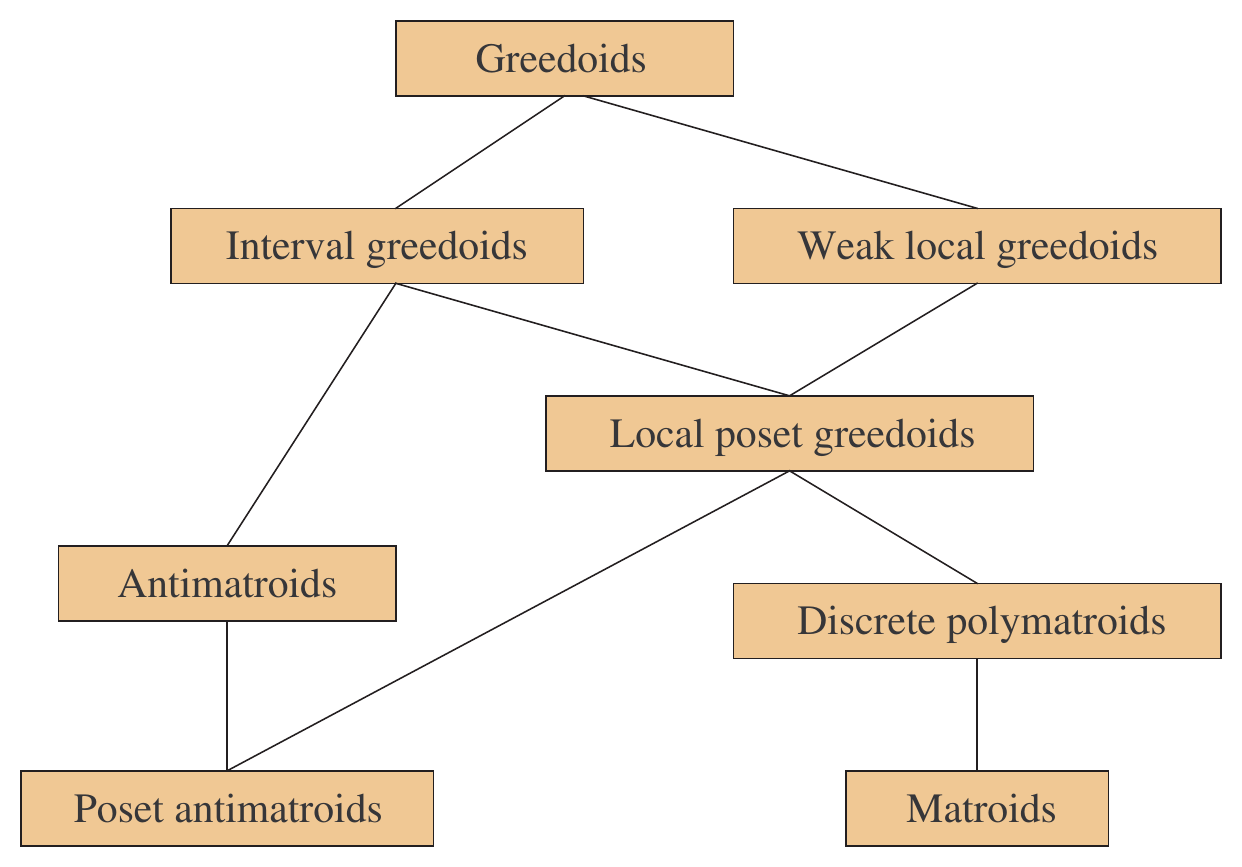}
		\vskip-.05cm
		\caption{Diagram of inclusions of greedoid classes.}
		\label{f:diagram}
\end{figure}

\subsection{}\label{ss:finrem-reversible}
Let \ts $\Gf=(X,\Lc)$ \ts be an interval  greedoid, and let \ts
$\Bc \subseteq \Lc$ \ts be the set of feasible words
\ts $\al = x_1 \cdots x_\ell$ \ts
of maximum length \ts $\ell=\rk(\Gf)$.
Denote by \ts $\Bc^{\op}$ \ts the set of words \ts
$\al^{\op} := x_\ell  \cdots x_1$.
An interval greedoid is called \defn{reversible} if \ts $\Bc^{\op}$ \ts
is the set of basis feasible words of an interval greedoid.
Note that matroids, polymatroids and poset antimatroids are
examples of reversible greedoids.

Let us note that our proof of Stanley's inequality~\eqref{eq:Sta}
can be generalized to
reversible interval greedoids.  Unfortunately, in the examples above
the corresponding generalization of Stanley's inequality is trivial.
It would be interesting to characterize reversible greedoids or at least find
new interesting examples.

\subsection{} \label{ss:finrem-complexity}
From the \emph{computational complexity} point of view, one can
distinguish ``easy inequalities'' from ``hard inequalities'',
depending whether the components (or their differences) are
computationally easy or hard.  For example, \emph{Hoffman’s bound}
(see e.g.~\cite[Thm~8.8]{Big}),
relates the \emph{independence number of a graph} which is \ts
$\NP$-hard, to the ratio of graph eigenvalues
which can be computed in polynomial time.  Assuming \ts $\poly \ne \NP$,
one would expect such bound not to be sharp in many natural cases.

By contrast, Alon's lower bound on the number of spanning trees
in regular graphs (see~\cite{Alon}) has both sides computable in
polynomial time.  This suggests that complexity approach may not
necessarily capture the mathematical difficulty of the result.

In this context, the inequalities in this paper are the ``hardest''
of all.  For the Mason's conjectures, even in the simplest case
of graphical matroids (Example~\ref{ex:intro-graphical}), the number
of $k$-forests is known to be $\SP$-complete,
see e.g.~\cite{Welsh1}.  Similarly, in Stanley's inequality
(Theorem~\ref{t:Sta}), the number of linear extensions of
a poset is \ts $\SP$-complete even for posets of height
two or dimension two, see~\cite{BW1,DP}.

\subsection{}\label{ss:finrem-complexity-GapP}
Another computational complexity approach to combinatorial inequalities
is to understand whether their difference of two sides is nonnegative
for combinatorial reason, i.e.\ whether it has a combinatorial
interpretation.  This is a natural question we previously
discussed in~\cite{Pak}.

For example, observe that both sides in Stanley's
inequality~\eqref{eq:Sta} are \ts $\SP$-functions, i.e., they have
a natural combinatorial interpretation.  The difference of LHS and RHS
is then a function in \ts $\GapP=\SP-\SP$.  Now the problem
whether it lies in~$\ts\SP$.  Although our proof is elementary, this
question remains unresolved.

Similarly, in the case of graphical matroid, the equation~\eqref{eq:matroid-LC}
also corresponds to a nonnegative function in \ts $\GapP$.  Again,
no combinatorial interpretation is known in this case.  This is in sharp
contrast, e.g., with the Heilmann--Lieb theorem (see~$\S$\ref{ss:hist-main})
on log-concavity of the matching polynomial, where a combinatorial
interpretation of the difference follows from Krattenthaler's
combinatorial proof~\cite{Kra}, see also~\cite{Pak}. We intend
to return to this problem in the future.\footnote{After this paper
appeared we continued our investigation in \cite{IP,Pak-what}.}

\subsection{}\label{ss:finrem-conj-Brenti}
Going back to the discussion in Foreword~$\S$\ref{ss:intro-for} and
Final Remark~$\S$\ref{ss:finrem-first} above, it seems, the importance of poset
log-concavity conjectures is yet to be settled.  Back in 1989, Francesco
Brenti wrote in this context:

\smallskip

\begin{center}\begin{minipage}{11cm}%
{\emph{``In this author's opinion, conjectures and open
problems in mathematics are not so much interesting and important `per se' but because
they are symptoms that our knowledge is not complete in some area. Their greatest value is
not whether they are true or false but that they stimulate and lead us into deeper
knowledge.''}~\cite[p.~6]{Bre}}
\end{minipage}\end{center}

\smallskip

\nin
One can disagree with these sentiments, but speaking for ourselves we
certainly owe these conjectures a debt of gratitude, as we find ourselves
in the midst of unexplored territory we neither sought nor expected to discover.

\vskip.7cm
	
\subsection*{Acknowledgements}
We are grateful to Marcelo Aguiar, Noga Alon, Nima Anari, Andy Berget,
Olivier Bernardi, June Huh, Jeff Kahn, Jonathan Leake, Alejandro Morales,
Yair Shenfeld, Hunter Spink and  Cynthia Vinzant for helpful discussions
and remarks on the subject.
Special thanks to Ramon van Handel for many insightful
comments on the draft of the paper, for his detailed suggestions on
how to streamline the proof of Theorem~\ref{t:equality Hyp}, and
for help with the references.
We are also thankful to Christian Ikenmeyer and Greta Panova
for numerous interesting conversations on poset inequalities,
and to Joseph Kung for the personal story~\cite{Kung-letter}
behind~$\S$\ref{ss:finrem-Rota}.
We give no thanks to COVID-19, which made writing this paper miserable.
The first author was partially supported by the Simons Foundation.
The second author was partially supported by the~NSF.

}

\vskip1.1cm


\end{document}